\documentclass[twoside, notables, nofigures, noexamples, nocopyright, nosignature]{brownthesis}
\usepackage{amssymb, latexsym, amsmath, amscd, amsthm, rotating}
\usepackage[mathscr]{eucal}
\usepackage[dvips]{xypic}
\xyoption{curve}
\xyoption{rotate}

\newtheorem{thm}{Theorem}[section]
\newtheorem{lemma}[thm]{Lemma}
\newtheorem{prop}[thm]{Proposition}
\newtheorem{cor}[thm]{Corollary}

\newtheorem*{thmI}{Theorem I}
\newtheorem*{thmII}{Theorem II}
\newtheorem*{thmIII}{Theorem III}
\newtheorem*{thmIII part2}{Lemma \ref{thmIII part2}}
\newtheorem*{thmIII part3}{Lemma \ref{thmIII part3}}

\theoremstyle{definition}
\newtheorem{defn}[thm]{Definition}
\newtheorem{ex}[thm]{Example}

\newtheorem{note}[thm]{Remark}

\newdir{|>}{%
 !/4.5pt/@{|}*:(1,-.2)@^{>}*:(1,+.2)@_{>}}

\def\G{{\scshape g}}
\def\DG{{\scshape dg}}
\def\GC{{\scshape gc}}
\def\GL{{\scshape gl}}
\def\DGC{{\scshape dgc}}
\def\DGA{{\scshape dga}}
\def\DGL{{\scshape dgl}}

\def\mDGL{\text{\scshape DGL}}

\def\mDGC{\text{\scshape DGC}}
\def\mDG{\text{\scshape DG}}
\def\mGL{\text{\scshape GL}}
\def\mGC{\text{\scshape GC}}
\def\mG{\text{\scshape G}}
\def\dgv{\mathcal{DG}}
\def\dg{\mathcal{DG}}
\def\g{\mathcal{G}}
\def\gl{\mathcal{GL}}
\def\gc{\mathcal{GC}}
\def\dgc{\mathcal{DGC}}
\def\dgl{\mathcal{DGL}}
\def\fdgl{\mathcal{FDGL}}
\def\fdgc{\mathcal{FDGC}}
\def\hfdgl{{{\widehat{\mathcal{F}}}\mathcal{DGL}}}
\def\hfdgc{{{\widehat{\mathcal{F}}}\mathcal{DGC}}}
\def\Fdgl{\mathit{f}\mathcal{DGL}}
\def\Fdgc{\mathit{f}\mathcal{DGC}}
\def\hFdgl{{{\widehat{\mathit{f}}}\mathcal{DGL}}}
\def\hFdgc{{{\widehat{\mathit{f}}}\mathcal{DGC}}}

\def\sp{{\mathcal{S}\!\mathit{pec}}}
\def\top{{\mathcal{T}\!\!\mathit{op}}}
\def\sets{{\mathcal{S}\!\mathit{ets}}}
\def\svect{{\mathit{s}\mathcal{V}\!\mathit{ect}}}
\def\sset{{\mathit{s}\mathcal{S}\!\mathit{et}}}
\def\sgp{{\mathit{s}\mathcal{G}\!\mathit{rp}}}
\def\scha{{\mathit{s}\mathcal{C}\!\mathcal{H}\!\mathcal{A}}}
\def\slie{{\mathit{s}\!\mathcal{L}\!\mathit{ie}}}

\def\catC{\mathcal{C}}
\def\catB{\mathcal{B}}
\def\catM{\mathcal{M}}
\def\catN{\mathcal{N}}
\def\catI{\mathcal{I}}
\def\catJ{\mathcal{J}}
\def\Funct{{\mathcal{F}\!unct}}

\def\Sinf{{\Sigma^\infty}}
\def\Linf{{\Omega^\infty}}
\def\hS{{h\Sigma}}
\def\weq{\simeq}

\def\Q{\mathbb{Q}}
\def\freeL{\mathbb{L}}
\def\fibr{\twoheadrightarrow}
\def\cofibr{\hookrightarrow}

\def\sur{\Sigma}

\def\Id{{\rm 1\kern -2.65pt {l}}}
\def\Idm{{\rm 1\kern -2.25pt {l}}}

\def\L{{\mathscr{L}}}
\def\C{{\mathscr{C}}}
\def\ab{{\mathrm{ab}}}
\def\pr{{\mathrm{pr}}}
\def\liesum{{\circledast}}
\def\B{{\mathrm{B}}}
\def\bi{{\mathrm{bi}}}

\def\thcof{\operatorname*{thcof}}

\def\thfib{\operatorname*{thfib}}

\def\hfib{\operatorname*{hofib}}
\def\hofib{\operatorname*{hofib}}

\def\hocof{\operatorname*{hocof}}

\def\tildeotimes{\,\tilde\otimes\,}
\def\tildeoplus{\,\tilde\oplus\,}
\def\coker{\operatorname{coker}}
\def\hocolim{\operatorname*{hocolim}}
\def\holim{\operatorname*{holim}}
\def\colim{\operatorname*{colim}}
\def\tot{\operatorname*{Tot}}

\def\trunc#1{\underline{#1}}
\def\red#1{\underleftrightarrow{#1}}

\begin{document}
\title{Rational Homotopy Calculus of Functors}
\author{Ben C. Walter}
\dept{Mathematics}
\degrees{B.~A., Rice University, 1998\\
         Sc.~M., Brown University, 2000}
\principaladvisor{Tom Goodwillie}
\submitdate{July 2005}
\copyrightyear{2005}
\reader{Bruno Harris}
\reader{Mark Behrens (MIT)}
\dean{Karen Newman}

\abstract{
We construct a homotopy calculus of functors in the sense of Goodwillie 
for the categories of rational
homotopy theory.  More precisely, given a homotopy functor between any
of the categories of differential graded vector spaces ($\dg$), reduced differential
graded vector spaces, differential graded Lie algebras ($\dgl$), and differential 
graded coalgebras ($\dgc$), we show that there is an associated approximating rational
Taylor tower of excisive functors.  The fibers in this tower are homogeneous
functors which factor as homogeneous endomorphisms of the category of differential
graded vector spaces.
\vskip 11pt
Furthermore, we develop very straightforward and simple models for all of 
the objects in this tower.  Constructing these models entails first building 
very simple models for homotopy pushouts and pullbacks in the  
categories $\dg$, $\dgl$, and $\dgc$.  Also, we point out that the category $\dg$ 
is equivalent to the stabilizations of the categories $\dgl$ 
and $\dgc$.  Derived from our models for homotopy pushouts and pullbacks in 
$\dgl$ and $\dgc$ are models for suspensions and loops in these categories.  These
functors in turn induce natural stabilization and infinite loop functors
between the categories $\dgl$ (and $\dgc$) and $\dg$. 
\vskip 11pt
We end with a short example of the usefulness of our computationally simple
models for rational Taylor towers, as well as a preview of some further 
results dealing with the structure of rational (and non-rational) Taylor 
towers.
}

\beforepreface

\abstractpage

\afterpreface

\chapter{Introduction}

\markboth{Ben Walter}{0.1  Introduction}

The superficial goal of this work is the construction of a homotopy calculus of functors in the algebraic setting of rational homotopy theory.  Numerous authors have already constructed such calculi on various different categories in various different settings (many of them even pleasantly algebraic -- see [JM03a], [JM03b], [JM04], and [M02]).  However, our motivation in the current work drastically differs from those of other authors.  In particular, homotopy calculus of functors has generally been constructed in different categories in order to better understand the structures of the categories in question.  Rational homotopy theory is already fairly well understood; our goal instead is to leverage our great understanding of rational homotopy theory in order to better understand homotopy calculus of functors itself.  The current work is meant to be the first in a series which goes on to analyse the structure of the rational homotopy calculus of functors constructed here and then to extend this analysis to the much more complicated case of Goodwillie's homotopy calculus of functors on topological spaces.

Since our eventual goal is to leverage our understanding of the structures of rational homotopy theory in order to better understand calculus of functors, the first half of the current work is largely devoted to building, organizing, and analysing these structures.  In particular we wish to have very simple, explicit constructions of homotopy pullbacks and pushouts in rational homotopy theory as well as knowing that certain homotopy limits and colimits commute.  Also we are interested in having a single associated stable category of ``rational spectra.''  As a result, the first half of our work is essentially a very long and extensive exercise in homological algebra and model category theory.  

In the second half, we outline a construction of homotopy calculus of functors in the realm of rational homotopy theory.  Using our nice models for loop, suspension, and homotopy pushouts and pullbacks from the first half, we are able to build very simple and pleasant models for approximating towers of rational homotopy functors.  Much of this is actually quite standard and formally follows from Goodwillie's methods of [GIII].  The main departure which we make from Goodwillie occurs in our sections dealing with rational spectra.  We do not define our category of rational spectra to be the immediate stabilization of either $\dgl$ or $\dgc$, in particular our ``$\Sinf$'' and ``$\Linf$'' functors are not the canonical functors between a model category and its stabilization.  A few slight modifications to Goodwillie's methods are required to get around this.  
We end with an easy application showing how having our particularly simple models can greatly simplify computations of rational Taylor towers.  Also we give a short preview of some further results in the continuation of our program of study.

As well as setting the stage for our later computations and results, our desire is for the current work to serve as both an introduction to the basics of rational homotopy theory for topologists as well as an introduction to the basics of calculus of functors for rational homotopy theorists.  To this end, we have attempted to add enough background to our discussions to enable ready understanding for the neophyte in either theory.  Our own experience has been that homotopy calculus becomes much more clear viewed through the lens of rational homotopy theory, and rational homotopy theory becomes much more clear when it is dealt with in a very topological way (with model categories, homotopy functors, products, coproducts, loops, suspensions, and the like).

\section{Homotopy Calculus of Functors}

Homotopy calculus of functors is a way of analyzing homotopy functors.
A functor $F:\catM\to\catN$ between two model categories is called a homotopy functor if it preserves weak equivalences.
Given a homotopy functor between two suitably nice pointed model categories, say $F:\top_* \to \top_*$, Goodwillie's homotopy calculus of functors (c.f. [GI], [GII], and [GIII]) constructs an inverse system (standardly called a ``tower'') of objectwise fibrations of homotopy functors 
$$\dots\to P_nF \to \dots \to P_1F \to P_0F$$
 along with natural maps $F\to P_nF$ which, under good circumstances, are of connectivity increasing with $n$ (in which case we say that the tower ``converges'' to $F$).  The fibers in this tower are generally written $D_nF$, and referred to as the {\it $n$th homogeneous layers}.\footnote{$D_nF$ is the fiber of the map $P_nF \to P_{n-1}F$}  The term {\it calculus of functors} comes in part from a curious analogy -- the fibers $D_nF$ (in the case of $F:\top_* \to \top_*$) have the form 
$$D_nF \simeq \Linf\bigl(\partial^n\!F(\ast) \wedge (\Sinf X)^{\wedge n}\bigr)_{\hS_n}$$
up to natural weak equivalence (where $\partial^n\!F(\ast)$ is a spectrum with $\Sigma_n$ action).  This is analogous to the $n^\mathrm{th}$ summand 
$$\bigl(\partial^n\!f(0) \cdot x^n\bigr)\,/\,n!$$
 of the Taylor series approximation of the function $f$ in ordinary calculus.  Furthermore, each functor $P_nF$ in the tower of $F$ is characterized by a property analogous to a property characterizing degree $n$ polynomial functions (and $P_nF$ is furthermore universal among all functors satisfying this property and having a (weak) map from $F$).

In keeping with this analogy, the functors $P_nF$ are referred to as the \emph{polynomial approximations} of the functor (and considered to be analogous to the $n^\mathrm{th}$ partial sums of the approximating Taylor series of a function),\footnote{The fibers $D_nF$ are called ``homogeneous'' because they are analogous to the homogeneous polynomial summands of the Taylor series approximation of a function.} and the spectra $\partial^n\!F(\ast)$ are referred to as the \emph{derivatives} of the functor (analogous to the derivatives of a function).

Of particular interest is the Taylor tower of the identity functor $\Id:\top_\ast \to \top_\ast$.  Unlike in ordinary calculus, the identity functor in homotopy calculus is not linear -- in fact, it has an infinite approximating tower nontrivial at every level.  The spectral sequence associated to the approximating tower of the identity functor evaluated at a space $X$ begins with the homotopy of a sequence of infinite loopspaces ($D_n\Id(X)$) and converges to the homotopy of $X$.  The tower of the identity functor thus becomes a useful tool for the understanding and computation of homotopy groups.
The derivatives of the identity functor were first computed by Johnson in her thesis [J95] and a reformulation of her answer was later given by Arone and Mahowald [AM99].  In particular, Arone and Mahowald were able to gain information about the groups, $\pi_\ast(D_\ast\Id(S^n))$ which make up the tower's spectral sequence in the case $X = S^n$.  Even without knowing the differentials in the spectral sequence, they were then able to draw a number of conclusions about the periodic homotopy of spheres.

Note that the constructions and universal properties of the approximating Taylor towers of functors $F$ are all ``derived'' constructions and properties.  That is, given a homootopy functor $F:\catM\to\catN$, the constructions of Goodwillie build an approximating tower of homotopy functors $\catM\to\catN$, but the universal properties which characterize this approximating tower are all properties in the homotopy category of homotopy functors $\catM\to\catN$ -- in fact, the tower only even approximates $F$ in the homotopy category.\footnote{The situation is similar to that of homotopy limits and colimits.}  In particular, any other tower of functors which is naturally weakly equivalent to the tower constructed by Goodwillie is also an approximating tower of the functor.  Also, any two towers related by zig-zags of natural weak equivalences (for example the approximating towers of two different functors which are related by a zig-zag of natural weak equivalences) are viewed as being the same.  

When doing computations, therefore, there comes a time when one must face the question of which models to choose.  For example, it matters which models one uses for homotopy limits and colimits when constructing approximating towers, because different models for $\holim$ and $\hocolim$ yield different models for the polynomial approximations in the tower.  There is no difference between the homotopy of the different models, but there may be a vast difference between the ease of computation with (and understanding of) the different models.  Similarly the properties characterizing the derivatives of a functor $\partial^n\!F(\ast)$ only determine them up to (weakly) equivariant weak equivalence.  Goodwillie gives a specific construction which yields $\Sigma_n$-spectra satisfying the universal properties of the derivatives of a functor; but again the result depends on the models for $\holim$ and $\hocolim$ which were chosen as well as (for more general model categories $\catM$) the method that has been chosen to stabilize $\catM$ in order to create a category of associated spectra, and the models chosen for the functors $\Sinf$ and $\Linf$ to and from this associated category of spectra.

The desire to use the simplest possible models in order to gain the computationally cleanest formulations drives most of our constructions in this work.

\section{Rational Homotopy Theory}

Rational homotopy theory is the study of spaces up to rational homotopy equivalence.  The theory's roots lie in a geometric construction by Sullivan in the 1960's
that, given a simply connected space $X$, builds a rationalization of that space $X_\Q$, whose homotopy and homology groups are those of $X$ modulo all torsion.  Furthermore, given a continuous map of simply connected spaces $f:X\to Y$, he builds a rationalization $f_\Q:X_\Q\to Y_\Q$.  In more modern terms, there is a localization functor $L_{H\Q}$ which localizes $\top$ with respect to the homology theory $H_\ast(-;\Q)$.  
The rational homotopy of spheres is very simple -- odd spheres have only one nonzero rational homotopy group, and even spheres have only two (a generator and its Whitehead square).  A result of this is that one is able to compute the rational homotopy of a {\scshape cw}-complex from its cellular chain complex.  Dually, the rational homology of the classifying spaces $K(\Q, n)$ is very simple -- for odd $n$ there is only trivial homology group, for even $n$ homology is trivial except in dimensions $kn$.  Due to this, one is able to compute the rational homology of a space from its rational Postnikov tower.

In [Q69] Quillen extends this theory by noting that, from the point of view of homotopy, the category $\top_\Q$ of rational, simply-conected spaces is naturally equivalent to either the category of differential graded coalgebras ($\dgc$) or the category of differential graded Lie algebras ($\dgl$).  The precise form of this equivalence is as a chain of Quillen equivalences between simply connected spaces, $\dgl$, and $\dgc$.  The \DGC\ corresponding to the rational simply-connected space $X_\Q$ has the property that its homology is equal to the homology of $X_\Q$, while the homology of the \DGL\ corresponding to $X_\Q$ is equal to the (shifted) homotopy of $X_\Q$ (with Lie brackets corresponding to Whitehead products).  Replacing $\top_\Q$ by $\dgc$ or $\dgl$ transforms the topological problem of classification of spaces up to homotopy into a purely algebraic problem.  In general, the very difficult problems of homotopy theory (such as homotopy of spheres) tend to become solvable in the realm of rational homotopy theory.

Our conceptual view is that rational spaces $X_\Q$ have two incarnations -- as a \DGL\ and as a \DGC.  The associated \DGL\ of a space is the ``homotopy friendly'' incarnation, the associated \DGC\ of a space is the ``homology friendly'' incarnation.

Computations in rational homotopy theory are generally not made using either \DGC\,s or \DGL\,s but rather using the differential graded algebras of piece-wise linear differential forms of Sullivan (see [FHT] or [GM]).  For finite complexes, \DGA\,s and \DGC\,s are precisely dual and the two approaches are equivalent.  However, we would like to couch our constructions in terms of model category structures and Quillen equivalences.  Also, for our constructions to be as natural and general as possible, we do not wish to be hampered by finiteness concerns.  For these reasons, our work will all occur in the setting of Quillen's model categories $\dgc$ and $\dgl$.

\section{Rational Homotopy Calculus}

We will construct a ``rational'' homotopy calculus for homotopy functors between the model categories $\dg$, $\dg_r$, $\dgl$, and $\dgc$.  The categories $\dgl$ and $\dgc$ are our categories of rational spaces.  The category $\dg$ is our category of ``rational spectra'' -- it is equivalent to the stabilizations of $\dgl$ and $\dgc$ as well as to the $H\Q$-localization of $\sp$ the category of spectra.  Objects in $\dg_r$ are ``connective rational spectra.''

In his overview of homotopy calculus of functors, Kuhn notes that the constructions of Goodwillie in [GIII] formally extend to any simplicial model category which is (left) proper and has very small homotopy limits commuting with filtered homotopy colimits (see [Kuhn \S3]).  The simplicial enrichment is used largely for the convenience of implying the existence of canonical homotopy limit and colimit functors (defined as ends and coends) as well as a canonical suspension functor (simplicial tensor with the simplicial $S^1$) yielding canonical stabilizations (see [Ho01] and [S97]). 

The categories $\dg$, $\dg_r$, $\dgl$, and $\dgc$ all satisfy the requirements Kuhn sets out.  In particular, Hinich constructs in [Hi97 4.1.1 and 4.8] and [Hi01 3.1] simplicial enrichments of the model categories $\dgl$ and $\dgc$.  Following Kuhn's arguments, we can therefore construct a homotopy calculus for functors between any two of the categories $\dg$, $\dg_r$, $\dgl$, and $\dgc$. That is, we may construct universal towers of approximating polynomial functors associated to any functor between $\dg$, $\dg_r$, $\dgl$ and $\dgc$.  Furthermore, if the domain and range of $F$ are the same (say, $F:\catM \to \catM$) then the homogeneous layers of its approximating tower will have the form 
$$D_nF(X) \,\weq\, \Omega^\infty_\catM\bigl(\partial_n\!F\otimes (\Sigma^\infty_\catM X)^{\otimes n}\bigr)_{\hS_n}$$
where $\partial_n\!F$ is a $\catM$-spectrum with $\Sigma_n$-action (i.e. a $\Sigma_n$-diagram in the stabilization of $\catM$ with respect to the canonical suspension), $\otimes$ is the symmetric monoidal product in the category of $\catM$-spectra induced by the smash product of simplicial sets, and $\Sigma^\infty_\catM$ and $\Omega^\infty_\catM$ are the associated canonical Quillen adjoint pair between $\catM$ and the stabilization of $\catM$.

For our purposes this approach suffers from a few deficiencies.  One minor deficiency is that the canonical homotopy limits and colimits in the categories $\dgl$ and $\dgc$ arising from the simplicial enrichment of Hinich are rather difficult to write out and compute with and also are larger than we would like.  In fact in each of the categories $\dg$, $\dg_r$, $\dgl$, and $\dgc$ there is an alternate model for homotopy pushouts and pullbacks in particular which is smaller and easier to work with than the canonical model.

Another deficiency is that using the canonical constructions of Kuhn, it is not clear what relationships will hold between the different calculi constructed.  Our conceptual view is that $\dgl$ and $\dgc$ are merely two faces of the same category, so we would like to have a unified theory which discriminates as little as possible between the functors whose domain is $\dgl$ and those whose domain is $\dgc$.  Similarly for functors whose range is $\dgl$ and those whose range is $\dgc$.  In fact, we would like a theory which is compatible with all of the natural Quillen functors between $\dg$, $\dg_r$, $\dgl$, and $\dgc$.

Finally it is unclear what to say about the structure of homogeneous functors which do not have the same domain and range, since the definitions of the canonical categories of spectra are different for each of $\dg_r$, $\dgl$, and $\dgc$.

Our work wil avoid these concerns by relying on our own explicit constructions of homotopy limits and colimits in the categories $\dg$, $\dg_r$, $\dgl$, and $\dgc$.  Also we will give our own explicit construction of a homotopy calculus of functors between these categories.

\section{History -- Jets and Rational Homotopy Calculus}

Recall that the homogeneous layers of the Taylor tower of a functor $F$ (at the space $\ast$) are each determined by a spectrum $\partial^n\!F(\ast)$ with $\Sigma_n$-action.  In general knowing the symmetric sequence of all of the derivative spectra of a functor $\{\partial^n\!F(\ast)\}_{n\ge 0}$ is not enough to determine the functor's Taylor tower since the tower could have nontrivial $k$-invariants.  This stands in sharp contrast to standard calculus, where (for good functions) knowing the collection of derivatives of a function is equivalent to knowing the Taylor series of the function.  Speaking metaphorically, the failure of homotopy calculus to simplify in this way is due to the existence of ``nontrivial ways of adding'' the homogeneous layers to make approximating polynomials.

It has been conjectured and expected for some time that the $k$-invariants of the Taylor tower of a functor were themselves determined by a series of natural, equivariant ``structure'' maps between (some models for) the derivative spectra of the functor.  While there were simple arguments displaying such maps in the homotopy category, it was not clear how to explicitly construct them.  Also, it was not clear what kind of equivariance properties the structure maps should satisfy or what coherence relations should exist between the structure maps.  

Our convention has been to say the ``jet of $F$''  to mean the (conjectural) natural object consisting of the symmetric sequence of derivatives of the functor $F$ along with all of the necessary extra structure (maps between the derivatives) required to recover the Taylor tower of $F$.  More generally we say the ``$n$-jet'' of $F$ for the jet of $P_nF$.  Note that the structure maps making up the jet of a functor induce all of the differentials in the associated spectral sequence of the approximating Taylor tower of the functor.
Furthermore, it is conjectured that a chain rule in homotopy calculus of functors will then be given in terms of jets by a statement of the form 
$$\text{``}\text{Jet of }(F\circ G) = (\text{Jet of }F) \bigotimes_{\substack{\text{structure}\\ \text{maps}}} (\text{Jet of }G)\text{''}$$

In fact, recent work by Michael Ching [C05] strongly suggests that the information encapsulated by a jet is precisely the information given by a symmetric sequence being a module over the \emph{Lie} operad [MSS 1.13 and 1.28] (along with a little more structure in the non-rational case).

The current work rose out of a long-running project to investigate the existence and properties of the structure maps making up the jet of a functor.  Our desire was to look first in rational homotopy calculus of functors where the objects and towers of homotopy calculus may be given a (relatively) simple, algebriac form.  There are particularly simple constructions of homotopy pullbacks in the category $\dgl$ and homotopy pushouts in the category $\dgc$.  Using these constructions we were able to make very simple models for homogeneous functors $\dgc \to \dgl$ and this allowed us to make great headway in the analysis of polynomial functors $\dgc \to \dgl$.

However, proving that our constructions of homogeneous functors and decompositions of the approximating towers to homogeneous parts were correct required us first to show that our simple models for homotopy pullbacks and pushouts were correct, as well as that certain maps and objects in the approximating tower could be given certain nice models.  To show that our homotopy pullbacks and pushouts were correct, it was easier to use a comparison to homotopy pullbacks and pushouts in $\dg$ than a comparison to the canonical homotopy pullbacks and pushouts in $\dgc$ and $\dgl$.  Furthermore, upon completing our construction and verification of homotopy pullbacks in $\dgl$ and homotopy pushouts in $\dgc$ it became apparent that a slight modification would also give homotopy pushouts in $\dgl$ and pushouts in $\dgc$.  

Proving that our models for the objects and maps in the approximating towers of functors $\dgc\to\dgl$ were correct amounted to stepping through Goodwillie's work, replacing proofs and constructions where necessary to show that our specific models satisfied the required properties.  Upon completing this, it became apparent that the same methods would work to give a pleasantly unified theory of approximating towers for functors $\dgc\to\dgc$, $\dgl\to\dgl$, and $\dgl\to\dgc$ as well.\footnote{Unified in the sense that homogeneous functors between all of the above categories factor simply and explicitly through homogeneous endomorphisms of a single stable category $\dg\to \dg$.}

The current document is composed of all of this preparatory work as well as introductions to the categories $\dg$, $\dgl$, and $\dgc$ and a (brief) introduction to rational homotopy theory.

\chapter{Prerequisites}\label{S:Model Cat}

\markboth{Ben Walter}{0.2  Prerequisites}

In this chapter we will give a brief review of the areas of model category theory most critical to our work.  The following material comes primarily from [Q69], [BK72], [Hir], and [DHKS].  For a more in-depth description of modern model category theory see [Hov], [Hir], and [DHKS].

\section{Limits and Colimits}

We begin by recalling the definition of limits and colimits in a category from [DHKS, \S19].  
Given a category $\catC$ and a small category $\catI$, an $\catI$-diagram in $\catC$ is a functor $\mathscr{D}:\catI \to \catC$.  
\begin{defn}[Limits and Colimits]\label{(co)lim adjoint}Let $\catI$ be a small category.  An {\em $\catI$-limit} functor on $\catC$ is a right adjoint $\lim^\catI_\catC$ to the constant diagram functor $c:\catC \to \catC^\catI$.\footnote{This is sometimes called the ``inverse limit'' $\varprojlim \mathscr{D}$.} 

Dually, an {\em $\catI$-colimit} functor  on $\catC$ is a left adjoint $\colim^\catI_\catC$ to the constant diagram functor $c:\catC \to \catC^\catI$.\footnote{This is sometimes called the ``direct limit'' $\varinjlim \mathscr{D}$.}
\end{defn}
A category is called {\em complete} if $\catI$-limit functors exist for all small categories $\catI$, 
and {\em cocomplete} if $\catI$-colimit functors exist for all small categories $\catI$. 
We will standardly suppress the $\catI$ and $\catC$ in our notation and write simply $\lim$ and $\colim$.  Note that if $\mathscr{D}:\catI\to \catC$ is an $\catI$-diagram in $\catC$ then $\lim(\mathscr{D})$ (if it exists) is the final object in the comma category $(c\downarrow \mathscr{D})$ of objects over the diagram $\mathscr{D}$.  Dually, $\colim(\mathscr{D})$ (if it exists) is the initial object in the comma category $(\mathscr{D}\downarrow c)$ of objects under $\mathscr{D}$. 

It is a standard result that if $\catB$ and $\catC$ are categories with adjoint functors $F:\catB\rightleftarrows\catC:U$ then colimits are preserved by the left adjoint  and limits are preserved by the right adjoint. That is, if $\mathscr{D}:\catI \to \catB$ is a diagram in $\catB$ such that $\colim_\catB(\mathscr{D})$ exists then $\colim_\catC(F\mathscr{D})$ exists and is isomorphic to $F\bigl(\colim_\catB(\mathscr{D})\bigr)$; dually for limits.

A common way of showing that a category $\catC$ has limits and colimits is to show that the category has products and coproducts of all small collections of objects as well as all equalizers and coequalizers.  Limits and colimits are then constructed as follows:
\begin{thm}[\lbrack ML98, V.2.2\rbrack]\label{(co)lim (co)equalizer}Suppose that $\catC$ is a category with small products and coproducts as well as all equalizers and coequalizers.  If $\mathscr{D}:\catI\to \catC$ is a small diagram in $\catC$ then 
\begin{itemize}
\item $\lim(\mathscr{D})$ exists and is given by the equalizer
 $$\xymatrix{
  \lim(\mathscr{D}) \ar[r] & 
  \displaystyle \phantom{\prod_{i\in \catI}\mathscr{D}(i)} \save[]-<0pt,4pt> 
   *{\displaystyle \prod_{i\in \catI} \mathscr{D}(i)}
   \restore \ar@<2pt>[r]^(.45)\phi \ar@<-2pt>[r]_(.45)\psi
  & \displaystyle \phantom{\prod_{(\sigma: j\to k)\in \catI} 
   \!\!\!\!\!\! \mathscr{D}(k)} \save[]-<0pt,5pt> 
   *{\displaystyle \prod_{(\sigma: j\to k)\in \catI} \!\!\!\!\!\!
   \mathscr{D}(k)} \restore
 }$$
where $\phi$ is given by the projections
$pr_\sigma\,\phi:\prod_{i}\mathscr{D}(i) \xrightarrow{}
  \mathscr{D}(j) \xrightarrow{\mathscr{D}(\sigma)} \mathscr{D}(k)$
and $\psi$ is given by the projections
$pr_\sigma\,\psi:\prod_{i}\mathscr{D}(i) \xrightarrow{}
  \mathscr{D}(k)$
\item $\colim(\mathscr{D})$ exists and is given by the coequalizer
 $$\xymatrix{
  \displaystyle \phantom{\coprod_{(\sigma: j \to  k)\in \catI} \!\!\!\!\!\!
   \mathscr{D}(j)} \save[]-<0pt,5pt>
   *{\displaystyle\coprod_{(\sigma: j \to  k)\in \catI} \!\!\!\!\!\!
   \mathscr{D}(j)} \restore \ar@<2pt>[r]^(.55)\phi \ar@<-2pt>[r]_(.55)\psi & 
  \displaystyle \phantom{\coprod_{ i\in \catI} \mathscr{D}(i)}
   \save[]-<0pt,4pt> *{\displaystyle\coprod_{ i\in \catI} \mathscr{D}(i)}
   \restore \ar[r] & 
  \colim(\mathscr{D})
 }$$
where $\phi$ is given by the components 
$\phi_\sigma:\mathscr{D}(j) \xrightarrow{\mathscr{D}(\sigma)} 
 \mathscr{D}(k)\xrightarrow{} \coprod_{i} \mathscr{D}(i)$
and $\psi$ is given by the components
$\psi_\sigma:\mathscr{D}(j) \xrightarrow{} \coprod_{i \in \catI} 
 \mathscr{D}(i)$
\end{itemize}
\end{thm}

\section{Model Categories}

The usual modern definition of a model category is as follows [Hir, 7.1.3]:

\begin{defn}[Model Category] 
A model category is a category $\catM$ along with three distinguished sub-classes of maps called weak equivalences ($\mathscr{W}$), fibrations ($\mathscr{F}$), and cofibrations ($\mathscr{C}$) satisfying the following axioms: 
\begin{enumerate}
\item[M1:] $\catM$ is complete and cocomplete.
\item[M2:] $\mathscr{W}$ satisfies the two out of three property (i.e. if two of $f$, $g$, and $fg$ are in $\mathscr{W}$ then so is the third).
\item[M3:] $\mathscr{W}$, $\mathscr{F}$, and $\mathscr{C}$ are closed under retracts.\footnote{
The map $f$ is a {\em retract} of $g$ if there is a commutative diagram
$$\xymatrix@C=10pt@R=10pt{
X \ar[d]^f \ar[r] \ar@/^7pt/[rr]^{\Idm} & A \ar[d]^g \ar[r] & X \ar[d]^f \\
Y \ar[r] \ar@/_7pt/[rr]_{\Idm} & B \ar[r] & Y
}$$}
\item[M4:] Cofibrations have the left lifting property with respect to trivial fibrations (fibrations which are weak equivalences), and fibrations have the right lifting property with respect to trivial cofibrations.\footnote{
We say $f$ has the {\em left lifting property} with respect to $g$ and $g$ has the {\em right lifting property} with respect to $f$ if for every
commutative square diagram as indicated by the solid arrows below
$$\xymatrix@C=10pt@R=10pt{
A \ar[r] \ar[d]_f & X \ar[d]^g \\
B \ar[r] \ar@{..>}[ur]|{h} & Y }$$
there exists a lift $h:B\to X$.}
\item[M5:] Maps in $\catM$ may be functorially factored both as cofibrations followed by trivial fibrations and as trivial cofibrations followed by fibrations.
\end{enumerate}
\end{defn}

This is a slight strengthening of the axioms originally presented by Quillen in [Q69] in that (M1) requires the existence of all {\em small} limits and colimits rather than just {\em finite} limits and colimits, and (M5) requires that factorizations are {\em functorial} rather than just asking that they exist.

We indicate that a map is a fibration by decorating its arrow as $\fibr$ and we indicate cofibrations with the decoration $\cofibr$.  The ``correct'' way to compare two model categories is with a Quillen adjoint pair of functors.
The following definitions are taken from [DHKS 14.1, 17.3]:

\begin{defn}[Quillen Adjoints]\label{Quillen pair}
An adjoint pair of functors between two model categories $F:\catM \rightleftarrows \catN: U$ is a {\em Quillen adjoint pair} if $F$ preserves cofibrations and trivial cofibrations and $U$ preserves fibrations and trivial fibrations.\footnote{By ``preserves cofibrations'' we mean that the functor takes cofibrations to cofibrations; similarly for ``preserves trivial cofibrations'', ``preserves fibrations'', etc.}
\end{defn}

\begin{note}
If the right (or left) adjoint of an adjoint pair of functors satisfies the above requirement, then the left (or right) adjoint will as well.
\end{note}

An object $X$ in a model category $\catM$ is called {\em fibrant} if the map $X \to 1$ to the final object is a fibration.  $X$ is called {\em cofibrant} if the map $0 \to X$ from the initial object is a cofibration.  The full subcategories of fibrant and cofibrant objects in $\catM$ play a critical role.  The standard notation is to write $\catM_\mathsf{f}$, $\catM_\mathsf{c}$, $\catM_\mathsf{fc}$ for the full subcategories of $\catM$ consisting of all fibrant objects, all cofibrant objects, and all fibrant and cofibrant objects respectively.  If $X$ is an object of $\catM$ we write $X_\mathsf{f}$ and $X_\mathsf{c}$ for the fibrant and cofibrant replacements of $X$ furnished by axiom (M5): $X \xrightarrow{\weq} X_\mathsf{f} \fibr 1$ and $0 \cofibr X_\mathsf{c} \xrightarrow{\weq} X$.

The general goal of imposing a model category structure on a category $\catM$ with weak equivalences is to aid us in discussing the homotopy category of $\catM$.  If $\catM$ is a category with weak equivalences, then its homotopy category is $\mathit{Ho}(\catM)$, the category obtained by formally inverting the weak equivalences.  Given $\catM$ and $\catN$ two model categories with a Quillen adjoint pair of functors between them, their homotopy categories are isomorphic $\mathit{Ho}(\catM)\cong\mathit{Ho}(\catN)$ if the Quillen pair is also a {\em Quillen equivalence}.

\begin{defn}[Quillen Equivalence]\label{Quillen equiv}
A Quillen adjoint pair $F:\catM\rightleftarrows \catN:U$ is a {\em Quillen equivalence} if for all $X\in \catM_\mathsf{c}$ and $Y \in \catN_\mathsf{f}$ we have $FX \to Y$ a weak equivalence if and only if $X\to UY$ is a weak equivalence.
\end{defn}

\section{Homotopy Limits and Colimits}

\subsection{Definition}

Throughout this work we will be very interested in functors which preserve weak equivalences.  We standardly refer such functors as ``homotopy functors.''\footnote{Another common description is to say the functor ``reflects'' weak equivalences.} 

Given a model category $\catM$, we cannot expect limits and colimits to be homotopy functors.  That is, given two $\catI$-diagrams $\mathscr{D}_1, \mathscr{D}_2: \catI \to \catM$ with a natural transformation $f:\mathscr{D}_1 \to \mathscr{D}_2$ which sends objects of $\catI$ to weak equivalences in $\catM$ (such a natural transformation is called a {\em natural weak equivalence}) or more generally a zig-zag of natural weak equivalences $\mathscr{D}_1 \xrightarrow{f_1} \cdots \xleftarrow{f_n}\mathscr{D}_2$ there is no guarantee of a weak equivalence or even a zig-zag of weak equivalences between $\lim(\mathscr{D}_1)$ and $\lim(\mathscr{D}_2)$ or $\colim(\mathscr{D}_1)$ and $\colim(\mathscr{D}_2)$.  

Homotopy limit and colimit functors are defined in such a way as to be the best remedy of this defect.  In particular they are determined by three properties
\begin{defn}[\lbrack DHKS, 19.2\rbrack]\label{D: ho(co)lim} 
A {\em homotopy $\catI$-limit} (and {\em $\catI$-colimit}) functor on a model category $\catM$ is a functor $\holim^\catI_\catM:\catM^\catI \to \catM$ (and $\hocolim^\catI_\catM:\catM^\catI \to \catM$) satisfying:
\begin{enumerate}
\item $\holim^\catI$ (and $\hocolim^\catI$) takes natural weak equivalences of diagrams $\catM^\catI$ to weak equivalences in $\catM$.
\item There is a natural transformation $e:\lim^\catI \to \holim^\catI$ (and $e:\hocolim^\catI \to \colim^\catI$).
\item $\holim^\catI$ is {\em homotopically initial} (and $\hocolim^\catI$ is {\em homotopically final}) among all functors satisfying (1) and (2).
\end{enumerate}
\end{defn}

By {\em homotopically initial} (dually homotopy final) in the above definition we essentially mean ``initial up to zig-zags of natural weak equivalences''.  More specifically:

Given a category $\catC$, recall that $c_0\in\catC$ is {\em initial} if there is a natural transformation of functors $\catC \to \catC$ 
$$\mathrm{cst}_{c_0} \xrightarrow{\,f\,} \Id_\catC$$
between the constant functor $\mathrm{cst}_{c_0}:c \mapsto c_0$ and the identity functor on $\catC$, such that $f(c_0):c_0 \to c_0$ is the identity map of $c_0$.  If $\catC$ has a good class of weak equivalences (see [DHKS] 26.2), then we say $c_0$ is {\em homotopy initial} if there is a zig-zag of natural transformations of functors $\catC\to\catC$ 
$$\mathrm{cst}_{c_0} \cdots F_0 \xrightarrow{\,f\,} F_1 \cdots \Id_\catC$$
where the ``$\cdots$'' stand for zig-zags of natural weak equivalences and $f$, while not a natural weak equivalece, is at least a weak equivalence on $c_0$ -- i.e. $f(c_0):F_0(c_0) \xrightarrow{\,\weq\,} F_1(c_0)$.  [This last requirement ensures that evaluating the zig-zag at $c_0$ gives a zig-zag of weak equivalences from $c_0$ to itself.]

Write $\bigl(\lim^\catI \downarrow \mathcal{F}^\mathit{h}(\catC^\catI, \catM)\bigr)$ for the category of functors satisfying (1) and (2) -- i.e. homotopy functors under $\lim^\catI$.  Morphisms in this category are natural transformation triangles and weak equivalences are natural weak equivalences after forgetting the transformations from $\lim^\catI$.  The functor $\holim^\catI:\catC^\catI \to \catM$ is homotopically initial if there is a zig-zag of natural transformations of functors 
$\bigl(\lim^\catI \downarrow \mathcal{F}^\mathit{h}(\catC^\catI, \catM)\bigr) 
\longrightarrow
\bigl(\lim^\catI \downarrow \mathcal{F}^\mathit{h}(\catC^\catI, \catM)\bigr)$: 
$$\mathrm{cst}_{\holim^\catI} \cdots F_0 \xrightarrow{\,f\,} F_1 \cdots \Id_{\left(\lim^\catI \downarrow \mathcal{F}^\mathit{h}(\catM^\catI, \catM)\right)}$$ 
where $\mathrm{cst}_{\holim^\catI}:\bigl(\lim^\catI \downarrow \mathcal{F}^\mathit{h}(\catC^\catI, \catM)\bigr) \longrightarrow \bigl(\lim^\catI \downarrow \mathcal{F}^\mathit{h}(\catC^\catI, \catM)\bigr)$ is the constant functor on $\holim^\catI$ and $\Id_{\left(\lim^\catI \downarrow \mathcal{F}^\mathit{h}(\catM^\catI, \catM)\right)}$ is the identity functor on the category $\bigl(\lim^\catI \downarrow \mathcal{F}^\mathit{h}(\catM^\catI, \catM)\bigr)$, the ``$\cdots$''s are zig-zags of natural weak equivalences, and $f$ may not be a natural weak equivalence but at least $f(\holim^\catI)$ is a natural weak equivalence.

The notion of homotopically terminal is dual.

\begin{note}
Given a diagram $\mathscr{D}$ the homotopy limit in general does not have a natural map $\holim\mathscr{D} \to \mathscr{D}$; unlike the limit which does come with a natural map $\lim\mathscr{D}\to\mathscr{D}$.  Similarly, the homotopy colimit does not in general have a map $\mathscr{D}\to\hocolim\mathscr{D}$.  
\end{note}

\subsection{Homotopy Limit and Colimit Functors}

\paragraph{Existences of Homotopy Limits and Colimits}

Dwyer, Hirschhorn, Kan, and Smith prove:

\begin{thm}[\lbrack DHKS, 20.2\rbrack]\label{t:hocomplete} 
All homotopy $\catI$-limits and $\catI$-colimits exist for all model categories $\catM$ -- i.e. model categories are both homotopically complete and homotopically cocomplete.
\end{thm}

This is proven by identifying certain full subcategories of the category of all $\catI$-diagram functors $\catM^\catI$.  These subcategories are called $(\catM^\catI)_\mathsf{vc}$ the {\em virtually-cofibrant diagrams} and $(\catM^\catI)_\mathsf{vf}$ the {\em virtually-fibrant diagrams}, and satisfy the following properties [DHKS 20.5]:
\begin{enumerate}
\item[(i)] There exist ``virtually-cofibrant replacement'' and ``virtually-fibrant replacement'' functors  $(-)_\mathsf{vc}:\catM^\catI\to(\catM^\catI)_\mathsf{vc}$ and $(-)_\mathsf{vf}:\catM^\catI\to(\catM^\catI)_\mathsf{vf}$ equipped with maps 
$$\mathscr{D}_\mathsf{vc} \xrightarrow{\,\weq\,} \mathscr{D}\qquad \text{and}\qquad 
  \mathscr{D} \xrightarrow{\,\weq\,} \mathscr{D}_\mathsf{vf}$$
which are natural in $\mathscr{D}$.
\item[(ii)] If $F:\catM\rightleftarrows\catN:U$ is a Quillen adjoint pair then the induced pair $F^\catI:\catM^\catI \rightleftarrows \catN^\catI:U^\catI$ has the property that 
\begin{itemize}
\item $F^\catI$ preserves virtual-cofibrancy and weak equivalences of virtually-cofibrant diagrams.
\item $U^\catI$ preserves virtual-fibrancy and weak equivalences of virtually-fibrant diagrams.
\end{itemize}
\item[(iii)] The colimit and limit functors on $\catM$ preserve weak equivalences of virtually-cofibrant diagrams and virtually-fibrant diagrams respectively.
\end{enumerate}
Homotopy limit and colimit functors are then given by the compositions: 
$$\holim{\!}_\catM^\catI(-) := \lim{\!}_\catM^\catI(-)_\mathsf{vf}\qquad \text{and}\qquad
  \hocolim{\!}_\catM^\catI(-) := \colim{\!}_\catM^\catI(-)_\mathsf{vc}$$

Note that in general Quillen adjoint pairs do not preserve homotopy limits and colimits.  Recall that Quillen adjoints are not required to preserve all weak equivalences, thus the composition of a right adjoint and homotopy limit will likely no longer be a homotopy functor and so has little hope of being a homotopy limit functor.  However, if a Quillen adjoint pair happens to have the property that each functor preserves all weak equivalences\footnote{Note that this is not enough to imply that the Quillen pair is a Quillen equivalence; though all Quillen equivalences have this property.} (and not just weak equivalences which are either fibrations or cofibrations) then it follows immediately from [DHKS, 20.4] that it also preserves homotopy limits and colimits:

\begin{lemma}\label{adjoint holim}
If $F:\catM\rightleftarrows\catN:U$ is a Quillen adjoint pair between model categories and $U$  preserves weak equivalences, then there is a zig-zag of natural weak equivalences between the compositions $U\holim^\catI_\catN:\catN^\catI\to \catM$ and $\holim^\catI_\catM U:\catN^\catI\to\catM$ where $\holim^\catI_\catN$ and $\holim^\catI_\catM$ are any $\catI$-$\holim$ functors on $\catN$ and $\catM$ respectively.  Dually for homotopy colimits.
\end{lemma}

Given functors such as the above, which commute (up to zig-zags of natural weak equivalences) with $\holim$ or $\hocolim$ functors, we say that they {\em preserve homotopy limits} or {\em preserve homotopy colimits}.

Again, we standardly omit the $\catI$ in the $\holim$ and $\hocolim$ notation.

\paragraph{Canonical Homotopy Limits and Colimits}
In our later work we require actual explicit models for homotopy limits and colimits in certain categories.  In particular, we need to make statements about the structure of homotopy pullbacks and pushouts; for example, loops and suspensions.  In order to construct homotopy limit and colimit functors, we use generalizations of the homotopy limit and colimit functors originally defined by Bousfield and Kan for the categories $s\sets$ and $\top$ of simplicial sets and topological spaces (based or unbased).  The following constructions essentially come from Bousfield and Kan [BK72] and Hirschhorn [Hir]:\footnote{Compare to \ref{(co)lim (co)equalizer} and \ref{(co)lim adjoint}.}

\begin{lemma}[\lbrack Hir, 18.1.8, 18.1.2\rbrack]\label{D:holim end} 
Let $\catC$ be one of $s\sets_\ast$, $\top_\ast$, $s\sets$, or $\top$ and $\mathscr{D}:\catI \to \catC$ be a diagram in $\catC$.
\begin{itemize}
\item A homotopy limit of $\mathscr{D}$ is given by the equalizer
$$\xymatrix{
\holim(\mathscr{D}) 
 \ar[r] & \displaystyle
\phantom{\prod_{ i\in \catI} \mathscr{D}( i)^{\B(\catI\downarrow i)}}
\save[]-<0pt,4pt> 
 *{\displaystyle \prod_{ i\in \catI} \mathscr{D}( i)^{\B(\catI\downarrow i)} }
\restore \ar@<2pt>[r]^(.45)\phi \ar@<-2pt>[r]_(.45)\psi
  & \displaystyle
\phantom{\prod_{(\sigma: i\to j)\in \catI} \!\!\!\!\!\!
 \mathscr{D}(j)^{\B(\catI\downarrow i)} }
\save[]-<0pt,5pt> *{\displaystyle \prod_{(\sigma: i\to j)\in \catI} \!\!\!\!\!\!
 \mathscr{D}( j)^{\B(\catI\downarrow i)} } \restore
}$$
where $\phi$ is given by the projections
$$pr_\sigma\,\phi:\prod_{ i\in\catI}(\cdots) \longrightarrow 
  \mathscr{D}( i)^{\B(\catI\downarrow i)} \longrightarrow
  \mathscr{D}( j)^{\B(\catI\downarrow i)}$$
and $\psi$ is given by the projections
$$pr_\sigma\,\psi:\prod_{ i\in\catI}(\cdots) \longrightarrow
  \mathscr{D}( j)^{\B(\catI\downarrow j)} \longrightarrow
  \mathscr{D}( j)^{\B(\catI\downarrow i)}$$

\item A homotopy colimit of $\mathscr{D}$ is given by the coequalizer
$$\xymatrix{\displaystyle
\phantom{\coprod_{(\sigma: i \to  j)\in \catI} \!\!\!\!\!\!
 \mathscr{D}( i)\otimes \B(\catI\downarrow j) } 
\save[]-<0pt,5pt>
*{\displaystyle\coprod_{(\sigma: i \to  j)\in \catI} \!\!\!\!\!\!
 \mathscr{D}( i)\otimes \B( j\downarrow\catI) } \restore
 \ar@<2pt>[r]^(.55)\phi \ar@<-2pt>[r]_(.55)\psi  & \displaystyle
\phantom{\coprod_{ i\in \catI} \mathscr{D}( i)\otimes \B(\catI\downarrow i) }
\save[]-<0pt,4pt>
 *{\displaystyle\coprod_{ i\in \catI} \mathscr{D}( i)\otimes \B( i\downarrow\catI) }
 \restore \ar[r] & 
\hocolim(\mathscr{D})
}$$
where $\phi$ is given by the components 
$$\phi_\sigma:\mathscr{D}( i)\otimes \B(j\downarrow \catI) \longrightarrow 
 \mathscr{D}( j)\otimes \B(j\downarrow \catI) \longrightarrow \coprod_{ i \in \catI} (\cdots)$$
and $\psi$ is given by the components
$$\psi_\sigma:\mathscr{D}( i)\otimes \B(j\downarrow \catI) \longrightarrow 
 \mathscr{D}( i)\otimes \B(i\downarrow \catI) \longrightarrow \coprod_{ i \in \catI} (\cdots)$$
\end{itemize}
\end{lemma}

\begin{lemma}[\lbrack BK72, XI 3.3, XII 2.2\rbrack]\label{D:holim adj}
Let $\catC$ be one of $s\sets_\ast$, $\top_\ast$, $s\sets$, or $\top$.
\begin{itemize}
\item 
A homotopy limit functor $\hocolim:\catC^\catI \to \catC$ is given by a right adjoint of 
$$-\otimes \B(\catI\downarrow -):\catC \to \catC^\catI$$
the functor $X\mapsto X\otimes \B(\catI\downarrow -) \in \catC^\catI$.

\item
A homotopy colimit functor $\hocolim:\catC^\catI \to \catC$ is given by a left adjoint of 
$$\mathrm {hom}_\catC\bigl(\B(\catI\downarrow -), -\bigr):\catC \to \catC^\catI$$ 
the functor $X\mapsto \mathrm {hom}(\B(\catI\downarrow -), X) \in \catC^\catI$.
\end{itemize}
\end{lemma}

In the above constructions $(\catI\downarrow i)$ means the comma category of objects over $i$ ($k\in \catI$ with a map $k\to i$), and $(i\downarrow \catI)$ means the comma category of objects under $i$ ($k\in \catI$ with a map $i\to k$); the functor $\B(-)$ takes the nerve of a category if $\catC = s\sets$ and the realization of the nerve if $\catC = \top$; and the operation ``${\otimes}$'' is $\times$ of simplicial sets or spaces or else half smash (i.e. $-\wedge (-)_+$) if $\catC$ is based; and $X^Y = \mathrm {hom}(Y, X)$ is an object of $\catC$ (the possible categories for $\catC$ are all enriched over themselves).

Note that we may make sense of the these constructions in any category which is simplicially enriched.\footnote{See [Hir Ch. 9] for a discussion of simplicial model categories.  Basically a simplicial enrichment on $\catC$ consists of the extra structure of:
\begin{itemize}
\item Extensions of morphism classes to simplicial sets $\mathrm{Map}(X,Y)_\bullet\in\sset$ with $\mathrm{Map}(X,Y)_0 = \mathrm{Mor}_\catC(X,Y)$.
\item A simplicial tensor functor $X\otimes K\in\catC$ for $X\in \catC$, $K\in\sset$.
\item A simplicial mapping functor $X^K\in\catC$ for $X\in \catC$, $K\in\sset$.
\end{itemize}
satisfying all of the desired composition and compatibility conditions.
}
Furthermore, it is standard that given such a category, these constructions define homotopy limit and colimit functors.  If $\catM$ is a model category enriched over $\mathit{s}\sets$ then these are called the {\em canonical homotopy limit and colimit functors} in $\catM$.

\paragraph{Creation of Homotopy Limits and Colimits}
We will give two theorems -- \ref{dgl holim plan} and \ref{dgc holim plan} -- which are essentially both applications of a general theorem about creating homotopy limits and colimits in one model category from those of another model category by using Quillen adjoint pairs between the two model categories:

\begin{thm}[Creation of Homotopy Limits and Colimits]\label{general holim plan}
Let $F:\catM\rightleftarrows \catN:U$ be a Quillen adjoint pair such that $U$ detects weak equivalences of maps between fibrant objects\footnote{That is if $f:X_\mathsf{f}\to Y_\mathsf{f}$ is a map between fibrant objects in $\catN$ such that $U\!f$ is a weak equivalence in $\catM$, then $f$ was a weak equivalence in $\catN$.} and let $\holim_\catM$ be any $\catI$-homotopy limit functor on $\catM$.    
Suppose $L:\catN^\catI \to \catN_\mathsf{f}$ is a homotopy functor from $\catI$-diagrams in $\catN$ to $\catN_\mathsf{f}$ such that for all $\catI$-diagrams $\mathscr{D}:\catI\to\catN$ it satisfies
\begin{itemize}
\item $U\bigl(L(\mathscr{D})\bigr) = \holim_\catM U(\mathscr{D})$.
\item $L$ is equipped with natural maps $e:\lim_\catN\mathscr{D} \to L(\mathscr{D})$.
\item $Ue$ is the canonical map $\lim_\catM U(\mathscr{D}) \to \holim_\catM U(\mathscr{D})$.
\end{itemize}
Then $L$ is an $\catI$-homotopy limit functor on $\catN$.  

Dually for homotopy colimits.
\end{thm}

We save the proof of this theorem for the specific cases \ref{dgl holim plan} and \ref{dgc holim plan} to be mentioned later, where it is more illuminating.

\subsection{Commuting Homotopy Limits and Colimits}

Given a diagram $\mathscr{D}:\catI \!\times\! \catJ \to \catM$ we may view $\mathscr{D}$ as a diagram of diagrams in two ways -- either as $\mathscr{D}:\catJ \to \catM^\catI$ (a $\catJ$-diagram of $\catI$-diagrams) or else $\mathscr{D}:\catI \to \catM^\catJ$ (an $\catI$-diagram of $\catJ$-diagrams).  Given $j\in \catJ$ and $i\in \catI$ consider the corresponding homotopy limit and colimit of the diagrams $\mathscr{D}(j):\catI\to\catM$ and $\mathscr{D}(i):\catJ\to\catM$.  These are called the ``homotopy limit over $\catI\times j$'' and the ``homotopy colimit over $i\times\catJ$'' and written:
\begin{align*}
 \holim_{\catI\times j}\mathscr{D} &:= 
               \holim \bigl(\mathscr{D}(j):\catI\to\catM\bigr) \\
 \hocolim_{i\times\catJ}\mathscr{D} &:= 
               \hocolim \bigl(\mathscr{D}(i):\catJ\to\catM\bigr)
\end{align*}

From the functoriality of homotopy limits and colimits, it follows that the above objects themselves define diagrams $\catJ\to\catM$ and $\catI\to\catM$ respectively.  These diagrams are the ``homotopy limit over $\catI$'' and the ``homotopy colimit over $\catJ$,'' written:
\begin{align*}
\left(\holim_\catI \mathscr{D}\right)(j) &:= \holim_{\catI\times j} \mathscr{D} \\
\left(\hocolim_\catJ \mathscr{D}\right)(i) &:= \hocolim_{i\times \catJ} \mathscr{D} \\
\end{align*}

It follows immediately from the definition of homotopy limits and colimits (\ref{D: ho(co)lim}) that both $\displaystyle \bigl(\holim_\catI\, \holim_\catJ\bigr)$ and $\displaystyle \bigl(\holim_\catJ\, \holim_\catI\bigr)$ define $(\catI\!\times\!\catJ)$-homotopy limit functors on $\catM$.\footnote{A key step is to note that 
$\displaystyle \lim_\catI \lim_\catJ \mathscr{D} \cong \lim_\catJ \lim_\catI \mathscr{D} \cong \lim \mathscr{D}$.}  Therefore there is a zig-zag of natural weak equivalences between them.  Similarly there is a zig-zag of natural weak equivalences between the functors $\displaystyle \bigl(\hocolim_\catI\, \hocolim_\catJ\bigr)$ and $\displaystyle \bigl(\hocolim_\catJ\, \hocolim_\catI\bigr)$.

It is natural to ask under what conditions there are also zig-zags of natural weak equivalences between the functors $\displaystyle \bigl(\holim_\catI\, \hocolim_\catJ\bigr)$ and $\displaystyle \bigl(\hocolim_\catJ\, \holim_\catI\bigr)$.  If such zig-zags exist in the model category $\catM$, then we say that ``$\catI$-homotopy limits and $\catJ$-homotopy colimits {\em commute} in $\catM$.''\footnote{For a discussion and some examples of commuting limits and colimits see [MacXL \S IX.2].}


As a trivial example, empty homotopy limits commute with empty homotopy colimits (because model categories are all pointed) in any model category.  Also all homotopy limits or colimits of singleton diagrams (one object, one morphism) commute with all homotopy colimits or limits.  In the category $\sp$ of spectra all very small homotopy limits commute with very small homotopy colimits (see the next section for definitions).  In the category $\top_*$ of based topological spaces all very small homotopy limits commute with all filtered homotopy colimits.  In fact, in all of the categories which we consider, all very small homotopy limits commute with all filtered homotopy colimits.

\section{Special Diagrams}

We will be particularly interested in the (homotopy) limits and colimits of certain very special classes of diagrams.

\subsection{Pullback and Pushout Diagrams}

If $\mathscr{D}$ is a diagram in $\catM$ of the form 
$$\mathscr{D}:
 \left(\begin{aligned}\xymatrix@R=10pt@C=10pt{ 
        & \bullet \ar[d] \\ \bullet \ar[r] & \bullet
 }\end{aligned}\right) \xrightarrow{\ \ \ \ \ }
 \catM$$
we call $\mathscr{D}$ a {\em pullback diagram} in $\catM$.  Dually, if $\mathscr{D}'$ is a diagram in $\catM$ of the form
$$\mathscr{D}':\left(\begin{aligned}\xymatrix@R=10pt@C=10pt{
         \bullet \ar[r] \ar[d] & \bullet \\ \bullet & 
 }\end{aligned}\right) \xrightarrow{\ \ \ \ \ }
 \catM$$
then we call $\mathscr{D}'$ a {\em pushout diagram} in $\catM$.  More generally, we define $n$-dimensional pushout and pullback diagrams as follows:

Given $S$ a set let $\mathcal{P}(S)$ be the poset of subests of $S$ and inclusion maps viewed as a category.  Also, write $\mathcal{P}_0(S)$ for the full subcategory of all nonempty subsets of $S$, and $\mathcal{P}_1(S)$ for the full subcategory of all proper subsets of $S$.  Finally, given an integer $n\ge 0$ we write $\underline{n}$ for the set $\underline{n} = \{1,\dots,n\}$ with the understanding that $\underline{0} = \emptyset$.
\begin{defn}\label{n-dim pullback diagram}
An $n$-dimensional pullback diagram in $\catM$ is a diagram of the form 
$$\mathscr{D}:\mathcal{P}_0(\underline{n}) \longrightarrow \catM$$

An $n$-dimensional pushout diagram in $\catM$ is a diagram of the form
$$\mathscr{D}':\mathcal{P}_1(\underline{n}) \longrightarrow \catM$$
\end{defn}

\begin{ex}\label{pullback ex}
 We standardly denote generic cubes by $\mathscr{D}:S\mapsto X_S$.
\begin{itemize}  
\item Two-dimensional pullback and pushout diagrams have the forms:
$$
\begin{aligned}
\xymatrix@C=20pt@R=20pt{
 & X_{\{2\}} \ar[d] \\
 X_{\{1\}} \ar[r] & X_{\{1,2\}} 
}\end{aligned} \qquad \text{and} \qquad
\begin{aligned}
\xymatrix@C=20pt@R=20pt{
  X_{\emptyset} \ar[r] \ar[d] & X_{\{2\}} \\ 
  X_{\{1\}} &
}\end{aligned}$$
\item Three-dimensional pullback and pushout diagrams have the forms:
$$
\begin{aligned}
\xymatrix@C=7pt@R=7pt{
  & & X_{\{3\}} \ar'[d][dd] \ar[dr] &   \\ 
 & X_{\{2\}} \ar[rr] \ar[dd] & & X_{\{2,3\}} \ar[dd] \\
 X_{\{1\}} \ar'[r][rr] \ar[dr] & & X_{\{1,3\}} \ar[dr] &  \\
 & X_{\{1,2\}} \ar[rr] & & X_{\{1,2,3\}}  
}\end{aligned}\qquad \text{and}\qquad
\begin{aligned}
\xymatrix@C=7pt@R=7pt{
  X_{\emptyset} \ar[rr] \ar[dd] \ar[dr] & & X_{\{3\}} \ar'[d][dd] \ar[dr] &  \\
  & X_{\{2\}} \ar[rr] \ar[dd] & & X_{\{2,3\}}   \\
  X_{\{1\}} \ar'[r][rr] \ar[dr] & & X_{\{1,3\}}  &  \\
  & X_{\{1,2\}}  & & 
}\end{aligned}$$
\end{itemize}
\end{ex}

If $\mathscr{D}$ is an $n$-dimensional pullback diagram, then we often refer to the limit and homotopy limit of $\mathscr{D}$ as the {\em pullback} and {\em homotopy pullback} respectively.  Dually, if $\mathscr{D}'$ is an $n$-dimensional pushout diagram, then we refer to the colimit and homotopy colimit of $\mathscr{D}'$ as the {\em pushout} and {\em homotopy pushout} respectively.

\subsection{Very Small Diagrams}

A category $\catI$ is {\em very small} if its nerve is a simplicial set with only finitely many non-degenerate simplicies.   
Such categories are characterized by having only finitely many morphisms and no nontrivial loops.  That is, there cannot be $a\neq b\in \catI$ with morphisms 
$$\xymatrix{ a\ar@/^4pt/[r]^f & b\ar@/^4pt/[l]^g }\qquad \text{or} \qquad 
 \xymatrix{a \ar@(ur,dr)[]^{h\neq \Idm_a}}$$
We say that a diagram is {\em very small} if its indexing category $\catI$ is very small; $\catI$-limit, colimit, homotopy limit, and homotopy colimit functors are {\em very small} if they are functors of very small diagrams.

\begin{ex}Some standard examples of very small diagrams are:
\begin{itemize}
\item
Diagrams which are a disjoint union of finitely many points are very small.
\item
All $n$-dimensional pushout and pullback diagrams are very small.
\end{itemize}
\end{ex}
In particular, note that if $G$ is a (nontrivial) group and $\mathbf{G}$ is the category consisting of one object with morphisms labelled by elements of $G$, then $\mathbf{G}$ is \underline{not} a very small category -- even if the group $G$ is finite.

\subsection{Filtered Diagrams}

A nonempty category $\catI$ is called {\em filtered}\footnote{Or maybe ``right filtered'' or ``cofiltered''.} (see [MacL \S IX.1]) if for every pair of objects $a,b\in\catI$ there is an object $c\in\catI$ above them both:
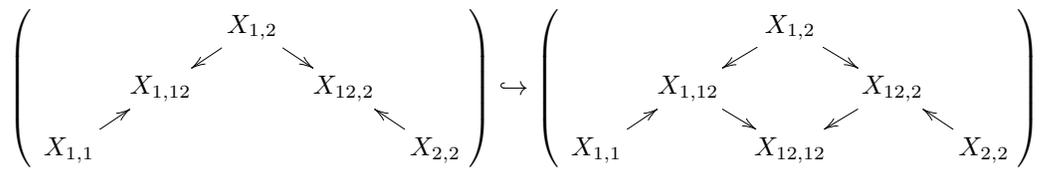
$$\xymatrix@R=0pt@C=20pt{a \ar[dr] & \\ & c \\ b\ar[ur] & }$$
and for every pair of parallel arrows $\xymatrix{a \ar@/^4pt/[r]^f \ar@/_4pt/[r]_g & b}$ there is an object $c\in\catI$ above them both:
$$\xymatrix@R=0pt@C=20pt{& b \ar[dr] & \\ a \ar[ur]^f \ar[dr]_g & & c \\ & b \ar[ur] & }$$
Diagrams $\catI\to\catM$ where $\catI$ is a filtered category are called {\em filtered diagrams}.  Similarly $\catI$-colimit and $\catI$-homotopy colimit functors are called {\em filtered} if $\catI$ is a filtered category.

\begin{ex}
Standard examples of filtered categories are:
\begin{itemize}
\item
If $\catI$ is an ordered set viewed as a category, then $\catI$ is filtered.

In particular, sequential homotopy colimits are filtered homotopy colimits.
\item
If $\catI$ has a final object, then $\catI$ is filtered.
\end{itemize}
\end{ex}

\subsection{Cofinal Diagrams}\label{S:cofinal}
Suppose $\catI$ is a category and $\catJ$ is some full subcategory of $\catI$. 
The subcategory $\catJ$ is {\em left cofinal} in $\catI$ (or more properly, the inclusion $\catJ\cofibr\catI$ is {\em left cofinal}) if the comma category $(\catJ\downarrow i)$ (is nonempty and) has contractible nerve for all $i\in\catI$.

Let $\mathscr{D}:\catI\to \catM$ be a diagram in the model category $\catM$. 
Bousfield and Kan show (c.f. [BK69 XI \S9]) that if $\catJ$ is left cofinal in $\catI$ and $\catM$ is one of $\mathit{s}\sets_*$, $\top_*$, $\mathit{s}\sets$, and $\top$ 
there is a weak equivalence between the canonical homotopy limit of $\mathscr{D}$ and the canonical homotopy limit of the restricted diagram $\mathscr{D}|_\catJ$:
$$\holim_\catI (\mathscr{D}) :=  \holim(\mathscr{D})
  \xrightarrow{\ \weq\ } \holim \bigl(\mathscr{D}|_\catJ\bigr)  =:  \holim_\catJ (\mathscr{D})$$ 
This weak equivalence is given by simplicial homotopies contracting the nerves of the comma categories $(\catJ \downarrow i)$ in the construction of the canonical homotopy limit (see \ref{D:holim end}).

More generally, if $\catM$ is any simplicially enriched model category so that the Bousfield-Kan construction yields canonical homotopy limits, then the same result holds.

\begin{ex}The inclusions of subcategories in the below diagrams are cofinal:
$$\left(\begin{aligned}\xymatrix@C=7.5pt@R=7.5pt{
 & X_{1} \ar[dr] & & X_{2} \ar[dl]  &  \\
 & & X_{1,2} &&
}\end{aligned}\right) \cofibr
 \left(\begin{aligned}\xymatrix@C=7.5pt@R=7.5pt{
 & X_{1} \ar[dl] \ar[dr] & & X_{2} \ar[dl] \ar[dr] &  \\
X_{1,1} & & X_{1,2} && X_{2,2}
}\end{aligned}\right)
$$
\item
$$\left(\begin{aligned}\xymatrix@C=7.5pt@R=7.5pt{
 & & X_{1,2} \ar[dl] \ar[dr] & & \\
 & X_{1,12} & & X_{12,2} & \\
 X_{1,1} \ar[ur] & &  & & X_{2,2} \ar[ul]
}\end{aligned}\right)\cofibr
 \left(\begin{aligned}\xymatrix@C=7.5pt@R=7.5pt{
 & & X_{1,2} \ar[dl] \ar[dr] & & \\
 & X_{1,12} \ar[dr] & & X_{12,2} \ar[dl] & \\
 X_{1,1} \ar[ur] & & X_{12,12} & & X_{2,2} \ar[ul]
}\end{aligned}\right)
$$
\end{ex}

\part{Rational Homotopy Theory}

\markboth{Ben Walter}{I.  Rational Homotopy Theory}

All of our constructions are over the ground field $\Q$.  As such we generally suppress the symbol $\Q$ in notation (e.g. we write $\otimes$ to mean $\otimes_\Q$, etc).  We will be concerned with a number of different categories whose objects are vector spaces over $\Q$ supporting some extra structures such as a grading, a differential, or a coproduct.  These categories will be denoted by identifying the extra structure that their objects carry; for example, $\g$ for graded vector spaces (over $\Q$) and $\dgl$ for differential graded Lie algebras (over $\Q$).  Script type will be used to refer to categories and small caps will be used for abbreviations of the category names: i.e. $C \in Obj(\dgc)$ means $C$ is a \DGC.  In general, the constructions of this chapter may be discussed over any field {\slshape{l\kern -1.75pt{k}}} of characteristic 0, but we will specialize sooner rather than later.  

The main players in our work are differential graded vector spaces (\DG\,s), differential graded Lie algebras (\DGL\,s), and differential graded (cocommutative, counital, coaugmented) coalgebras (\DGC\,s).  We also refer at times to graded vector spaces (\G\,s), graded Lie algebras (\GL\,s), and graded coalgebras (\GC\,s).  We will be concerned with many different functors between these categories.  Given categories $\catB$ and $\catC$ our convention is to abusively use the notation $[-]_\catC:\catB \to \catC$ to denote either a forgetful functor $\catB \to \catC$ (e.g. $[-]_\mDG:\dgl \to \dg$ by forgetting the Lie bracket) or a trivial section of a forgetful functor $\catC \to \catB$ (e.g. $[-]_\mDGL:\dg \to \dgl$ by equipping a \DG\ with the trivial Lie bracket $[-,-]=0$).  Other functors are given more creative and unique notation.  Finally, our convention when referring to adjoint pairs of functors is to write $F:\catB \rightleftarrows \catC : U$ to mean that $F:\catB \to \catC$ is {\em left adjoint} to $U: \catC \to \catB$ (that is, $hom(Fb,c) \cong hom(b,Uc)$ for $b\in\catB$ and $c\in \catC$).

In Chapters~\ref{S:DG}, \ref{S:DGL}, and \ref{S:DGC} we remind the reader of the definitions and develop some of the properties and structures of the categories $\dg$, $\dgl$, and $\dgc$.  The work that we present is a mixture of previously known results, new results, and previously known results placed in a new framework.  We have endeavored to indicate all statements already proven in the literature by supplying references appropriately -- furthermore we standardly omit the proofs of these results.

The first section of each of Chapters~\ref{S:DG}-\ref{S:DGC} is devoted to basic definitions -- in particular, limits, colimits, as well as cylinder and path objects.  In Chapter~\ref{S:DGL} (and Chapter~\ref{S:DGC}) we also give the definitions of free maps and free objects (and cofree maps and cofree objects) which are the cofibrations and cofibrant objects (and fibrations and fibrant objects) under our model category structure.  We also note the very tight connection between $\dg$ and $\dgl$ (and $\dgc$).
Much of the material in the first section of Chapter~\ref{S:DG} is standard homotopical algebra (from, for example, [Weib]).   Similarly almost all of the first section of Chapters~\ref{S:DGL} and \ref{S:DGC} consists of either standard or simple extensions of standard facts (from, for example, [FHT] and [Q69]).  The first section of each of these chapters ends with a definition of loop and suspension functors $\Omega$ and $\Sigma$.  These functors are Quillen adjoint pairs compatible with the standard maps between $\dg$, $\dgl$, and $\dgc$.

The second section of these chapters is largely devoted to homotopy limits and colimits in the categories $\dg$, $\dgl$, and $\dgc$ (after beginning with a small note about the standard model category structure of the categories in question).  The primary goal of these sections is to note two things -- (1) we may construct certain very simple models for homotopy pushouts and pullbacks so that the functors $\Omega$ and $\Sigma$ described at the end of the previous section are indeed given by homotopy pullbacks and homotopy pushouts; and (2) homotopy pullbacks commute with sequential homotopy limits.  

We do not have a specific reference to the literature duplicating our precise construction of homotopy limits and colimits in $\dg$, but we suspect that it is standard.  The constructions which we make in $\dgl$ and $\dgc$, however, are new -- in particular even knowing that these constructions deserve to be called homotopy limits and colimits relies on relatively recent work of Dwyer, Hirschhorn, Kan, and Smith [DHKS].  Our homotopy limits and colimits in $\dgl$ and $\dgc$ are essentially made by ``lifting'' the homotopy limits and colimits which we give in $\dg$.  Note that our models for homotopy limits and colimits in $\dg$, $\dgc$, and $\dgl$ are \underline{not} the canonical homotopy limits and colimits arising from the simplicial enrichment of $\dg$, $\dgl$, and $\dgc$ -- they are in general much smaller.

  In Chapter~\ref{S:Q homotopy} we begin by recalling the standard framework of rational homotopy theory from Quillen [Q69].  We then go on to discuss rational spectra as well as rational $\Sinf$ and $\Linf$ functors.  
There are standard methods to stabilize a model category to construct an associated category of spectra by inverting either a generic ``suspension endofunctor'' [Ho01] or the simplicial suspension functor [S97].  We do not explicitly use any of these methods; though we do show that our category of ``rational spectra'' is stable, is equivalent to the rational localization of the cateogory of spectra, and that our $\Sinf$ and $\Linf$ functors satisfy the desired properties.

Note that we are interested in the rational homotopy of highly connected spaces, where by ``highly connected'' we mean $n$-connected for $n\ge 1$, possibly $n\gg1$.

\chapter{Differential Graded Vector Spaces}\label{S:DG}
\markboth{Ben Walter}{I.3  Differential Graded Vector Spaces}

Differential graded vector spaces will serve as our rational spectra.  Reduced differential graded vector spaces are connective rational spectra.  The primary goal of this chapter is the construction in \ref{l:dg pushout} of specific models for homotopy pullbacks and pushouts which are simpler than those given by the canonical homotopy limits and colimits (see \ref{D:holim end} and \ref{D:holim adj}) built using the simplicial enrichment of $\dg$ and $\dg_r$.  The simple models which we construct for homotopy pushouts and pullbacks in $\dg$ and $\dg_r$ will later be used to make simple models for homotopy pushouts and pullbacks in $\dgl$ and $\dgc$ as well as in our eventual construction of rational homotopy calculus of functors.

\section{Category Structure}

Most of the proofs in this section are omitted since they are either well-known or trivial.  A good reference for much of this material is [Weib \S1].

We write $\g$ to denote the category of graded $\Q$-vector spaces (and degree 0 maps) and $\dg$ to denote the category of differential graded $\Q$-vector spaces (i.e. rational chain complexes).  Objects of $\g$ are written $V_\bullet = \{V_i\}_{i\in \mathbb{Z}}$.  Given an element $v \in V_n$ we say that the degree of $v$ is $n$ and write $|v|= n$.  
Note that we do not require objects of $\g$ (or $\dg$) to be bounded below.  More explicitly:

\begin{defn}[DG]
A differential graded vector space (\DG) is $V = (V_\bullet,\, d_V)$ where $V_\bullet$ is a graded vector space and $d_V$ is a degree $-1$ endomorphism of $V_\bullet$ \footnote{Our notation is $d_V=\{d_i\}$ where $d_i:V_i \to V_{i-1}$.} (i.e. $|d_Vv| = |v| - 1$) with $d_V\circ d_V = 0$.
\end{defn}  

In general we use $(-)_\bullet$ for the forgetful functor $\dg\to\g$; given $V$ a \DG\ when we write $V_\bullet$ we mean the underlying \G\ of $V$.  The category $\dg$ has objects all rational chain complexes and arrows all degree 0 chain maps (i.e. degree 0 maps respecting the differential: $fd = df$).  Injections, surjections, and isomorphisms of \DG\,s are \DG-maps which are degree-wise injections, surjections, or isomorphisms of vector spaces.  {\em Quasi-isomorphisms} of \DG\,s are \DG-maps which induce isomorphisms on homology.  We say that a \DG\ is {\em contractible} if it has trivial homology $H_*(V) = 0$.

A \DG\ is called $r$-reduced if it is trivial below grading $r$.  We write $\dg_r$ for the full subcategory of $\dg$ consisting of all $r$-reduced \DG\,s.
There are two particularly useful functors $\dg \to \dg_r$ called the {\em $r$-truncation} and {\em $r$-reduction} functors (we may leave out the $r$ when it is clear from context).  
\begin{defn}[Truncation and Reduction]\footnote{
These are called the ``brutal'' and ``good'' truncations by [Weib].}
Given $V = (V_\bullet,\, d)$ a \DG, its $r$-truncation $\mathrm{trunc}_r(V) = \trunc{V}\in\dg_r$ is the quotient-\DG\ $\trunc{V}$ given by $\trunc{V}_k = V_k$ if $k \ge r$ and $\trunc{V}_k = 0$ if $k < r$.

Its $r$-reduction $\mathrm{red}_r(V) = \red{V}\in\dg_r$ is the sub-\DG\  $\red{V}$ given by where $\red{V}_k = V_k$ for $k > r$, $\red{V}_k = 0$ for $k < r$, and $\red{V}_r = {\rm ker}(d_r:V_r \to V_{r-1})$.\footnote{Reduction of a \DG\ corresponds to taking the universal $(r-1)$-connected cover of a topological space.}
\end{defn}

We use the same symbols $\mathrm{red}_r$ and $\mathrm{trunc}_r$ to denote the reduction and truncation functors $\dg \to \dg_r$ as well as the reduction and truncation functors $\dg_t\to\dg_r$ (for $t<r$).  Also, we use the same symbol to denote the inclusion of categories functor $\dg_r\to\dg$ and $\dg_r \to \dg_t$ (for $t<r$).

\begin{lemma}[Adjointness of Truncation and Reduction]\label{dg red}
The following are adjoint pairs
\begin{alignat*}{7}
\mathrm{trunc}_r:&\, \dg &\,\rightleftarrows\, &\dg_r&\,:&\,\mathrm{incl} \\
\mathrm{incl}:&\,\dg_r &\,\rightleftarrows\, &\dg&\,:&\,\mathrm{red}_r 
\end{alignat*} 
between $\dg$ and $\dg_r$; and 
\begin{alignat*}{7}
\mathrm{trunc}_r:&\, \dg_t &\,\rightleftarrows\, &\dg_r&\,:&\,\mathrm{incl} \\
\mathrm{incl}:&\,\dg_r &\,\rightleftarrows\, &\dg_t&\,:&\,\mathrm{red}_r \\
\end{alignat*}
between $\dg_r$ and $\dg_t$ (for $t<r$).

Furthermore, $\mathrm{incl}$ the inclusion of categories functor is a section of both the truncation and reduction functors.
\end{lemma}

\begin{note}\label{H_* red} 
The reduction functor has the additional useful property that it is compatible with homology.  By construction, it is clear that $H_i(\red{V}) = H_i(V)$ for $i\ge r$ and $H_i(\red{V}) = 0$ for $i < r$. In particular, $\mathrm{red}_r$ preserves quasi-isomorphisms.
For this reason, we are primarily interested only in reductions and truncations will play very little part in the remainder.
\end{note}

\subsection{Limits and Colimits}

Recall that the category $\dg$ (and $\dg_r$) supports  natural symmetric monoidal operations $\oplus$, $\times$, and $\otimes$ defined by

\begin{itemize}
\item $V \oplus W := (V_\bullet,\, d_V) \oplus (W_\bullet,\, d_W) = ((V\oplus W)_\bullet,\, d_{\oplus})$, \newline where $(V\oplus W)_n = V_n \oplus W_n$ with differential $d_{\oplus}(v \oplus w) = d_V(v) + d_W(w)$
\item $V \times W := (V_\bullet,\, d_V) \times (W_\bullet,\, d_W) = ((V\times W)_\bullet,\, d_{\times})$, \newline where $(V\times W)_n = V_n \times W_n$ with differential $d_{\times}(v ,\, w) = (d_V(v),\,d_W(w))$
\item $V \otimes W := (V_\bullet,\, d_V) \otimes (W_\bullet,\, d_W) = ((V\otimes W)_\bullet,\, d_{\otimes})$, \newline where $(V\otimes W)_n = \oplus_{i+j=n} V_i \otimes W_j$ with differential $d_{\otimes}(v\otimes w) = d_V(v)\otimes w + (-1)^{|v|}v \otimes d_W(w)$
\end{itemize}

The operations $\oplus$ and $\times$ give isomorphic answers on finite collections of \DG\,s but differ on infinite collections\footnote{Our practice we is to use $\times$ and $\oplus$ to indicate whether we are thinking of the categorical product or coproduct in a given situation -- even when the two give isomorphic answers.} -- in general $\oplus$ and $\times$ extend to define operations on all small collections of objects.  Given $\{V_i\}_{i\in \catI}$ a collection of \DG\,s, $\oplus_{i\in\catI} V_i$ is their categorical coproduct and their categorical product is given by $\times_{i\in\catI} V_i$ (in $\dg_r$ as well as $\dg$).
Note that the trivial \DG\ $0_\mDG := (0_\bullet,\, d=0)$ (where $0_i = \Q^0=0$) is both initial and final in the category $\dg$ (and $\dg_r$), thus $\dg$ (and $\dg_r$) is a pointed category.  Also, the \DG\ $0_\mDG$ is the unit for $\oplus$, and the \DG\ $1_\mDG := (1_\bullet,\, d=0)$ (where $1_0 = \Q$, $1_i = 0$ all other $i$) is the unit for $\otimes$.     
Note that $V$ a \DG\ is contractible if and only if the map to the final object $V \to 0_\mDG$ is a quasi-isomorphism. 

The \DG\,s $0_\mDG$ and $1_\mDG$ are part of a more general framework.  There is an adjoint pair of functors between the category of sets (and all set maps) and $\dg$ given by $\Q- : \sets \rightleftarrows \dg : -_0$ where given $S$ a set, $\Q S$ is the free \DG\ concentrated on degree 0 generated by $S$ (i.e. $\Q S = (V_\bullet, d=0)$ with $V_0 = \oplus_{s\in S}\, \Q s$, $V_i = 0$ for $i\neq 0$); and $-_0:V \mapsto V_0$, its degree 0 vector space (viewed as a set).  Even more generally there is an adjoint pair $-_\mDG:\sset \rightleftarrows \dg_0:N^*$ between the categories of simplicial sets and $\dg$, given by the composition of adjoint pairs
\begin{equation}\label{adj:sset dg}
\xymatrix{
-_\mDG:\sset \ar@<2pt>[r]^(.6){\Q-} & \svect \ar@<2pt>[l]^(.4)U \ar@<2pt>[r]^(.52)n & \dg_0 \ar@<2pt>[l]^(.45){n^*} \ar@<2pt>[r]^(.43){\mathrm{incl}} & \dg \ar@<2pt>[l]^(.56){\mathrm{red}_0} : N^*
}
\end{equation}
where $\Q-$ makes the free simplicial vector space on a simplicial set, $U$ forgets vector space structure, $n$ is the normalization functor, and $n^*$ is the adjoint of normalization.  

Two other \DG's of note are $s := (s_\bullet,\, d)$ where $s_\bullet$ has a $\Q$ in degree 1 and zero elsewhere (abusing notation\footnote{So $s\in s$ but hopefully context will keep this notational circularity from becoming confusing.} we write $s_\bullet = \Q s$ where $s$ is in degree 1), and $s^{-1} := (s^{-1}_\bullet,\, d)$ where $s^{-1}_\bullet$ has a $\Q$ in degree $-1$ and zero elsewhere (again, abusing notation we write $s^{-1}_\bullet = \Q s^{-1}$ where $s^{-1}$ is in degree $-1$).  Given a \DG, V, we write $sV$ and $s^{-1}V$ for $s\otimes V$ and $s^{-1} \otimes V$; similarly given $v\in V$ we write $sv$ and $s^{-1}v$ for $s\otimes v\in sV$ and $s^{-1}\otimes v \in s^{-1}V$.  These have the properties that $(sV)_i \cong V_{i-1}$ and $(s^{-1}V)_i \cong V_{i+1}$; furthermore $d_{sV}(sv) = -sd_V(v)$ and $d_{s^{-1}V}(s^{-1}v) = - s^{-1}d_V(v)$.  At times, we may write $s^nV$.  By this we mean $\underbrace{s\otimes \cdots \otimes s}_{n} \otimes V$; similarly, $s^{-n}V = \underbrace{s^{-1}\otimes \cdots \otimes s^{-1}}_{n} \otimes V$.

\begin{lemma}
The pushout in $\dg$ (and $\dg_r$) of the \DG-diagram $U\xleftarrow{f} V \xrightarrow{g} W$ is the quotient-\DG: $$U\oplus_V W := (U\oplus W)/\langle f(v)+g(v)\rangle$$  

The pullback in $\dg$ (and $\dg_r$) of the \DG-diagram $U\xrightarrow{f} V \xleftarrow{g} W$ is the sub-\DG:
$$U\times_V W := \{(u,w) \in U\times W\ \ | \ f(u)+g(w)=0\}$$

General small limits and colimits in $\dg$ (and $\dg_r$) are defined analogously (using the categorical product and coproduct along with equalizers and coequalizers as in Theorem~\ref{(co)lim (co)equalizer}).
\end{lemma}


Recall that fibers and cofibers of maps in pointed categories may be defined as pushouts and pullbacks:
\begin{itemize} 
\item The fiber of $f:V \to W$ is the pullback of $V \xrightarrow{f} W \xleftarrow{} 0_\mDG$.
\item The cofiber of $f:V \to W$ is the pushout of $0_\mDG\xleftarrow{}V \xrightarrow{f} W$.
\end{itemize}

\begin{lemma}
Let $f$ be a map $f:V \to W$.  The fiber of $f$ is $\ker(f)$; the cofiber of $f$ is $\coker(f)$.
\end{lemma}

\subsection{Cones, Suspensions, Paths, and Loops}

\begin{defn}[Cone]
Given $V$ a \DG, the cone on $V$ is the \DG\ 
$\mathrm{c}V := c\otimes V$ where $c$ is the \DG\ given by $c=(\Q v_0\oplus \Q v_1,\, d_cv_1 = v_0)$, $|v_i| = i$.
\end{defn}

More explicitly, $\mathrm{c}V = (V\oplus sV,\, d_cv=d_Vv, d_c(sv)=-sd_Vv + v)$.  Just as for topological spaces, the cone on $V$ is a contractible \DG\ (i.e. quasi-isomorphic to 0) with an injection $V \to \mathrm{c}V$.

For more on mapping cones (and cylinders) see [Weib \S1.2].

\begin{defn}[Suspension]
Given $V$ a \DG, the {\em suspension of $V$} is the \DG\
$\Sigma V := sV$ defined earlier.\end{defn}

At times we wish to use a slightly larger model for the suspension:  $\hat\Sigma V := S\otimes V$ where $S$ is the \DG\ given by $S=(\Q v_0\oplus \Q w_1 \oplus \Q v_1,\, d_Sw_1 = v_0, d_Sv_1 = v_0)$, $|v_i|=i$, $|w_1|=1$.

More explicitly, $\hat\Sigma V$ is given by:
\begin{itemize}
\item $(\hat\Sigma V)_\bullet = (sV_\bullet\oplus V_\bullet\oplus sV_\bullet)$ 
\item $d_{\hat\Sigma}(sv_1 \oplus v_2 \oplus sv_3) = -sd_Vv_1 \oplus (d_Vv_2 + v_1 + v_3) \oplus-sd_Vv_3$ for $sv_1 \oplus v_2 \oplus sv_3 \in (sV\oplus V\oplus sV)_\bullet$
\end{itemize}

\begin{note} 
The \DG\,s $\Sigma V$ and $\hat\Sigma V$ are both pushouts:
\begin{itemize}
\item $\Sigma V$ is the pushout of $0_\mDG \xleftarrow{} V \xrightarrow{} \mathrm{c}V$.
\item $\hat \Sigma V$ is the pushout of $\mathrm{c}V \xleftarrow{} V \xrightarrow{} \mathrm{c}V$.
\end{itemize}
\end{note}
The \DG\ ${\hat\Sigma} V$ admits a map $V \to {\hat\Sigma} V$.  This corresponds to the map of topological spaces sending $X$ to the equator $X \to {\Sigma} X$.  Also, the two projection maps (onto the first and second copies of $sV$) induce maps ${\hat\Sigma}V \to \Sigma V$ which is are both quasi-isomorphism (however, the compositions $V\to \hat\Sigma V \to \Sigma V$ are each zero).
  
We define loops similarly:

\begin{defn}[Paths]
Given $V$ a \DG, the paths on $V$ is the \DG\
$\mathrm{p}V := p\otimes V$ where $p$ is the \DG\ given by $p=(\Q v_{-1}\oplus \Q v_0,\, d_p v_0 = v_{-1})$, $|v_i|=i$.
\end{defn}

More explicitly, $\mathrm{p}V = (V\oplus s^{-1}V,\, d_pv=d_Vv + s^{-1}v, d_p(s^{-1}v)=-s^{-1}d_Vv)$.  Just as for topological spaces, the paths on $V$ is a contractible \DG\ (i.e. quasi-isomorphic to 0) with a surjection $\mathrm{p}V \to V$.

\begin{defn}[Loops]
Given $V$ a \DG, the {\em loops on $V$} is the \DG\ $\Omega V := s^{-1}V$ defined earlier.\end{defn}

At times we also wish to use a slightly larger model for loops:  ${\hat\Omega} V := P\otimes V$ where $P$ is the \DG\ given by $P=(\Q v_0\times \Q w_{-1} \times \Q v_{-1},\, d_Pv_0 = v_{-1}+w_{-1})$, $|v_i|=i$, $|w_{-1}|=-1$.

More explicitly, ${\hat\Omega} V$ is given by:
\begin{itemize}
\item $({\hat\Omega} V)_\bullet = (s^{-1}V \times V \times s^{-1}V)_\bullet$
\item $d_{\hat\Omega}(s^{-1}v_1,\,v_2,\,s^{-1}v_3) = \big(s^{-1}( - d_Vv_1 + v_2),\,d_Vv_2,\,s^{-1}( -d_Vv_3+v_2)\big)$
\end{itemize}

\begin{note} 
The \DG\,s $\Omega V$ and $\hat\Omega V$ are both pullbacks:
\begin{itemize}
\item $\Omega V$ is the pullback of $0_\mDG \xrightarrow{} V \xleftarrow{} \mathrm{p}V$.
\item $\hat\Omega V$ is the pullback of $\mathrm{p}V \xrightarrow{} V\xleftarrow{} \mathrm{p}V$.
\end{itemize}
\end{note} 

The \DG\ ${\hat\Omega} V$ admits a map ${\hat\Omega} V \to V$.  This corresponds to the map of topological spaces evaluating at the midpoint of loops ${\Omega} X \to X$.  Also the two injection maps (into the first and second copy of $s^{-1}V$) induce maps $\Omega V \to {\hat\Omega} V$ which are both quasi-isomorphisms (however, each composition $\Omega V \to \hat\Omega V \to V$ is zero). 

Since $\mathrm{c}V$, $\mathrm{p}V$, $\Sigma V$, $\Omega V$, $\hat \Sigma V$, and $\hat \Omega V$ are each given by tensoring with a \DG, they give functors $c, p, \Sigma, \Omega, \hat \Sigma, \hat \Omega:\dg \to \dg$.  Also, there is also a clear 1-1 correspondence between maps $\Sigma V \to W$ and maps $V\to \Omega W$.  In fact,

\begin{lemma}\label{dg S L adjoint}
The following functors are adjoint pairs:
\begin{itemize}
\item $c:\dg \rightleftarrows \dg : p$
\item $\Sigma:\dg \rightleftarrows \dg :\Omega$
\item $\hat\Sigma:\dg \rightleftarrows \dg: \hat\Omega$
\end{itemize}
\end{lemma}
\begin{proof}[Sketch of proof]
This follows from the fact that for a fixed, degree-wise finite dimensional $V \in \dg$, the functor $-\otimes V: W \mapsto W\otimes V$ is adjoint to the functor $-\otimes V^*$ (where $V^* = \mathcal{M}\!\mathit{ap}_\dg(V, 1_\mDG)$ is the dual of $V$).  It is not hard to show that $c^* = p$, $s^* = s^{-1}$, and $S^* = P$. [Note that $\dg$ becomes enriched over itself if we allow chain maps of arbitrary degree.]
\end{proof}
In fact, the above pairs are also adjoint in the opposite direction as well (e.g. $p:\dg \rightleftarrows \dg:c$).

Note that $\Sigma$ restricts to a functor $\dg_r \to \dg_r$ for any $r$.  However, $\Omega$ at best restricts to a functor $\dg_{r+1} \to \dg_r$.  We define the {\em reduced loop} functor to be the ``correct'' loop functor on $\dg_r$:

\begin{defn}[Reduced Loops]
If $V\in \dg_r$ is an $r$-reduced \DG\ then the reduced loops on $V$ is $\red{\Omega}V := \red{\Omega V} = \mathrm{red}_r(\Omega V)$.
\end{defn}

The following lemma is immediately implied by \ref{dg S L adjoint} and \ref{dg red}:

\begin{lemma}
Reduced loops is right adjoint to suspension on $\dg_r$:
$$\Sigma:\dg_r \rightleftarrows \dg_r: \red{\Omega}$$
\end{lemma}

The reduced loops $\red\Omega V$ is the pullback of the diagram in $\dg_r$ given by $0_\mDG \xrightarrow{} V \xleftarrow{} \red{\mathrm{p}V}$.  Note that the reduced paths $\red{\mathrm{p}V}$ is still contractible (by \ref{H_* red}).  Also, although the map $V \xleftarrow{} \red{\mathrm{p}V}$ is no longer necessarily a surjection of \DG\,s (since it may fail to be surjective in degree $r$), we will place a model category structure on $\dg_r$ such that it will be a fibration.  If we wished, we could similarly define $\red{\hat\Omega}V$ as the pullback of the diagram in $\dg_r$ given by $\red{\mathrm{p}V}\xrightarrow{} V \xleftarrow{} \red{\mathrm{p}V}$.

\section{Model Category Structure of $\dg$}

Recall that a {\em quasi-isomorphism} of \DG\,s is a (degree 0) chain map $f:V\to W$ which induces an isomorphism on homology $f_\ast: H_\ast(V) \xrightarrow{\,\cong\,} H_\ast(W)$. 
The standard model category structure on $\dg$ (see, e.g. [Hov] \S2.3) is to take quasi-isomorphisms as weak equivalences, degree-wise surjections as fibrations, and degree-wise injections as cofibrations.  Under this structure, all objects are both cofibrant and fibrant. 

Given $V$ a \DG\ we may view its homology $H_\ast(V)$ as a \DG\ with trivial differential.  By choosing basis elements and generators it is possible to construct non-canonical \DG-maps $V \to H_\ast(V)$ and $H_\ast(V) \to V$ which are quasi-isomorphisms since $H_\ast\bigl(H_\ast (V)\bigr) = H_\ast(V)$.  In particular, every \DG\ is quasi-isomorphic to its own homology.  This is a very special occurence:  In general, the special class of objects in the categories $\dgl_r$ and $\dgc_r$ which are quasi-isomorphic to their own homology are called {\em formal}, and are the subject of much interest.

\subsection{Homotopy Limits and Colimits}\label{dg holim}

By [DHKS] all homotopy limits and colimits exist (see \ref{t:hocomplete}).  We give explicit nice models for most homotopy limits and colimits similar to the standard constructions of homotopy limits and colimits in $\top$ as ends, coends, and right and left adjoints as stated in Lemmas~\ref{D:holim end} and \ref{D:holim adj} in the previous chapter.

The primary result which we make use of in later chapters is the existence of a convenient, model for homotopy pushouts and pullbacks which is easier to use in computations than the canonical Bousfield-Kan homotopy pushouts and pullbacks:

\begin{thm}\label{l:dg pushout} Homotopy pushouts and pullbacks
 in $\dg$ may be given as follows:
\begin{enumerate}
\item\label{l:dg pushout 1} Given $\mathscr{D}$ the diagram $U\xrightarrow{f} V \xleftarrow{g} W$, 
a homotopy pullback of $\mathscr{D}$ is given by the path \DG:
$$\mathscr{P}_\mathscr{D} := \bigl((U\times s^{-1}V \times W)_\bullet,\
d_{{\mathscr{P}}} = d_{\times} + d_{fg}\bigr)$$
where $d_{\times}$ is the differential on the \DG\ $(U\times s^{-1}V \times W)$ and $d_{fg}$ is determined by $f$ and $g$:
$d_{fg}(u,\, v ,\, w) = \bigl(0,\, s^{-1}(f(u) + g(w)),0\bigr)$.
\item\label{l:dg pushout 2} Given $\mathscr{D}'$ the diagram $U\xleftarrow{f} V \xrightarrow{g} W$ a homotopy pushout of $\mathscr{D}'$ is given by the cylinder \DG:
$$\mathscr{C}_{\mathscr{D}'} := \bigl((U \oplus sV \oplus W)_\bullet,\ d_\mathscr{C} = d_{\oplus} + d^{fg}\bigr)$$
where $d_{\oplus}$ is the differential on the \DG\ $(U \oplus sV \oplus W)$ and $d^{fg}$ is determined by $f$ and $g$:
$d^{fg}(sv) = f(v) + g(v)$.
\end{enumerate}
\end{thm}

We prove this in steps. 

\begin{lemma}
The objects $\mathscr{P}_\mathscr{D}$ and $\mathscr{C}_{\mathscr{D}'}$ defined above are indeed \DG\,s.
\end{lemma}
\begin{proof}
We already have that $(d_{\times}\, d_{\times}) = 0$, $(d_{fg}\, d_{fg}) = 0$, $(d_{\oplus}\, d_{\oplus}) = 0$ and $(d^{fg}\, d^{fg}) = 0$.  It remains to show only that $(d_{\times}\, d_{fg}) + (d_{fg}\, d_{\times}) = 0$ and $(d_{\oplus}\, d^{fg}) + (d^{fg}\, d_{\oplus}) = 0$.  These follow from the commutativity of $f$ and $g$ with the differentials $d_U$, $d_V$, and $d_W$:
\begin{align*}
d_{\times}\, d_{fg}(u,\, s^{-1}v,\, w) &= 
 d_{\times}\bigl(0,\,s^{-1}\bigl((f(u) + g(w)\bigr),\,0\bigr) \\ 
 &= \bigl(0,\ -s^{-1}\bigl(d_Vf(u) + d_Vg(w)\bigr),\ 0\bigr) \\
d_{fg}\,d_{\times} (u,\, s^{-1}v,\, w) &= 
 d_{fg}\bigl(d_Uu,\, -s^{-1}d_Vv,\, d_Ww\bigr) \\
 & = \bigl(0,\  s^{-1}\bigl(f(d_Uu) + g(d_Ww)\bigr),\  0\bigr)
\end{align*}
Similarly,
\begin{align*}
d_{\oplus}\,d^{fg}(u + sv + w) &=
 d_{\oplus} \bigl(f(v) + g(v)\bigr) \\
 & = d_Vf(v) + d_V g(v) \\
d^{fg}\,d_{\oplus}(u + sv + w) &=
 d^{fg}\bigl(d_Uu - sd_Vv + d_Ww\bigr) \\
 & = -f(d_Vv) - g(d_Ww)
\end{align*}
\end{proof}

\begin{lemma}\label{dg ho(co)lim functors}
The objects $\mathscr{P}_\mathscr{D}$ and $\mathscr{C}_{\mathscr{D}'}$ determine functors of diagrams in $\dg$
\begin{align*}
\mathscr{P}_{(-)}&:\dg^{(\bullet\xrightarrow{}\bullet\xleftarrow{}\bullet)}
 \xrightarrow{\ \ \ \ }\dg \\
\mathscr{C}_{(-)}&:\dg^{(\bullet\xleftarrow{}\bullet\xrightarrow{}\bullet)}
 \xrightarrow{\ \ \ \ }\dg
\end{align*}
\end{lemma}
\begin{proof}
We show this only for $\mathscr{P}_{(-)}$ since the proof for $\mathscr{C}_{(-)}$ is similar.

Suppose $t:\mathscr{D}_1 \to \mathscr{D}_2$ is a map of diagrams given by
\begin{equation}\label{E1:dg ho(co)lim functor}
\begin{aligned}
\xymatrix@R=7.5pt@C=13pt{
\mathscr{D}_1 \ar[dd]^t & & 
 U_1 \ar[r]^{f_1} \ar[dd]^{t_U} & V_1 \ar[dd]^{t_V} & W_1 \ar[l]_{g_1} \ar[dd]^{t_W} \\
{}\phantom{XXX}\ar@{|~|>}[rr] & & {}\phantom{XXX} & {}\phantom{XXX} & {}\phantom{XXX} \\
\mathscr{D}_2 & & 
 U_2 \ar[r]^{f_2} & V_2 & W_2  \ar[l]_{g_2}
}\end{aligned}
\end{equation}
This induces a \DG-map
\begin{equation}\label{E2:dg ho(co)lim functor}
\begin{aligned}
\xymatrix{
U_1\times s^{-1}V_1 \times W_1 \ar[d]^{(t_U\times s^{-1}t_V \times t_W)} \\ 
U_2\times s^{-1}V_2 \times W_2
}\end{aligned}
\end{equation}
We will show this extends to a \DG-map 
$(t_U\times s^{-1}t_V \times t_W):
 \mathscr{P}_{\mathscr{D}_1} \longrightarrow \mathscr{P}_{\mathscr{D}_2}$:

Since (\ref{E2:dg ho(co)lim functor}) is a \DG-map, we already have that $(t_U\times s^{-1}t_V \times t_W)$ commutes with $d_\times$.  It remains to show only that it commutes with $d_{fg}$.  However, this follows from the commutativity of (\ref{E1:dg ho(co)lim functor}) because
\begin{align*}
d_{f_2g_2}\,(t_U\times s^{-1}t_V \times t_W) &= 
  \bigl(0\,\times\, s^{-1}(f_2t_U + g_2t_W)\,\times\, 0\bigr) \\
(t_U\times s^{-1}t_V \times t_W)\,d_{f_1g_1} &=
  \bigl(0\,\times\, s^{-1}t_V(f_1 + g_1) \,\times\, 0 \bigr)
\end{align*}
\end{proof}

\begin{lemma}\label{dg ho(co)lim natural maps}
There are natural maps $e_\mathscr{P}:\lim_\dg\mathscr{D} \to \mathscr{P}_\mathscr{D}$ and $e_\mathscr{C}:\mathscr{C}_{\mathscr{D}'} \to \colim_\dg\mathscr{D}'$ given by $e_\mathscr{P}(u,\,w) = (u,\,0,\,w)$ and $e_\mathscr{C}(u+sv+w) = [u+w]$.
\end{lemma}
\begin{proof} The given maps are clearly natural so long as they are indeed \DG-maps.  The following shows commutativity of $e_\mathscr{P}$ with differentials:
\begin{align*} 
d_\mathscr{P}\, e_\mathscr{P}(u,\, w) &= d_\mathscr{P} (u,\, 0,\, w) \\
 &= \bigl(d_Uu,\ s^{-1}\bigl(f(u) + g(w)\bigr),\ d_Ww\bigr) \\
e_\mathscr{P}\, d_{\lim} (u,\,w) &= e_\mathscr{P} \bigl(d_Uu,\,d_Ww\bigr) \\
 &= \bigl(d_Uu,\ 0,\ d_Ww\bigr)
\end{align*}
These are equal since for $(u,\,w)\in U\times_V W$ we have $f(u) + g(w) = 0$.

Similarly,
\begin{align*}
d_{\colim}\, e_\mathscr{C}(u + sv + w) &= d_{\colim}[u+w] \\
 &= [d_Uu + d_Ww] \\
e_\mathscr{C}\, d_{\mathscr{C}} (u+sv+w) 
 &= e_\mathscr{C} \bigl(d_Uu + f(v) - sd_Vv + d_Ww + g(v)\bigr) \\
 &= \bigl[d_Uu + f(v) + d_Ww + g(v)\bigr]
\end{align*}
These are equal since $\bigl[f(v) + g(v)\bigr] = [0]$ in $U\oplus_V W$.
\end{proof}

To finish the proof of Theorem~\ref{l:dg pushout}, we will note that the above functors are related by zig-zags of natural weak equivalences to another pair of functors which we know (by other means) are homotopy limit and colimit functors.  Before doing this, however, we would like to make a short detour in order to give a generalization of the above construction to $n$-dimensional homotopy pullbacks and pushouts.

\subsubsection{Differential Bigraded Vector Spaces and $n$-Dimensional Pullbacks and Pushouts}

Recall that an $n$-dimensional pullback diagram in $\dg$ is a functor $\mathscr{D}:\mathcal{P}_0(\underline{n})\to \dg$ where $\mathcal{P}_0(\underline{n})$ is the poset of non-empty subsets of $\{1,\dots,n\}$ and inclusion maps.  Dually, an $n$-dimensional pushout diagram is a functor $\mathscr{D}':\mathcal{P}_{1}(\underline{n}) \to \dg$ where $\mathcal{P}_{1}(\underline{n})$ is the poset of proper subsets of $\{1,\dots,n\}$ and inclusion maps.  We will define the $n$-dimensional path or cylinder \DG\ associated to an $n$-dimensional pullback or pushout diagram to be the totalization of a certain differential bigraded vector space constructed from the diagram.

\begin{defn}[biDG]\label{biDG}
A differential bigraded vector space (bi\DG) is $V = (V_{\bullet, \bullet},\, d^h,\,d^v)$ where $V_{\bullet,\bullet}$ is a bigraded vector space and $d^h$ and $d^v$ are endomorphisms of $V_{\bullet,\bullet}$ such that $d^h$ has bidegree $(-1,0)$ and $d^v$ has bidegree $(0,-1)$ and the two satisfy $d^hd^h = 0$, $d^vd^v=0$, and $d^hd^v + d^vd^h=0$.

A bi\DG-map is a bidegree $(0,0)$ bigraded vector space map commuting with $d^h$ and $d^v$.  Write $\mathit{bi}\dg$ for the category of bi\DG\,s and bi\DG-maps.
\end{defn}

We write $V_{\bullet,\bullet}$ for the underlying bigraded vector space of the bi\DG\ $V$.  Also we refer to $(i,\bullet)$ as {\em horizontal grading $i$} and $(\bullet, i)$ as {\em vertical grading $i$}.  Similarly, we call $d^h:V_{i,\bullet}\to V_{i-1,\bullet}$  the {\em horizontal differential} and $d^v:V_{\bullet,i}\to V_{\bullet,i-1}$ the {\em vertical differential}.  And we say {\em total degree $i$} for $(\bullet,i-\bullet)$.

There is a standard functor $\tot:\mathit{bi}\dg \to \dg$ given by taking the total \DG\ of a bi\DG.   If $V=(V_{\bullet,\bullet},\, d^h,\,d^v)$ is a bi\DG, then the total \DG\ of $V$ is given by $(\tot V)_n = \oplus_{i+j=n} V_{i,j}$ with differential $d_{\tot}=d^h+d^v$.\footnote{[Weib] calls this $\tot^{\oplus}$ to differentiate it from $\tot^{\prod}$.  Since we will only consider finite bi\DG\,s these two will always be equal and we will not need to make this distinction.}  It follows immediately from the definition of bi\DG\,s that $d_{\tot}$ is indeed a differential on $(\tot V)_\bullet$.

It is clear that the symmetric monoidal operations $\otimes$ and $\oplus$ on graded vector spaces induce symmetric monoidal operations $\otimes$ and $\oplus$ on bigraded vector spaces.  Note that there are two suspension functors on bi\DG\,s.  Write $s^h$ for the horizontal suspension, given by tensoring with the bi\DG\ $s^h$ which is $\Q$ in bidegree $(1,0)$ and 0 elsewhere; and $s^v$ for the vertical suspension, given by tensoring with the bi\DG\ $s^v$ which is $\Q$ in bidegree $(0,1)$ and 0 elsewhere.

Given a \DG\ $V$ write $\bi(V)$ for the associated bi\DG\ with $\bi(V)_{i,0} = V_i$ and $\bi(V)_{i,j} = 0$ for $j\neq 0$, equipped with differentials $d^h_{\bi(V)} = d_V$, $d^v_{\bi(V)} = 0$.  This clearly defines a functor $\bi:\dg\to\mathit{bi}\dg$.  

Furthermore $\bi(-)$ extends to a functor $\bi:\mathcal{M}\!\mathit{or}(\dg)\to\mathit{bi}\dg$ as follows:  For $f:V\to W$ a \DG-map, let $\bi(f)$ be the bi\DG\ given by $\bigl(\bi(f)\bigr)_{\bullet,\bullet} = \bigl(s^v\,\bi(V)\oplus\bi(W)\bigr)_{\bullet,\bullet}$ with horizontal differential $d^h_{\bi(f)} = d^h_{s^v\bi(V)\oplus\bi(W)}$ and vertical differential $d^v_{\bi(f)} = f$.  [To show that this is a bi\DG\ we need only show that $d^h_{\bi(f)}d^v_{\bi(f)} + d^v_{\bi(f)}d^h_{\bi(f)} = 0$.  But this is equal to $d^h_{\bi(W)}f + fd^h_{s^v\bi(V)} = d_Wf + f(-d_V)$.]  That is, $\bi(f)$ is given by $V$ in vertical grading 1 and $W$ in vertical grading 0 equipped with a vertical differential given by $f$.

Even more generally, $\bi(-)$ extends to functors of $n$-dimensional pullback and pushout diagrams 
$\hat{\mathscr{P}}_{(-)}:\dg^{\mathcal{P}_0(\underline{n})}\longrightarrow\mathit{bi}\dg$
and 
$\hat{\mathscr{C}}_{(-)}:\dg^{\mathcal{P}_1(\underline{n})}\longrightarrow\mathit{bi}\dg$ as follows:
Given $n$-dimensional pullback and pushout diagrams $\mathscr{D}\in\dg^{\mathcal{P}_0(\underline{n})}$ and $\mathscr{D}'\in\dg^{\mathcal{P}_1(\underline{n})}$, define $\hat{\mathscr{P}}_\mathscr{D}$ and $\hat{\mathscr{C}}_{\mathscr{D}'}$ by 
\begin{align*}
\bigl(\hat{\mathscr{P}}_\mathscr{D}\bigr)_{\bullet,\bullet} 
  &= \biggl(\,\bigoplus_{T\in\mathcal{P}_0(\underline{n})} 
     (s^v)^{1-|T|}\,bi\bigl(\mathscr{D}(T)\bigr)\biggr)_{\bullet,\bullet} \\
\bigl(\hat{\mathscr{C}}_{\mathscr{D}'}\bigr)_{\bullet,\bullet} 
  &= \biggl(\,\bigoplus_{T\in\mathcal{P}_1(\underline{n})} 
     (s^v)^{n-1-|T|}\,bi\bigl(\mathscr{D}(T)\bigr)\biggr)_{\bullet,\bullet}
\end{align*}
with horizontal differentials given by the differentials of $\mathscr{D}(T)$ and $\mathscr{D}'(T)$ and vertical differentials given by the maps of $\mathscr{D}$ and $\mathscr{D}'$ with alternating signs.  More explicitly, the horizontal differentials are $d^h_{\hat{\mathscr{P}}} = d^h_{\oplus}$ and $d^h_{\hat{\mathscr{C}}} = d^h_{\oplus}$ and the vertical differentials are given by 
$$d^v_{\mathscr{P}} = 
  \!\!\!\!\underset{\sigma:S\to S\cup\{t\}}{\bigoplus}
  \!\!{\mathrm{sgn}(\sigma)}\,\mathscr{D}(\sigma)\qquad
\text{and} \qquad
d^v_{\mathscr{C}} = 
  \!\!\!\!\underset{\sigma:S\to S\cup\{t\}}{\bigoplus}
  \!\!{\mathrm{sgn}(\sigma)}\,\mathscr{D}'(\sigma)$$
where $\mathrm{sgn}(\sigma) = (-1)^{\bigl|\{s\in S\ |\ s<t\}\bigr|}$ for $\sigma:S\to S\cup\{t\}$.

\begin{ex} 
Let $\mathscr{D}\in\dg^{\mathcal{P}_0(\underline{3})}$ and $\mathscr{D}'\in\dg^{\mathcal{P}_1(\underline{3})}$ be three dimensional pullback and pushout diagrams.  Then $\hat{\mathscr{P}}_\mathscr{D}$ and $\hat{\mathscr{C}}_{\mathscr{D}'}$ are as follows:\footnote{See Example~\ref{pullback ex} for notation.}
\begin{itemize}
\item $\hat{\mathscr{P}}_\mathscr{D}$ is given by 
$$\xymatrix@C=14pt@R=7pt{
{}\vphantom{|}  &  X_{\{1\}} \ar[dd]_(.40){-} \ar[ddrr]_(.25){-} 
 &{}\save[] *{\oplus} \restore & X_{\{2\}} \ar'[dl]^(.5){+}[ddll] \ar'[dr]_(.5){-}[ddrr] 
 &{}\save[] *{\oplus} \restore & X_{\{3\}} \ar[ddll]^(.25){+} \ar[dd]^(.4){+} & {}\vphantom{|} &
 {\text{$\xleftarrow{\ \ }$ vertical grading $0$}}   \\
& && && && \\
{}\vphantom{|}  &  X_{\{1,2\}} \ar[ddrr]_(.4){+}      
 &{}\save[] *{\oplus} \restore & X_{\{1,3\}} \ar[dd]_(.45){-}        
 &{}\save[] *{\oplus} \restore & X_{\{2,3\}} \ar[ddll]^(.4){+}   &  {}\vphantom{|} &
 {\text{$\xleftarrow{\ \ }$ vertical grading $-1$}}  \\
& && && && \\
{}\vphantom{|}  &  && X_{\{1,2,3\}}  && & {}\vphantom{|} &
 {\text{$\xleftarrow{\ \ }$ vertical grading $-2$}}
}$$
\item $\hat{\mathscr{C}}_{\mathscr{D}'}$ is given by
$$\xymatrix@C=14pt@R=7pt{
{}\vphantom{|}  &  && X_\emptyset \ar[ddll]_(.55){+} \ar[dd]^(.45){+} \ar[ddrr]^(.55){+} && & {}\vphantom{|} &
 {\text{$\xleftarrow{\ \ }$ vertical grading $2$}} \\
& && && &&\\
{}\vphantom{|}  &  X_{\{1\}} \ar[dd]_(.4){-} \ar[ddrr]_(.25){-} 
 &{}\save[] *{\oplus} \restore & X_{\{2\}} \ar'[dl]^(.5){+}[ddll] \ar'[dr]_(.5){-}[ddrr] 
 &{}\save[] *{\oplus} \restore & X_{\{3\}} \ar[ddll]^(.25){+} \ar[dd]^(.4){+} & {}\vphantom{|} &
 {\text{$\xleftarrow{\ \ }$ vertical grading $1$}} \\
& && && &&\\
{}\vphantom{|}  &  X_{\{1,2\}}   
 &{}\save[] *{\oplus} \restore & X_{\{1,3\}}      
 &{}\save[] *{\oplus} \restore & X_{\{2,3\}} & {}\vphantom{|} &
 {\text{$\xleftarrow{\ \ }$ vertical grading $0$}}
}$$ 
\end{itemize}
where the marked arrows with the indicated signs give the vertical differentials.
\end{ex}


\begin{defn}[$n$-Dimensional Homotopy Pullback and Pushouts]
For  $\mathscr{D}\in\dg^{\mathcal{P}_0(\underline{n})}$ an $n$-dimensional pullback diagram, define $\mathscr{P}_\mathscr{D} := \tot\bigl(\hat{\mathscr{P}}_\mathscr{D}\bigr)$.

Similarly, for $\mathscr{D}'\in\dg^{\mathcal{P}_1(\underline{n})}$ an $n$-dimensional pushout diagram, define  
$\mathscr{C}_{\mathscr{D}'} := \tot\bigl(\hat{\mathscr{C}}_{\mathscr{D}'}\bigr)$.
\end{defn}

\subsubsection{Return to Homotopy Limits and Colimits}

In general we may construct small homotopy limits and colimits in $\dg$ in the vein of Bousfield and Kan [BK72].  Let $B$ be the functor from small categories to \DG\,s given by taking the nerve of a category and then applying the functor $-_\mDG:\sset \to \dg$ described in (\ref{adj:sset dg}).  We may use this to define homotopy limits and colimits according to the end and coend construction of Bousfield and Kan given in Lemma~\ref{D:holim end} in the previous section.  

Bousfield and Kan's proof that the end and coend construction is weakly equivalent to the adjoint functor definition formally applies to this setting as well.  Furthermore their proof that their adjoint functor construction satisfies the properties of a homotopy limit and colimit also formally applies to this setting.

Note that the Bousfield-Kan construction yields slightly larger models for homotopy pushouts and pullbacks than the construction mentioned above.  In particular for $\mathscr{D}$ and $\mathscr{D}'$ the 2-dimensional pullback and pushout diagrams in Lemma~\ref{l:dg pushout}, the Bousfield-Kan construction gives 
\begin{itemize}
\item Bousfield-Kan homotopy limit of $\mathscr{D}$ is $\holim(\mathscr{D}) = \bigl((U\times s^{-1} V \times V \times s^{-1}V \times W)_\bullet,\, d_{\holim}\bigr)$ 
\item Bousfield-Kan homotopy colimit of $\mathscr{D}'$ is $\hocolim(\mathscr{D}') = \bigl((U\oplus sV \oplus V \oplus sV \oplus W)_\bullet,\, d_{\hocolim}\bigr)$ 
\end{itemize}
where $d_{\holim}$ and $d_{\hocolim}$ are the obvious extensions of $d_{{\mathscr{P}}}$ and $d_{{\mathscr{C}}}$:
\begin{itemize}
\item $d_{\holim} = d_\times + d_{f\Idm\Idm g}$ where 
 $$d_{f\Idm\Idm g}\bigl(u,\,s^{-1}v_0,\,v_1,\,s^{-1}v_2,\,w\bigr) 
 = \Bigl(0,\ s^{-1}\bigl(f(u)+v_1\bigr),\ 0,\ s^{-1}\bigl(v_1 + g(w)\bigr),\ 0\Bigr)$$
\item $d_{\hocolim} = d_\oplus + d^{f\Idm\Idm g}$ where
 $$d^{f\Idm\Idm g}\bigl(u\oplus s^{-1}v_0\oplus v_1\oplus s^{-1}v_2\oplus w\bigr)
 = \bigl(f(v_0)\oplus 0 \oplus (v_0+v_2) \oplus 0 \oplus g(v_2)\bigl)$$
\end{itemize}

However, (just as in the cases of $\hat \Omega$ versus $\Omega$ and $\hat \Sigma$ versus $\Sigma$) there are clear quasi-isomorphisms ${\mathscr{P}}_\mathscr{D} \xrightarrow{\weq} \holim(\mathscr{D})$ and  $\hocolim(\mathscr{D}')  \xrightarrow{\weq} {\mathscr{C}}_{\mathscr{D}'}$.  In fact, there are two quasi-isomorphisms ${\mathscr{P}}_\mathscr{D} \xrightarrow{\weq} \holim(\mathscr{D})$ given by 
\begin{itemize}
\item $F(u,\,s^{-1}v,\,w) = \bigl(-u,\ s^{-1}d_Vf(u),\ f(u),\ s^{-1}v,\ w)$ 
\item $G(u,\,s^{-1}v,\,w) = \bigl(u,\ s^{-1}v,\ g(w),\ s^{-1}d_Vg(w),\ w)$
\end{itemize}
and similarly for $\hocolim(\mathscr{D}')  \xrightarrow{\weq} {\mathscr{C}}_{\mathscr{D}'}$.\footnote{These are clearly quasi-isomorphisms if they are DG-maps.  We leave this computation as an exercies.}

In the category of topological spaces the functors $\Sigma X$ and $\Omega X$ may be defined as the homotopy pushout of $\ast \xleftarrow{} X \xrightarrow{} \ast$ and the homotopy pullback of $\ast \xrightarrow{} X \xleftarrow{} \ast$ respectively.  We have defined $\Omega$ and $\Sigma$ in $\dg$ so that this is again the case.

\begin{cor}
$\Omega V$ is the homotopy pullback in $\dg$ of the diagram $0_\mDG \xrightarrow{\ } V \xleftarrow{\ } 0_\mDG$.

$\Sigma V$ is the homotopy pushout in $\dg$ of the diagram $0_\mDG \xleftarrow{\ } V \xrightarrow{\ } 0_\mDG$.
\end{cor}

Recall that homotopy fibers and cofibers of maps may be defined as homotopy pushouts and pullbacks: 
\begin{itemize}
\item The homotopy fiber of $f:V\to W$ is the homotopy pullback of $V \xrightarrow{f} W \xleftarrow{} 0_\mDG$
\item The homotopy cofiber of $f:V\to W$ is the homotopy pushout of $0_\mDG \xleftarrow{} V \xrightarrow{f} W$
\end{itemize}

\begin{cor}\label{l:dg hofib}
If $f:V\to W$ is a \DG-map, then its homotopy fiber and cofiber are:
\begin{itemize}
\item $\hofib(f) \weq (V_\bullet \oplus s^{-1}W_\bullet,\, d_{\hofib})$ where $d_{\hofib}$ is defined by $d_{\hofib}(v) = d_V(v) - s^{-1}f(v)$, and $d_{\hofib}(s^{-1}w) = -s^{-1}d_W(w)$ for $v\in V$ and $w \in W$.
\item $\hocof(f) \weq (W_\bullet \oplus sV_\bullet,\, d_{\hocof})$ where $d_{\hocof}$ is defined by $d_{\hocof}(sv) =  f(v) - sd_V(v)$, and $d_{\hocof}(w) = d_W(w)$ for $v\in V$ and $w \in W$.
\end{itemize}
\end{cor} 

Recall that a commutative square given by 
$$\xymatrix@=15pt{
U \ar[r]^{t} \ar[d]^{s} & V \ar[d]^{g} \\ 
W \ar[r]^{f} & X
}$$ 
is called {\em homotopy cartesian} if the composition of canonical maps 
$$U \longrightarrow \lim\bigl(W\xrightarrow{f}X\xleftarrow{g}V\bigr) 
 \longrightarrow \holim\bigl(W\xrightarrow{f}X\xleftarrow{g}V\bigr)$$ 
is a weak equivalence, and {\em homotopy cocartesian} if the composition of canonical maps
$$\hocolim\bigl(W\xleftarrow{s}U\xrightarrow{t}V\bigr) \longrightarrow 
 \colim\bigl(W\xleftarrow{s}U\xrightarrow{t}V\bigr)\longrightarrow X$$
is a weak equivalence.  In the model category $\dg$ we are in the very pleasant situation where

\begin{lemma}\label{dg cart=cocart} 
A commutative square in $\dg$ is homotopy cartesian if and only if it is homotopy cocartesian (i.e. $\dg$ is a stable model category).
\end{lemma}
\begin{proof}
This follows from the fact that the statement that the square is homotopy cartesian and the statement that the square is homotopy cocartesian are each equivalent to the statement that the following bi\DG\ has trivial homology:
$$\xymatrix@C=14pt@R=7pt{
{}\vphantom{|}  &  && U \ar[ddl]_{s} \ar[ddr]^{t} && & {}\vphantom{|} &
 {\text{$\xleftarrow{\ \ }$ vertical grading $1$}} \\
& && && &&\\
{}\vphantom{|}  &  
 & W  \ar[ddr]_{-f}  &  {}\save[] *{\oplus} \restore
 & V \ar[ddl]^{g} & 
 & {}\vphantom{|} &
 {\text{$\xleftarrow{\ \ }$ vertical grading $0$}} \\
& && && &&\\
{}\vphantom{|}  &  
 & & X     
 & & & {}\vphantom{|} &
 {\text{$\xleftarrow{\ \ }$ vertical grading $-1$}}
}$$ 
\end{proof}

\begin{note}In particular the squares
$$\begin{aligned} \xymatrix@R=10pt@C=10pt{
 V \ar[r] \ar[d]    & \mathrm{c}V \ar[d] \\ 
 \mathrm{c}V \ar[r] &  \Sigma V} \end{aligned}
\qquad \mathrm{and} \qquad
\begin{aligned} \xymatrix@R=10pt@C=10pt{
 \Omega V \ar[r] \ar[d] & \mathrm{p}V \ar[d] \\
 \mathrm{p}V \ar[r]     &  V} \end{aligned}
$$
are each both homotopy cartesian and homotopy cocartesian.
\end{note}

It is a standard property of stable model categories that:

\begin{lemma}\label{dg holim commute}
Very small homotopy limits commute with homotopy colimits in $\dg$.

Very small homotopy colimits commute with homotopy limits in $\dg$.
\end{lemma}

\begin{cor}
In particular, all sequential homotopy colimits commute with homotopy pullbacks in $\dg$.
\end{cor}

\section{Model Category Structure of $\dg_r$}\label{dg_r model category}

The standard model category structure on $\dg_r$ is to take quasi-isomorphisms to be weak equivalences, maps surjective in degree $>r$ to be fibrations, and degree-wise injections to be cofibrations.   Just as in $\dg$ all objects in $\dg_r$ are both fibrant and cofibrant.  Also, all objects in $\dg_r$ are formal as well.



Under this model category structure, the inclusion and reduction functors defined in \ref{dg red} as well as the suspension and desuspension functors ($s$ and $s^{-1}$) are Quillen adjoint pairs.

\begin{lemma} 
The adjoint pair given in Lemma~\ref{dg red}
$$\mathrm{incl}:  \dg_r  \rightleftarrows \dg : \mathrm{red}_r $$
is a Quillen adjoint pair.  
Also the adjoint pair 
$$s:\dg_r \rightleftarrows \dg_{r+1}:s^{-1}$$ 
is a Quillen adjoint pair.

Furthermore, both of the pairs above satisfy the stronger condition each adjoint preserves all weak equivalences.
\end{lemma}
\begin{proof}
The inclusion of categories functor preserves cofibrations, since all cofibrations of $\dg_r$ are also cofibrations of $\dg$.  Similarly the suspension functor $s$ preserves cofibrations.  Also, all of the above functors clearly preserve weak equivalences.

Since the left adjoints preserve cofibrations and trivial cofibrations, we must have that the right adjoints preserve fibrations and trivial fibrations and the adjoint pairs are each Quillen pairs.
\end{proof}

It follows that $\Sigma:\dg_r\rightleftarrows \dg_r:\red{\Omega}$ is also a Quillen adjoint pair preserving all weak equivalences.
Furthermore, we are supplied with the following corollary by Lemma~\ref{adjoint holim}:

\begin{cor}\label{dg_r holim plan}
The right adjoints $\mathrm{red}_r(-):\dg\to\dg_r$ and $s^{-1}:\dg_{r+1}\to\dg_r$ preserve all homotopy limits.

The left adjoints $\mathrm{incl}:\dg_r\to\dg$ and $s:\dg_r\to\dg_{r+1}$ preserve all homotopy colimits.
\end{cor}

\subsection{Homotopy Limits and Colimits}

Homotopy limits and colimits in $\dg_r$ are not much different than in $\dg$.  In particular, homotopy limits and colimits in $\dg_r$ are ``created'' in $\dg$ in the following sense:

\begin{lemma}\label{dg_r holim}
Let $\holim_\dg$ and $\hocolim_\dg$ be the homotopy limit and colimit functors described in Section~\ref{dg holim}.  Homotopy limit and colimit functors on $\dg_r$ are given by: 
\begin{align*}
\holim{\!}_{\dg_r}\,\mathscr{D} &= 
  \mathrm{red}_r\,\bigl(\holim{\!}_\dg\,\mathrm{incl}(\mathscr{D})\bigr) \\
\hocolim{\!}_{\dg_r}\,\mathscr{D} &= 
  \mathrm{red}_r\,\bigl(\hocolim{\!}_\dg\,\mathrm{incl}(\mathscr{D})\bigr)
\end{align*}
Where $\mathrm{incl}$ is the functor on diagrams given by objectwise inclusion $\dg_r\to\dg$, and $\mathrm{red}_r$ is the $r$-reduction functor $\dg\to\dg_r$.
\end{lemma}

\begin{proof}
Let $\widehat{\holim}_{\dg_r}$ and $\widehat{\hocolim}_{\dg_r}$ be any homotopy limit and colimit functors on $\dg_r$.  We will show that there are zig-zags of natural weak equivalences giving 
$\widehat{\holim}_{\dg_r} \weq 
 \mathrm{red}_r\,\bigl(\holim{\!}_\dg\,\mathrm{incl}(\mathscr{D})\bigr)$ 
and 
$\widehat{\hocolim}_{\dg_r} \weq \
 \mathrm{red}_r\,\bigl(\hocolim{\!}_\dg\,\mathrm{incl}(\mathscr{D})\bigr)$:

By \ref{dg_r holim plan}, there are zig-zags of natural weak equivalences exhibiting
\begin{align}
\mathrm{red}_r\,\bigl(\holim{\!}_\dg(-)\bigr) &\weq 
 \widehat{\holim}_{\dg_r}\,\mathrm{red}_r(-) \label{E:dg_r holim}\\ 
\hocolim{\!}_\dg\,\mathrm{incl}(-) &\weq
 \mathrm{incl}\,\bigl(\widehat{\hocolim}_{\dg_r}(-)\bigr) \label{E:dg_r hocolim}
\end{align}
Recall that $\mathrm{incl}$ is a section of the reduction functor $\mathrm{red}_r:\dg\to\dg_r$.  Thus (\ref{E:dg_r holim}) and (\ref{E:dg_r hocolim}) yield zig-zags of natural weak equivalences showing 
\begin{align*}
\mathrm{red}_r\,\bigl(\holim{\!}_\dg\,\mathrm{incl}(-)\bigr) &\weq
 \,\widehat{\holim}_{\dg_r}\,\mathrm{red}_r\,\mathrm{incl}(-) = \widehat{\holim}_{\dg_r}(-)
\\
\mathrm{red}_r\,\bigl(\hocolim{\!}_\dg\,\mathrm{incl}(-)\bigr) &\weq
 \,\mathrm{red}_r\, \mathrm{incl}\,\bigl(\widehat{\hocolim}_{\dg_r}(-)\bigr) =
 \widehat{\hocolim}_{\dg_r}(-)
\end{align*}
\end{proof}

Note that the homotopy colimit functor defined in Section~\ref{dg holim}, takes diagrams of $r$-reduced \DG\,s to $r$-reduced \DG\,s.  Thus for $\mathscr{D}:\catI\to\dg_r$ a diagram of $r$-reduced \DG\,s $\hocolim_\dg\,\mathrm{incl}(\mathscr{D})$ is {\em already} $r$-reduced.  In particular, $\mathrm{red}_r$ does not change it.  In other words

\begin{cor}[Homotopy Limits and Colimits in $\dg_r$]
Let $\mathscr{D}:\catI\to\dg_r$ be an $\catI$-diagram in $\dg_r$.  Then 
\begin{itemize}
\item The homotopy limit of $\mathscr{D}$ in $\dg_r$ is given by the reduction of the homotopy limit of $\mathscr{D}$ in $\dg$.
\item The homotopy colimit of $\mathscr{D}$ in $\dg_r$ is given by the homotopy colimit of $\mathscr{D}$ in $\dg$.
\end{itemize}
\end{cor}

\begin{cor}
$\red{\Omega} V$ is the homotopy pullback in $\dg_r$ of the diagram $0_\mDG \xrightarrow{\ } V \xleftarrow{\ } 0_\mDG$.

$\Sigma V$ is the homotopy pushout in $\dg_r$ of the diagram $0_\mDG \xleftarrow{\ } V \xrightarrow{\ } 0_\mDG$.
\end{cor}

\begin{note} 
The category $\dg_r$ is no longer stable.  For example, the square
$$\xymatrix@=15pt{
0_\mDG \ar[r] \ar[d] & 0_\mDG \ar[d] \\
0_\mDG \ar[r] & s^r
}$$
is homotopy cartesian in $\dg_r$, but not homotopy cocartesian.
\end{note}

\begin{lemma}\label{dg_r holim commute}
In $\dg_r$, sequential homotopy colimits commute with homotopy pullbacks.
\end{lemma}
\begin{proof}This may be shown by an explicit calculation.  We give the calculation for 2-dimensional pullbacks.  In general the the calculation is the same, but bigger.

First, note that the homotopy colimit of the sequential diagram $\catI$
 $$V_1 \xrightarrow{\,g_1\,} V_2 \xrightarrow{\,g_2\,} V_3 \xrightarrow{\,g_3\,} V_4 \xrightarrow{\,g_4\,} \cdots$$
is given by  the mapping telescope $\Bigl(\bigl(V_1 \oplus sV_1 \oplus V_2 \oplus sV_2 \oplus V_3 \oplus sV_3 \oplus \cdots\bigr)_\bullet,\ d = d_\oplus + d^\catI\Bigr)$  where $d^\catI(sv_i) = v_i + g_i(v_i)$, for $sv_i \in sV_i$ (and $d^\catI(v_i) = 0$ for $v_i\in V_i$).

Let $\catJ$ be the category $(\bullet\xrightarrow{\ \ }\bullet\xleftarrow{\ \ }\bullet)$ and $\mathscr{D}$ be the $\catI\!\times\!\catJ$-diagram in $\dg_r$ given by
$$\xymatrix{
W_1 \ar[r]^{f_1} \ar[d]_{r_1} & W_2 \ar[r]^{f_2} \ar[d]_{r_2} & 
  W_3 \ar[r]^{f_3} \ar[d]_{r_3} & W_4 \ar[r]^{f_4} \ar[d]_{r_4} & \cdots \\
V_1 \ar[r]^{g_1} & V_2 \ar[r]^{g_2} & V_3 \ar[r]^{g_3} & V_4 \ar[r]^{g_4} & \cdots \\
U_1 \ar[r]^{h_1} \ar[u]^{t_1} & U_2 \ar[r]^{h_2} \ar[u]^{t_2} & 
  U_3 \ar[r]^{h_3} \ar[u]^{t_3} & U_4 \ar[r]^{h_4} \ar[u]^{t_4} & \cdots 
}$$ 
A short computation shows that $\displaystyle [\hocolim_\catJ \holim_\catI \mathscr{D}]_\mG$ is
$$
 \red{(U_1\times s^{-1} V_1 \times W_1)} 
 \oplus s\red{(U_1\times s^{-1}V_1\times W_1)} 
 \oplus \red{(U_2\times s^{-1}V_2\times W_2)} 
 \oplus s\red{(U_2 \times s^{-1}V_2\times W_2)} \oplus \cdots
$$
Collecting infinte $\oplus$ from the finite $\times$'s gives a \G-isomorphism to
$$
 \red{(U_1 \oplus sU_1 \oplus U_2 \oplus sU_2\oplus \cdots) \times 
 (s^{-1}V_1 \oplus V_1 \oplus s^{-1}V_2 \oplus V_2 \oplus \cdots) \times  
 (W_1 \oplus sW_1 \oplus W_2 \oplus sW_2\oplus \cdots)} 
$$
which is $\displaystyle [\holim_\catI \hocolim_\catJ \mathscr{D}]_\mG$.  Also, the differentials of both $\displaystyle \hocolim_\catJ \holim_\catI \mathscr{D}$ and 
$\displaystyle \holim_\catI \hocolim_\catJ \mathscr{D}$ are given by $d = d_{\times\oplus} + d^\catJ + d^\catI$ where 
\begin{itemize}
\item $d^\catJ$ is zero on $s^{-1}V_i$ and on $V_i$ and elsewhere is  $d^\catJ(w_i) = s^{-1}r_i(w_i)$, $d^\catJ(u_i) = s^{-1}t_i(u_i)$, $d^\catJ(sw_i) = -r_i(w_i)$, and $d^\catJ(su_i) = -t_i(u_i)$.
\item $d^\catI$ is zero on $U_i$, $W_i$, and $s^{-1}V$ and elsewhere is  $d^\catI(sw_i) = w_i + f_i(w_i)$, $d^\catI(su_i) = u_i + h_i(u_i)$ and $d^\catI(v_i) = -s^{-1}v_i - s^{-1}g(v_i)$. 
\end{itemize}

\end{proof}

More generally, all filtered homotopy colimits commute with small homotopy limits in $\dg_r$, but the previous result is enough for our purposes.

\section{Symmetric DG\,s}\label{symm dg}

Write $\Sigma_n$ for the symmetric group on $n$ elements.  A \DG\ with $\Sigma_n$ action is a diagram $V:{\boldsymbol{\Sigma}}_n \to \dg$ where ${\boldsymbol{\Sigma}}_n$ is the group $\Sigma_n$ viewed as a category in the standard way.\footnote{The category ${\boldsymbol{\Sigma}}_n$ has only one object and it has morphisms labelled by elements of $\Sigma_n$.}  The usual notation is to write $V$ for both the diagram and for the value of the diagram on the single object of its indexing category.  If $V\in \dg^{{\boldsymbol{\Sigma}}_n}$ then the $\Sigma_n$-fixed points and $\Sigma_n$-orbits of $V$ are 
$$V^{\Sigma_n} = \lim_{{\boldsymbol{\Sigma}}_n} V\qquad \text{and}\qquad
 V_{\Sigma_n} = \colim_{{\boldsymbol{\Sigma}}_n} V$$
The canonical map given by the composition $V^{\Sigma_n} \to V \to V_{\Sigma_n}$ is called the trace map.  Since we are working over $\Q$, the trace map is an isomorphism -- its inverse is given by the norm map $V_{\Sigma_n}\to V^{\Sigma_n}$  
$$\mathrm{tr}:[v] \mapsto \frac{1}{|\Sigma_n|}\,\sum_{\sigma \in \Sigma_n} \sigma(v)$$
(This is a \DG-map since the maps $\sigma:v\mapsto \sigma(v)$ are each \DG-maps.)

For $V \in \dg^{{\boldsymbol{\Sigma}}_n}$, $\Sigma_n$-homotopy fixed points and $\Sigma_n$-homotopy orbits of $V$ are similarly defined as 
$$V^{h\Sigma_n} = \holim_{{\boldsymbol{\Sigma}}_n} V\qquad \text{and}\qquad
 V_{h\Sigma_n} = \hocolim_{{\boldsymbol{\Sigma}}_n} V$$
However in this case, the fixed points and orbits are already homotopy functors, so in particular they give ${\boldsymbol{\Sigma}}_n$-homotopy limit and ${\boldsymbol{\Sigma}}_n$-homotopy colimit functors.  That is, the orbits and fixed points {\em are} homotopy orbits and homotopy fixed points.

\begin{lemma}
The functors $\displaystyle \lim_{{\boldsymbol{\Sigma}}_n}$ and $\displaystyle \colim_{{\boldsymbol{\Sigma}}_n}$ are homotopy functors $\dg^{{\boldsymbol{\Sigma}}_n} \to \dg$.
\end{lemma}
\begin{proof}
It is enough to show this only for $\displaystyle \lim_{{\boldsymbol{\Sigma}}_n}$ since the two are isomorphic.

Note that any $\Sigma_n$-equivariant map $f:V\to W$ restricts to a map $V^{\Sigma_n} \to W^{\Sigma_n}$.  Thus the natural map $\displaystyle \lim_{{\boldsymbol{\Sigma}}_n} f = f|_{V^{\Sigma_n}}$.  
Suppose that $f:V\xrightarrow{\,\weq\,}W$ is a quasi-isomorphism.  By brute force, we show that $f|_{V^{\Sigma_n}}:V^{\Sigma_n} \to W^{\Sigma_n}$ is also a quasi-isomorphism:  Write $Avg_{\Sigma_n}$ for the average over the group maps in $V$ and $W$
$$Avg_{\Sigma_n}(v) = \frac{1}{|\Sigma_n|}\sum_{\sigma\in\Sigma_n} \sigma(v) \qquad \text{and}\qquad
  Avg_{\Sigma_n}(w) = \frac{1}{|\Sigma_n|}\sum_{\sigma\in\Sigma_n} \sigma(w)$$
Note that $Avg_{\Sigma_n}$ gives \DG-maps to $V^{\Sigma_n}$ and $W^{\Sigma_n}$ which are the identity on $V^{\Sigma_n}\subset V$ and $W^{\Sigma_n}\subset W$.

First we show $H_*(f|_{V^{\Sigma_n}})$ is a surjection.  Suppose $w\in W^{\Sigma_n}$ is a cycle.  The map $f$ is a quasi-isomorphism, so there are $v\in V$ and $w_0\in W$ so that $f(v) = w + d_W w_0$.  Then $Avg_{\Sigma_n}(v)\in V^{\Sigma_n}$ and $Avg_{\Sigma_n}(w_0)\in W^{\Sigma_n}$ with $f\bigl(Avg_{\Sigma_n}(v)\bigr) = w + d_W Avg_{\Sigma_n}(w_0)$.

Now we show $H_*(f|_{V^{\Sigma_n}})$ is an injection.  Suppose $v_1,v_2 \in V^{\Sigma_n}$ such that there is $w\in W^{\Sigma_n}$ with $f(v_1) - f(v_2) = d_W w$.  Since $f$ is a quasi-isomorphism there is $v\in V$ with $v_1 - v_2 = d_V v$.  Then $Avg_{\Sigma_n}(v)\in V^{\Sigma_n}$ with $v_1 - v_2 = d_V Avg_{\Sigma_n}(v)$.
\end{proof}

This same proof also shows that orbits and fixed points are homotopy functors on $\dg_r$ (and $\dgl_r$ and $\dgc_r$) as well; though, we will not need this fact.

\chapter{Differential Graded Lie Algebras}\label{S:DGL} 
\markboth{Ben Walter}{I.4  Differential Graded Lie Algebras}

Differential graded Lie algebras are models for rational spaces.  The primary goal of this chapter is the construction in \S\ref{S:dgl holim}  of specific simple models for homotopy pullbacks and pushouts in $\dgl_r$.  The models which we give are essentially lifted from the category $\dg_r$.  In \S\ref{s:adjoints dgl} we note the close connection between $\dg_r$ and $\dgl_r$, which we use to ``lift'' homotopy limits and colimits from $\dg_r$ to those of $\dgl_r$.

\section{Category Structure}\label{dgl S1}

\begin{defn}[DGL] 
A {\em differential graded Lie algebra} (\DGL) $L = (L_\bullet,\, d_L,\, [-,-]_L)$ consists of a differential graded vector space $(L_\bullet, d_L)$ equipped with a bilinear, degree zero, graded Lie bracket map $[-,-]_L:(L_\bullet,\, d_L)\otimes (L_\bullet,\, d_L) \to (L_\bullet,\, d_L)$ which  satisfies:
\begin{enumerate}
\item[(i)] $[x,y] = -(-1)^{|x|\cdot |y|}[y,x]$ 
   \hfill {\em (graded anti-symmetry)}
\item[(ii)]$\big[x,[y,z]\big] + (-1)^{|x|(|y|+|z|)}\big[y,[z,x]\big] + 
  (-1)^{|z|(|x| + |y|)}\big[z,[x,y]\big] = 0$ 
   \hfill {\em (graded Jacobi identity)}
\end{enumerate}
\end{defn}

\begin{note} 
Defining the graded Lie bracket map as a \DG-map imposes the following compatibility with the differential on the bracket:  $d_L[x,y] = [d_Lx, y] + (-1)^{|x|}[x, d_Ly]$.
\end{note}

Maps of \DGL\,s are given by graded vector space maps $f:V_\bullet \to W_\bullet$ which are degree 0 and preserve differentials and Lie brackets (i.e. $f(d_Vx) = d_Wf(x)$ and $f([x,y]) = [f(x), f(y)]$).  Isomorphisms, injections, and surjections of \DGL\,s are given by \DGL-maps which induce isomorphisms, injections, and surjections of underlying graded vector spaces.  Similarly, quasi-isomorphisms of \DGL\,s are given by \DGL-maps which are quasi-isomorphisms of underlying \DG\,s.

A \DGL\ is called $r$-reduced if its underlying \DG\ is $r$-reduced.  Recall that this means the \DGL\ is trivial below grading $r$.  All of the \DGL\,s which we consider are $r$-reduced for some $r$:  our convention is to write $\dgl_r$ for the category of all $r$-reduced \DGL\,s and \DGL-maps between them.  At times we may write $\dgl$ without any subscript.  By this we either mean that the $r$ is implicitly present and determined by context or else we mean the category $\dgl_0$ of 0-reduced \DGL\,s.  Note that $\dgl_r$ is a full subcategory of $\dgl_t$ for all $t<r$.

Just as in $\dg$, we may define a reduction functor mapping to more highly reduced \DGL\,s: 
\begin{defn}[Reduction]\footnote{
There is also a truncation functor, but we will be uninterested in it.}
Given a $t$-reduced \DGL\ $L\in\dgl_t$ for $t<r$ 
 the {\em $r$-reduction} $\mathrm{red}_r(L) = \red{L} \in \dgl_r$ is the sub-\DGL\ $\red{L}$ given by $\red{L}_i = L_i$ for $i>r$, $\red{L}_i = 0$ for $i<r$, and $\red{L}_r = \ker(d_L:L_r \to L_{r-1})$.  
\end{defn}
\begin{lemma}\label{dgl red}

Reduction is right adjoint to the inclusion of categories functor $\dgl_r \to \dgl_t$
\begin{align*}
[-]_{\mDGL_t}:\dgl_r &\rightleftarrows \dgl_t:\mathrm{red}_r 
\end{align*}
\end{lemma}

Note that the categories $\dgl_r$ are all pointed categories:  The \DGL\ $0_\mDGL$ on the trivial differential graded vector space $0_\mDGL = (0_\mDG,\, [-,-]=0)$ equipped with the trivial Lie bracket is both initial and final.  If a \DGL\ has trivial homology $H_*(L) = 0$ (for example, if $L = 0_\mDGL$) then we say that $L$ is {\em contractible}.  A \DGL\ $L$ is contractible if and only if the map $L\to 0_\mDGL$ is a quasi-isomorphism.  In general of course, the homology of a \DGL\ is a graded Lie algebra, and taking homology gives a functor $H_\ast:\dgl_r \to \gl_r$.

\subsection{Free DGLs}

We could just as well have defined \DGL\,s as graded Lie algebras equipped with a differential compatible with their Lie bracket (rather than as differential graded vector spaces equipped with a Lie bracket compatible with the differential).  We use this point of view in defining free \DGL\,s.  Given a graded vector space $V = \{V_i\}_{i\ge0}$ starting in degree 0, we may consider its tensor algebra $TV = \bigoplus_{k} V^{\otimes k}$.  This can be given a Lie algebra structure by defining the bracket of two elements to be the graded commutator, $[x,y] = x\otimes y - (-1)^{|x|\cdot|y|} y\otimes x$.  The sub-Lie algebra of this generated by $V$ is called the {\em free graded Lie algebra on $V$} and denoted by $\freeL_V$.  The functor $\freeL_{(-)}:\g_r \to \gl_r$ ($n\ge 0$) by $V \mapsto \freeL_V$ is the left adjoint of the forgetful functor $[-]_{\mG_r}:\gl_r \to \g_r$.

\begin{defn}[Free DGL]\label{free DGL} 
A {\em free differential graded Lie algebra} is a \DGL\ $L$ which, as a graded Lie algebra is isomorphic to a free graded Lie algebra:
$$[L]_\mGL \cong \freeL_V\qquad \text{where}\qquad V\in\g$$

Write $\Fdgl_r$ for the full subcategory of all free \DGL\,s in $\dgl_r$.
\end{defn}
Writing $L$ as $L=([L]_\mGL,\, d)$ and pushing the differential of $L$ across the isomorphism $[L]_\mGL \cong \freeL_V$, we get the equivalent statement that $L$ is free if and only if $L$ is isomorphic to a \DGL\ of the form $L\cong (\freeL_V,\, d)$ for some graded vector space $V$ and differential $d$.

\begin{note}\label{truly free} 
What we have just defined would more precisely be called a ``differential free graded Lie algebra'' since it is not left adjoint to the forgetful functor $[-]_{\mDG_r}:\dgl_r \to \dg_r$.  However, it is standard practice to call objects isomorphic to $(\freeL_V,\, d)$ ``free \DGL\,s'' nonetheless -- see e.g. [FHT].  The forgetful functor $[-]_{\mDG_r}:\dgl_r \to \dg_r$ does have a left adjoint, which we abusively denote by $\freeL_{(-)}:(V,\, d_V) \mapsto \freeL_{(V,\,d_V)}$.  The \DGL\ $\freeL_{(V,\,d_V)}$ is defined in the same manner as the \GL\ $\freeL_V$ was above, starting with a \DG\ instead of a \G.  At times we will wish to refer to \DGL\,s which are free in this true sense: we will signify these with the notation $L\cong\freeL_{(V,\,d_V)}$ or by calling them {\em truly free} \DGL\,s.
\end{note}

Closely related to $\Fdgl$ is the category $\hFdgl$ of all \DGL\,s of the form $(\freeL_V,\,d)$ and \DGL-maps between them.  Note that an object of $\hFdgl$ consists of the data of a free \DGL\ $L$ along with 
an isomorphism $L\cong (\freeL_V,\,d)$ (however the maps are not required to preserve these isomorphisms).  We may lazily write $\freeL_V = (\freeL_V,\,d)$ for such \DGL\,s with the differential unwritten (but not forgotten).  The categories $\hFdgl$ and $\Fdgl$ are equivalent but not isomorphic.  This allows us to prove many facts about free \DGL\,s by considering only those which specifically have the form $(\freeL_V,\,d)$.  In general our intuition about free \DGL\,s comes from our intuition about objects $(\freeL_V,\,d)$, and most statements which we make about free \DGL\,s we first make about objects $(\freeL_V,\,d)$.  

\DGL\,s of the form $\freeL_V = (\freeL_V,\, d)$ have an extra grading by bracket length (inherited from the word-length grading on $TV$).  Using this grading, we may write their differential as $d=d_0 + d_1 + \cdots$ where $d_i$ increases bracket length by $i$.  [To prove this it is enough to show it on the generating graded vector space $V$, since the differential respects the bracket.  However, there it is trivially true.]  $d\circ d = 0$ implies in particular that $d_0\circ d_0 = 0$.  Thus, given $(\freeL_V,\, d)$ the differential $d$ restricts to $d_V = d_0|_V = d|_V$ a differential on the generating graded vector space $V$.  We write $(V,\,d_V\!)$ for 
the resulting \DG\ consisting of the generating \G\ of $\freeL_V$ equipped with the restriction of the differential of $\freeL_V$ to this generating graded vector space.
Note that $(\freeL_V,\, d)$ is truly free in the sense of \ref{truly free} if $d=d_0$ -- in this case $(\freeL_V,\,d) = \freeL_{(V,\,d_V)}$. 

If $L$ is any \DGL, we write $(L)^{\rm ab}$ for the graded\footnote{Graded abelian means multiplication is graded commutative -- i.e. $ab = (-1)^{|a|\cdot |b|}ba$.  Note that though Lie brackets are already abelian on odd degree elements of $L$ they are not {\em graded} abelian.} abelianization $(L)^{\rm ab} := L/[-,-]$.    Although, strictly speaking, $(L)^{\rm ab}$ is an object of $\dgl$, since it has trivial Lie bracket we standardly abuse notation and treat it as an object of $\dg$ -- i.e. we write $(L)^{\rm ab}$ to mean $\big[(L)^\ab\big]_\mDG$.  If $L\cong (\freeL_V, d)$, then there is a clear isomorphism $[(L)^{\mathrm{ab}}]_\mG \cong V$.  Note that $L$ is truly free in the sense of \ref{truly free} if and only if $L \cong \freeL_{(L)^{ab}}$ (i.e. if there is some $(V,\,d_V)$ so that $L\cong \freeL_{(V,d_V)}$).  We will show that $(L)^{\mathrm{ab}}$ captures much of the homology information of $L$ even when $L$ is merely free in the weak sense of Definition~\ref{free DGL}.

We are now in a position to state the rational Hurewicz Theorem (for free \DGL\,s):
\begin{thm}[Rational Hurewicz (free DGLs)]\label{Hurewicz dgl} 
If $L\in\Fdgl_1$ is a free \DGL\, then the \DG\ $(L)^\ab$ is contractible if and only if $L$ is contractible.

Moreover, if $L_1,L_2\in\Fdgl_1$ are 1-reduced, free \DGL\,s then a \DGL-map $L_1 \xrightarrow{\,f\,} L_2$ induces a \DG-quasi-isomorphism $(L_1)^\ab \xrightarrow{\,(f)^\ab\,} (L_2)^\ab$ if and only if the map $f$ itself is a quasi-isomorphism.\footnote{That is $(-)^\ab$ detects and reflects quasi-isomorphisms of free \DGL\,s.}
\end{thm}

\begin{proof}[Proof Sketch]
Note that the first statement above is implied by the second.

Felix, Halperin, and Thomas give a proof of the second statement for free \DGL\,s in $\hFdgl$ in [FHT 22.2], which is enough to imply this statement.  Their proof uses the functor $\C:\dgl_{r-1}\to \dgc_r$ which is described in Chapter~\ref{S:Q homotopy} -- in particular this functor both detects and reflects quasi-isomorphisms.  Furthermore there is a natural weak equivalence $\bigl[\C(L)\bigr]_\mDG \xrightarrow{\,\weq\,}s(L)^\ab$ (see \ref{free homotopy}).  The result is therefore implied by the commutative diagram
$$\xymatrix@R=7pt@C=4pt{
 L_1\ar[dd]^f & &&& 
   & \bigl[\C(L_1)\bigr]_\mDG \ar[r]^{\weq} \ar[dd]^{[\C(f)]_\mDG} 
   & s(L_1)^\ab \ar[dd]^{s(f)^\ab} \\
 & & \ar@{|~|>}[rr] & {}\phantom{XXX} & &{}\phantom{XXXX}&{}\phantom{XXXXXXXX} \\
 L_2        & &&& & \bigl[\C(L_2)\bigr]_\mDG \ar[r]^{\weq}         & s(L_2)^\ab
}$$
and the fact that $[-]_\mDG$ and $s$ each both reflect and detect weak equivalences.

A more direct proof may be given by introducing the induced map on the bracket-length filtrations (see \S\ref{DGL to DGL}) of $L_1$ and $L_2$.  If $(f)^\ab$ is a quasi-isomorphism then the induced map of filtrations is a quasi-isomorphism on the fibers $H_n(f)$ (see \S\ref{DGL to DGL}) and so a quasi-isomorphism of the entire tower.

In $\Fdgl_2$ the theorem also follows immediately from the Hurewicz Theorem for topological spaces along with Theorem~\ref{Quillen Thm I} and Lemma~\ref{free homotopy} which is stated later.
\end{proof}

Note that the corresponding statement about non-free \DGL\,s is not true.  In particular, the abelienization functor on non-free \DGL\,s, $(-)^\ab:\dgl\to\dg$ does not preserve quasi-isomorphisms.  For example consider the quasi-isomorphism  
$$f:\bigl(\Q v \oplus \Q [v,v] \oplus \Q w,\ d(w)=[v,v]\bigr) \xrightarrow{\ \weq\ } 
 \bigl(\Q v,\ d=0,\ [-,-]=0\bigr)$$
sending $v$ to $v$ and the other generators to 0 (where $|v|=2k+1$).\footnote{In our description of the left-hand \DGL\ above, we mean that the differential of the other genrators is 0 and the brackets of other generators are 0.}  The abelianization of this is 
$$(f)^\ab:\bigl(\Q v \oplus \Q w,\ d=0\bigr) \longrightarrow \bigl(\Q v,\ d=0\bigr)$$
which is clearly not a quasi-isomorphism.







We are particularly interested in a special class of maps between free \DGL\,s which we call ``freely generated maps''.  Before defining these, however, we first describe another class of maps called ``free maps''.

 Recall that on $\gl$ the categorical coproduct is the ``free product'' written $\liesum$.  The free product of \GL\,s $L$ and $K$ is given by allowing $L$ and $K$ to freely generate all brackets between them (we give a more explicit definition later in \ref{dgl sum}).  Quillen defines free maps of \DGL\,s as follows ([Q69, II.5]):

\begin{defn}[Free Map]\label{free map}
A map of \DGL\,s $f:L\to K$ is a {\em free map} if as a map of graded Lie algebras, it is  isomorphic to an inclusion of $[L]_\mGL$ with free cokernel, as expressed by the  diagram
$$\xymatrix@R=7pt@C=25pt{
\bigl[L\bigr]_\mGL \ar[r]^{[f]_\mGL} \ar@{_(->}[dr]
   &  \bigl[K\bigr]_\mGL \\
   & \bigl[L\bigr]_\mGL \liesum \freeL_V
    \ar@{{}{}{}}[u]^*[right]{\cong}
}$$

\end{defn}

Note in particular that:

\begin{lemma}
Let $f:L\to K$ be a free map.
If $L$ is a free \DGL, then $K$ is also free.
\end{lemma}
\begin{proof}

Suppose $L$ is free.  Coproducts are given by colimits, so they commute with left adjoint functors.  In particular, coproducts commute with $\freeL_{(-)}:\g \to \gl$; thus $\freeL_V \liesum \freeL_W = \freeL_{V\oplus W}$.  Therefore if $L \cong (\freeL_W,\,d)$ and $f:L\to K$ is a free map, then 
$$[K]_\mGL \ \cong\  [L]_\mGL\liesum \freeL_V \ \cong\  \freeL_W \liesum \freeL_V \ =\  
  \freeL_{V\oplus W}$$
In particular, $K$ is free.  

[This also follows from the later result that the free maps are precisely the cofibrations in $\dgl$ and the free \DGL\,s are precisely the cofibrant \DGL\,s.]
\end{proof}

One useful way of constructing free maps in $\hFdgl$ is given by the following lemma:

\begin{lemma}\label{hfdgl}
Let $(\freeL_V,\,d),\, (\freeL_{W},\,d') \in \hFdgl_r$ and $f:(\freeL_V,\,d)\to(\freeL_{W},\,d')$ be a \DGL-map such that $[f]_\mGL = \freeL_{\hat f}$ where $\hat f:V \cofibr W$ is an injection.  Then $f$ is a free map.
\end{lemma}
\begin{proof}
Recall that injections in $\g$ always split.  Thus $\hat f$ may be written as $\hat f:V \cofibr V\oplus U$ where $W\cong V\oplus U$.  Then
$[f]_\mGL : \freeL_V \cofibr \freeL_{U} = \freeL_V \liesum \freeL_U$, is the inclusion map $\freeL_V \cofibr \freeL_V\liesum \freeL_U$.
\end{proof}

The actual class of maps which we are interested in are a little weaker than those in \ref{hfdgl}.  We call these maps ``freely generated maps''.  

\begin{defn}[Freely Generated Map]
Given $f:(\freeL_V,\,d)\to(\freeL_W,\,d')$ a \DGL-map.  The map $f$ is {\em freely generated} if $[f]_\mGL = \freeL_{\hat f}$ where $\hat f: V \to W$ is any \G-map.

Write $\hfdgl_r$ for the subcategory of $\hFdgl_r$ consisting of free \DGL\,s of the form $(\freeL_V,\, d)$ and all freely generated maps between them.
\end{defn}

\subsection{Adjoints between $\dgl$ and $\dg$}\label{s:adjoints dgl}

Recall the functor $[-]_{\mDG_r}:\dgl_r\to \dg_r$ which forgets the bracket structure of a \DGL.  Our notation is $\big[(L_\bullet,\, d,\, [-,-])\big]_{\mDG_r} = (L_\bullet,\, d)$.  The forgetful functor has a section the functor $[-]_{\mDGL_r}\!:\dg_r\to \dgl_r$ which equips a \DG\ with a trivial Lie bracket (i.e. $[V]_{\mDGL_r}$ is given by $[V]_{\mDGL_r}\!\! = (V,\, [-,-]=0)$).  We are particularly interested in associated functors to and from $\dg$ given by composing the above with the inclusion of categories functor and reduction functor.  Abusing notation, we write $[-]_{\mDG}$  and $[-]_{\mDGL_r}$ for the compositions 
\begin{align*}
 [-]_{\mDG}&:\dgl_r \xrightarrow{[-]_{\mDG_r}} \dg_r \xrightarrow{\mathrm{incl}} \dg \\
 [-]_{\mDGL_r}&:\dg \xrightarrow{\mathrm{red}_r} \dg_r \xrightarrow{[-]_{\mDGL_r}} \dgl_r
\end{align*} 
Oftentimes we will sloppily write just $[-]_\mDGL$, $[-]_\mDG$, and $[-]_\mG$ for the functors $[-]_{\mDGL_r}$, $[-]_{\mDG_r}$, and $[-]_{\mG_r}$ with the understanding that $r$ is implicitly present and determined by context.

The above compositions are part of adjoint pairs between $\dgl$ and $\dg$ given by:
\begin{lemma}[Adjoints]\label{dgl adjoints} 
The following give adjoints between $\dgl_r$ and $\dg$:\footnote{The truncation functor $\dg \to \dg_r$ may be used to extend the top adjoint pair to $\dg \rightleftarrows \dgl_r$, but since truncation doesn't preserve weak equivalences, we do not wish to go that far.}
$$\xymatrix@C=30pt@R=15pt{
& \freeL_{(-)}:\dg_r \ar@<2pt>[r]^(.43){\freeL_{(-)}} &
  \dgl_r:[-]_\mDG \phantom{XXx} \ar@<2pt>[l]^(.54){[-]_{\mDG_r}} }$$
$$\xymatrix@C=30pt@R=15pt{
(-)^\mathrm{ab}:\dgl_r \ar@<2pt>[r]^(.63){(-)^\mathrm{ab}} &
  \dg_r \ar@<2pt>[l]^(.34){[-]_{\mDGL_r}} \ar@<2pt>[r]^(.4){\mathrm{incl}} &
  \dg:[-]_{\mDGL_r} \ar@<2pt>[l]^(.6){\mathrm{red}_r}
}$$          
where $\mathrm{incl}$ is the inclusion of categories functor.

Furthermore these adjoints are sections in the following sense:
\begin{itemize}
\item $[-]_{\mDGL_r}:\dg_r \to \dgl_r$ is a  section of $[-]_{\mDG_r}:\dgl_r\to \dg_r$.
\item $\freeL_{(-)}:\dg_r \to \dgl_r$ is a  section of $(-)^\mathrm{ab}:\dgl_r \to \dg_r$.
\end{itemize}
\end{lemma}

Note that related to the rational Hurewicz Theorem is the much more obvious lemma:
\begin{lemma}\label{weq dgl equals weq dg}
If $L\in\dgl_0$ is a \DGL, then the \DG\ $[L]_\mDG$ is contractible if and only if $L$ is contractible.

Moreover if $L_1, L_2 \in \dgl_0$ are \DGL\,s then a \DGL-map $L_1\xrightarrow{f} L_2$ induces a \DG-quasi-isomorphism $[L_1]_\mDG \xrightarrow{[f]_\mDG} [L_2]_\mDG$ if and only if $f$ itself is a quasi-isomorphism.
\end{lemma}

\begin{note}\label{dgl factorization} 
Given $L$ a \DGL, a composition of the right adjoints $\bigl[[L]_\mDG\bigr]_{\mDGL}$ recovers $L$ up to bracket information: 
$$L \ =\  \Bigl(\Bigl[\bigl[[L]_\mDG\bigr]_\mDGL\Bigr]_\mDG,\, [-,-]_L\Bigr)$$  
Similarly, given $L \cong (\freeL_V,\, d)$ a free \DGL, a composition of the left adjoints $\freeL_{(L)^{\mathrm{ab}}}$ recovers $L$ up to higher order differential information: 
$$L \ \cong\ \bigl(\freeL_V,\ d = d_0 + d_{>0}\bigr) \ \cong\ \left(\bigl[\freeL_{(L)^{\mathrm{ab}}}\bigr]_\mGL,\ d=d_{\freeL_{(L)^{\mathrm{ab}}}}\!\!+d_{>0}\right)$$ 
where $d_0 = d_V$ is the degree 0 part of the differential of $(\freeL_V,\,d)$; $d_{>0} = d -  d_0$ is the part of the differential of $(\freeL_V,\, d)$ which increases bracket length by at least 1; and $d_{\freeL_{(L)^{\mathrm{ab}}}}$ is the differential on the truly free \DGL\ $\freeL_{(L)^{\mathrm{ab}}}$.
\end{note}

\subsection{Limits and Colimits}

The category $\dgl_r$ has both products and coproducts.  We denote the coproduct of two \DGL\,s by $\liesum$, and the product by $\times$.  Below, we neglect to mention the precise category where products, coproducts, limits and colimits are taken.  The understanding should be that if objects or diagrams are in $\dgl_r$ then the products, coproducts, limits and colimits are those of the category $\dgl_r$ as well.  We implicitly assume that $r\ge 0$ (though the definitions of products and limits given below are still correct for $r$ negative).  Note that the operations $\liesum$ and $\times$ are defined in such a way that for $X, Y$ simply-connected spaces we have:
\begin{itemize}
\item $\pi_\ast(X\vee Y) \,\cong\, \bigl(\pi_\ast(X)\otimes\Q\bigr)\,\liesum\,\bigl(\pi_\ast(Y)\otimes \Q\bigr)$
\item $\pi_\ast(X\times Y) \,\cong\, \bigl(\pi_\ast(X)\otimes\Q\bigr)\times\bigl(\pi_\ast(Y)\otimes \Q\bigr)$
\end{itemize}

\begin{lemma}[Sums]\label{dgl sum} 
Given two \DGL\,s $L_1 = (V_\bullet,\, d_V,\, [-,-]_V)$ and $L_2 = (W_\bullet,\, d_W,\, [-,-]_W)$, their categorical coproduct is given by their free product 
$$L_1\liesum L_2 := (\freeL_{(V\oplus W)}/I,\, d_\liesum)$$ 
where $I$ is the ideal of $\freeL_{(V\oplus W)}$ generated by elements of the form $i([v_1,v_2]_V) - [i(v_1), i(v_2)]$ and $j([w_1, w_2]_W) - [j(w_1), j(w_2)]$, $i:V \hookrightarrow V\oplus W$, $j:W\hookrightarrow V\oplus W$; and $d_{\liesum}$ is the free differential induced on $\freeL_{V\oplus W}$ by the differentials $d_V$ and $d_W$.

This generalizes to give all small coproducts.
\end{lemma}

\begin{lemma}[Products]\label{dgl prod} 
Given \DGL\,s $L_1 = (V_\bullet,\, d_V,\, [-,-]_V)$ and $L_2 = (W_\bullet,\, d_W,\, [-,-]_W)$, their categorical product is the \DGL\ 
$$L_1\times L_2 := (V_\bullet \times W_\bullet,\, d_{\times},\, [-,-]_V \times [-,-]_W)$$ 
where $d_{\times}$ is the differential induced on $V_\bullet \times W_\bullet$ by the differentials $d_V$ and $d_W$
(it is defined by $d_{\otimes}(v,w) = (d_Vv,d_Ww)$).

This generalizes to give all small products.
\end{lemma}

In general, if $L_1$ and $L_2$ are \GL\,s then $L_1\liesum L_2$ is their categorical sum in $\gl$.  Recall that left adjointness of $\freeL_{(-)}$ implies that $\freeL_V \liesum \freeL_W = \freeL_{V\oplus W}$.  From this it follows that on free \DGL\,s the categorical product and coproduct (of $\dgl$) may be given a slightly simpler form.  We make critical use of this when defining cones, suspensions, and general homotopy colimits.

\begin{cor}[Sums on $\Fdgl$]\label{fdgl sum}
Given the two free \DGL\,s (with chosen isomorphisms) $L_1 \cong (\freeL_V,\, d_1)$ and $L_2 \cong (\freeL_W,\, d_2)$ their coproduct is isomorphic to 
$$L_1 \liesum L_2 \cong (\freeL_{V\oplus W},\, d)$$
where $d$ is the differential induced by $d_1$ and $d_2$.
\end{cor}

\begin{cor}[Products on $\Fdgl$]\label{products on fdgl} 
Given the two free \DGL\,s (with chosen isomorphisms) $L_1 \cong (\freeL_V,\, d_1)$ and $L_2 \cong (\freeL_W,\, d_2)$ there is a quasi-isomorphism  of \DGL\,s
$$L_1 \times L_2 \cong (\freeL_V,\, d_1) \times (\freeL_W,\, d_2) \xleftarrow{\,\weq\,} 
  (\freeL_V\liesum\freeL_W\liesum\freeL_{s(V\otimes W)},\,d) 
  = \big(\freeL_{V\oplus W \oplus s(V\otimes W)},\, d\big)$$
where $d$ is $d_1$ on $\freeL_V$, $d_2$ on $\freeL_W$, $d(s(v\otimes w)) = [v, w] - s(d_1(v)\otimes w) - (-1)^{|v|}s(v\otimes d_2(w))$ on the generators of $\freeL_{s(V\otimes W)}$, and $d$ is extended freely to brackets of these elements.
\end{cor}
\begin{proof}[Proof Sketch]
The quasi-isomorphism above is 
$f:(\freeL_V\liesum\freeL_W\liesum\freeL_{s(V\otimes W)},\,d) \longrightarrow (\freeL_V,\, d_1) \times (\freeL_W,\, d_2)$ given by the identity on $\freeL_V$ and $\freeL_W$, 0 on $\freeL_{s(V\otimes W)}$, and 0 on all mixed brackets.
\end{proof}

Finally, note that the category $\dgl_r$ has pullbacks and pushouts:

\begin{lemma}[Pullbacks and Pushouts]\label{dgl (co)lim} 
The pullback of the \DGL-diagram $L_1\xrightarrow{f_1} K\xleftarrow{f_2} L_2$ is given by:
$$L_1\times_{K} L_2 := \{(l_1,l_2)\in L_1\times L_2\ \ |\ f_1(l_1) + f_2(l_2) = 0\}$$ 

The pushout of the \DGL-diagram $L_1 \xleftarrow{f_1} K\xrightarrow{f_2} L_2$ is given by:
$$L_1\liesum_{K} L_2 := (L_1\liesum L_2)/\langle f_1(k) + f_2(k)\rangle$$
where $\langle f_1(k) + f_2(k)\rangle$ is the two-sided Lie-ideal generated by the elements $f_1(k) + f_2(k)$.

In general all small limits and colimits are given by the appropriate equalizers and coequalizers (kernels and cokernels) of products and coproducts as in Theorem~\ref{(co)lim (co)equalizer}.
\end{lemma}

\begin{cor}
Let $f$ be a map $f:K \to L$.  The fiber of $f$ is $\ker(f)$; the cofiber of $f$ is $\coker(f)$.
\end{cor}

\begin{note}\label{dgl lim plan} 
We can also recover the formula for the categorical product and more generally all limits in $\dgl_r$ using the adjoint functors from the previous section along with the comments in Remark~\ref{dgl factorization}.  Since the functors $[-]_\mDG$ and $[-]_\mDGL$ are both right adjoints, they commute with limits.  
In particular, if $\mathscr{D}:\catI\to\dgl_r$ is a diagram in $\dgl_r$ then we have
$$
\bigl[[\lim{\!}_\dgl\mathscr{D}]_\mDG\bigr]_\mDGL =
 \bigl[\lim{\!}_\dg[\mathscr{D}]_\mDG\bigr]_\mDGL 
$$

According to Remark~\ref{dgl factorization}, this means that in order to compute a limit in $\dgl$, we may instead compute the corresponding limit in $\dg$ and then figure out the correct Lie bracket structure on the result.  We leave it to the interested reader to prove that in general there is only one Lie bracket structure possible such that there are natural maps $\lim{\!}_\dgl\mathscr{D} \longrightarrow \mathscr{D}$ which descend to the existing natural maps $[\lim{\!}_\dgl \mathscr{D}]_\mDG = \lim{\!}_\dg[\mathscr{D}]_\mDG \longrightarrow [\mathscr{D}]_\mDG$.  This is what is meant when people say that ``limits in $\dgl$ are {\em created} in $\dg$.''\footnote{See [MacL, p122] for a precise definition.}

Similarly, since the functors $(-)^\ab$ and $\freeL_{(-)}$ are both left adjoints they commute with colimits.
For $\mathscr{D}':\catJ\to\hfdgl_r$ a diagram in $\hfdgl_r$ we have
$$\freeL_{(\colim{\!}_\hfdgl \mathscr{D}')^\ab} =
 \freeL_{\colim{\!}_\dg (\mathscr{D}')^\ab} $$
[It is an easy exercise that $\colim_\hfdgl\mathscr{D}' = \colim_\dgl\mathscr{D}'$.]  Now \ref{dgl factorization} suggests that colimits in $\fdgl$ may also be created in $\dg$.  This is in fact the case.  Again we leave it to the interested reader that there is a unique derivation which recovers the higher order differential information on $\colim{\!}_\hfdgl\mathscr{D}'$ so that there are natural maps $\mathscr{D}' \longrightarrow \colim{\!}_\hfdgl\mathscr{D}'$ which descend to the existing natural maps $(\mathscr{D}')^\ab \longrightarrow (\colim{\!}_\hfdgl\mathscr{D}')^\ab = \colim{\!}_\dg(\mathscr{D}')^\ab$. 
\end{note}

\subsection{Cones, Suspensions, Paths, and Loops}

\begin{defn}[Unreduced Paths]\label{d:path dgl} 
Given $L$ an $r$-reduced \DGL\ define the {\em unreduced paths on $L$} to be the $(r-1)$-reduced \DGL\ given by:
$${\tilde{\mathrm{p}}}L := (L_\bullet\oplus s^{-1}L_\bullet,\, d_{\tilde{\mathrm{p}}},\, [-,-]_{\tilde{\mathrm{p}}})$$
\begin{itemize}
\item $d_{\tilde{\mathrm{p}}}(l_1 + s^{-1}l_2) = d_Ll_1 - s^{-1}(d_Ll_2 + l_1)$
\item $[-,-]_{\tilde{\mathrm{p}}}$ is defined by:
\begin{itemize}
 \item[\textbullet] $[l_1,l_2]_{\tilde{\mathrm{p}}} = [l_1,l_2]_L$ 
 \item[\textbullet] $[l_1,s^{-1}l_2]_{\tilde{\mathrm{p}}} = \frac{1}{2}s^{-1}[l_1,l_2]_L$ (anti-commutativity forces 
 $[s^{-1}l_1,l_2]_{\tilde{\mathrm{p}}} = (-1)^{|l_2|}\frac{1}{2}s^{-1}[l_1,l_2]_L$)
 \item[\textbullet] $[s^{-1}l_1,s^{-1}l_2]_{\tilde{\mathrm{p}}} = 0$
\end{itemize}
\end{itemize} 
\end{defn}

Note that this is just the paths on the underlying \DG\ of $L$ --- $\mathrm{p}[L]_\mDG$ --- equipped with a Lie bracket which is compatible with the differential of $\mathrm{p}[L]_\mDG$.
Again, just like paths on a topological space $X$, the unreduced paths on $L$ is a contractible \DGL\ which comes equipped with a (degree-wise) surjection ${\tilde{\mathrm{p}}}L \to L$.  

\begin{defn}[Unreduced Loops] 
Given $L$ an $r$-reduced \DGL\ define the {\em unreduced loops on $L$} to be the $(r-1)$-reduced \DGL\ given by $s^{-1} L := (s^{-1}[L]_\mDG,\, [-,-]=0) = \big[s^{-1}[L]_\mDG\big]_\mDGL$.
\end{defn}

At times we wish to use a slightly larger model for unreduced loops:
$$\hat s^{-1} L := \big(s^{-1}L_\bullet \oplus L_\bullet\oplus s^{-1}L_\bullet,\, d_{\hat s},\,[-,-]_{\hat s}\big)$$
\begin{itemize}
\item $d_{\hat s}(s^{-1}l_1,\,k,\,s^{-1}l_2) =\big( s^{-1}( - d_Ll_1 + k),\ d_Lk,\ s^{-1}( -d_Ll_2+k)\big)$ 
\item $[-,-]_{\hat s}$ is 0 on elements both from one of the two copies of $s^{-1}L$, the inherited bracket on $L\otimes L$, and on mixed brackets is given by 
\begin{itemize}
\item[\textbullet] $[k,\, s^{-1}l_i]_{\hat s} = \frac{1}{2} s^{-1}[k,\,l_i]_L$ (anti-commutativity forces $[s^{-1}l_i,\, k]_{\hat s} = (-1)^{|k|} \frac{1}{2} s^{-1}[l_i,\, k]_L$)
\item[\textbullet] $[s^{-1}l_1, s^{-1}l_2]_{\hat s} = 0$
\end{itemize}
where $(s^{-1}l_1,\, k,\, s^{-1}l_2) \in s^{-1}L_\bullet \times L_\bullet \times s^{-1}L_\bullet$
\end{itemize}

Again, $\hat s^{-1} L$ is just ${\hat \Omega}[L]_\mDG$ equipped with a Lie bracket which is compatible with the differential of ${\hat\Omega}[L]_\mDG$.
Also, there is a surjection ${\hat s^{-1}} L \to L$ which plays the role of the midpoint evaluation map in $\top$ sending $\Omega X \to X$.  Furthermore, two injection maps from $s^{-1}L$ each give quasi-isomorphisms $s^{-1}L \to \hat s^{-1}L$.  

\begin{note}
The \DGL\,s $s^{-1} L$ and ${\hat s^{-1}}L$ are both pullbacks:
\begin{itemize}
\item $s^{-1} L$ is the pullback of $0_\mDGL \xleftarrow{} L \xrightarrow{} {\tilde{\mathrm{p}}}L$
\item $\hat s^{-1} L$ is the pullback of ${\tilde{\mathrm{p}}}L \xleftarrow{} L \xrightarrow{} {\tilde{\mathrm{p}}}L$
\end{itemize}
\end{note}

In general, if $L \in \dgl_r$ then $\tilde{\mathrm{p}}L,  s^{-1} L, {\hat s^{-1}} L \in \dgl_{r-1}$.  In fact, from their definitions it is clear that these define functors $\dgl_r \to \dgl_{r-1}$.  In general we wish to instead have functors $\dgl_r \to \dgl_r$.  [Note in particular that if $L$ is 0-reduced then the above objects are all $(-1)$-reduced -- which is unfortunate, since for much of the section so far we have been assuming that all \DGL\,s are trivial in negative gradings; but not terrible since the structures above still make sense for negatively graded \DGL\,s.]  To remedy this, we compose with the $r$-reduction functor:

\begin{defn}[Paths and Loops]
Define the paths and loops functors $\dgl_r \to \dgl_r$ by
\begin{itemize}
\item $p:\dgl_r \to \dgl_r$ by $\mathrm{p}L = \red{\tilde{\mathrm{p}} L}$.
\item ${\Omega}:\dgl_r \to \dgl_r$ by ${\Omega}L = \red{s^{-1} L} = \big[s^{-1}[L]_{\mDG}\big]_{\mDGL_r}$.
\end{itemize}
\end{defn}

From \ref{H_* red} we have that $\mathrm{p}L$ is still a contractible \DGL.  Also there is a map $\mathrm{p}L\to L$ which, although it is no longer necessarily a surjection in degree $r$, is still a fibration (see the following subsection). 
Since reduction is a right adjoint, $\Omega L$ is the pullback of $0_\mDGL \xrightarrow{} L \xleftarrow{} \mathrm{p}L$.  We may similarly define $\hat \Omega L = \red{\hat s^{-1} L}$ and we get that $\hat \Omega L$ is the pullback of $\mathrm{p}L \xrightarrow{} L \xleftarrow{} \mathrm{p}L$.

Given a free \DGL\ of the form $\freeL_V = (\freeL_V, d)$ we we write $\freeL_{sV}$ for the truly free \DGL\ given by $\freeL_{sV} = \freeL_{s(V,\,d_V)}$ (recall that $d_V$ is the restriction of $d$ to $V$ the generating \G).  This \DGL\ plays the same role that $sV$ played for \DG\,s in defining cones and suspensions. 

\begin{defn}[Cone of $\freeL_V$]\label{cone free dgl} 
Given $\freeL_V = (\freeL_V,\, d)$  define the {\em cone of $\freeL_V$} to be $\mathrm{c}\freeL_V$ the free \DGL\ given by taking the free product of \DGL\,s $\freeL_V\liesum \freeL_{sV}$ and then modifying the differential:
$$\mathrm{c}\freeL_V := ([\freeL_V \liesum \freeL_{sV}]_{\mathrm GL},\, d_\mathrm{c} = d_{\liesum} + \freeL_{d'})$$
where $d_{\liesum}$ is the differential on $\freeL_V \liesum \freeL_{sV}$ and $\freeL_{d'}$ is the derivation freely generated by the differential on $(V\oplus sV)$ given by $d'(sv) = v$ for $v\in {V}$ (recall that $[\freeL_V \liesum \freeL_{sV}]_{\mathrm GL} = \freeL_{V\oplus sV}$).
\end{defn}

By Corollary~\ref{fdgl sum}, $\mathrm{c}\freeL_V$ is merely $\freeL_{V\oplus sV}$ with a modified differential.  Also, by the rational Hurewicz theorem, $\mathrm{c}\freeL_V$ is contractible.  And of course, there is a \DGL-map $\freeL_V \to \mathrm{c}\freeL_V$ which is as we desire.  Note that the map $\freeL_V \to \mathrm{c}\freeL_V$ is in fact a free map of free \DGL\,s.

\begin{defn}[Suspension of $\freeL_V$]\label{sigma free dgl} 
Given $\freeL_V = (\freeL_V,\, d)$ let the {\em suspension of $\freeL_V$} be the truly free \DGL\ $\Sigma\freeL_V := \freeL_{sV}$ defined above.
\end{defn}

At times we desire a slightly different model.  Define ${\hat \Sigma}\freeL_V$ to be the free product of \DGL\,s $\freeL_{sV}\liesum \freeL_V\liesum \freeL_{sV}$ with a modified differential:
$$\hat\Sigma \freeL_V := ([\freeL_{sV}\liesum \freeL_V\liesum \freeL_{sV}]_{\rm GL},\, d_{\hat \Sigma} = d_{\liesum} + \freeL_{d'})$$
where $\freeL_{d'}$ is the derivation freely generated by the differential $d'(sv_1\oplus v_2 \oplus sv_3) = 0\oplus (v_1 + v_3)\oplus 0$ on $sV\oplus V \oplus sV$ (recall that $[\freeL_{sV}\liesum\freeL_V\liesum\freeL_{sV}]_\mGL = \freeL_{sV\oplus V\oplus sV}$).

By Corollary~\ref{fdgl sum}, $\hat\Sigma \freeL_V$ is merely $\freeL_{sV\oplus V\oplus sV}$ with a modified differential. 

\begin{note} 
By construction, we have that 
\begin{itemize}
\item $\Sigma \freeL_V$ is the pushout of $0_\mDGL \xleftarrow{} \freeL_V \xrightarrow{} \mathrm{c}\freeL_V$.
\item $\hat \Sigma \freeL_V$ is the pushout of $\mathrm{c}\freeL_V \xleftarrow{} \freeL_V \xrightarrow{} \mathrm{c}\freeL_V$.
\end{itemize}
\end{note}

More generally, if $L$ is any free \DGL, then we may define the suspension of $L$ to be $\Sigma L := \freeL_{s(L)^\mathrm{ab}}$. 
Note that this construction extends naturally to all of $\dgl_r$ and furthermore that: 

\begin{lemma}[Suspension]
A left adjoint of the functor $\Omega:\dgl_r \to \dgl_{r}$ is given by the functor $\Sigma:L\mapsto \freeL_{s(L)^{\rm ab}}$.
\end{lemma} 

\begin{proof}
By definition $\Omega$ and $\Sigma$ are the following compositions of adjoints:
$$\xymatrix{
\Sigma:\dgl_r \ar@<2pt>[r]^(.51){\mathrm{incl}} &
  \dgl_{r-1} \ar@<2pt>[r]^(.55){(-)^{\mathrm{ab}}} \ar@<2pt>[l]^(.46){\mathrm{red}_r} &
  \dg_{r-1} \ar@<2pt>[r]^(.52){s} \ar@<2pt>[l]^(.43){[-]_\mDGL} &
  \dg_r \ar@<2pt>[r]^(.44){\freeL_{(-)}} \ar@<2pt>[l]^(.45){s^{-1}} &
  \dgl_r: \Omega \ar@<2pt>[l]^(.52){[-]_\mDG}
}$$
\end{proof}

\section{Model Category Structure}

Recall that quasi-isomorphisms of \DGL\,s are \DGL-maps which are quasi-isomorphisms on the underlying \DG\,s.  The usual model category structure on $\dgl_r$ (from e.g. [Q69 \S5]) takes quasi-isomorphisms to be weak equivalences, degree-wise surjections in degree $> r$ to be fibrations, and allows cofibrations to be determined by left lifting with respect to acyclic fibrations.  Under this model category structure, all objects are fibrant.

\begin{thm}
This gives a model category structure on $\dgl_r$ ($r\ge1$).\footnote{It is possible to extend this model category structure to define a viable model category structure on $\dgl_0$; however we do not need this, so we stick with the easily defined and proven case.}
\end{thm}

\begin{proof}
Quillen shows in [Q69, B.5.1], that this structure satisfies all of the modern model category axioms except for M1 (small limits and colimits) and M5 (functorial factorizations).  We have already constructed small limits and colimits in $\dgl_r$ in Lemma~\ref{dgl (co)lim}.  So it remains only to show that $\dgl_r$ has functorial factorizations.  However, this implicitly follows from the proof of factorizations given by Quillen, since he uses the small object argument.
\end{proof}

The following proposition is also proven by Quillen:

\begin{prop}[\lbrack Q69, 5.5\rbrack]
Cofibrations in $\dgl$ are the free maps.
\end{prop}

\begin{cor}
The cofibrant objects in $\dgl$ are precisely the free objects $\Fdgl$.
\end{cor}

The identity functor serves as a fibrant replacement functor.  $\dgl$ also has a (fibrant) cofibrant replacement functor (given by $\L\C$) which is described in Chapter~\ref{S:Q homotopy} and in particular by Corollary~\ref{funct approx}.  This cofibrant replacement functor is particularly nice.  The following lemma immediately follows from the definitons of $\L$ and $\C$ which is given in \ref{L functor} and \ref{C functor}:

\begin{lemma}\label{LC is free}
The cofibrant replacement functor in $\dgl_r$ is a functor $\L\C:\dgl_r\to\hfdgl_r$.

That is, cofibrant replacement takes \DGL\,s to free \DGL\,s of the form $(\freeL_V,\, d)$ and \DGL-maps to freely generated maps of the form $\freeL_{\hat f}:(\freeL_V,\,d) \to (\freeL_W,\, d')$.
\end{lemma}

Note that under the above model category structure, the adjoints which we gave in Lemmas~\ref{dgl red} and \ref{dgl adjoints} become Quillen adjoint pairs (\ref{Quillen pair}): 

\begin{lemma}
The adjoint pair given in Lemma~\ref{dgl red}
$$\xymatrix@C=20pt@R=5pt{
[-]_{\mDGL_t}:\dgl_r \ar@<2pt>[r] & \dgl_t:\mathrm{red}_r(-) \ar@<2pt>[l] 
}$$
($t<r$) is a Quillen adjoint pair.
Also both of the adjoint pairs from Lemma~\ref{dgl adjoints}
$$\xymatrix@C=10pt@R=1pt{
\freeL_{(-)}:\dg_r \ar@<2pt>[r] &
  \dgl_r:[-]_{\mDG_r} \phantom{xxx} \ar@<2pt>[l]  \\
{}\kern 2pt(-)^\ab:\dgl_r \ar@<2pt>[r] &
  \dg:[-]_{\mDGL_r} \ar@<2pt>[l]
}$$      
are Quillen adjoint pairs.
\end{lemma}

\begin{proof}
The right adjoints $\mathrm{red}_r$(-), $[-]_{\mDG_r}$, and $[-]_{\mDGL_r}$ all preserve fibrations, since these are merely degree-wise surjections in $\dg$ and degree-wise surjections except for in degree $r$ in $\dg_r$ and $\dgl_r$.  
The functors $[-]_{\mDGL_t}$ and $\mathrm{red}_r$ clearly both preserve all weak equivalences.
And, since weak equivalences are defined on the level of underlying \DG\,s, they are also preserved by $[-]_\mDG$ and $[-]_\mDGL$. 
Since the right adjoints in the above pairs each preserve all fibrations and trivial fibrations, the left adjoints must also preserve all cofibrations and trivial cofibrations and the pairs are Quillen adjoints as claimed.
\end{proof}  

Note that the top adjoint from Lemma~\ref{dgl adjoints} does not extend as a Quillen adjoint all the way to $\dg$.  This is because the adjoint pair $\mathrm{trunc}_r:\dg\rightleftarrows\dg_r:\mathrm{incl}$ is not a Quillen adjoint pair.  Lemma~\ref{adjoint holim} now provides us with the following corollary:

\begin{cor}\label{C:dgl holim plan}The right adjoint functors $\mathrm{red}_r(-):\dgl_t\to\dgl_r$, $[-]_{\mDG_r}:\dgl_r\to\dg_r$, and $[-]_\mDGL:\dg\to\dgl_r$ preserve all weak equivalences and therefore all homotopy limits.

The left adjoints $[-]_{\mDGL_t}:\dgl_r\to\dgl_t$, $\freeL_{(-)}:\dg_r\to\Fdgl_r$, and $(-)^\ab:\Fdgl_r\to\dg$ induced by the above preserve all weak equivalences and therefore all homotopy colimits.  
\end{cor}
\begin{proof}
It remains to show only that the left adjoints $\freeL_{(-)}$ and $(-)^\ab$ preserve weak equivalences; however this follows from the rational Hurewicz Theorem (\ref{Hurewicz dgl}).
\end{proof}

\subsection{Homotopy Limits and Colimits}\label{S:dgl holim}

By [DHKS] all homotopy limits and colimits in $\dgl_r$ exist (see \ref{t:hocomplete}).  We use \ref{C:dgl holim plan} and \ref{dgl factorization} in order to construct nice models for homotopy limits and colimits in $\dgl_r$ in much the same way as we commented in \ref{dgl lim plan} that we could have created limits and colimits in $\dgl_r$.    Essentially, our construction is as follows:  If $\mathscr{D}:\catI\to \dgl_r$ is a diagram in $\dgl_r$ then up to bracket, the homotopy limit of $\mathscr{D}$ in $\dgl_r$ is given by  
$$\bigl[[\holim{\!}_\dgl\mathscr{D}]_\mDG\bigr]_\mDGL = \bigl[\holim{\!}_\dg[\mathscr{D}]_\mDG\bigr]_\mDGL$$
Similarly (if $\mathscr{D}$ is a diagram in $\hfdgl$) up to higher order differential information, the homotopy colimit of $\mathscr{D}$ is given by 
$$\freeL_{\left(\hocolim_\dgl\mathscr{D}\right)^\mathrm{ab}} = \freeL_{\hocolim_\dg(\mathscr{D})^\mathrm{ab}}$$
To construct models for homotopy limits or colimits in $\dgl_r$ we insert our models for homotopy limits or colimits in $\dg_r$ (given in Section~\ref{dg_r model category}) into the above and then supply the missing Lie bracket or higher order differential information. 

More precisely, note that if $\holim_{\dgl}^\catI:(\dgl_r)^\catI\longrightarrow\dgl_r$ is any $\catI$-homotopy limit functor on $\dgl_r$ then it satisfies the following properties:
\begin{itemize}
\item $\holim_{\dgl}^\catI$ is a homotopy functor (i.e. it preserves weak equivalences\footnote{Recall that weak equivalences of diagrams are natural transformations which give weak equivalences objectwise.}).
\item $\holim_\dg^\catI(-) = \bigl[\holim_{\dgl}^\catI[-]_\mDGL\bigr]_\mDG$ is an $\catI$-homotopy limit functor on $\dg$ such that 
$$\bigl[\holim{\!}_\dgl^\catI\mathscr{D}\bigr]_\mDG = \holim{\!}_\dg^\catI[\mathscr{D}]_\mDG$$ 
for all $\catI$-diagrams $\mathscr{D}:\catI\to\dgl_r$.
\item There are natural maps $\lim_{\dgl}\mathscr{D} \to \holim_{\dgl}^\catI\mathscr{D}$.
\item These maps descend to the natural maps $\lim_\dg[\mathscr{D}]_\mDG \to \holim_\dg^\catI[\mathscr{D}]_\mDG$.
\end{itemize}
According to the following theorem, these properties are enough to characterize $\catI$-homotopy limit functors on $\dgl_r$:

\begin{thm}[Creation of Homotopy Limits and Colimits]\label{dgl holim plan}
Let $\holim_\dg$ and $\hocolim_\dg$ be any $\catI$-homotopy limit and colimit functors on $\dg_r$.
\begin{enumerate}
\item\label{dgl holim plan 1} Suppose $F:(\dgl_r)^{\catI} \longrightarrow \dgl_r$ is a functor from $\catI$-diagrams in $\dgl_r$ to $\dgl_r$ such that for all $\catI$-diagrams $\mathscr{D}:\catI\to\dgl_r$, we have 
\begin{itemize}
\item $\bigl[F(\mathscr{D})\bigr]_\mDG = \holim_\dg[\mathscr{D}]_\mDG$.\footnote{Note that this is stronger than the statement $\bigl[F([\mathscr{D}]_\mDGL)\bigr]_\mDG = \holim_\dg(\mathscr{D})$ for all $\mathscr{D}:\catI\to\dg$.}
\item $F$ is equipped with natural maps $e_F:\lim_{\dgl}\mathscr{D} \to F(\mathscr{D})$. 
\item $[e_F]_\mDG$ is the canonical map 
 $\lim_\dg[\mathscr{D}]_\mDG \to \holim_\dg[\mathscr{D}]_\mDG$.  
\end{itemize}
Then $F$ is an $\catI$-homotopy limit functor on $\dgl_r$.

\item\label{dgl holim plan 2} Dually, suppose $G:(\dgl_r)^\catI\longrightarrow\Fdgl_r$ is a homotopy functor such that for all $\catI$-diagrams $\mathscr{D}:\catI\to\dgl_r$ we have 
\begin{itemize}
\item $\bigl(G(\mathscr{D})\bigr)^\ab = \hocolim_\dg(\mathscr{D})^\ab$.
\item $G$ is equipped with natural maps $e_G:G(\mathscr{D}) \to \colim_{\dgl}(\mathscr{D})$.
\item $(e_G)^\ab$ is the canonical map $\hocolim_\dg(\mathscr{D})^\ab \longrightarrow \colim_\dg(\mathscr{D})^\ab$.
\end{itemize}
Then $G$ is an $\catI$-homotopy colimit functor on $\dgl_r$.
\end{enumerate}
\end{thm}

Recall that [DHKS] construct $\catI$-homotopy limit and colimit functors on the model category $\catM$ using virtually-fibrant and virtually-cofibrant diagrams in $\catM$ -- setting $\holim_\catM\mathscr{D} = \lim_\catM\mathscr{D}_\mathsf{vf}$ and $\hocolim_\catM\mathscr{D} = \colim_\catM\mathscr{D}_\mathsf{vc}$.  We will show that there is a zig-zag of natural weak equivalences between $F$ and the functor $\holim_\dgl$ constructed in this way as well as between $G$ and the functor $\hocolim_\dgl$ constructed in this way.

The proof of \ref{dgl holim plan} relies on the following two technical lemmas:

\begin{lemma}\label{L1:dgl holim plan}
The functors $[-]_\mDG:(\dgl_r)^\catI \to (\dg_r)^\catI$ and $(-)^\ab:(\dgl_r)^\catI \to (\dg_r)^\catI$ preserve virtually-fibrant diagrams and virtually-cofibrant diagrams respectively.
\end{lemma}

\begin{lemma}\label{L2:dgl holim plan}
Let $\mathscr{D}:\catI\to\dgl_r$ be an $\catI$-diagram in $\dgl_r$.  Then $\colim_\dgl(\mathscr{D}_\mathsf{vc})$ is cofibrant (i.e. $\colim_\dgl(\mathscr{D}_\mathsf{vc})$ is a free \DGL).
\end{lemma}

These lemmas follow directly from the properties of virtually-fibrant and virtually-cofibrant diagrams (see [DHKS 20.5] or else our discussion of the proof of Theorem~\ref{t:hocomplete}).

\begin{proof}[Proof of \ref{dgl holim plan}]
We begin with (\ref{dgl holim plan 1}):

First we show that $F$ is a homotopy functor.  Suppose $i:\mathscr{D}_1\xrightarrow{\,\weq\,}\mathscr{D}_2$ is a weak equivalence of diagrams.  Then $[i]_\mDG:[\mathscr{D}_1]_\mDG \xrightarrow{\,\weq\,}[\mathscr{D}_2]_\mDG$ which induces a weak equivalence
$$\bigl[F(\mathscr{D}_1)\bigr]_\mDG = \holim{\!}_\dg [\mathscr{D}_1]_\mDG \xrightarrow[\,\holim\lbrack i\rbrack_\mDG\,]{\,\weq\,}
 \holim{\!}_\dg[\mathscr{D}_2]_\mDG = \bigl[F(\mathscr{D}_2)\bigr]_\mDG$$
However, this map is $[F(i)]_\mDG$.  So $F(i)$ is a weak equivalence by \ref{weq dgl equals weq dg}. 

Since $F$ is a homotopy functor there is a natural weak equivalence $F(\mathscr{D}) \xrightarrow{\,\weq\,} F(\mathscr{D}_\mathsf{vf})$ induced by the virtually-fibrant replacement functor in $(\dgl_r)^\catI$.  
Consider the map $e_F:\lim_\dgl(\mathscr{D}_\mathsf{vf}) \to F(\mathscr{D}_\mathsf{vf})$ descending to $[e_F]_\mDG:\lim_\dg\bigl([\mathscr{D}_\mathsf{vf}]_\mDG\bigr) \to \holim_\dg[\mathscr{D}_\mathsf{vf}]_\mDG$.  
Since $[-]_\mDG:(\dgl_r)^\catI \to (\dg_r)^\catI$ preserves fibrations,  the diagram $[\mathscr{D}_\mathsf{vf}]_\mDG$ is virtually-fibrant in $(\dg_r)^\catI$.  
Thus $\lim_\dg\bigl([\mathscr{D}_\mathsf{vf}]_\mDG\bigr)$ is a homotopy limit functor on $\dg$ which means $[e_F]_\mDG$ is a weak equivalence.  Therefore $e_F$ is a weak equivalence by~\ref{weq dgl equals weq dg}.  
We have now completed a zig-zag of natural weak equivalences 
$$F(\mathscr{D}) \xrightarrow{\ \weq\ } F(\mathscr{D}_\mathsf{vf}) \xleftarrow{\ \weq\ } 
  \lim{\!}_\dgl\mathscr{D}_\mathsf{vf} = \holim{\!}_\dgl\mathscr{D}$$

Part (\ref{dgl holim plan 2}) is proven similarly.  
Virtual-cofibrant replacement in $(\dgl_r)^\catI$ yields a natural weak equivalence $G(\mathscr{D}_\mathsf{vc}) \xrightarrow{\,\weq\,} G(\mathscr{D})$.  
The map $(e_G)^\ab:\hocolim_\dg(\mathscr{D}_\mathsf{vc})^\ab \longrightarrow \colim_\dg(\mathscr{D}_\mathsf{vc})^\ab$ is a weak equivalence because $(\mathscr{D}_\mathsf{vc})^\ab$ is a virtually-cofibrant diagram in $\dg$ which means $\colim_\dg(\mathscr{D}_\mathsf{vc})^\ab$ is a homotopy colimit functor on $\dg$.
Since $\colim_\dgl(\mathscr{D}_\mathsf{vc})$ and $G(\mathscr{D}_\mathsf{vc})$ are both free, the rational Hurewicz Theorem (\ref{Hurewicz dgl}) tells us that $(e_G)^\ab$ is a weak equivalence if and only if $e_G$ is a weak equivalence.  Thus we have a zig-zag of natural weak equivalences
$$G(\mathscr{D}) \xleftarrow{\ \weq\ } G(\mathscr{D}_\mathsf{vc}) \xrightarrow{\ \weq\ }
  \colim{\!}_\dgl\mathscr{D}_\mathsf{vc} = \hocolim{\!}_\dgl\mathscr{D}$$
\end{proof}

\subsubsection{Homotopy Limits in $\dgl_r$}

Given $\mathscr{D}$ the \DGL-diagram $L_1\xrightarrow{f_1} K \xleftarrow{f_2} L_2$ with each $L_i$ $r$-reduced, let $\tilde {\mathscr{P}}_\mathscr{D}$ denote the $(r-1)$-reduced \DGL\ given by:
$$\tilde {\mathscr{P}}_\mathscr{D} := \big((L_1 \times s^{-1}K \times L_2)_\bullet,\  d_{\tilde{\mathscr{P}}},\  [-,-]_{\tilde{\mathscr{P}}})$$ 
\begin{itemize}
\item $d_{\tilde{\mathscr{P}}}(l_1,\  s^{-1}k,\  l_2) = \Bigl(d_{L_1}l_1,\  s^{-1}\bigl(f_1(l_1) - f_2(l_2) - d_Kk\bigr),\  d_{L_2}l_2)\Bigr)$
\item $[-,-]_{\tilde{\mathscr{P}}}$ is 0 on $s^{-1}K\otimes s^{-1}K$, the inherited bracket ($[-,-]_{L_i}$) on $L_i\otimes L_i$,  and on mixed brackets is given by: 
\begin{itemize}
\item[\textbullet] $[s^{-1}k,\, l_i]_{\tilde{\mathscr{P}}} = \frac{1}{2}\, s^{-1}\bigl[k,\,f_i(l_i)\bigr]_K$ (forcing $[l_i,\, s^{-1}k]_{\tilde{\mathscr{P}}} = (-1)^{|l_i|}\, \frac{1}{2}\, s^{-1}\bigl[f_i(l_i),\, k\bigr]_K$) 
\item[\textbullet] $[l_1, l_2]_{\tilde{\mathscr{P}}} = 0$
\end{itemize}
where $l_i \in L_i$ and $s^{-1}k \in s^{-1}K$.
\end{itemize}
Note that $\tilde {\mathscr{P}}_\mathscr{D}$ is just ${\mathscr{P}}_{[\mathscr{D}]_\mDG}$ equipped with a compatible Lie bracket.

\begin{proof}[Proof that this is a \DGL]
Since $\bigl[\tilde {\mathscr{P}}_\mathscr{D}\bigr]_\mDG = {\mathscr{P}}_{[\mathscr{D}]_\mDG}$ the fact that $d_{\tilde{\mathscr{P}}}\circ d_{\tilde{\mathscr{P}}} = 0$ follows from the corresponding fact about the \DG\ ${\mathscr{P}}_{[\mathscr{D}]_\mDG}$ proven in \ref{l:dg pushout}.
Also, we have defined the bracket of $\tilde{\mathscr{P}}_\mathscr{D}$ specifically so that it is anti-commutative.  Thus it remains to show only that the bracket and differential are compatible.  This is trivial for brackets of two elements of $s^{-1}K$.  
For brackets $[l_1,l_2]_{\tilde{\mathscr{P}}}$ of elements $l_i \in L_i$ the following shows compatibility:
\begin{align*}
\bigl[d_{\tilde{\mathscr{P}}}l_1,l_2\bigr]_{\tilde{\mathscr{P}}} 
  + (-1)^{|l_1|}\bigl[l_1, d_{\tilde{\mathscr{P}}}l_2\bigr]_{\tilde{\mathscr{P}}}
 &= \bigl[d_{L_1}l_1, l_2\bigr]_{\tilde{\mathscr{P}}} + \bigl[s^{-1}f_1(l_1), l_2\bigr]_{\tilde{\mathscr{P}}} \\ 
 &\qquad + (-1)^{|l_1|} \bigl[l_1, d_{L_2}l_2\bigr]_{\tilde{\mathscr{P}}}  - 
   (-1)^{|l_1|} \bigl[l_1, s^{-1}f_2(l_2)\bigr]_{\tilde{\mathscr{P}}} \\
 &= {\textstyle \frac{1}{2}}\,s^{-1}\bigl[f_1(l_1), f_2(l_2)\bigr]_K - (-1)^{|l_1| + |l_1|}\,{\textstyle \frac{1}{2}}\,s^{-1}[f_1(l_1), f_2(l_2)]_K \\
 &= 0 \\ &= d_{\tilde{\mathscr{P}}} \bigl[l_1, l_2\bigr]_{\tilde{\mathscr{P}}}
\end{align*}
For brackets $[v, w]_{\tilde{\mathscr{P}}}$ of elements $v,w \in L_i$ the following shows compatibility: 
\begin{align*}
d_{\tilde{\mathscr{P}}}[v,\, w]_{\tilde{\mathscr{P}}} &= d_{\tilde{\mathscr{P}}}[v,\, w]_{L_i} \\
&= d_{L_i} [v,\,w]_{L_i} - (-1)^i s^{-1}f_i\bigl([v,\, w]_{L_i}\bigr)  \\
&= \bigl[d_{L_i} v,\, w\bigr]_{L_i} + (-1)^{|v|}\bigl[v,\, d_{L_i} w\bigr]_{L_i} - 
   (-1)^i s^{-1}\bigl[f_i(v),\, f_i(w)\bigr]_K \\
&= \bigl[d_{L_i}v,\, w\bigr]_{L_i} - (-1)^i\,{\textstyle \frac{1}{2}}\,s^{-1}\bigl[f_i(v),\, f_i(w)\bigr]_K  \\
&\qquad + (-1)^{|v|}\bigl[v,\, d_{L_i}w\bigr]_{L_i} - 
   (-1)^i\,{\textstyle \frac{1}{2}}\,s^{-1}\bigl[f_i(v),\, f_i(w)\bigr]_K  \\
&= \bigl[d_{L_i}v,\, w\bigr]_{\tilde{\mathscr{P}}} - \bigl[(-1)^is^{-1}f_i(v),\, w\bigr]_{\tilde{\mathscr{P}}} \\
&\qquad + (-1)^{|v|}\bigl[v,\, d_{\tilde{\mathscr{P}}}w\bigr]_{\tilde{\mathscr{P}}}
   - (-1)^{|v|}\bigl[v,\, (-1)^is^{-1}f_i(w)\bigr]_{\tilde{\mathscr{P}}} \\
&= \bigl[d_{\tilde{\mathscr{P}}}v,\, w\bigr]_{\tilde{\mathscr{P}}} 
   + (-1)^{|v|}\bigl[v,\, d_{\tilde{\mathscr{P}}}w\bigr]_{\tilde{\mathscr{P}}}
\end{align*}
For brackets $[s^{-1}k,\, l_i]_{\tilde{\mathscr{P}}}$ of elements $k\in K$ and $l_i \in L_i$ the following shows compatibility:
\begin{align*}
\bigl[d_{\tilde{\mathscr{P}}}s^{-1}k,\, l_i\bigr]_{\tilde{\mathscr{P}}} - (-1)^{|k|}\bigl[s^{-1}k,\, d_{\tilde{\mathscr{P}}}l_i\bigr]_{\tilde{\mathscr{P}}} 
 &= \bigl[-s^{-1}d_Kk,\, l_i\bigr]_{\tilde{\mathscr{P}}}  \\
 &\qquad - (-1)^{|k|}\bigl(\bigl[s^{-1}k,\, d_{L_i} l_i\bigr]_{\tilde{\mathscr{P}}}  - 
      \bigl[s^{-1}k,\, (-1)^{i}s^{-1}f_i(l_i)\bigr]_{\tilde{\mathscr{P}}}\bigr)  \\ 
 &= - {\textstyle \frac{1}{2}}\, s^{-1}\bigl[d_Kk,\, f_i(l_i)\bigr]_K - 
    (-1)^{|k|}{\textstyle \frac{1}{2}}\,s^{-1}\bigl[k,\, f_i(d_{L_i} l_i)\bigr]_K \\
 &= - {\textstyle \frac{1}{2}}\, \bigl(s^{-1}\bigl[d_Kk,\, f_i(l_i)\bigr]_K + 
    (-1)^{|k|}s^{-1}\bigl[k,\, d_Kf_i(l_i)\bigr]_K\bigr)  \\
 &= - {\textstyle \frac{1}{2}}\, s^{-1}d_K\bigl[k,\, f_i(l_i)\bigr]_K  
  = {\textstyle \frac{1}{2}}\, d_{\tilde{\mathscr{P}}} s^{-1}\bigl[k,\, f_i(l_i)\bigr]_K  \\
 &= d_{\tilde{\mathscr{P}}}[s^{-1}k,\, l_i]_{\tilde{\mathscr{P}}} 
\end{align*}
\end{proof}

This \DGL\ essentially gives us homotopy pullbacks:

\begin{lemma}[Homotopy Pullback]\label{L:dgl pullback} 
Given $\mathscr{D}$ the \DGL-diagram $L_1\xrightarrow{f_1} K \xleftarrow{f_2} L_2$ with each $L_i$ $r$-reduced, the homotopy pullback of $\mathscr{D}$ in $\dgl_{r-1}$ is $\tilde {\mathscr{P}}_\mathscr{D}$.
\end{lemma}

\begin{proof}
By \ref{dgl holim plan}, all that is required is to construct a map $\lim{\!}_{\dgl_{r-1}}(\mathscr{D}) \to \tilde{\mathscr{P}}_\mathscr{D}$.  Note that $\lim{\!}_{\dgl_{r-1}}(\mathscr{D}) = \lim{\!}_{\dgl_r}(\mathscr{D}) = L_1 \times_K L_2$.  The desired map is given by $(l_1,\,l_2) \mapsto (l_1,\,0,\,l_2)$ (by construction, this map respects the bracket and differential).
\end{proof}

\begin{cor} Given $\mathscr{D}$ the \DGL-diagram $L_1\xrightarrow{f_1} K \xleftarrow{f_2} L_2$ with each $L_i$ $r$-reduced, the homotopy pullback of $\mathscr{D}$ in $\dgl_{r}$ is ${\mathscr{P}}_\mathscr{D} := \red{\tilde {\mathscr{P}}_\mathscr{D}} = \mathrm{red}_r(\tilde{\mathscr{P}}_\mathscr{D})$.
\end{cor}

\begin{cor} $\Omega L$ is the homotopy pullback in $\dgl_r$ of the diagram $0_\mDGL \xrightarrow{\ } L \xleftarrow{\ } 0_\mDGL$.
\end{cor}

\begin{cor}[Homotopy Fiber]\label{dgl hofib}
The homotopy fiber of the map $f:L\to K$ in $\dgl_r$ is given by
$$\hofib(f) = \mathrm{red}_r\bigl((s^{-1}K\times L)_\bullet,\ d_{\hofib},\ [-,-]_{\hofib}\bigr)$$
where 
\begin{itemize}
\item $d_{\hofib}(s^{-1}k,\ l) = \bigl(s^{-1}(f(l) - d_Kk),\ d_Ll\bigr)$
\item $[s^{-1}k,\,l]_{\hofib} = \frac{1}{2}\,s^{-1}\bigl[k,\,f(l)\bigr]_K$ 
  (forcing $[l,\,s^{-1}k]_{\hofib} = (-1)^{|l|}\,\frac{1}{2}\,\bigl[f(l),\,k\bigr]_K$)
\item $[l_1,\,l_2]_{\hofib} = [l_1,\,l_2]_L$
\item $[s^{-1}k_1,\,s^{-1}k_2]_{\hofib} = 0$
\end{itemize}
\end{cor}

In general these constructions may be extended in the obvious manner in order to define higher dimensional homotopy pullbacks or indeed any homotopy limit in $\dgl_n$.  However, for our purposes all that we will explicitly require is the above homotopy pullback.

\subsubsection{Homotopy Colimits in $\dgl_r$}

Our strategy is to first construct homotopy colimits (in $\dgl_r$) of diagrams in $\hfdgl_r$ and then use these to get homotopy colimits of diagrams in $\dgl_r$.  We call diagrams in the image of the inclusion $(\hfdgl_r)^\catI \cofibr (\dgl_r)^\catI$ {\em freely generated diagrams}.  Just as a free $\dgl$ is given by the \DGL\ freely generated by a \DG\ along with some higher order differential information, a freely generated diagram is given by the diagram in $\dgl$ freely generated by a diagram in $\dg$ plus some higher order differential information.  The homotopy colimit of a freely generated diagram is the \DGL\ freely generated by the homotopy colimit of the associated diagram in $\dg$ plus some higher order differential information.
It follows from Lemma~\ref{LC is free} that our cofibrant replacement functor ($\L \C$) on $\dgl_r$ takes all diagrams to freely generated diagrams.  

Writing $\widehat{\hocolim}_\hfdgl:(\hfdgl_r)^\catI\to \Fdgl$ for our functor giving homotopy colimits of freely generated diagrams in $\dgl_r$, an $\catI$-homotopy colimit functor in $\dgl_r$ is then given by the composition
$$\hocolim{\!}_\dgl(-) := \widehat{\hocolim}_\hfdgl\bigl(\L \C (-)\bigr)$$
This is our model for general homotopy colimits in $\dgl_r$;
however, in practice almost all of the diagrams which we are interested in are already freely generated diagrams.  In these cases, we may use the simpler model $\widehat{\hocolim}_\hfdgl$ for their homotopy colimit, since for freely generated diagrams $\mathscr{D}$ there is a natural weak equivalence 
$$\widehat{\hocolim}_\hfdgl\mathscr{D} \xleftarrow{\ \weq\ } 
  \widehat{\hocolim}_\hfdgl\L \C (\mathscr{D}) = \hocolim{\!}_\dgl\mathscr{D}$$ 
induced by the natural weak equivalence $\mathscr{D} \xleftarrow{\,\weq\,} \L\C(\mathscr{D})$.

It follows from the definition of $\hfdgl$ that an $\catI$-diagram $\mathscr{D}:\catI\to\dgl_r$ in $\dgl_r$ is a freely generated diagram if, as a diagram of graded Lie algebras, it is the free extension of a diagram of graded vector spaces:
$$[\mathscr{D}]_\mGL = \freeL_{\hat{\mathscr{D}}}\qquad \text{where} \qquad
 \hat{\mathscr{D}}:\catI\to\g_r$$

\begin{ex}
Let $(\freeL_V,\,d)$ be a free \DGL.  The following diagrams are both freely generated:
\begin{enumerate}
\item $0_\mDGL\xleftarrow{\ \ }(\freeL_V,\,d)\xrightarrow{\ \ }0_\mDGL$
\item $0_\mDGL\xleftarrow{\ \ }(\freeL_V,\,d)\xrightarrow{\,\Idm\,}(\freeL_V,\,d)$
\end{enumerate}
We define homotopy colimits of freely generated diagrams so that the homotopy colimit of (1) is $\Sigma(\freeL_V,\,d)$ as defined in \ref{sigma free dgl} and the homotopy colimit of (2) is $\mathrm{c}(\freeL_V,\,d)$ as defined in \ref{cone free dgl}.
\end{ex}

We begin with simple diagrams:

Suppose $\mathscr{D}$ the freely generated pushout diagram $(\freeL_{U},\, d') \xleftarrow{\,\freeL_{f}\,} (\freeL_{V}, d'') \xrightarrow{\,\freeL_{g}\,} (\freeL_{W},\,d''')$ in $\dgl_r$.  This diagram has $[\mathscr{D}]_\mGL = \freeL_{\hat{\mathscr{D}}}$ where $\hat{\mathscr{D}}$ is the diagram $U\xleftarrow{f} V \xrightarrow{g} W$ in $\g_r$.  Note furthermore that $(\mathscr{D})^\ab$ is the diagram $(U,\,d'_U)\xleftarrow{\,f\,}(V,\,d''_V)\xrightarrow{\,g\,}(W,\,d'''_W)$ in $\dg_r$ and $\bigl[(\mathscr{D})^\ab\bigr]_\mG = \hat{\mathscr{D}}$.  

Let $\mathscr{C}_{\mathscr{D}}$ be the free cylinder \DGL\ given by:
$$\mathscr{C}_{\mathscr{D}} := (\freeL_{U}\liesum\freeL_{sV}\liesum\freeL_{W},\  
  d_\mathscr{C}=d_{\!\liesum} + \freeL_{d^{fg}})$$
where
\begin{itemize}
\item $d_{\!\liesum}$ is the differential on $(\freeL_{U},\, d') \liesum \freeL_{s(V,d''_V)} 
 \liesum (\freeL_{W},\,d''')$.
\item $\freeL_{d^{fg}}$ is the differential on the \GL\ $\freeL_{U}\liesum\freeL_{sV}\liesum\freeL_{W} = \freeL_{U\oplus sV \oplus W}$ freely generated by the differential on $U\oplus sV \oplus W$ given by $d^{fg}(sv)=f(v)+g(v)$ for $v\in V$.
\end{itemize}
Since the functor $\freeL_{(-)}$ is a section of $(-)^\ab$, it follows that 
$$(\mathscr{C}_\mathscr{D})^\ab 
 = \Bigl(\bigl[(U,d'_U)\oplus s(V,d''_V) \oplus (W,d'''_W)\bigr]_\mG,\ d=d_\oplus + d^{fg}\Bigr) 
 = \hocolim{\!}_\dg (\mathscr{D})^\ab$$
Furthermore there is a map $\mathscr{C}_\mathscr{D} \to \colim\mathscr{D}$ given by the composition
$$\mathscr{C}_\mathscr{D} \longrightarrow 
 (\freeL_{U},\, d') \liesum (\freeL_{W},\,d''') \longrightarrow \colim{\!}_\dgl\mathscr{D}$$
And this composition descends to
$$\hocolim{\!}_\dg(\mathscr{D})^\ab \longrightarrow
 (U,\,d'_U)\oplus (W,\,d'''_W) \longrightarrow \colim{\!}_\dg(\mathscr{D})^\ab$$
which is the natural map $\hocolim_\dg \to \colim_\dg$ described in \ref{dg ho(co)lim natural maps}.  
Also, $\mathscr{C}_{(-)}$ is a homotopy functor on diagrams in $\hfdgl$ because $i:\mathscr{D}_1\xrightarrow{\,\weq\,}\mathscr{D}_2$ a weak equivalence implies $\hocolim_\dg(i)^\ab = \bigl(\mathscr{C}_{(i)}\bigr)^\ab$ is a weak equivalence and so $\mathscr{C}_{(i)}$ is a weak  equivalence by the rational Hurewicz theorem.

Thus Theorem~\ref{dgl holim plan} gives us
\begin{lemma} 
The composition $\mathscr{C}_\mathscr{\L\C (-)}$ is a homotopy colimit functor.  

In particular, if $\mathscr{D}$ is a 
freely generated diagram in $\dgl_r$ then the \DGL\ $\mathscr{C}_\mathscr{D} \weq \mathscr{C}_\mathscr{\L\C \mathscr{D}}$ is a model for the homotopy colimit (in $\dgl_r$) of $\mathscr{D}$.
\end{lemma}

\begin{cor}
The following functors (defined earlier) are models for homotopy colimits
\begin{itemize}
\item $\Sigma (\freeL_V,\,d) = \hocolim_\dgl\Bigl(0_\mDGL\xleftarrow{\ \ }(\freeL_V,\,d)\xrightarrow{\ \ }0_\mDGL\Bigr)$.
\item $\mathrm{c}(\freeL_V,\,d) = \hocolim_\dgl\Bigl(0_\mDGL\xleftarrow{\ \ }(\freeL_V,\,d)\xrightarrow{\,\Idm\,}(\freeL_V,\,d)\Bigr)$.
\end{itemize}
\end{cor}

\begin{note}  
$\mathscr{C}_\mathscr{D}$ could also be written as 
$$\mathscr{C}_\mathscr{D} = \Bigl(\bigl[\freeL_{
  (U,d'_U)\oplus s(V,d''_V) \oplus (W,d'''_W)}\bigr]_\mGL,\ 
  d_\mathscr{C} = d_\freeL + d_{> 0}\Bigr)$$ 
\begin{itemize}
\item $d_\freeL$ is the differential on the truly free \DGL\ $\freeL_{(U,d'_U)\oplus s(V,d''_V) \oplus (W,d'''_W)}$.
\item $d_{>0}$ is the degree $-1$ map on $\freeL_{U\oplus sV\oplus W} = \freeL_U\liesum \freeL_{sV} \liesum \freeL_W$ generated by $d'-d'_U$ on $\freeL_U$ and $d'''-d'''_W$ on $\freeL_W$.
\end{itemize}
\end{note}

More generally,
given $\mathscr{D}:\catI\to\hfdgl_r$ a freely generated diagram, we could define the homotopy colimit of $\mathscr{D}$ along the lines of
$$\widehat{\hocolim}_\hfdgl\mathscr{D} := \Bigl(\bigl[\freeL_{\hocolim{\!}_\dg(\mathscr{D})^\ab}\bigr]_\mGL,\ d = d_\freeL + d_{> 0}\Bigr)$$
where $d_\freeL$ is the differential on the truly free \DGL\ $\freeL_{\hocolim{\!}_\dg(\mathscr{D})^\ab}$ and $d_{> 0}$ is the degree $-1$ map adding back in all of the higher order differential information (as above) so that there is a natural map $\widehat{\hocolim}_\hfdgl\mathscr{D} \to \colim_\dgl\mathscr{D}$ (recall that $\hocolim_\dg$ of a diagram is a large direct sum (with modified differential) of the \DG\,s in the diagram and their iterated suspensions).
We will not be more specific here because we will not need explicit models for the homotopy colimits of complicated diagrams.  The only remaining fact which we need is:

\begin{thm}
In $\dgl_r$, sequential homotopy colimits commute with homotopy pullbacks.
\end{thm}
\begin{proof}[Proof Sketch]

It is enough to show this for freely generated diagrams $\mathscr{D}:\catI\times \catJ \to \hfdgl$, since cofibrant replacement gives a weak equivalence $\L\C(\mathscr{D}) \xrightarrow{\,\weq\,}\mathscr{D}$ natural in $\mathscr{D}$ and $\L\C(\mathscr{D})$ is freely generated.

For freely generated diagrams, however, this may be proven by a large computation just as the corresponding theorem for $\dg_r$ (\ref{dg_r holim commute}).  The critical fact is that infinite $\liesum$'s of free \DGL\,s commute with finite $\times$.
 
Kuhn also comments in [Kuhn p6] that this may be shown using fact that the sequential small object argument applies in $\dgl$.

\end{proof}

\chapter{Differential Graded Coalgebras}\label{S:DGC}
\markboth{Ben Walter}{I.5 Differential Graded Coalgebras}

Differential graded coalgebras are also models for rational spaces.  In this chapter we mirror the constructions of the previous chapter.
Many of the proofs in this chapter are omitted since they are precisely dual to the proofs of corresponding statements for $\dgl$ in the previous chapter.

\section{Category Structure}

As well as differential graded Lie algebras, we are also particularly interested in differential graded coalgebras.  By ``differential graded coalgebra'' we mean what is more precisely called a ``differential graded, coassociative, cocommutative, coaugmented, counital coalgebra'':

\begin{defn}[DGC]
A differential graded coalgebra (\DGC) $C=(V_\bullet,\, d,\, \Delta,\, \varepsilon, \varrho)$ consists of 
 a differential graded vector space $V = (V_\bullet,\, d)$ 
equipped with
 a coassociative, graded cocommutative comultiplication map (of \DG\,s) $\Delta:(V_\bullet,\, d) \to (V_\bullet,\, d)\otimes (V_\bullet,\, d)$ 
as well as
 a counit\footnote{Recall: $1_\mDG = (\Q v_0, d=0)$ the \DG\ consisting of only a $\Q$ in degree 0.} $\varepsilon:(V_\bullet,\, d) \to 1_\mDG$ and coaugmentation $\varrho:1_\mDG \to (V_\bullet,\,d)$ properly compatible with $\Delta$.\footnote{The compatibility condition is that if we give $1_\mDG$ a trivial coproduct then these maps commute with coproducts as well as differentials.}
\end{defn}

\begin{note} 
Defining the comultiplication map as a \DG-map imposes the following compatibility condition with the differential:  If $v\in V$ with $\Delta(v) = \sum_i a_i\otimes b_i$ then $\Delta (dv) = \sum_i da_i\otimes b_i + (-1)^{|a_i|}a_i\otimes db_i$. 
Also, our \DGC\,s are {\em graded cocommutative} meaning that for $\Delta(v) = \sum_i a_i \otimes b_i$ we have $\sum_i a_i \otimes b_i = \sum_i (-1)^{|a_i|\cdot |b_i|} b_i \otimes a_i$.
\end{note}

Our general convention is that our \DGC\,s have $\varepsilon|_{V_0}:V_0 \to \Q$ an isomorphism, and nothing in negative degrees (though some of our constructions are still valid even if this is not the case); thus we standardly define and write only $\varepsilon$ with the understanding that $\varrho = \varepsilon^{-1}$.  Such \DGC\,s are also known as ``1-reduced'' \DGC\,s.  A map of \DGC\,s is a graded vector space map $f:V_\bullet \to W_\bullet$ which is degree 0 and preserves differentials ($fd = df$), coproducts ($(f\otimes f) \circ \Delta = \Delta \circ f$), augmentation ($\varepsilon f = \varepsilon$), and coaugmentation ($f \varrho = \varrho$).  Surjections, injections, isomorphisms, and quasi-isomorphisms of \DGC\,s are \DGC-maps which induce a surjections, injections, isomorphisms, and quasi-isomorphisms on underlying graded vector spaces.  

A \DGC\ is called $r$-reduced if, aside from the coaugmentation, it is 0 below grading $r$ (i.e. $C_\bullet =\Q\oplus\{C_i\}_{i \ge r}$).  All of the \DGC\,s which we consider are $r$-reduced for some $r$:  our convention is to write $\dgc_r$ for the category of all $r$-reduced \DGC\,s and \DGC-maps.  The notation $\dgc$ without a subscript means either that the $r$ is implicitly present and determined by context or else $\dgc_1$.
Note that any differential graded coalgebra which in grading 0 has only $\Q$  (which is not in the image of $d$) is canonically counital.  Also, note that $\dgc_r$ is a full subcategory of $\dgc_t$ for all $t<r$.

At times we find it desirable to work with non-coaugmented \DGC\,s (given by $(V_\bullet,\, d,\, \Delta)$).  A non-coaugmented \DGC\ is $r$-reduced if its underlying \DG\ is $r$-reduced.  We write $\widetilde{\dgc}_r$ for the category of $r$-reduced non-coaugmented \DGC\,s.  Given a \DGC, $C = (C_\bullet,\, d_C,\, \Delta_C,\, \varepsilon) \in \dgc_r$,  by the de-augmentation\footnote{We would like to call this ``reduction,'' but that word is already taken.} $\widetilde{C}$ we mean the non-coaugmented \DGC\ $\widetilde{C} = (\ker(\varepsilon),\, \widetilde{d},\, \widetilde{\Delta}) \in \widetilde{\dgc}_r$, given by $\widetilde{d} = d|_{\ker(\varepsilon)}$ and $\widetilde{\Delta}(c) = \Delta(c) - 1\otimes c  - c\otimes 1$.  For $r>1$ there is a 1-1 correspondence between $\dgc_r$ and $\widetilde{\dgc}_r$ given by taking a \DGC\ to its de-augmentation $C\mapsto \widetilde{C}$.  The inverse functor is given by adding a disjoint coaugmentation $C \mapsto C_+$ given by $C_+ = ([C]_\mDG \oplus 1_\mDG,\,\Delta_+,\, \varepsilon:1_\mDG \xrightarrow{\cong} 1_\mDG)$ where $\Delta_+(c) = \Delta_C(c) + 1\otimes c + c\otimes 1$ for $C = (C,\, d_C,\, \Delta_C)$.
  At times we find it convenient to describe \DGC\,s by describing, instead, their de-augmentation (see for example Definition~\ref{dgc trunc red}).  When describing the de-augmentation of a coalgebra, our convention is to indicate this by putting a tilde over the coalgebra's name (e.g. $\widetilde{C}$).  

The difference between coaugmented and non-coaugmented \DGC\,s is similar to the difference between based and unbased topological spaces.  This is the reason for our choice of the suggestive notation $(-)_+$ and the terminology ``adding a disjoint coaugmentation'' for the map $\widetilde{\dgc} \to \dgc$.  

Reduction is defined essentially by reducing the de-augmentations of a \DGC\,s and then adding back a coaugmentation:

\begin{defn}[Reduction]\label{dgc trunc red}\footnote{
There is also a truncation functor, but we are uninterested in it.}
Given a $t$-reduced \DGC\ $C\in\dgc_t$  with $1<t<r$
 the {\em $r$-reduction} $\mathrm{red}_r(C) = \red{C} \in \dgc_r$ is the sub-\DGC\ $\red{C}$ given by 
the following procedure:   Write $\mathrm{red}^\mDG_r\widetilde{C}$ for the reduction of $\widetilde{C}$ as a \DG.  Define $$\red{\widetilde{C}} = \Bigl\{c\in \mathrm{red}^\mDG_r\widetilde{C}\ \ |\ \ \widetilde{\Delta}(c) \in \bigl(\mathrm{red}^\mDG_r\widetilde{C}) \otimes (\mathrm{red}^\mDG_r\widetilde{C}\bigr)\Bigr\}$$
and let $\red{C} = \red{\widetilde{C}}_+$.
\end{defn}

\begin{lemma}\label{dgc red}
Reduction is right adjoint to the inclusion of categories functor $\dgc_r \to \dgc_t$
\begin{align*}
[-]_{\mDGC_t}:\dgc_r &\rightleftarrows \dgc_t:\mathrm{red}_r 
\end{align*}
\end{lemma}

The categories $\dgc_r$ are all pointed:  The \DGC\ consisting of only the counit with trivial differential and coproduct $0_\mDGC = (1_\mDG,\, \Delta(c) = 1\otimes c + c\otimes 1,\,\varepsilon:1_\mDG \xrightarrow{\cong}1_\mDG)$ is both initial and final.\footnote{Note: $0_\mDGC$ is given by $0_\mDGC = (0_\mDG,\, \Delta=0)_+$.}  If a \DGC\ has trivial reduced homology $\widetilde {H}_*(C) = H_*(\widetilde{C})= 0$ (for example if $C=0_\mDGC$) then we say that $C$ is {\em contractible}.  Note that $C$ is contractible if and only if the map $C\to 0_\mDGC$ is a quasi-isomorphism.  In general of course, the homology of a \DGC\ is a graded coalgebra, and taking homology gives a functor $H_\ast:\dgc \to \gc$.

\subsection{Cofree DGCs}

We define cofree \DGC\,s in a manner similar to the way that we defined free \DGL\,s.  A similar note about heinous abuse of notation applies in this case as well -- what we define should really be called a ``differential cofree graded coalgebra.''   We nonetheless persist in calling it a ``cofree \DGC''. 
 
Given a graded vector space $V = \{V_i\}_{i\ge 2}$ starting in grading 2, we may consider the tensor coalgebra $TV_\bullet = \bigoplus_k\, T^kV_\bullet = \bigoplus_k\, \{v_1|\cdots|v_k \text{ where } v_i \in V_\bullet\}$.  This is the cofree, counital, coagumented (the unit is an isomorphism on degree 0 $\varepsilon:TV \to T^0V \cong \Q$) coassociative coalgebra (with reduced comultiplication given by $\widetilde{ \Delta}(v_1|\cdots|v_n) = \sum_{i=2}^n (v_1|\cdots|v_{i-1})\otimes (v_{i}|\cdots|v_n)$) primitively cogenerated by $V$, but is not cocommutative.  Note that the symmetric group $\Sigma_n$ acts on the left of $T^nV$ ({\em with signs according to the Koszul convention})\footnote{That is, $\Sigma_n$ acts by a permutation of indices twisted by the sign of the permutation:
$$\pi(v_1\otimes\cdots\otimes v_n) = \mathrm{sgn}(\pi)\,(v_{\pi(1)}\otimes\cdots\otimes v_{\pi(n)})$$}.   By $\Lambda V$ we mean the invariants (or equivalently coinvariants) of these actions $\Lambda V = \bigoplus_k (T^kV)^{\Sigma_k} = \bigoplus_k \Lambda^k V$: the cofree, coaugmented, counital, graded cocommutative coalgebra primitively cogenerated by the graded vector space $V = V_\bullet$.  The functor $\Lambda -:\g_1 \to \gc_1$ is the right adjoint of the functor $[-]_\mG:\gc_1 \to \g_1$ (see e.g. [Q69 B.4.1] or [FHT 22.1]).\footnote{As a graded vector space, $\Lambda V$ is given by the graded-symmetric powers of $V$, so maybe $SV$ would be better notation; however, $\Lambda V$ appears to already be entrenched in the literature as the conventional notation for this object.  Also $S$ is already quite overused.}  We write elements of $\Lambda^k V$ by $v_1\wedge \cdots \wedge v_k$.

\begin{ex} 
If $V_\bullet$ is given by $V_i = \Q v$ for $i$ some fixed integer and $V_j = 0$ for all $j \neq i$, then
\begin{enumerate}
\item[(i)]  if $i$ is even, $\Lambda V$ is the polynomial coalgebra on $v$ over $\Q$:
$$\Lambda V = (\Q[v],\, \Delta(v^n) = \sum_{i+j=n} v^i\otimes v^j)$$
\item[(ii)] if $i$ is odd, then $\Lambda V$ is the exterior coalgebra on $v$ over $\Q$:
$$\Lambda V = \left(\Q v,\, \Delta(v)=1\otimes v + v\otimes 1 \right)$$ 
\end{enumerate}
\end{ex}

\begin{defn}[Co-Free DGC] 
A {\em cofree differential graded coalgebra} is a \DGC\ $C$ which, as a graded coalgebra, is isomorphic to a cofree graded coalgebra:
$$[C]_\mGL\cong\Lambda V\qquad \text{where}\qquad V\in\g$$

Write $\Fdgc_r$ for the full subcategory of all cofree \DGC\,s in $\dgc_r$.
\end{defn}

Writing $C$ as $C=([C]_\mGC,\,d)$ and pushing the differential of $C$ across the isomorphism $[C]_\mGL \cong \Lambda V$, we get the equivalent statement that $C$ is cofree if and only if $C$ is isomorphic to a \DGC\ of the form $C\cong(\Lambda V,\, d)$ for some graded vector space $V$ and differential $d$.

\begin{note}\label{truly cofree}
Again the {\em truly cofree} objects in $\dgc_r$ are those isomorphic to $\Lambda(V,\,d_V)$ (for some $(V,\,d_V)$) in the image of the obvious functor $\Lambda-:\dg_r\to \dgc_r$ which is right adjoint to the functor $[-]_\mDG:\dgc_r\to \dg_r$.  Truly cofree objects are denoted by $C\cong\Lambda(V,\,d_V)$ as opposed to regular cofree \DGC\,s which are written $C\cong \Lambda V = (\Lambda V,\, d)$.
\end{note}

Closely related to $\Fdgc$ is the category $\hFdgc$ of all \DGC\,s of the form $(\Lambda V,\, d)$ and \DGC-maps between them.  Objects of $\hFdgc$ consist of a cofree \DGC\ $C$ along with 
an isomorphism $C\cong (\Lambda V,\, d)$ (however, maps are not required to respect these isomorphisms).  We may lazily write $\Lambda V = (\Lambda V,\, d)$ for such \DGC\,s (leaving the differential unwritten).  The cateogries $\Fdgc$ and $\hFdgc$ are equivalent, though not isomorphic.  This allows us to construct most proofs about cofree \DGC\,s by considering only cofree \DGC\,s of the form $(\Lambda V,\, d)$.  In general, our intuition about $\Fdgc$ follows from our intuition about objects $(\Lambda V,\, d)$.

Dual to the case with free \DGL\,s, free \DGC\,s of the form $\Lambda V = (\Lambda V,\, d)$ inherit an extra grading from the word-length grading on $TV$, which allows us to write the differential in the form $d=d_0+d_{-1}+\cdots$ where $d_{-i}:\Lambda^k V \to \Lambda^{k-i} V$.  Again $d_0 \circ d_0 = 0$, so $d$ restricts to a differential $d_V = d|_V$ on $V$.  We write $(V,\,d_V)$ for the \DG\ which this restriction defines.  The \DGC\ $(\Lambda V,\, d)$ is truly cofree (in the sense of~\ref{truly cofree}) if  $d=d_0$ -- in this case $(\Lambda V,\, d) = \Lambda(V,\, d_V)$.
 
More generally, if $C$ is any \DGC\ we write $(C)^\pr$ for the primitive elements of $C$.  This is a \DG\ given by $(C)^\pr := \ker(\widetilde{\Delta}_C)$, where $\Delta_C$ is the coproduct of $C$ (and $\widetilde{\Delta}_C$ is the coproduct of the de-augmentation $\widetilde{C}$).  [Since the differential must commute with the coproduct, it follows that $d$ cannot increase word length.  Note in particular that the differential on a \DGC\ must preserve primitives.]
If $C\cong(\Lambda V,\,d)$ then there is an isomorphism $[(C)^\pr]_\mG \cong V$, and furthermore $C$ is truly cofree if and only if $C\cong \Lambda(C)^\pr$.  The \DGC\ version of the rational Hurewicz theorem states that $(C)^\pr$ captures much of the homology information of cofree \DGC\,s:

\begin{thm}[Rational Hurewicz (cofree DGCs)]\label{Hurewicz dgc} 
If $C\in\Fdgc_2$ is a cofree \DGC\ then the \DG\ $(C)^\pr$ is contractible if and only if $C$ is contractible.

Moreover if $C_1,C_2\in\Fdgc_2$ are cofree \DGC\,s then a \DGC-map $C_1 \xrightarrow{\,f\,} C_2$ induces a \DG-quasi-isomorphism $(C_1)^\pr \xrightarrow{\,(f)^\pr\,} (C_2)^\pr$ if and only if $f$ itself is a quasi-isomorphism.\footnote{That is, $(-)^\pr$ detects and reflects quasi-isomorphisms of cofree \DGC\,s.}
\end{thm}

\begin{proof}

In $\Fdgc_2$ this follows immediately from the Hurewicz Theorem for topological spaces along with Theorem~\ref{Quillen Thm I} and Lemma~\ref{free homotopy} which are given later.

\end{proof}

Note that the corresponding statement about non-cofree \DGC\,s is not true.  In particular, the primitives functor $(-)^\pr:\dgc\to\dg$ does not preserve quasi-isomorphisms.

Just as with free \DGL\,s we are interested in two special classes of maps called ``cofree maps'' (in $\dgl$) and ``cofreely generated maps'' (in $\hFdgl$).  A cofree map of \DGC\,s is as follows:

\begin{defn}[Cofree Map]\label{cofree map}
A map of \DGC\,s $f:B\to C$ is a {\em cofree map} if as a map of graded coalgebras, it is isomorphic to a projection to $[C]_\mGC$ with cofree kernel, as expressed by the diagram:
$$\xymatrix@R=7pt@C=25pt{
\bigl[C\bigr]_\mGC \ar[r]^{[f]_\mGC}
  \ar@{{}{}{}}[d]_*[right]{\cong}
   &  \bigl[D\bigr]_\mGC \\
\bigl[D\bigr]_\mGC \otimes \Lambda V \ar@{->>}[ur] &
}$$
\end{defn}

\begin{lemma} Let $f:B\to C$ be a cofree map.
If $C$ is a cofree \DGC, then $B$ is also cofree.
\end{lemma}
\begin{proof}

Suppose $C$ is cofree.  Products commute with left adjoint functors, so $\Lambda V \times \Lambda W = \Lambda (V \oplus W)$.  Therefore, if $C \cong (\Lambda W,\, d)$ and $f:C\to D$ is a cofree map, then
$$[B]_\mGC \ \cong\  [C]_\mGC \otimes \Lambda V \ \cong\ \Lambda W \otimes \Lambda W
 \ = \ \Lambda(V\oplus W)$$
So $B$ is also cofree.

[This also follows from the later result that cofree maps are precisely the fibrations in $\dgc$ and cofree \DGC\,s are precisely the fibrant \DGC\,s.]
\end{proof}

One useful way of constructing cofree maps in $\hFdgc$ is given by the following lemma:

\begin{lemma}
Let $(\Lambda V,\, d),\, (\Lambda W,\, d') \in \hFdgc_r$ and 
$f:(\Lambda V,\, d) \to (\Lambda W,\, d')$ be a \DGC-map such that $[f]_\mGC = \Lambda(\hat f)$ where $\hat f:V \fibr W$ is a surjection.  Then $f$ is a cofree map.
\end{lemma}
\begin{proof}[Proof Sketch]
This follows from the fact that surjections in \G\ split, so $\hat f:V=U\oplus W \cofibr W$ (recall that $\Lambda (U\oplus W) = \Lambda U \otimes \Lambda W$). 
\end{proof}

The actual class of maps which we are interested in are a little weaker than those of the previous lemma. 

\begin{defn}[Cofreely Generated Map] 
Given $f:(\Lambda V,\, d) \to (\Lambda W,\, d')$ a \DGC-map.  The map $f$ is {\em cofreely generated} if $[f]_\mGC = \Lambda(\hat f)$ where $\hat f:V \to W$ is any \G-map.

Write $\hfdgc_r$ for the subcategory of $\dgc_r$ consisting of all cofree \DGC\,s of the form $(\Lambda V,\, d)$ and all cofreely generated maps between them.
\end{defn}

\subsection{Adjoints between $\dgc$ and $\dg$}\label{s:adjoints dgc}

Recall the functor $[-]_\mDG:\dgc_r\to \dg_r$ which {\em de-augments and then} forgets the coproduct structure of a \DGC.  Our notation is $[(C_\bullet,\, d,\, \Delta,\, \varepsilon)]_\mDG = (\widetilde{C}_\bullet,\, \widetilde{d})$.  This functor has a section the functor $[-]_{\mDGC_r}\!:\dg_r\to \dgc_r$ which equips a \DG\ with a trivial coproduct {\em and then adds a disjoint augmentation} ($[V]_{\mDGC_r}\!$ is given by $[V]_{\mDGC_r}\! = (V,\, \Delta=0)_+$).  Again, we are particularly interested in the functors to and from \DG\ given by composing the above with the inclusion and reduction functors.  Abusing notation we write:
\begin{align*}
 [-]_{\mDG}&:\dgc_r \xrightarrow{[-]_{\mDG_r}} \dg_r \xrightarrow{\mathrm{incl}} \dg \\
 [-]_{\mDGC_r}&:\dg \xrightarrow{\mathrm{red}_r} \dg_r \xrightarrow{[-]_{\mDGC_r}} \dgc_r
\end{align*} 
There are obvious corresponding maps to and from $\widetilde{\dgc}_r$.

Oftentimes we write just $[-]_\mDGC$ for the functor $[-]_{\mDGC_r}$ with the understanding that the $r$ is implicit and determined by context.  We are also interested in the forgetful functor $[-]_\mGC$ from $\dgc$ to $\gc$ as well as the functor $[-]_\mG = (\widetilde{-})_\bullet$ from $\dgc$ to $\g$ which de-augments a coalgebra and then forgets the differential and coproduct.

Note that we have two adjoint pairs between $\dgc$ and $\dg$ given by:
\begin{lemma}[Adjoints]\label{dgc adjoints} 
The following give adjoints between $\dgc_r$ and $\dg$:\footnote{The truncation functor $\dg\to\dg_r$ may be used to extend the top adjoint to $\dg\rightleftarrows\dgc_r$, but since truncation doesn't preserve weak equivalences, we do not wish to go that far.}
$$\xymatrix@C=30pt@R=15pt{
& [-]_{\mDGC_r}:\dg_r \ar@<2pt>[r]^(.54){[-]_{\mDGC_r}} &
  \dgc_r:(-)^\pr \ar@<2pt>[l]^(.45){(-)^\pr} } $$
$$\xymatrix@C=30pt@R=15pt{
{}\phantom{XX}[-]_\mDG:\dgc_r \ar@<2pt>[r]^(.67){[-]_\mDG} &
  \dg_r \ar@<2pt>[l]^(.32){\Lambda-} \ar@<2pt>[r]^(.46){\mathrm{incl}} &
  \dg:\Lambda- \ar@<2pt>[l]^(.52){\mathrm{red}_r}
}$$          
where $\mathrm{incl}$ is the inclusion of categories functor.

Furthermore these adjoints are sections in the following sense:
\begin{itemize}
\item $\Lambda-:\dg_r \to \dgc_r$ is a section of $(-)^\pr:\dgc_n\to \dg_n$.
\item $[-]_\mDGC:\dg_r \to \dgc_r$ is a section of $[-]_\mDG:\dgc_r\to \dg_r$.
\end{itemize}
\end{lemma}

Again related to the rational Hurewicz Theorem is the more simple lemma:
\begin{lemma}\label{weq dgc equals weq dg}
If $C\in\dgc$ is a \DGC, then the \DG\ $[C]_\mDG$ is contractible if and only if $C$ is contractible.

Moreover if $C_1, C_2 \in \dgc$ are \DGC\,s then a \DGC-map $C_1\xrightarrow{f} C_2$ induces a \DG-quasi-isomorphism $[C_1]_\mDG \xrightarrow{[f]_\mDG} [C_2]_\mDG$ if and only if $f$ itself is a quasi-isomorphism.
\end{lemma}

\begin{note}\label{dgc factorization} 
As in the previous chapter, we may use the above sections to recover \DGC\,s modulo some structure.  In particular, given $C\in\dgc_r$, a composition of the left adjoints $\bigl[[C]_\mDG\bigr]_\mDGC$ recovers $C$ modulo coproduct structure:
$$C = \Bigl(\Bigl[\bigl[[C]_\mDG\bigr]_\mDGC\Bigr]_\mDG,\ \Delta_C\Bigr)$$
Similarly, given $C\cong(\Lambda V,\, d)$ a cofree \DGC, a composition of the right adjoints $\Lambda(C)^\pr$ recovers $C$ up to higher order differential information:
$$C\ \cong\ \bigl(\Lambda V,\ d=d_0 + d_{>0}\bigr) \ \cong\  \Bigl(\bigl[\Lambda (C)^\pr\bigr]_\mGC,\ d = d_{\Lambda (C)^\pr} + d_{>0}\Bigr)$$
where $d_0 = d_V$ is the degree 0 part of the differential of $(\Lambda V,\, d)$; $d_{>0} = d - d_0$ is the part of the differential of $(\Lambda V,\, d)$ which increases word-length by at least 1; and $d_{\Lambda(C)^\pr}$ is the differential on the truly cofree \DGC\ $\Lambda(C)^\pr$.
\end{note}

\subsection{Limits and Colimits}

The category $\dgc$ has supports four interesting operations -- non-coaugmented and coaugmented sums ($\oplus$ and $\tildeoplus$), products ($\otimes$), and coaugmented products ($\tildeotimes$).  The operations $\oplus$ and $\tildeoplus$ are analogous to the disjoint union ($\coprod$) and wedge ($\vee$) of based spaces -- $\oplus$ is the categorical coproduct in the category $\widetilde{\dgc}$ just as $\coprod$ is the categorical coproduct in $\top$; similarly $\tildeoplus$ is the categorical coproduct in $\dgc$ just as $\vee$ is the categorical coproduct in $\top_*$.  The operations $\otimes$ and $\tildeotimes$ are analogous to the cross product ($\times$) and smash ($\wedge$) of $\top_*$.  In fact, if we view the homology of the based spaces $X$ and $Y$ as \DGC\,s, these four operations are defined specifically so that the following hold:
\begin{itemize}
\item
$H_*(X\coprod Y;\, \Q) \cong H_*(X;\, \Q)\oplus H_*(Y;\, \Q)$
\item
$H_*(X\vee Y;\, \Q) \cong H_*(X;\, \Q)\tildeoplus H_*(Y;\, \Q)$
\item
$H_*(X\times Y;\, \Q) \cong H_*(X;\, \Q)\otimes H_*(Y;\, \Q)$
\item
$H_*(X\wedge Y;\, \Q) \cong H_*(X;\, \Q)\tildeotimes H_*(Y;\, \Q)$
\end{itemize}

We neglect below to mention the precise category where products, coproducts, limits, and colimits are being taken with the understanding that if the objects or diagrams are all in $\dgc_r$ (or $\widetilde{\dgc}_r$) then so are the products, coproducts, limits, and colimits.  We implicitly assume that $r\ge 1$ (though the definitions of coproducts and colimits below are still correct for $r<1$).

\begin{lemma}[Non-coaugmented Sums] 
Given non-coaugmented \DGC\,s \, $C = (C_\bullet,\, d_C,\, \Delta_C)$ and $D = (D_\bullet,\, d_D,\, \Delta_D)$, their categorical coproduct in $\widetilde{\dgc}$ is given by taking $\oplus$ of the underlying \DG\,s, and equipping it with a coproduct map:
$$C\oplus D := \big((C_\bullet,\, d_C)\oplus (D_\bullet,\, d_D),\, \Delta_\oplus\big)$$
where $\Delta_\oplus = \Delta_C \oplus \Delta_D$.

This generalizes to give small coproducts of non-coaugmented \DGC\,s.
\end{lemma}

\begin{lemma}[(Coaugmented) Sums] 
Given 1-reduced \DGC\,s \, $C = (C_\bullet,\, d_C,\, \Delta_C,\, \varepsilon_C)$ and $D = (D_\bullet,\, d_D,\, \Delta_D,\, \varepsilon_D)$, their categorical coproduct in $\dgc$ is 
$C\tildeoplus D$ given by 
$${C\tildeoplus D} = (\tilde C \oplus \tilde D)_+$$

This generalizes to give small coproducts.
\end{lemma}

The relationship between non-coaugmented sums and coaugmented sums is precisely the relationship between the disjoint union of two based topological spaces and the wedge of the two based topological spaces.  

\begin{lemma}[Products] 
Given two \DGC\,s $C = (C_\bullet,\, d_C,\, \Delta_C,\, \varepsilon_C)$ and $D = (D_\bullet,\, d_D,\, \Delta_D,\, \varepsilon_D)$, their categorical product is given by taking $\otimes$ of the underlying \DG\,s and equipping it with a coproduct map:
$$C\otimes D := ((C_\bullet,\, d_C)\otimes (D_\bullet,\, d_D),\, \Delta_{\otimes},\,\varepsilon_{\otimes} = \varepsilon_C \otimes \varepsilon_D)$$
where $\Delta_{\otimes}$ is given by 
$$\Delta_{\otimes}(v\otimes a) = \frac{1}{2}\sum_{i,j} (v_i \otimes a_j) \otimes (w_i\otimes b_j) + 
(v_i\otimes b_j) \otimes (w_i \otimes a_j)$$ 
for $v\in C$, $a\in D$ with $\Delta_C(v) = \sum_i v_i\otimes w_i$ and $\Delta_D(a)= \sum_j a_j \otimes b_j$.

This generalizes to give finite products.
\end{lemma}

Note that this descends to a natural definition of $\otimes$ for non-coaugmented \DGC\,s.  However, in the category $\widetilde{\dgc}$ it does not define a categorical product, since the projection maps $C_1 \otimes C_2 \to C_i$ require the counits in order to be defined.

\begin{defn}[Coaugmented Products]\label{def: tildeotimes} 
Given 1-reduced \DGC\,s \, $C$ and $D$  define their coaugmented product to be $\tildeotimes$ where $C\tildeotimes D$ is the \DGC\ defined by
${C\tildeotimes D} := (\tilde C\otimes \tilde D)_+$.
\end{defn}

\begin{note}The coaugmented product is {\em NOT} the product in the category $\dgc$ of (coaugmented) \DGC\,s.  This is just as smash is not the product on the category of based topological spaces.
\end{note}

\begin{cor}[Sums on $\Fdgc$]\label{fdgc sum} 
The coproduct of cofree \DGC\,s has
$\Lambda V \tildeoplus \Lambda W \xrightarrow{\weq} \Lambda(V\oplus W)$.
\end{cor}

\begin{cor}[Products on $\Fdgc$]\label{fdgc prod} 
The product of cofree \DGC\,s has $\Lambda V \otimes \Lambda W \cong \Lambda(V \otimes W)$.
\end{cor}

Finally, note that the category $\dgc$ has pushouts and pullbacks given by:

\begin{lemma}[Pullbacks and Pushouts] 
The pullback of the \DGC-diagram $C_1\xrightarrow{f_1} B\xleftarrow{f_2} C_2$ is given by:
$$C_1\otimes_B C_2 := \{c_1\otimes c_2 \in C_1\otimes C_2\ \ |\ \ \ f_1(c_1)+f_2(c_2)=0\}$$

The pushout of the \DGC-diagram $C_1\xleftarrow{f_1} B\xrightarrow{f_2} C_2$ is given by:
\begin{align*}C_1\tildeoplus_B C_2 &:= (C_1\tildeoplus C_2)/\langle f_1(b) + f_2(b)\rangle \\
&= (C_1\oplus C_2)/\langle f_1(b) + f_2(b)\rangle
\end{align*} 

General finite limits and small colimits are defined analogously using products, coproducts, equalizers (kernels) and coequalizers (cokernels) as in Theorem~\ref{(co)lim (co)equalizer}. 
\end{lemma}

\begin{cor} 
Let $f$ be a map $f:B \to C$.  The fiber of $f$ is $\ker(f)$; the cofiber of $f$ is $\coker(f)$.
\end{cor}

In fact, a slight modification of an argument by Getzler and Goerss proves:

\begin{thm}[\lbrack GG99, 1.8\rbrack]
The category $\dgc$ has small limits.
\end{thm}
\begin{proof}[Sketch of Proof]
In [GG99] Getzler and Goerss construct small limits in the category of counital, coassociative, {\em non-cocommutative} differential graded coalgebras.  Merely adding ``cocommutative'' and ``commutative'' throughout their construction transports it to our case.
The argument goes roughly as follows:

First, show that every homogeneous element of a \DGC\ is contained in a finite dimensional sub-\DGC.  This implies that every \DGC\ is a filtered colimit of its finite dimensional sub-\DGC\,s.  Thus \DGC\,s are equivalent to {\em ind-finite}-\DGC\,s.  Dualizing an {\em ind-finite}-\DGC\ gives a {\em pro-finite}-\textsc{dga}.  The continuous dual of a {\em pro-finite}-\textsc{dga} recovers the original \DGC, since all of our dualization is of finite objects.

To compute the limit of a diagram of \DGC\,s we may therefore dualize the diagram to one of {\em pro-finite}-\textsc{dga}\,s, take the colimit of that diagram, and then dualize the answer back to $\dgc$.  The colimit of a diagram of {\em pro-finite}-\textsc{dga}\,s is given by a completion of the colimit of the diagram in the category of \textsc{dga}\,s. 
\end{proof}

All of the above operations ($\oplus$, $\tildeoplus$, $\otimes$, $\tildeotimes$) clearly descend to the category $\gc_r$ of graded (cocommutative, counital) coalgebras where the same operations again give products, coproducts, limits, and colimits.  Furthermore note that if $V,\, W\in \g_r$ are graded vector spaces, then  $\Lambda V \tildeoplus \Lambda W$ maps to $\Lambda (V\oplus W)$ by an inclusion of a sub-\GC:  $\Lambda V\tildeoplus \Lambda W \cofibr \Lambda(V\oplus W)$.  Later on, we will wish to define differentials on the \GC\ $\Lambda V\tildeoplus \Lambda W$ by defining them first on $\Lambda(V\oplus W)$ and then restricting them to $\Lambda V\tildeoplus \Lambda W$.

\begin{note}\label{dgc colim plan} 
We can also recover the formula for the categorical coproduct and more generally all colimits in $\dgc_r$ using the adjoint functors from the previous section along with the comments in Remark~\ref{dgc factorization}.  Since the functors $[-]_\mDG$ and $[-]_\mDGC$ are both left adjoints, they commute with colimits.  
In particular, if $\mathscr{D}:\catI\to\dgc_r$ is a diagram in $\dgc_r$ then we have
$$
\bigl[[\colim{\!}_\dgc\mathscr{D}]_\mDG\bigr]_\mDGC =
 \bigl[\colim{\!}_\dg[\mathscr{D}]_\mDG\bigr]_\mDGC 
$$

According to Remark~\ref{dgc factorization}, this means that in order to compute a colimit in $\dgc$, we may instead compute the corresponding colimit in $\dg$ and then figure out the correct coproduct structure on the result.  We leave it to the interested reader that in general, there is only one coproduct structure possible so that there are natural maps $\mathscr{D}\longrightarrow\colim_\dgl\mathscr{D}$ which descend to the existing natural maps $[\mathscr{D}]_\mDG\longrightarrow \colim{\!}_\dg[\mathscr{D}]_\mDG = [\colim_\dgc \mathscr{D}]_\mDG$.  This is what is meant when people say that ``colimits in $\dgc$ are {\em created} in $\dg$.''

Similarly, since the functors $(-)^\pr$ and $\Lambda(-)$ are both right adjoints, they commute with limits.  For $\mathscr{D}':\catJ\to \fdgc_r$ a diagram in $\hfdgc_r$ we have
$$\Lambda(\lim{\!}_\hfdgc\mathscr{D}')^\pr = \Lambda\bigl(\lim{\!}_\dg(\mathscr{D}')^\pr\bigr)$$
[It is an easy exercise that $\lim_\hfdgc\mathscr{D}' = \lim_\dgc\mathscr{D}'$.]  Now \ref{dgc factorization} suggests that limits in $\fdgc$ may also be created in $\dg$.  This is in fact the case.  Again we leave it to the interested reader that there is a unique differential which recovers the higher order differential information on $\lim{\!}_\hfdgc\mathscr{D}'$ so that there are natural maps $\lim{\!}_\hfdgc\mathscr{D}' \longrightarrow \mathscr{D}'$ which descend to the existing natural maps 
$\lim{\!}_\dg(\mathscr{D}')^\pr = (\lim{\!}_\hfdgc\mathscr{D}')^\pr  \longrightarrow (\mathscr{D}')^\pr$. 
\end{note}

\subsection{Cones, Suspensions, Paths, and Loops}

\begin{defn}[Cones] 
If $C$ is a $r$-reduced \DGC, then its cone is the $(r+1)$-reduced \DGC\ $\mathrm{c}C$ defined by ${\mathrm{c}C} := (c \otimes  \tilde C)_+$, where $c$ is the non-coaugmented \DGC\ on $c_\bullet = \Q v_0 \oplus \Q v_1$, $|v_i| = i$ with $d(v_1) = v_0$ and $\Delta(v_0) = v_0\otimes v_0$, $\Delta(v_1) = v_0\otimes v_1 + v_1 \otimes v_0$.
\end{defn}

Note that this is just the cone on $[C]_\mDG$ equipped with a coproduct which is compatible with the differential of $\mathrm{c}[C]_\mDG$.  Recall that this means $(\widetilde{\mathrm{c}C})_\bullet = \tilde C_\bullet \oplus s\tilde C_\bullet$.
Again, just like cone on a topological space $X$, the cone on a $C$ is a contractible \DGC\ (i.e. quasi-isomorphic to 0) which comes equipped with an injection $C\to \mathrm{c}C$.

\begin{defn}[Suspensions] 
If $C$ is a \DGC, then its suspension  is the \DGC\ $\Sigma C$ defined by 
$\Sigma C := (s\otimes \widetilde{C})_+$ where $s$ is the non-coaugmented \DGC\ with trivial differential and bracket on $s_\bullet = \Q s$, $|s|=1$. 
\end{defn}

\begin{note}
$\Sigma C$ is equal to $\Sigma C =  \big[s[C]_\mDG\big]_\mDGC$.
\end{note}

Again, we can define a slightly larger model for $\Sigma C$ by introducing a  coproduct structure on $\hat\Sigma [C]_\mDG$.  Define
$\hat\Sigma C := (S\otimes \widetilde{C})_+$
where $S$ is the non-coaugmented \DGC\ $S=(S_\bullet,\, d_S,\, \Delta_S)$ with:
\begin{itemize}
\item $S_\bullet = (\Q v_0 \oplus \Q v_1\oplus \Q w_1)$ (for $|v_0| = 0$, and $|v_1|=1=|w_1|$)
\item $d_S(v_0) = 0$, $d_S(w_1) = -v_0$, $d_S(v_1)=v_0$
\item $\Delta_S(v_0)=v_0\otimes v_0$, $\Delta_S(v_1)=\frac{1}{2}(v_1\otimes v_0 + v_0\otimes v_1)$, $\Delta_S(w_1) = \frac{1}{2}(w_1\otimes v_0 + v_0\otimes w_1)$.
\end{itemize}

More explicitly, $\hat \Sigma C$ is given by
$$\hat {\Sigma} (C,\, d_C, \Delta_C) = (s\tilde C_\bullet\oplus \tilde C_\bullet\oplus s\tilde C_\bullet,\, d_{\hat\Sigma},\, \Delta_{\hat\Sigma})_+$$
where, for $sc_1+c_2+sc_3  \in s\tilde C_\bullet\oplus \tilde C_\bullet\oplus s\tilde C_\bullet$,
\begin{itemize}
\item $d_{\hat\Sigma}$ is defined by $d_{\hat\Sigma}(sc_1 + c_2 + sc_3) = -s d_Cc_1 +(d_Cc_2 + c_1 + c_3) - sd_Cc_3$ and 
\item $\Delta_{\hat\Sigma}$ is defined by 
\begin{itemize}
\item[\textbullet] $\Delta_{\hat\Sigma}(c_2) = \tilde \Delta_C(c_2)$ 
\item[\textbullet] $\Delta_{\hat\Sigma}(sc_1) = \frac{1}{2}\sum_i sa_i\otimes b_i + (-1)^{|a_i|}a_i\otimes b_i$, where $\tilde\Delta_C(c_1) = \sum_i a_i\otimes b_i$, 
\item[\textbullet] $\Delta_{\hat\Sigma}(sc_3)$ is defined similar to $\Delta_{\hat\Sigma}(sc_1)$
\end{itemize}
\end{itemize}

There is an injection $C\to \hat \Sigma C$ which plays the role of the map of spaces $X\to \Sigma X$ sending $X$ to the equator.  Also, either of the two projection maps to $s\widetilde{C}_\bullet$ induce quasi-isomorphisms $\hat\Sigma C \xrightarrow{\weq} \Sigma C$.

\begin{note}
The \DGC\,s $\Sigma C$ and $\hat\Sigma C$ are both pushouts:
\begin{itemize}
\item $\Sigma C$ is the pushout of $0_\mDGC \xleftarrow{} C \xrightarrow{} \mathrm{c}C$
\item $\hat \Sigma C$ is the pushout of $\mathrm{c}C \xleftarrow{} C \xrightarrow{} \mathrm{c}C$
\end{itemize}
\end{note}
 
Given a cofree 2-reduced \DGC\ $\Lambda V = (\Lambda V, d)$ we write $\Lambda s^{-1}V$ for the truly cofree \DGC\ given by $\Lambda s^{-1}V = \Lambda s^{-1}(V,\,d_V)$ (recall that $d_V$ is the restriction of $d$ to the cogenerating \G).  This \DGC\ plays the dual role to that which $\freeL_{sV}$ played for \DGL\,s in defining paths and loops.

\begin{defn}[Unreduced Paths of $\Lambda V$] 
Given $\Lambda V = (\Lambda V,\, d)$  a cofree $r$-reduced \DGC\ where $r\ge 2$ define the {\em unreduced paths of $\Lambda V$} to be $\tilde{\mathrm{p}}\Lambda V$ the $(r-1)$-reduced cofree \DGC\ given by taking the coproduct $\Lambda V \tildeoplus \Lambda s^{-1}V$ and modifying the differential:
$$\tilde{\mathrm{p}}\Lambda V := \bigl([\Lambda V \tildeoplus \Lambda s^{-1}V]_\mGC,\, d_c = d_{\tilde\oplus} + \Lambda(d')\bigr)$$
where $d_{\tilde\oplus}$ is the differential on $\Lambda V \tildeoplus \Lambda s^{-1}V$ and $\Lambda(d')$ is the differential on $[\Lambda V \tildeoplus \Lambda s^{-1}V]_\mGC$ given by the restriction of the differential on $\Lambda(V\oplus s^{-1}V)$ cofreely generated by the differential $d'(v) = s^{-1}v$ on $(V\oplus s^{-1}V)$.
\end{defn}

Note that by Corollary~\ref{fdgc sum}, $\tilde{\mathrm{p}}\Lambda V$ is merely $\Lambda(V \oplus s^{-1} V)$ with a modified differential.  By the Rational Hurewicz Theorem, $\tilde{\mathrm{p}}\Lambda V$ is contractible.  Also there is a \DGC-map $\tilde{\mathrm{p}}\Lambda V \to \Lambda V$.  Furthermore, the map $\tilde{\mathrm{p}}\Lambda V \to \Lambda V$ is a cofree map of cofree \DGC\,s.

\begin{defn}[Unreduced Loops on $\Lambda V$] 
Given $\Lambda V = (\Lambda V,\, d)$  a cofree $r$-reduced \DGC\ with $r\ge 2$ define the {\em unreduced loops on $\Lambda V$} to be the truly cofree $(r-1)$-reduced \DGC\ $\tilde{\Omega} \Lambda V := \Lambda s^{-1}V$ defined earlier.
\end{defn}

\begin{defn}[Paths and Loops of $\Lambda V$]\label{dgc paths and loops}
Define the paths and loops of the cofree \DGC\ $\Lambda V = (\Lambda V,\, d)$ to be $\mathrm{p}\Lambda V := \mathrm{red}_r(\tilde{\mathrm{p}}\Lambda V)$ and $\Omega \Lambda V := \mathrm{red}_r(\tilde{\Omega}\Lambda V)$.
\end{defn}

Note that since the reduction functor preserves weak equivalences, $\mathrm{p}\Lambda V$ is still a contractible \DGC.  Furthermore the map $\mathrm{p}\Lambda V \to \Lambda V$ is still a cofree map of cofree \DGC\,s.

At times we desire a slightly larger model.  Dual to the construction we made in $\dgl$ we define $\hat \Omega \Lambda V$ by taking the coproduct $\Lambda s^{-1}V \tildeoplus \Lambda V \tildeoplus \Lambda s^{-1}V$ and modifying the differential:
$$\hat\Omega \Lambda V := \mathrm{red}_r\bigl([\Lambda s^{-1}V \tildeoplus \Lambda V \tildeoplus \Lambda s^{-1}V]_{\rm GC},\, d_{\hat \Omega} = d_{\tildeoplus} + \Lambda(d')\bigr)$$
where $\Lambda(d')$ is the differential on $[\Lambda s^{-1}V \tildeoplus \Lambda V \tildeoplus \Lambda s^{-1}V]_{\rm GC} \cofibr \Lambda (s^{-1}V\oplus V \oplus s^{-1}V)$ given by the restriction of the differential on $\Lambda (s^{-1}V\oplus V \oplus s^{-1}V)$ cofreely generated by $d'(0,v_1,0) = (s^{-1}v_1,0, s^{-1}v_1)$ on $(s^{-1}V\oplus V \oplus s^{-1}V)$.

More generally, if $C$ is any cofree \DGC, then we may define the unreduced loops and loops on $C$ to be the truly cofree \DGC\,s $\hat{\Omega} C := \Lambda s^{-1}(C)^\pr$ and $\Omega C := \Lambda\bigl(\mathrm{red}_rs^{-1}(C)^\pr\bigr)$.  

\begin{lemma}[Loops] 
A right adjoint of the functor $\Sigma:\dgc_n \to \dgc_{n+1}$ is given by the functor $$\tilde{\Omega}:C \mapsto \Lambda s^{-1}(C)^\pr$$

A right adjoint of the functor $\Sigma:\dgc_n \to \dgc_n$ is given by the functor $$\Omega:C\mapsto \Lambda\bigl(\mathrm{red}_rs^{-1}(C)^\pr\bigr)$$
\end{lemma}

\section{Model Category Structure}

Quasi-isomorphisms in $\dgc$ are \DGC-maps which are quasi-isomorphisms on the level of underlying \DG\,s.  The standard model category structure on $\dgc$ (see [Q69] or [GG] for a more general version than we need) is to take quasi-isomorphisms to be weak equivalences, degree-wise injections to be cofibrations, and allow fibrations to determined by right lifting with respect to acyclic cofibrations.  Under this model category structure, all objects are cofibrant.

\begin{thm} 
This gives a model catgegory structure on $\dgc_r$ $(r\ge 2)$.\footnote{It is possible to extend this model category structure to one on $\dgc_1$ and even $\dgc_r$ for $r\le 0$; however we do not need this, so we stick with the easily defined and proven case.}
\end{thm}
\begin{proof}[Proof sketch]
This is just as the proof of the corresponding theorem for the category $\dgl$,
the proof of this follows from Quillen's proof of [Q69, II.5.2] as well as the results of the previous section.  
\end{proof}

\begin{prop}
Fibrations in $\dgc$ are the cofree maps.
\end{prop}

\begin{cor}
The fibrant objects in $\dgc$ are precisely the cofree objects $\Fdgc$.
\end{cor}

The identity functor serves as a functorial fibrant cofibrant approximation.  $\dgc$ also has a functorial cofibrant fibrant approximation (given by $\C\L$) which is described in Chapter~\ref{S:Q homotopy} and in particular by Corollary~\ref{funct approx}.  The following lemma immediately follows from the definitions of $\C$ and $\L$ which is given in \ref{C functor} and \ref{L functor}:

\begin{lemma}\label{dgc fibr replacement is cofree}
The fibrant replacement functor in $\dgc_r$ is a functor $\C\L:\dgc_r\to\hfdgc_r$.

That is, fibrant replacement takes \DGC\,s to cofree \DGC\,s of the form $(\Lambda V,\, d)$ and \DGC-maps to cofreely generated maps of the form $\Lambda(\hat f):(\Lambda V,\, d) \to (\Lambda W,\, d')$.
\end{lemma}

Under this model category structure, adjoints which we gave in Lemmas~\ref{dgc red} and \ref{dgc adjoints} become Quillen adjoint pairs: 

\begin{lemma}
The adjoint pair given in Lemma~\ref{dgc red}
$$\xymatrix@C=20pt@R=5pt{
[-]_{\mDGC_t}:\dgc_r \ar@<2pt>[r] & \dgc_t:\mathrm{red}_r \ar@<2pt>[l] 
}$$
($t<r$) is a Quillen adjoint pair.
Also both of the adjoint pairs from Lemma~\ref{dgc adjoints}
$$\xymatrix@C=10pt@R=1pt{
[-]_\mDGC:\dg_r \ar@<2pt>[r] &
  \dgc_r:(-)^\pr\phantom{XXx}  \ar@<2pt>[l]  \\
{}\phantom{xx}[-]_\mDG :\dgc_r \ar@<2pt>[r] &
  \dg:\Lambda(-) \ar@<2pt>[l]
}$$  
are Quillen adjoint pairs. 
\end{lemma}
\begin{proof}
The left adjoints $[-]_{\mDGC_t}$, $[-]_\mDGC$, and $[-]_\mDG$ each preserve all cofibrations, since cofibrations are degree-wise injections in each of $\dg$, $\dg_r$ and $\dgc_r$.
Also, the left adjoints each preserve all weak equivalences because quasi-isomorphism in $\dgc_r$ is defined on the level of \DG\,s.  
Since the left adjoints in the above pairs each preserve all cofibrations and trivial cofibrations, the right adjoints must also preserve all fibrations and trivial fibrations and the pairs are Quillen adjoints as claimed
\end{proof}

Note that the top adjoint from Lemma~\ref{dgc adjoints} does not extend as a Quillen adjoint all the way to $\dg$.  This is because the adjoint pair $\mathrm{trunc}_r:\dg\rightleftarrows\dg_r:\mathrm{incl}$ is not a Quillen adjoint pair.  Lemma~\ref{adjoint holim} now provides us with the following corollary:

\begin{cor}\label{C:dgc holim plan}
The left adjoint functors $[-]_{\mDGC_t}:\dgc_r\to\dgc_t$, $[-]_\mDGC:\dg_r\to\dgc_r$, and $[-]_\mDG:\dgc_r\to\dg$ above preserve all weak equivalences and therefore all homotopy colimits.

The right adjoints $(-)^\pr:\Fdgc_r\to\dg_r$ and $\Lambda(-):\dg\to\Fdgc_r$ induced by the above preserve all weak equivalences and therefore all homotopy limits.
\end{cor}

\begin{proof}
It remains to show only that the right adjoints preserve all weak equivalences.  This is trivially true for $\mathrm{red}_r$.  For $\Lambda(-)$ and $(-)^\pr$ this follows from the rational Hurewicz Theorem (\ref{Hurewicz dgc}).
\end{proof}  

As in the previous chapter, we use these adjoints to construct models for homotopy limits and colimits in $\dgc_r$:

\subsection{Homotopy Limits and Colimits}

By [DHKS] all homotopy limits and colimits exist (see \ref{t:hocomplete}).  We use \ref{dgc factorization} and \ref{C:dgc holim plan} in order to construct nice models for homotopy colimits and limits in $\dgc_r$, dual to the way in which we created homotopy limits and colimits in $\dgl_r$.  Essentially the construction is as follows:  If $\mathscr{D}:\catI\to\dgc_r$ is a diagram in $\dgc_r$ then up to coproduct, the homotopy colimit of $\mathscr{D}$ in $\dgc_r$ is given by
$$\bigl[[\hocolim{\!}_\dgc\mathscr{D}]_\mDG\bigr]_\mDGC = \bigl[\hocolim{\!}_\dg[\mathscr{D}]_\mDG\bigr]_\mDGC$$
Similarly up to higher order differential information, the homotopy limit of $\mathscr{D}$ is given by
$$\Lambda(\holim{\!}_\dgc\mathscr{D})^\pr = \Lambda\bigl(\holim{\!}_\dg(\mathscr{D})^\pr\bigr)$$
To construct models for homotopy colimits or limits in $\dgc_r$ we insert our models for homotopy colimits and limits in $\dg_r$ (given in Section~\ref{dg_r model category}) into the above and then supply the missing coproduct or higher order differential information.

More precisely, we rely on a characterization of homotopy colimit and limit functors on $\dgc$ as in the previous chapter: 

\begin{thm}[Creation of Homotopy Limits and Colimits]\label{dgc holim plan}
Let $\hocolim_\dg$ and $\holim_\dg$ be any homotopy colimit and limit functors on $\dg$.
\begin{enumerate}
\item Suppose $F:(\dgc_r)^{\catI} \longrightarrow \dgc_r$ is a functor from $\catI$-diagrams in $\dgc_r$ to $\dgc_r$ such that for all $\catI$-diagrams $\mathscr{D}:\catI\to\dgc_r$, we have 
\begin{itemize}
\item $\bigl[F(\mathscr{D})\bigr]_\mDG = \hocolim_\dg[\mathscr{D}]_\mDG$.
\item $F$ is equipped with natural maps $e_F:F(\mathscr{D}) \to \colim_{\dgc}\mathscr{D}$. 
\item $[e_F]_\mDG$ is the canonical map $\hocolim_\dg[\mathscr{D}]_\mDG \to \colim_\dg[\mathscr{D}]_\mDG$.
\end{itemize}
Then $F$ is an $\catI$-homotopy colimit functor on $\dgc_r$.

\item Dually, suppose $G:(\dgc_r)^{\catI} \longrightarrow \Fdgc_r$ is a homotopy functor such that for all $\catI$-diagrams $\mathscr{D}:\catI\to\dgc_r$, we have 
\begin{itemize}
\item $\bigl(G(\mathscr{D})\bigr)^\pr = \holim_\dg(\mathscr{D})^\pr$. 
\item $G$ is equipped with natural maps $e_G:\lim_{\dgc}\mathscr{D} \to G(\mathscr{D})$.
\item $(e_G)^\pr$ is the canonical map $\lim_\dg(\mathscr{D})^\pr \to \holim_\dg(\mathscr{D})^\pr$.  
\end{itemize}
Then $G$ is an $\catI$-homotopy limit functor on $\dgc_r$.
\end{enumerate}
\end{thm}

\begin{proof}[Proof Sketch]
The proof of Theorem~\ref{dgl holim plan} from the previous chapter translates directly to this setting.  Our strategy is to make natural zig-zags of weak equivalences (dual to those of \ref{dgl holim plan}):
$$F(\mathscr{D})\xleftarrow{\ \weq\ } F(\mathscr{D}_\mathsf{vc}) \xrightarrow{\ \weq\ }
  \colim{\!}_\dgc\mathscr{D}_\mathsf{vc} = \hocolim{\!}_\dgc\mathscr{D}$$
$$G(\mathscr{D})\xrightarrow{\ \weq\ } G(\mathscr{D}_\mathsf{vf}) \xleftarrow{\ \weq\ }
  \lim{\!}_\dgc\mathscr{D}_\mathsf{vf} = \holim{\!}_\dgc\mathscr{D}$$
Where $\hocolim_\dgc$ and $\holim_\dgc$ are the homotopy colimit and limit functors constructed by [DHKS].

The first arrows in each of the above zig-zags are weak equivalences because they are induced by the virtually-cofibrant and virtually-fibrant replacement maps $\mathscr{D}_\mathsf{vc}\xrightarrow{\,\weq\,}\mathscr{D}$ and $\mathscr{D}\xrightarrow{\,\weq\,}\mathscr{D}_\mathsf{vf}$.  The second arrows of the above zig-zags are $e_F:F(\mathscr{D}_\mathsf{vc}) \xrightarrow{\ \ } \colim_\dgc\mathscr{D}_\mathsf{vc}$ and $e_G:\lim_\dgc\mathscr{D}_\mathsf{vf} \xrightarrow{\ \ } G(\mathscr{D}_\mathsf{vf})$.   The map $[e_F]_\mDG$ is a weak equivalence because $\colim_\dg[\mathscr{D}_\mathsf{vc}]_\mDG$ is a homotopy colimit functor; thus $e_F$ is a weak equivalence by \ref{weq dgc equals weq dg}.  The map $(e_G)^\pr$ is a weak equivalence because $\lim_\dg[\mathscr{D}_\mathsf{vf}]_\mDG$ is a homotopy limit functor.  Since both $G(\mathscr{D}_\mathsf{vf})$ and $\lim_\dgl\mathscr{D}_\mathsf{vf}$ are cofree, the rational Hurewicz Theorem (\ref{Hurewicz dgc}) forces $e_G$ to be a weak equivalence as well.
\end{proof}

\subsubsection{Homotopy Colimits in $\dgc_r$}

\begin{lemma} 
If $\mathcal{D}$ is the \DGC-diagram $B_1\xleftarrow{f_1}C\xrightarrow{f_2}B_2$ then its homotopy pushout is given by the two-sided mapping cylinder $\mathcal{C}_\mathcal{D}$ defined by
$$\widetilde{\mathcal{C}}_\mathcal{D} := \bigl((\widetilde{B_1}\oplus s\widetilde{C}\oplus \widetilde{B_2})_\bullet,\, d_\mathcal{C},\, \widetilde{\Delta}_\mathcal{C})$$
where $d_\mathcal{C}(b_1+sc+b_2) = \big(d_{B_1}b_1 + f_1(c)\big) - sd_{C}c + \big(d_{B_2}b_2 + f_2(c)\big)$ and $\widetilde{\Delta}_\mathcal{C}(b_i) = \widetilde{\Delta}_{B_i}(b_i)$ and
\begin{align*}
\widetilde{\Delta}_\mathcal{C}(sc) &= \frac{1}{2}\sum_j \big(s\alpha_j\otimes f_1(\beta_j) + (-1)^{|a_j|}f_1(\alpha_j)\otimes s\beta_j\big) \\
&\quad + \frac{1}{2}\sum_j \big(s\alpha_j\otimes f_2(\beta_j) + (-1)^{|\alpha_j|}f_2(\alpha_j)\otimes s\beta_j\big)
\end{align*}
for $b_i\in\widetilde{B_i}$, $c\in \widetilde{C}$, and $\widetilde{\Delta}_{C}(c) = \sum_j \alpha_j\otimes \beta_j$.
\end{lemma}
Note that $\mathcal{C}_\mathcal{D}$ is merely the homotopy pushout of $\mathcal{D}$ as a \DG-diagram ($\mathcal{C}_{[\mathcal{D}]_\mDG}$) equipped with a certain coproduct.  Thus by~\ref{dgc holim plan} to show that it is the homotopy pushout in the category $\dgc$, all that remains is to show that $\mathcal{C}_\mathcal{D}$ is indeed a \DGC\ and that it maps correctly to the colimit $B_1\oplus_C B_2$.  However, the composition
$\mathcal{C}_\mathcal{D} \to B_1 \oplus B_2 \to B_1\oplus_C B_2$ clearly gives the desired natural map.  

\begin{proof}[Proof that this is a \DGC] 
From Lemma~\ref{l:dg pushout} we know that $d_\mathcal{C} \circ d_\mathcal{C} = 0$.  It remains to show that the coproduct is cocommutative and compatible with the differential.  By definition, $\widetilde{\Delta}_\mathcal{C}$ is cocommutative and compatible with $d_\mathcal{C}$ on ${B_1}$ and ${B_2}$.  Note that for $c\in \widetilde{\mathcal{C}}_\mathcal{D}$, 
$$\widetilde{\Delta}_\mathcal{C}(c)\in \big((s{C}\otimes {B_1}) \oplus ({B_1}\otimes {B_1}) \oplus ({B_1}\otimes s{C})\big) \oplus \big((s{C}\otimes {B_2})\oplus ({B_2}\otimes {B_2}) \oplus ({B_2}\otimes s{C})\big)$$  
To simplify our computations we write $\widetilde{\Delta}_\mathcal{C}(c)$ as $\widetilde{\Delta}_\mathcal{C}^{B_1}(c) + \widetilde{\Delta}_\mathcal{C}^{B_2}(c)$ where 
$\widetilde{\Delta}_\mathcal{C}^{B_1}(c)\in (s{C}\otimes {B_1}) \oplus ({B_1}\otimes s{C})\oplus({B_1}\otimes {B_1})$ and $\widetilde{\Delta}_\mathcal{C}^{B_2}(c)\in (s{C}\otimes {B_2}) \oplus ({B_2}\otimes s{C})\oplus ({B_2}\otimes {B_2})$.  To show that $\widetilde{\Delta}_\mathcal{C}$ is cocommutative on $s{C}$, we show that the $\widetilde{\Delta}_\mathcal{C}^{B_i}$ are cocommutative on $s{C}$:  

Let $c \in \bar C$ with $\bar \Delta_{C} c = \sum_j \alpha_j\otimes \beta_j$.  Then
\begin{align*}T\widetilde{\Delta}_\mathcal{C}^{B_i}(sc) 
&= \frac{1}{2}T\sum_j\, \Bigl[ s\alpha_j \otimes f_i(\beta_j) + (-1)^{|\alpha_j|}f_i(\alpha_j)\otimes s\beta_j \Bigr] \\
&= \frac{1}{2}\sum_j\, \Bigl[(-1)^{|\beta_j|(|\alpha_j| +1)}f_i(\beta_j)\otimes s\alpha_j 
  + (-1)^{|\alpha_j|}(-1)^{|\alpha_j|(|\beta_j| +1)} s\beta_j\otimes f_i(\alpha_j)\Bigr] 
\end{align*}
Also, $\widetilde{\Delta}_{C}$ is cocommutative, so $\widetilde{\Delta}_{C} c = T\widetilde{\Delta}_{C} c = \sum_j (-1)^{|\alpha_j|\cdot|\beta_j|}\beta_j\otimes \alpha_j$.  Thus
\begin{align*}
 \widetilde{\Delta}_\mathcal{C}^{B_i}(sc) 
&= \frac{1}{2}\sum_j\, \Bigl[ (-1)^{|\alpha_j|\cdot|\beta_j|}s\beta_j\otimes f_i(\alpha_j) 
  + (-1)^{|\beta_j|}(-1)^{|\alpha_j|\cdot|\beta_j|}f_i(\beta_j)\otimes s\alpha_j\Bigr] 
\end{align*}
Therefore the $\widetilde{\Delta}_\mathcal{C}^{B_i}$ are cocommutative as claimed.

In order to show compatibility of the differential with the coproduct we display compatibility with $\widetilde{\Delta}_\mathcal{C}^{B_i}$:
Let $c \in \widetilde{C}$ with $\widetilde{\Delta}_{C} c = \sum_j \alpha_j\otimes \beta_j$.  Then 
$\widetilde{\Delta}_{C} d_{C}c = \sum_j d_{C}\alpha_j \otimes \beta_j + (-1)^{|\alpha_j|} \alpha_j \otimes d_{C} \beta_j$ and so
\begin{align}{}
\widetilde{\Delta}_\mathcal{C}^{B_i} d_\mathcal{C}(sc) &= \widetilde{\Delta}_\mathcal{C}^{B_i}\bigl(f_1(c) + f_2(c) - sd_{C}c\bigr) \notag \\
 &= \sum_j  f_i(\alpha_j) \otimes f_i(\beta_j)  - \frac{1}{2} \sum_j\, \Bigl[
    sd_{C}\alpha_j \otimes f_i(\beta_j) + (-1)^{|\alpha_j|}s\alpha_j\otimes f_i(d_{C}\beta_j) \notag \\
 &\hskip 130pt
   + f_i(\alpha_j)\otimes sd_{C}\beta_j  - (-1)^{|\alpha_j|}f_i(d_{C}\alpha_j)\otimes s\beta_j\Bigr]  \notag
\end{align}
\begin{align}
d_{\otimes}\widetilde{\Delta}_\mathcal{C}^{B_i}(sc) 
 &= d_{\otimes}\,\frac{1}{2}\sum_j\, s\alpha_j \otimes f_i(\beta_j) + (-1)^{|\alpha_j|}f_i(\alpha_j)\otimes s\beta_j \notag \\
 &= \frac{1}{2}\sum_j\, \Bigl[f_i(\alpha_j)\otimes f_i(\beta_j) - sd_{C}\alpha_j \otimes f_i(\beta_j) 
   - (-1)^{|\alpha_j|}s\alpha_j\otimes d_{B_i}f_i(\beta_j) \notag \\
 & \hskip 40pt + (-1)^{|\alpha_j|}d_{B_i}f_i(\alpha_j)\otimes s\beta_j + f_i(\alpha_j)\otimes f_i(\beta_j) 
   - f_i(\alpha_j)\otimes sd_{C}\beta_j\Bigr] \notag \\
 &= \sum_j\,  f_i(\alpha_j) \otimes f_i(\beta_j)  - \frac{1}{2} \sum_j\, \Bigl[
    sd_{C}\alpha_j \otimes f_i(\beta_j) + (-1)^{|\alpha_j|}s\alpha_j\otimes f_i(d_{C}\beta_j) \notag \\
 &\hskip 130pt
   + f_i(\alpha_j)\otimes sd_{C}\beta_j  - (-1)^{|\alpha_j|}f_i(d_{C}\alpha_j)\otimes s\beta_j\Bigr]  \notag
\end{align}

\end{proof}

\begin{cor} $\Sigma C$ is the homotopy pushout in $\dgc_r$ of the diagram 
$0_\mDGC \xleftarrow{\ } C \xrightarrow{\ } 0_\mDGC$.
\end{cor}

\begin{cor}\label{p: dgc hofib}
If $f:C \to B$ is a \DGC-map then its homotopy cofiber may be modelled by the mapping cone: 
$$\hocof(f) = \bigl((s\widetilde{C} \oplus B)_\bullet,\, d_f,\, \Delta_f,\, \varepsilon_B\bigr)$$
where $d_f$ is defined by $d_f(b) = d_B(b)$, and $d_f(sc) = f(c) - sd_C(c)$; and $\Delta_f$ is defined by $$\widetilde{\Delta}_f(sc) = \frac{1}{2}\sum_j \big(s\alpha_j\otimes f(\beta_j) + (-1)^{|\alpha_j|}f(\alpha_j)\otimes s\beta_j\big)$$ and $\Delta_f(b) = \Delta_B(b)$ for $c \in \widetilde{C}$ with $\widetilde{\Delta}_C(c) = \sum_j \alpha_j\otimes \beta_j$ and $b \in B$.  
\end{cor}

These constructions generalize in the obvious manner to define higher dimensional homotopy pushouts or indeed any homotopy colimit in $\dgc_r$.  However, for our purposes all that we explicitly require a model for is the above homotopy pushout.

\subsubsection{Homotopy Limits in $\dgc_r$}

Dual to our construction of homotopy colimits in $\dgl_r$, in order to build models for homotopy limits in $\dgc_r$, we first construct homotopy limits (in $\dgc_r$) of diagrams in $\hfdgc_r$ and then use these to build homotopy limits of general diagrams.  We call diagrams in the image of $(\fdgc_r)^\catI \cofibr (\dgc_r)^\catI$ {\em cofreely generated diagrams}.  
By Lemma~\ref{dgc fibr replacement is cofree} the fibrant replacement functor $(\C\L)$ on $\dgc_r$ takes all diagrams to cofreely generated diagrams.  The homotopy limit of a cofreely generated diagram is given by the \DGC\ cofreely generated by the homotopy limit (in $\dg_r$) of the underlying diagram in $\dg_r$, with some extra higher order differential structure.

Writing $\widehat{\holim}_\hfdgc:(\fdgc_r)^\catI\to\Fdgc$ for our functor giving homotopy colimits of cofreely generated diagrams in $\dgc_r$, an $\catI$-homotopy limit functor in $\dgc_r$ is given by the composition
$$\holim{\!}_\dgc(-) := \widehat{\holim}_\hfdgc\bigl(\C\L(-)\bigr)$$
This is our model for general homotopy limits in $\dgc_r$; however, in practice almost all of the diagrams which we are concerned with are already be cofreely generated.  In these cases, we use the simpler model $\widehat{\holim}_\hfdgl$ for their homotopy limit -- note that for cofreely generated diagrams $\mathscr{D}$ there is a natural weak equivalence 
$$\widehat{\holim}_\hfdgc \mathscr{D} \xrightarrow{\,\weq\,}
  \widehat{\holim}_\hfdgc \C\L(\mathscr{D}) = 
  \holim{\!}_\dgc \mathscr{D}$$
induced by the natural weak equivalence $\mathscr{D} \xrightarrow{\,\weq\,} \C\L\mathscr{D}$.

From the definition of $\hfdgc$ it follows that an $\catI$-diagram $\mathscr{D}:\catI\to\dgc_r$ is cofreely generated if and only if, as a diagram of graded coalgebras, it is given by the cofree extension of a diagram of graded vector spaces:
$$[\mathscr{D}]_\mGC = \Lambda(\hat{\mathscr{D}})\qquad \text{where}\qquad
  \hat{\mathscr{D}}:\catI\to \g_r$$
The most basic and important diagrams which we consider are thus cofreely generated:
\begin{ex}Let $(\Lambda V,\, d)$ be a cofree \DGC.  The following diagrams are both cofreely generated.
\begin{enumerate}
\item $0_\mDGC \xrightarrow{\ \ } (\Lambda V,\, d) \xleftarrow{\ \ } 0_\mDGC$
\item $(\Lambda V,\, d) \xrightarrow{\,\Idm\,} (\Lambda V,\, d) \xleftarrow{\ \ } 0_\mDGC$
\end{enumerate}
\end{ex}
We define homotopy limits of cofreely generated diagrams so that the homotopy limit of (1) is $\Omega(\Lambda V,\, d)$ and the homotopy limit of (2) is $\mathrm{p}(\Lambda V\, d)$ as defined in \ref{dgc paths and loops}. 

We begin with simple diagrams:

Suppose $\mathscr{D}$ is the cofreely generated pushout diagram $(\Lambda U,\, d') \xrightarrow{\,\Lambda(f)\,} (\Lambda V,\, d'') \xleftarrow{\,\Lambda(g)\,} (\Lambda W,\, d''')$ in $\dgc_r$ so that $\hat{\mathscr{D}}$ is the diagram $U\xrightarrow{\,f\,}V\xleftarrow{\,g\,}W$ in $\g_r$.  Note that $(\mathscr{D})^\pr$ is the diagram in $\dg_r$ given by $(U,\,d'_U) \xrightarrow{\,f\,} (V,\,d''_V)\xleftarrow{\,g\,}(W,\,d'''_W)$ and $\bigl[(\mathscr{D})^\pr\bigr]_\mG = \hat{\mathscr{D}}$.

Let $\mathscr{P}_\mathscr{D}$ be the cofree path \DGC\ given by:
$$\mathscr{P}_\mathscr{D} := 
 \mathrm{red}_r\bigl(\Lambda U \tildeoplus \Lambda s^{-1} V \tildeoplus \Lambda W,\
 d_\mathscr{P} = d_{\!\tildeoplus\!} + \Lambda(d_{fg})\bigr)$$
where 
\begin{itemize}
\item $d_{\!\tildeoplus\!}$ is the differential on 
  $(\Lambda U,\,d')\tildeoplus \Lambda s^{-1}(V,\,d''_V) \tildeoplus (\Lambda W,\, d''')$.
\item $\Lambda(d_{fg})$ is the differential on the \GC\ 
 $\Lambda U \tildeoplus \Lambda s^{-1} V \tildeoplus \Lambda W \cofibr \Lambda (U \oplus s^{-1}V \oplus W)$ given by the restriction of the differential on $\Lambda (U \oplus s^{-1}V \oplus W)$ cofreely generated by the differential on $(U \oplus s^{-1}V\oplus W)$ given by $d_{fg}(u + s^{-1}v + w) = s^{-1}\bigl(f(u) + g(w)\bigr)$ for $u\in U$, $v\in V$, $w\in W$. 
\end{itemize}
Since the functor $\Lambda(-)$ gives a section of $(-)^\pr$, we get that
$$(\mathscr{P}_\mathscr{D})^\pr = \mathrm{red}_r
 \Bigl(\bigl[(U,\, d'_U) \oplus s^{-1}(V,\, d''_V) \oplus (W,\, d'''_W)\bigr]_\mG,\ d = d_{\!\tildeoplus\!} + d_{fg}\Bigr) = \holim{\!}_{\dg_r}(\mathscr{D})^\pr$$
Furthermore we have a natural map $\lim\mathscr{D} \to \mathscr{P}_\mathscr{D}$ given by the composition
$$\lim{\!}_\dgc \mathscr{D} \longrightarrow 
 (\Lambda U,\, d')\tildeoplus (\Lambda W,\, d''') \longrightarrow \mathscr{P}_\mathscr{D}$$
which descends to the standard natural maps between $\lim$ and $\holim$ on $\dg_r$:
$$\lim{\!}_{\dg_r}(\mathscr{D})^\pr \longrightarrow (U,\,d'_U)\oplus(W,\,d'''_W) 
 \longrightarrow \colim{\!}_{\dg_r}(\mathscr{D})^\pr$$
as described in \ref{dg ho(co)lim natural maps}.  Also $\mathscr{P}_{(-)}$ is a homotopy functor on diagrams in $\fdgc$ because $i:\mathscr{D}_1\xrightarrow{\,\weq\,}\mathscr{D}_2$ a weak equivalence implies $\holim_{\dg_r}(i)^\pr = \bigl(\mathscr{P}_{(i)}\bigr)^\pr$ is a weak equivalence and so $\mathscr{P}_{(i)}$ is a weak equivalence by the rational Hurewicz theorem.

By Theorem~\ref{dgc holim plan} we therefore have
\begin{lemma}
The composition $\mathscr{P}_{\C\L(-)}$ is a homotopy limit functor on $\dgc_r$.

In particular, if $\mathscr{D}$ is a cofreely generated diagram in $\dgc_r$ then the \DGC\ 
$\mathscr{P}_\mathscr{D} \weq \mathscr{P}_{\C\L\mathscr{D}}$ is a model for the homotopy limit (in $\dgc_r$) of $\mathscr{D}$.
\end{lemma}

\begin{cor}
The following functors (defined earlier) are models for hoomtopy limits
\begin{itemize}
\item $\Omega(\Lambda V,\, d) = \holim_\dgc\Bigl(0_\mDGC \xleftarrow{\ \ } 
 (\Lambda V,\, d) \xrightarrow{\ \ } 0_\mDGC\Bigr)$.
\item $\mathrm{p}(\Lambda V,\, d) = \holim_\dgc\Bigl(
 (\Lambda V,\, d) \xrightarrow{\,\Idm\,} (\Lambda V,\, d) \xleftarrow{\ \ }
 0_\mDGC\Bigr)$.
\end{itemize}
\end{cor}

\begin{note}
$\mathscr{P}_\mathscr{D}$ could also be written as 
$$\mathscr{P}_\mathscr{D} = \Bigl(\bigl[\Lambda\, \mathrm{red}_r\bigr(
 (U,\,d'_U)\oplus s^{-1}(V,\,d''_V)\oplus (W,\,d'''_W)\bigr)\bigr]_\mGC,\ 
 d_\mathscr{P} = d_{\Lambda} + d_{>0}\Bigr)$$
where
\begin{itemize}
\item $d_{\Lambda}$ is the differential on the truly cofree \DGC\ 
 $\Lambda\,\mathrm{red}_r\bigr(
  (U,\,d'_U)\oplus s^{-1}(V,\,d''_V)\oplus (W,\,d'''_W)\bigr)$.
\item $d_{>0}$ is the degree $-1$ map on $\Lambda(U\oplus s^{-1}V\oplus W) = \Lambda U \tildeoplus \Lambda s^{-1}V \tildeoplus \Lambda W$ generated by $d'-d'_U$ on $\Lambda U$ and $d'''-d'''_W$ on $\Lambda W$.
\end{itemize}
\end{note}

More generally, given $\mathscr{D}:\catI\to \fdgc_r$ a cofreely generated diagram, we could define the homotopy limit of $\mathscr{D}$ along the lines of
$$\widehat{\holim}_\hfdgc \mathscr{D} := 
 \Bigl(\bigl[\Lambda\bigl(\holim{\!}_{\dg_r}(\mathscr{D})^\pr\bigr)\bigr]_\mGC,\ 
 d = d_{\Lambda} + d_{>0}\Bigr)$$
where $d_\Lambda$ is the differential on the truly cofree \DGC\ $\Lambda\bigl(\holim_{\dg_r}(\mathscr{D})^\pr\bigr)$ and $d_{>0}$ is the degree $-1$ map adding back in all of the higher order differential information (as in the remark above) so that there is a natural map $\lim_\dgc\mathscr{D} \to \widehat{\holim}_\hfdgc\mathscr{D}$ (recall that $\holim_{\dg_r}$ of a diagram is a large direct sum (with modified differential) of the \DG\,s in the diagram and their iterated desuspensions).  We will not be more specific here because we will not need explicit models for the homotopy limits of complicated diagrams.  The only remaining fact which we need is:

\begin{thm} 
In $\dgc_r$, sequential homotopy colimits commute with homotopy pullbacks.
\end{thm}
\begin{proof}[Proof Sketch]

It is enough to show this for cofreely generated diagrams $\mathscr{D}:\catI\times \catJ \to \hfdgc$, since fibrant replacement gives a weak equivalence $\mathscr{D} \xleftarrow{\ \weq\ }\C\L(\mathscr{D})$ natural in $\mathscr{D}$ and $\C\L(\mathscr{D})$ is cofreely generated.

For cofreely generated diagrams, however, this may be proven by a large computation just as the corresponding theorem for $\dg_r$ (\ref{dg_r holim commute}).  The critical fact is that infinite $\tildeoplus$'s of cofree \DGC\,s commute with finite $\otimes$.
 
Kuhn also comments in [Kuhn p6] that this may be shown using fact that the sequential small object argument applies in $\dgc$.

\end{proof}

\chapter{Rational Homotopy Theory}\label{S:Q homotopy}
\markboth{Ben Walter}{I.6 Rational Homotopy Theory}

\section{Framework}

We give a brief summary of the objects and structures of rational homotopy theory which we use.  For a more detailed discussion see [FHT], [GM], or [Q69].

In [Q69] Quillen constructed a chain of equivalences from simply connected topological spaces ($\top_2$), to 2-reduced simplicial sets ($\sset_2$), to 1-reduced simplicial groups ($\sgp_1$), to 1-reduced simplicial complete Hopf algebras (over $\Q$) ($\scha_1$), to 1-reduced simplicial Lie algebras (over $\Q$) ($\slie_1$), to 1-reduced \DGL\,s ($\dgl_1$), to 2-reduced \DGC\,s ($\dgc_2$):\footnote{In order to write this chain -- and for only this equation -- we break from our standard convention and adopt Quillen's convention of writing left adjoints above right adjoints -- i.e. ``$\C:\dgl_1 \leftrightarrows \dgc_2:\L$'' means $\C$ is {\em right adjoint} to $\L$ (even though it is written on the left).}
\begin{equation}
\xymatrix@=20pt{
 \top_2 \ar@<-.5ex>[r]_{\mathit{Sing}_2} 
& \sset_2 \ar@<-.5ex>[l]_{|\cdot|} \ar@<.5ex>[r]^{\Omega} 
& \sgp_1 \ar@<.5ex>[l]^{W} \ar@<-.5ex>[r]_{\widehat{\Q}}
& \scha_1 \ar@<-.5ex>[l]_{\mathcal{G}} \ar@<.5ex>[r]^{\mathcal{P}}
& \slie_1 \ar@<.5ex>[l]^{\widehat{U}} \ar@<-.5ex>[r]_N
& \dgl_1 \ar@<-.5ex>[l]_{K} \ar@<-.5ex>[r]_{\C}
& \dgc_2 \ar@<-.5ex>[l]_{\L} }\label{D:Quillen}
\end{equation}
where $|\cdot|$ denotes the geometric realization of a simplicial set, $\mathit{Sing}_2$ denotes the singular complex of a 1-connected space with trivial 1-skeleton; $\Omega$ denotes Kan loop-group functor, $W$ is Kan's classifying simplicial set functor; $\mathcal{G}$ is the grouplike elements of a $\scha$, $\widehat{\Q}$ is the completion of the group ring over $\Q$ at the powers of its augmentation ideal; $\mathcal{P}$ is the primitive elements of a Hopf algebra, $\widehat{U}$ is the completion of the universal enveloping algebra of a Lie algebra at powers of its augmentation ideal; $K$ is (Dold-Kan's) inverse to the normalization map $N$ taking a simplicial object to a differential graded object; and $\C$ and $\L$ are the maps defined below.

Given a space, we write $\mathfrak{C}(X)$ for ``the \DGC\ corresponding to $X$'' and $\mathfrak{L}(X)$ for ``the \DGL\ corresponding to $X$''.  By this we mean the compositions of the functors $\top_2$ to $\dgl_1$ and $\dgc_2$ given in (\ref{D:Quillen}).

The following definitions and lemmas are (with slight modification) taken from [FHT]:
 
\begin{defn}[\lbrack FHT, \S 22.b\rbrack]\label{C functor}
The functor $\C:\dgl_{r-1} \to \dgc_r$ is given by the Cartan-Eilenberg-Chevalley construction
 $$\C:(L,\, d_L,\, [-,-]_L) \longmapsto \bigl(\bigl[\Lambda s[L]_{\mDG}\bigr]_{\mGC},\, d=d_\Lambda + d_{[-,-]}\bigr)$$
where $d_\Lambda$ is the differential on $\Lambda s[L]_\mDG$, and $d_{[-,-]}$ is a differential coming from the Lie bracket of $L$.\footnote{$d_\Lambda$ does not change word-length; $d_{[-,-]}$ decreases word-length by 1.}  More specifically:
\begin{align}
d_\Lambda(sx_1\wedge \dots \wedge sx_k) &= -\sum_i (-1)^{n_i} sx_1 \wedge \dots \wedge sd_Lx_i \wedge \dots \wedge sx_k \notag \\
d_{[-,-]}(sx_1\wedge \dots \wedge sx_k) &= \sum_{i<j} (-1)^{n_{ij} + |x_i|} s[x_i,x_j]\wedge sx_1\wedge \dots \widehat{sx_i} \dots \widehat{sx_j} \dots \wedge sx_k \notag
\end{align}
(In the formulas above, $n_i = \sum_{j<i} |sx_j|$ the Koszul sign incurred by moving $d_L$ to the $sx_j$ term and past the $s$; and $n_{ij}$ is the sign change incurred by moving $sx_i\wedge sx_j$ to the front of $sx_1\wedge\dots\wedge sx_k$.)
\end{defn}

\begin{defn}[\lbrack FHT, \S 22.e\rbrack]\label{L functor}
The functor $\L:\dgc_r \to \dgl_{r-1}$ is given by 
 $$\L:(C,\, d_C,\, \Delta_C,\,\varepsilon_C) \longmapsto \bigl(\bigl[\freeL_{s^{-1}[C]_\mDG}\bigr]_\mGL,\, d=d_\freeL + d_\Delta)$$  where $d_\freeL$ is the differential on $\freeL_{s^{-1}[C]_\mDG}$ and $d_\Delta$ is a differential coming from the coproduct of $C$.\footnote{$d_\freeL$ does not change bracket-length; $d_\Delta$ increases bracket-length by 1.}   (Recall that our convention is for the functor $[-]_\mDG:\dgc_r \to \dg$ to {\em de-augment} and then forget coproduct structure.)  More specifically $d_\freeL$ and $d_\Delta$ are the free extensions of the following differentials:
\begin{alignat*}{5}
d_\freeL(s^{-1}\widetilde{x}) &= - s^{-1}d_C\widetilde{x}, && \text{where } \widetilde{x} \in \widetilde{C} \\
d_\Delta(s^{-1}\widetilde{x}) &= \frac{1}{2}\sum_i (-1)^{|\alpha_i|}[s^{-1}\alpha_i, s^{-1}\beta_i],  \qquad && \text{where } \ \widetilde{\Delta}_C\widetilde{x} = \sum_i \alpha_i \otimes \beta_i
\end{alignat*}
\end{defn}

Note that $\L$ maps to $\hfdgl$ and $\C$ maps to $\hfdgc$.

\begin{thm}[\lbrack Q69, II.5.3, B.7.5\rbrack]
The functors $$\L:\dgc_r \rightleftarrows \dgl_{r-1}:\C$$ define a Quillen equivalence for $r\ge 2$.  

Furthermore, $\L$ and $\C$ also satisify the stronger property that they preserve all weak equivalences.\end{thm}

This theorem has a number of immediate corollaries.

\begin{cor}
The functors $\C$ and $\L$ both detect and reflect weak equivalences of free and cofree objects respectively.

i.e. $f:K\xrightarrow{\,\weq\,} L$ iff $\C(f):\C(K)\xrightarrow{\,\weq\,} \C(L)$ and 
     $g:B\xrightarrow{\,\weq\,} C$ iff $\L(g):\L(B)\xrightarrow{\,\weq\,} \L(C)$ where $K,L\in\Fdgl_{r-1}$ and $B,C\in \Fdgc_r$ ($r\ge 2$).
\end{cor}

\begin{cor} 
The left adjoint $\L:\dgc_r \to \dgl_{r-1}$ preserves all homotopy colimits.  

The right adjoint $\C:\dgl_{r-1} \to \dgc_r$ preserves all homotopy limits.
\end{cor}

\begin{cor}
There are natural weak equivalences $C\xrightarrow{\weq}\C\L C$ and $\L \C L \xrightarrow{\weq} L$ for $C \in \dgc_r$ and $L\in \dgl_{r-1}$.
\end{cor}

Since $\L$ and $\C$ map to cofree \DGC\,s and free \DGL\,s respectively, the above natural weak equivalences serve as fibrant and cofibrant replacement functors in $\dgc_r$ and $\dgl_{r-1}$.  Also note that $\L$ (or $\C$) of any \DGC\ (or \DGL)-map is a (co)freely generated map in $\hfdgl$ (or $\hfdgc$).

\begin{cor}\label{funct approx} $\C\L$ gives functorial fibrant replacements in $\dgc_r$.  $\L\C$ gives functorial cofibrant replacements in $\dgl_{r-1}$.  In fact they respectively give cofibrant-fibrant and fibrant-cofibrant replacements.
\end{cor}

On free \DGL\,s and cofree \DGC\,s the functors $\C$ and $\L$ have a simpler form (up to Lie bracket and coproduct information):

\begin{lemma}[\lbrack FHT, 22.8\rbrack]\label{free homotopy} There are natural weak equivalences of \DG\,s 
\begin{align*}
[\C(\freeL_V, d)]_\mDG &\xrightarrow{\,\weq\,} s(V,\,d_V) \\
[\L(\Lambda V,\, d)]_{\mDG} &\xleftarrow{\,\weq\,} s^{-1}(V,\, d_V)
\end{align*}

More generally if $L$ is any free \DGL\ and $C$ is any cofree \DGC\ then there are natural weak equivalences: 
\begin{align*}
[\C L]_\mDG &\xrightarrow{\,\weq\,} s(L)^\ab \\
[\L C]_\mDG &\xleftarrow{\,\weq\,} s^{-1}(C)^\pr
\end{align*}
\end{lemma}


\begin{lemma}[\lbrack FHT 22.5\rbrack]
Let $f:L\to K$ be a map of free \DGL\,s.  Then $f$ is a quasi-isomorphism if and only if $\C(f)$ is a quasi-isomorphism -- i.e. $\C$ both detects and reflects quasi-isomorphisms of free \DGL\,s.
\end{lemma}

\begin{cor}
Let $g:B\to C$ be a map of cofree \DGC\,s.  Then $g$ is a quasi-isomorphism if and only if $\L(g)$ is a quasi-isomorphism -- i.e. $\L$ both detects and reflects quasi-isomorphisms of cofree \DGC\,s.
\end{cor}

Finally, note that the functors $\L$ and $\C$ are essentially nothing more than the bar and cobar constructions:

\begin{lemma}[\lbrack FHT, 22.7\rbrack] Given a \DGL\ $L$, there is a quasi-isomorphism of \DGC\,s 
$$\C L \xrightarrow{\,\weq\,} BUL$$ 
(where $BUL$ is the bar construction on the universal enveloping algebra of L).
\end{lemma}

\begin{lemma}[\lbrack Q69,  B.6.6\rbrack] Given a \DGC\ $C$, $\L C$ is the \DGL\ of primitive elements of the cobar construction\footnote{The cobar construction may be written as $\Omega(C,\, d_C,\, \Delta_C) = (Ts^{-1}\widetilde{C},\, d_0+d_1)$.  Recall that $\freeL_V$ is defined as a sub-Lie algebra of $TV$.} on $C$. 
\end{lemma}

Given a simply connected space $X$ we may consider  $\mathfrak{L}(X)$ the \DGL\ corresponding to $X$ under (\ref{D:Quillen}), and $\mathfrak{C}(X)$ the \DGC\ corresponding to $X$ under (\ref{D:Quillen}).  The central result of rational homotopy theory is that $\mathfrak{L}(X)$ gives the rational homotopy of $X$ (with Lie brackets corresponding to Whitehead products) and $\mathfrak{C}(X)$ gives the rational homology:

\begin{thm}[\lbrack Q69, Theorem I\rbrack]\label{Quillen Thm I}
There are natural isomorphisms $\pi_\ast(\Omega X)\otimes \Q \cong H_\ast\bigl(\mathfrak{L}(X)\bigr)$ (as graded Lie algebras) and 
$H_\ast(X;\, \Q) \cong H_\ast\bigl(\mathfrak{C}(X)\bigr)$ (as graded coalgebras).
\end{thm}

For this reason we sometimes say ``homotopy of'' $L$ (a \DGL) or $C$ (a \DGC) to mean $H_\ast(L)$ or $H_\ast(\L C)$ respectively and ``homology of'' $L$ or $C$ to mean $H_\ast(\C L)$ or $H_\ast(C)$ respectively.  Along with Lemma~\ref{free homotopy}, this is why Theorems~\ref{Hurewicz dgl} and~\ref{Hurewicz dgc} in the previous chapters were called ``rational Hurewicz'' theorems.

\section{Rational Spectra}



There is already a rich theory regarding the construction of spectra in general model categories (see for example [Ho01] and [S97]).
Following these, we could construct rational spectra by considering the stabilization of the categories $\dgl_{r-1}$ and $\dgc_r$.  Definitions and constructions would follow analogous to the classical definitions and constructions of $\sp$ as the stabilization of $\top$.  
To begin, recall the classical definition of a spectrum -- a spectrum $E$ is a sequence of spaces $\{E_i\}_{i\ge 0}$ equipped with maps $\Sigma E_i \to E_{i+1}$.  Analogously we may make the definitions:

\begin{defn}A {\em $\dgl$-spectrum} $E$ is a sequence of \DGL\,s $\{L_i\}_{i\ge 0}$ equipped with maps $$L_i \xrightarrow{} \Omega L_{i+1} = \Bigl[\mathrm{red}_r\bigl({ s^{-1} [L_{i+1}]_\mDG}\bigr)\Bigr]_\mDGL$$
\end{defn}

\begin{defn}A {\em $\dgc$-spectrum} $E$ is a sequence of \DGC\,s $\{C_i\}_{i \ge 0}$ equipped with maps $$\bigl[s [C_i]_\mDG\bigr]_\mDGC = \Sigma C_i \xrightarrow{} C_{i+1}$$
\end{defn}

Maps of $\dgl$-spectra and $\dgc$-spectra would then be defined in a more complicated manner.  

By the work of Schwede and Shipley the stabilizations of good simplicially enriched model categories which are Quillen equivalent are isomorphic.  So there is an isomorphism between the categories of $\dgl$-spectra and $\dgc$-spectra.  In fact, examining the categories of $\dgl$-spectra and $\dgc$-spectra defined above, it soon becomes apparent that they are both isomorphic to the category $\dg$.  The proof of this, however, is a bit unpleasant, and it seems unnecessarily complicated to carefully introduce and work seriously with either of the categories of $\dgl$-spectra or $\dgc$-spectra at all; when we would then immediately prove they are isomorphic to $\dg$ and subsequently ignore the more complicated categories.

\begin{note}In fact, Tom Goodwillie points out that there are functors between the categories of $\dg$ and rational spectra which are inverses up to zig-zags of natural weak equivalences.  By ``the category of rational spectra''  we mean the Bousfield localization of the category spectra with respect to the homology theory $H\Q = H_*(-;\,\Q)$.  
More precisely there are functors
$$\dg \xrightarrow{\,K\,} L_{H\Q}\,\sp \qquad \text{and}\qquad
 L_{H\Q}\,\sp \xrightarrow{\,C\,} \dg$$
where $C$ gives the stable singular chains of a rational spectrum, and $K$ creates an Eilenberg-MacLane spectrum from the given \DG.  It is clear that the composition of these two in either order is weakly equivalent to the identity.
\end{note}

Instead of working with categories of $\dgl$-spectra and $\dgc$-spectra which have complicated definitions, we merely fiat that the category of rational spectra is the category $\dg$.  We then exhibit Quillen adjoint functors $\Sinf$ and $\Linf$ between the categories $\dgl$ and $\dgc$ and the category $\dg$ which preserve all weak equivalences and are compatible with $\C$ and $\L$.  This is enough for us to blindly do (using $\dg$) everything which we would like to do with $\dgl$-spectra or $\dgc$-spectra without any further troubles or worries about isomorphisms of stablilizations.

\begin{note}
We grade $\dg$ the consistently with $\dgc$ rather than $\dgl$ -- this inserts a grading shift into the $\Linf$ and $\Sinf$ functors when going to/from $\dgl$.
\end{note}

\begin{defn} 
Define the functors $\Sigma^\infty_\dgl$, $\Omega^\infty_\dgl$, $\Sigma^\infty_\dgc$, and $\Omega^\infty_\dgc$ as follows:
\begin{itemize}
\item
$\Sigma^\infty_\dgl:\dgl_{r-1}\to\dg_r$ by $\Sigma^\infty_\dgl(L) = s(L)^\ab$
\item 
$\Omega^\infty_\dgl:\dg_r\to\dgl_{r-1}$ by $\Omega^\infty_\dgl(V) = [s^{-1}\mathrm{red}_rV]_\mDGL$
\item
$\Sigma^\infty_\dgc:\dgc_r \to \dg_r$ by $\Sigma^\infty_\dgc(C) = [C]_\mDG$
\item
$\Omega^\infty_\dgc:\dg_r \to \dgc_r$ by $\Omega^\infty_\dgc(V) = \Lambda \, \mathrm{red}_r V$
\end{itemize}
\end{defn}

Note that these functors are given by compositions of Quillen adjoints, thus they are also Quillen adjoints.

\begin{lemma}The above functors are Quillen adjoint pairs as follows:
\begin{align*}
\Sigma^\infty_\dgl:\dgl_{r-1} &\rightleftarrows \dg_r:\Omega^\infty_\dgl \\
\Sigma^\infty_\dgc:\dgc_r &\rightleftarrows \dg_r:\Omega^\infty_\dgc
\end{align*}
\end{lemma}

Furthermore they are compositions of functors preserving weak equivalences.

\begin{lemma}The right adjoint functors $\Omega^\infty_\dgl$ and $\Omega^\infty_\dgc$ preserve all weak equivalences and therefore homotopy pullbacks.  

The left adjoint functor $\Sigma^\infty_\dgc$ preserves all weak equivalences and therefore homotopy pushouts.

The left adjoint functor $\Sigma^\infty_\dgl$ preserves weak equivalences of cofibrant (free) objects, and therefore homotopy pushouts of cofibrant objects.
\end{lemma}

We generally treat $\Sigma^\infty_\dgl$ as though it preserves all homotopy pushouts since it at least preserves homotopy pushouts after cofibrant replacement.

\begin{note}
The functors $\Sigma^\infty_\dgc:\dgc_n \to \dg$ and $\Omega^\infty_\dgl:\dg\to \dgl_{n-1}$ and have the property that $\Sigma^\infty_\dgc$ takes a \DGC\ to a \DG\ with the same homology and $\Omega^\infty_\dgl$ takes a \DG\ to a \DGL\ with the same (higher) homology.  This is analogous to the non-rational situation where $\Sinf$ takes a space to a spectrum whose homotopy is the same as the (reduced) homology of the space and $\Linf$ takes a spectrum to a space whose (higher) homotopy is the same as the (higher) homotopy of the spectrum.
\end{note}

Note that the stabilization of $\dg_r$ is also given by $\dg$.  The correct Quillen adjoint pair between $\dg_r$ and $\dg$ is
\begin{defn}
Define $\Sigma^\infty_{\mDG_r}$ and $\Omega^\infty_{\mDG_r}$ as follows 
\begin{itemize}
\item
$\Sigma^\infty_{\mDG_r}:\dg_r \to \dg$ by $\Sigma^\infty_{\mDG_r}(V) = V$
\item
$\Omega^\infty_{\mDG_r}:\dg \to \dg_r$ by $\Omega^\infty_{\mDG_r}(V) = \mathrm{red}_rV$
\end{itemize}
\end{defn}
We have already commented that this is a Quillen adjoint pair preserving all weak equivalences.

\begin{note}
Gathered below are the definitions of the functors $\Sigma$, $\Omega$, $\Sinf$, and $\Linf$ in our various categories of interest:  
\begin{enumerate}
\item On $\dg_r$:
\begin{enumerate}
\item $\Sigma:\dg_r \to \dg_r$ by $V\mapsto sV$.
\item $\Omega:\dg_r \to \dg_r$ by $V\mapsto \red{s^{-1}V}$.
\item $\Sinf:\dg_r \to \dg$ by inclusion of categories.
\item $\Linf:\dg \to \dg_r$ by reduction $V\mapsto \red{V}$.
\end{enumerate}
\item On $\dgl_{r-1}$:
\begin{enumerate}
\item $\Sigma:\dgl_{r-1}\to\dgl_{r-1}$ by $L\mapsto \freeL_{s(L)^\ab}$ (the truly free \DGL\ on the shifted abelienization).
\item $\Omega:\dgl_{r-1}\to\dgl_{r-1}$ by $L\mapsto \bigl[\red{s^{-1}[L]_\mDG}\bigr]_\mDGL$ (a \DGL\ with trivial Lie bracket).
\item $\Sinf:\Fdgl_{r-1}\to\dg$ by $L \mapsto s(L)^\ab$ (abelienizing and shifting).
\item $\Linf:\dg\to\dgl_{r-1}$ by $V\mapsto [\red{s^{-1}V}]_{\mDGL}$ (reducing and giving trivial Lie bracket).
\end{enumerate}
\item On $\dgc_r$:
\begin{enumerate}
\item $\Sigma:\dgc_r\to\dgc_r$ by $C\mapsto \bigl[s[C]_\mDG\bigr]_\mDGC$ (a \DGC\ with trivial coproduct).
\item $\Omega:\dgc_r\to\dgc_r$ by $C\mapsto \Lambda \red{s^{-1}(C)^\pr}$ (the truely cofree \DGC\ on the shifted primitives).
\item $\Sinf:\dgc_r\to\dg$ by $C\mapsto [C]_\mDG$ (de-augmenting and forgetting coproduct).
\item $\Linf:\dg\to\Fdgc_r$ by $V\mapsto \Lambda\red{V}$ (the truely cofree \DGC\ on $V$).
\end{enumerate}
\end{enumerate}
\end{note}

\begin{lemma}There are natural zig-zags of weak equivalences exhibiting the following:
\begin{enumerate}
\item\label{L Linf} $\L \Omega^\infty_\dgc \weq \Omega^\infty_\dgl : \dg\to \dgl_{r-1}$
\item\label{C Linf} $\C \Omega^\infty_\dgl \weq \Omega^\infty_\dgc : \dg\to \dgc_r$
\item\label{Sinf C} $\Sigma^\infty_\dgc \C \weq \Sigma^\infty_\dgl : \dgl_{r-1} \to \dg$
\item\label{Sinf L} $\Sigma^\infty_\dgl \L \weq \Sigma^\infty_\dgc : \dgc_r \to \dg$
\end{enumerate}
where $\Omega^\infty_\catM$ and $\Sigma^\infty_\catM$ denote the functors $\Linf:\dg\to\catM$ and $\Sinf:\catM\to\dg$.
\end{lemma}
\begin{proof}
Let $V\in \dg$ and $C \in \dgc_r$.

For (\ref{C Linf}) note that $\Omega^\infty_\dgl V = [s^{-1}\red{V}]_\mDGL$ has trivial Lie bracket, so $\C\Omega^\infty_\dgl V$ is the truly cofree \DGC\ on $V$
$$(\C\Omega^\infty_\dgl) (V) = \bigl(\Lambda\, s(s^{-1}\red{V}),\, d_\Lambda\bigr) =  
 \Lambda\bigl(\mathrm{red}_r(V,\,d_V)\bigr) = \Omega^\infty_\dgc V$$  
Similarly for (\ref{Sinf L}), since $\L C$ is a free \DGC\ and $d_\Delta$ increases word-length by 1 we have that
$$(\Sigma^\infty_\dgl \L) (C) = s\bigl(\freeL_{s^{-1}C},\, d_\freeL + d_\Delta\bigr)^\ab =
 s\bigl(s^{-1}C,\, d_{s^{-1}C}\bigr) = (C,\,d_C) = \Sigma^\infty_\dgc C$$
(\ref{L Linf}) and (\ref{Sinf C}) now follow using the adjointness of $\C$ and $\L$.  Recall that there is a natural weak equivalence $\L\C \xrightarrow{\,\weq\,} \Id$.  This induces natural weak equivalences
\begin{align*}
\L\Omega^\infty_\dgc = \L\bigl(\C\Omega^\infty_\dgl\bigr) 
 &\xrightarrow{\,\weq\,} \Omega^\infty_\dgl \\
\Sigma^\infty_\dgc\C = \bigl(\Sigma^\infty_\dgl \L\bigr) \C 
 &\xrightarrow{\,\weq\,} \Sigma^\infty_\dgl
\end{align*}
\end{proof}

\begin{cor}[Rational Snaith Splitting]\label{Q Snaith} Rationally, $\Sigma^\infty_\dgc\Omega^\infty_\dgc:\dg\to\dg$ is given by $V \mapsto [\Lambda V]_\mDG$ the underlying module of the graded-symmetric coalgebra on $V$ with induced differential.
\end{cor}


\part{Rational Homotopy Calculus of Functors}

\markboth{Ben Walter}{II.  Rational Homotopy Calculus}

We mimic Goodwillie's construction of homotopy calculus of functors in the categories $\top$ and $\sp$ ([GIII]), this time working in the categories $\dgl$, $\dgc$, and $\dg$ (our omission to show $\dg$ is equivalent to $\dgl$-spectra and $\dgc$-spectra adds a slight wrinkle, but not cause too much trouble).  
The fundamental objects in Goodwillie's calculus of functors are $n$-excisive homotopy functors.  Recall that a functor is a homotopy functor if it takes weak equivalences to weak equivalences.  Such a functor is $n$-excisive if it takes strongly homotopy cocartesian $(n+1)$-cubes to homotopy cartesian $(n+1)$-cubes.  

An $n$-cube in $\catC$ is a diagram 
$$\mathcal{X}:\mathcal{P}(\underline{n}) \longrightarrow \catC$$ 
where $\mathcal{P}(\underline{n})$ is the poset of subsets of $\underline{n} = \{1,\dots, n\}$ and inclusion maps viewed as a category.
Recall that we write $\mathcal{P}_0(\underline{n})$ for the full subcategory of nonempty subsets of $\underline{n}$ and $\mathcal{P}_1(\underline{n})$ for the full subcategory of proper subsets of $\underline{n}$.
An $n$-cube $\mathcal{X}$ is {\em homotopy cartesian} if the natural composition 
$$\mathcal{X}(\emptyset) \xrightarrow{\ \ \ \ \ } 
  \!\lim_{S\in\mathcal{P}_0(\underline{n})}\!\!\mathcal{X}(S) \xrightarrow{\ \ \ \ \ } 
  \holim_{S\in\mathcal{P}_0(\underline{n})}\mathcal{X}(S)$$ 
is a weak equivalence.  Dually, we say that an $n$-cube $\mathcal{X}$ is {\em homotopy cocartesian} if the natural composition 
$$\hocolim_{S\in\mathcal{P}_1(\underline{n})} \mathcal{X}(S) \xrightarrow{\ \ \ \ \ }  
    \colim_{S\in\mathcal{P}_1(\underline{n})} \mathcal{X}(S) \xrightarrow{\ \ \ \ \ } 
    \mathcal{X}(\underline{n})$$ 
is a weak equivalence.  A {\em strongly homotopy cocartesian} $n$-cube is one in which every 2-dimensional face is cocartesian (it is an easy exercise to show this implies the cube itself as well as every face of dimension greater than 1 is also cocartesian).

Throughout this part, we often refer to ``towers of functors'' and ``towers of fibrations''.  By a ``tower'' of objects, we merely mean an inverse system 
$$\xymatrix@R=10pt{
 \vdots \ar[d] \\
 X_n \ar[d] \\
 X_{n-1} \ar[d] \\
 \vdots \ar[d] \\
 X_1 \ar[d] \\
 X_0
}$$
A ``tower of functors'' is an inverse limit diagram of functors.  A ``tower of fibrations'' is an inverse limit diagram where each map is a fibration.  Similarly ``tower of fibrations of functors'', ``tower of graded Lie algebras'', etc. all have the obvious meanings.

Our towers come equipped with an object $X$ which has maps $X\to X_i$ for all $i\ge 0$, and we are interested in comparing $X$ with the (homotopy) limit of the tower.  Good objects $X$ are recovered (up to homotopy) by this (homotopy) limit.  In order to analyse the (homotopy) limit of the tower, we examine the (homotopy) fibers of the maps $X_n \to X_{n-1}$.  Much effort is put into understanding these (homotopy) fibers as well as possible.

\chapter{Construction}\label{S:Q calc const}
\markboth{Ben Walter}{II.7 Construction}

Rational homotopy calculus of functors is an algebraic theory analogous to Goodwillie's homotopy calculus of functors as defined in [GIII] and described in [Kuhn].  Our objects of interest are rational homotopy functors between any of the categories $\dg$, $\dg_r$, $\dgl_{r-1}$, $\dgc_r$ -- i.e. functors between these categories which preserve rational homotopy equivalences (quasi-isomorphisms).  Following the outline of Kuhn, we divide our analysis into three main theorems which we call Theorems I-III.

Let $\catM$ and $\catN$ be any of the model categories $\dg$, $\dg_r$, $\dgl_{r-1}$, or $\dgc_r$.  Given a rational homotopy functor $F:\catM \to \catN$, Theorem I associates to it a natural tower of fibrations with maps from $F$
$$F\to \ \ \cdots\fibr P_nF \fibr P_{n-1}F \fibr\cdots\fibr P_0F$$
such that each $P_nF$ is $n$-excisive.  Ideally the tower converges to $F$ -- that is the maps $F\to P_nF$ become more and more highly connected as $n$ increases.  However, we do note concern ourselves overly much with this detail in the current monograph.

The next two main theorems are used to analyze the fibers $D_nF \to P_nF \fibr P_{n-1}F$ in this tower.  Essentially,\footnote{Actually, Theorems II and III are stronger that what is alluded to in this introduction.} they say that in order to understand the fibers $D_nF:\catM\to\catN$ it is enough to understand certain assiciated functors $\mathbb{D}_nF:\dg \to \dg$.  In particular, there is a zig-zag of natural weak equivalences of functors displaying 
$$D_nF \,\weq\, \Omega^\infty_\catN(\mathbb{D}_nF)\Sigma^\infty_\catM$$ 
where $\mathbb{D}_nF$ is an $n$-homogeneous functor $\dg \to \dg$.
Furthermore, the functors $\mathbb{D}_nF$ all have a particulary nice form:  up to a zig-zag of natural weak equivalences they are given by 
$$\mathbb{D}_nF(X) \weq (\partial_n\!F \otimes X^{\otimes n})_{\Sigma_n}$$
where $\partial_n\!F$ is a \DG\ with $\Sigma_n$-action.  Thus the fibers in tower of $F$ are determined (up to weak equivalence) by $\{\partial_n\!F\}$ a symmetric sequence of \DG\,s.  It is natural to ask what extra structure is required to determine the entire tower of $F$ and not just the fibers.  This question is left for a sequel.

As we noted in Chapter 1, the construction of rational homotopy calculus is largely formal, thus many proofs are left out or given only in sketch form.   The few proofs which we need to give are for those results dealing specifically with the category of rational spectra.  Since we did not define rational spectra as either the stabilization of $\dgl$ or $\dgc$, Goodwillie's methods do not extend to give proofs of these results.  Fortunately, rational homotopy theory is so pleasant that we are able to give alternate proofs which are still extremely straightforward.  We also make a slight detour in Section~\ref{calculus comparison} to show that homotopy calculus is consistent with good Quillen equivalences in the sense that approximating towers are preserved by Quillen equivalences where the left and right adjoints preserve all weak equivalences.

\section{Excisive Approximations and Theorem I}

As in the introduction, let $\catM$ and $\catN$ be any of the model categories $\dg$, $\dg_r$, $\dgl_{r-1}$, or $\dgc_r$, and let $F:\catM \to \catN$ be a rational homotopy functor from $\catM$ to $\catN$.  The construction of a tower of excisive approximations to $F$ follows along the lines of Goodwillie's work [GIII] and Kuhn's survey [Kuhn].  We give a brief outline of the construction along with some notes on our specific case.

Excisive approximations are built by making for every $X\in\catM$ a natural ``test strongly cocartesian $n$-cube''\footnote{This terminology is meant to bring to mind ``test cofibrations'' or ``test fibrations'' which one may use when building a model category structures.} $\mathcal{X}_n(X):\mathcal{P}(\underline{n})\to\catM$ and then examining the resulting $n$-cubes $F(\mathcal{X}_n(X))$.  The cubes $\mathcal{X}_n(X)$ are built in such a way that if all $F(\mathcal{X}_n(X))$ are cartesian, then $F$ is $(n-1)$-excisive.  If $F(\mathcal{X}_n(X))$ is not cartesian, then $F(\mathcal{X}_n)$ provides us with a new functor $T_{n-1}F$ which is equipped with a natural map $F \to T_{n-1}F$ and is (hopefully) slightly closer to being cartesian on the $\mathcal{X}_n(X)$.  We now check $T_{n-1}F$ on our test cubes $\mathcal{X}_n(X)$ and the process iterates.  Our $(n-1)$-excisive approximation to $F$ is then given by the sequential homotopy colimit of the diagram
$$P_{n-1}F = \hocolim(F\to T_{n-1}F \to T_{n-1}^2F \to \cdots)$$
Showing that $P_{n-1}F$ is indeed $(n-1)$-excisive consists essentially of noting that homtopy limits in $\catN$ commute with sequential homotopy colimits in $\catN$.

We define test cocartesian $n$-cubes following the method of Kuhn.  The first step is to define the join of an object of $\catM$ and a set:

\begin{defn}[\lbrack Kuhn, 4.6\rbrack] For $X\in \catM$ and $T$ a finite set, define $X\ast T$ -- the {\em join of $X$ and $T$} -- to be the homotopy cofiber of the folding map
$$X\ast T = \hocof\bigl(\coprod_{T} X \xrightarrow{\,\triangledown\,} X\bigr)$$

More generally, for $X\xrightarrow{f} Y \in \catM$ and $T$ a finite set, define $X\ast_Y T$ -- the {\em join of $X$ and $T$ over $Y$} -- to be the homotopy colimit
$$X\ast_Y T = \hocolim\bigl(\coprod_{T} Y \xleftarrow{\,f\,} \coprod_T X \xrightarrow{\,\triangledown\,} X\bigr)$$
\end{defn}

\begin{ex}
For $X$ either a \DG\ or a \DGC,
\begin{itemize}
\item $X \ast \underline{0} = X\ast \emptyset = X$.
\item $X \ast \underline{1} = cX$.
\item $X \ast \underline{2} = \hat \Sigma X$.
\item In general, for $T = \{t_1,\dots,t_n\}$ let 
\begin{align*}
&[T]_\mDG = (\Q v_0 \oplus \Q t_1 \oplus \cdots \Q t_n,\, d(t_i) =  v_0) \quad \mathrm{where}\ \ |v_0| = 0,\ \ |t_i| = 1\\
&[T]_\mDGC = \bigl([T]_\mDG,\, \Delta(v_0) = v_0\otimes v_0,\, \Delta(t_i) = {\textstyle \frac{1}{2}}(v_0 \otimes t_i + t_i\otimes v_0)\bigr)
\end{align*}
Then $V\ast T = [T]_\mDG\otimes V$ and $C\ast T = \bigl([T]_\mDGC \otimes \widetilde{C}\bigr)_+$ (for $V\in \dg$ and $C\in\dgc_r$).
\end{itemize}
\end{ex}

\begin{ex}
For $L$ a \DGL,
\begin{itemize}
\item $L \ast \underline{0} = \L \C L$.
\item $L \ast \underline{1} = c\L \C L$.
\item $L \ast \underline{2} = \hat\Sigma \L \C L$.
\end{itemize}
\end{ex}

Our test cocartesian cubes are formed by using joins with finite sets.  For $X\in\catM$ define $$\mathcal{X}_n(X):\mathcal{P}(\underline{n}) \to \catM\quad \mathrm{by}\quad \mathcal{X}_n(X):S \mapsto X\ast S$$

\begin{lemma}For every $X\in \catM$ the $n$-cube $\mathcal{X}_n(X)$ is strongly cocartesian.
\end{lemma} 
\begin{proof}
This follows from the fact that homotopy colimits commute with themselves.

We show that $\mathcal{X}_2(X)$ is cocartesian -- proving that each 2-dimensional face of $\mathcal{X}_n(X)$ is cocartesian follows the same (but is notationally more ugly).  The cube $\mathcal{X}_2(X)$ is given by:
$$\xymatrix@C=10pt@R=15pt{
\hocolim(0\leftarrow \coprod_{\emptyset} X \rightarrow X) \ar[r] \ar[d] &
 \hocolim(0\leftarrow \coprod_{\{1\}} X \rightarrow X) \ar[d] \\
\hocolim(0\leftarrow \coprod_{\{2\}} X \rightarrow X) \ar[r] &
 \hocolim(0\leftarrow \coprod_{\{1,2\}} X \rightarrow X)
}$$
We are interested in the homotopy colimit
$$\hocolim\left(
\begin{aligned}\xymatrix@R=15pt{ 
\hocolim(0 \longleftarrow \coprod_{\{1\}} X \longrightarrow X) \\ 
\hocolim(0 \longleftarrow \coprod_\emptyset X \longrightarrow X) \ar[d]\ar[u]\\
\hocolim(0 \longleftarrow \coprod_{\{2\}} X \longrightarrow X) 
}\end{aligned}\right)$$ 
However, this is equal to the homotopy colimit 
$$\hocolim\left(
\begin{aligned}\xymatrix@C=20pt@R=15pt{
0 & \coprod_{\{1\}} X \ar[l]\ar[r] & X  \\
0 \ar[d]\ar[u] & \coprod_\emptyset X \ar[l]\ar[r]\ar[d]\ar[u] & X \ar[d]\ar[u]\\
0 & \coprod_{\{2\}} X \ar[l]\ar[r] & X 
}\end{aligned}\right)$$
which is equal to the homotopy colimit 
$$\hocolim\left(\ \begin{aligned}
  \hocolim\left(\begin{aligned}
   \xymatrix@R=15pt{0\vphantom{\coprod_\{} \\ 0\vphantom{\coprod_\{} \ar[d]\ar[u] \\ 0\vphantom{\coprod_\{}}
   \end{aligned}\right) \longleftarrow
  \hocolim\left(\begin{aligned}
   \xymatrix@R=15pt{\coprod_{\{1\}} X \\ \coprod_\emptyset X \ar[d]\ar[u] \\ \coprod_{\{2\}} X}
   \end{aligned}\right) \longrightarrow
  \hocolim\left(\begin{aligned}
   \xymatrix@R=15pt{X\vphantom{\coprod_\{} \\ X\vphantom{\coprod_\{} \ar[d]\ar[u] \\ X\vphantom{\coprod_\{}}
   \end{aligned}\right)
 \end{aligned}\ \right)
$$
This homotopy colimit maps by a weak equivalence to $\hocolim(0 \leftarrow \coprod_{\{1,2\}} X \rightarrow X)$ because the diagram itself maps by a weak equivalence to the diagram $0 \leftarrow \coprod_{\{1,2\}} X \rightarrow X$.
\end{proof}

More generally for $X\xrightarrow{f}Y \in \catM$ define the $n$-cube $\mathcal{X}_n(X\xrightarrow{f}Y)$ to be
$$\mathcal{X}_n(X\xrightarrow{f}Y):\mathcal{P}(\underline{n}) \to \catM \quad \mathrm{by}\quad 
\mathcal{X}_n(X\xrightarrow{f}Y):S\mapsto X\ast_Y S$$
Using the commutativity of homotopy colimits as in the proof of the previous lemma it follows that:

\begin{lemma}For every $X\xrightarrow{f}Y \in \catM$ the $n$-cube $\mathcal{X}_n(X\xrightarrow{f}Y)$ is strongly cocartesian.
\end{lemma}

Consider the $n$-cubes given by evaluating $F$ objectwise on the cubes $\mathcal{X}_n(X)$.  If $F$ is $(n-1)$-excisive, then by definition the $n$-cube $F\bigl(\mathcal{X}_n(X)\bigr)$ is cartesian -- i.e. the natural map
$$\xymatrix@C=30pt{
F(X) = F\bigl(\mathcal{X}_n(X)(\emptyset)\bigr) \ar[r]^{t_{n-1}} & 
 {}\phantom{\displaystyle \holim_{S\in \mathcal{P}_0(\underline{n})} F\bigl(\mathcal{X}_n(X)(S)\bigr)}
 \save[]-<0pt,.05in>
   *{\displaystyle \holim_{S\in \mathcal{P}_0(\underline{n})} F\bigl(\mathcal{X}_n(X)(S)\bigr)} \restore
}$$ 
is a weak equivalence.
Define the functor $T_{n-1}F$ be the homotopy limit 
$$T_{n-1}F(X) = \holim_{S\in \mathcal{P}_0(\underline{n})} F(\mathcal{X}_n(X)(S))$$
By definition if $F$ is $(n-1)$-excisive, then the map $F(X) \to T_{n-1}F(X)$ is a quasi-isomorphism.  If $F$ is not $(n-1)$-excisive, then our goal is to use $T_{n-1}$ to attempt to make an $(n-1)$-excisive approximation for $F$. 

Write $T_{n-1}^iF$ for the iterated construction $\overbrace{T_{n-1}(T_{n-1}(\dots T_{n-1}F))}^{i}$ and note that the natural map $F\to T_{n-1}F$ gives maps $T_{n-1}^iF \to T_{n-1}^{i+1}F$.  Let $P_{n-1}F$ be the sequential homotopy colimit $$P_{n-1}F = \hocolim_i T_{n-1}^iF$$ 
Since $F = T_{n-1}^0F$ the functor $P_{n-1}F$ comes equipped with a map $F\xrightarrow{p_{n-1}\!F}P_{n-1}F$.  Furthermore, there are clear maps $T_nF \to T_{n-1}F$ induced by the inclusion of categories $\mathcal{P}_0(\underline{n}) \to \mathcal{P}_0(\underline{n+1})$.  As discussed in [GIII p664] this formally extends to give a commutative diagram
$$\xymatrix@R=20pt@C=55pt{
F\ar[r]^{t_n\!F}\ar[d] & T_nF \ar[r]^{t_n\!T_nF}\ar[d]^{q_{n,1}} & T^2_nF \ar[r]^{t_n\!T^2_nF}\ar[d]^{q_{n,2}} & \dots \\
F\ar[r]^{t_{n-1}\!F} & T_{n-1}F \ar[r]^{t_{n-1}\!T_{n-1}\!F} & T^2_{n-1}F \ar[r]^{t_{n-1}\!T^2_{n-1}\!F} & \dots 
}$$
and thus defines a natural map between the homotopy colimits $q_nF:P_nF \to P_{n-1}F$.  Since the $q_{n,i}F$ are the natural maps from the homotopy limit of a diagram to the homotopy limit of a restriction of the diagram, they are objectwise fibrations.  By inspection, sequential homotopy colimits of fibrations are again fibrations in $\dg$, $\dgc$, and $\dgl$.

\begin{thmI}[\lbrack GIII, 1.13\rbrack] A homotopy functor $F:\catM\to\catN$ determines a tower of functors $P_nF:\catM\to\catN$ with maps from $F$:
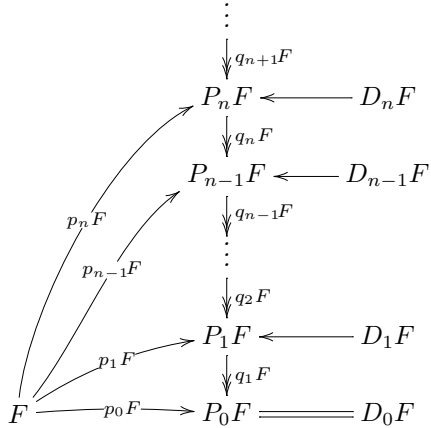
$$\xymatrix@R=15pt{
&& {}\phantom{|}\save[]+<0pt,4pt>*{\vdots}\restore \ar@{->>}[d]^{q_{n+1}\!F}  & \\
&& P_nF \ar@{->>}[d]^{q_nF} & D_nF \ar[l] \\
&& P_{n-1}F \ar@{->>}[d]^(.45){q_{n-1}\!F} & D_{n-1}F \ar[l] \\
&& {}\phantom{|}\save[]+<0pt,4pt>*{\vdots}\restore \ar@{->>}[d]^{q_2F} &  \\
&& P_1F \ar@{->>}[d]^{q_1F} & D_1F \ar[l] \\
F \ar@/^2pt/[rr]|{p_0F} \ar@/^/[urr]|{p_1F} \ar@(ur,l)[uuurr]|{p_{n-1}\!F} \ar@(u,l)[uuuurr]|{p_{n}F}
 && P_0F & D_0F \ar@{=}[l]
}$$
where the functors $P_nF$ are $n$-excisive, the maps $q_nF:P_nF(X)\fibr P_{n-1}F(X)$ are fibrations, the functors $D_nF = \hofib(q_nF)$ are $n$-homogeneous,\footnote{Recall that a functor is $n$-homogeneous if it is both $n$-excisive and $n$-reduced -- i.e. $D_nF$ is $n$-excisive and $P_{n-1}D_nF \weq 0$.} the maps $p_n$ and $q_n$ are natural in $F$, and $p_n$ is homotopy initial among all natural transformations from $F$ to an $n$-excisive homotopy functor.
\end{thmI}

\begin{note}
More generally, if $\catM$ and $\catN$ are any pointed model categories such that Theorem I may be proven, then we say that ``a calculus of functors may be constructed'' for functors $\catM\to\catN$.
\end{note}

We do not give a full proof of this here, since Goodwillie's proof applies almost immediately to the present situation.  Instead we note the main ingredients taken from [GIII] along with enough detail of their proofs to show that they carry over into the current context.  

Key to the proof of Theorem I is the following observation:

\begin{lemma}[\lbrack GIII, 1.7\rbrack]\label{G1.7} Up to a zig-zag of natural weak equivalences,
\begin{enumerate}
\item\label{1.7.1} $T_n$ commutes with $\holim_\catN$.
\item\label{1.7.2} $P_n$ commutes with finite $\holim_\catN$.
\item\label{1.7.3} $T_n$ and $P_n$ commute with $\hofib_\catN$.
\item\label{1.7.4} $T_n$ and $P_n$ commute with sequential $\hocolim_\catN$.
\item\label{1.7.5} for $\dg$-valued functors, $T_n$ and $P_n$ commute with any $\hocolim_\dg$.
\end{enumerate}
\end{lemma}
\begin{proof}
Note that $T_n$ is a homotopy pullback and $P_n$ is a sequential $\hocolim_\catN$ of homotopy pullbacks.  Critical to us is the commutativity of sequential homotopy colimits with homotopy pullbacks, which we have already pointed out in each of our categories of interest.

Statement~(\ref{1.7.1}) follows from the commutativity of $\holim_\catN$ with $\holim_\catN$.  (\ref{1.7.2}) follows from the commutativity (up to a weak equivalence) of homotopy pullbacks with a sequential $\hocolim_\catN$.  Together (\ref{1.7.1}) and (\ref{1.7.2}) imply (\ref{1.7.3}).  Statement (\ref{1.7.4}) also follows from the commutativity of homotopy pullbacks with a sequential $\hocolim_\catN$.  (\ref{1.7.5}) follows from the fact that homtopy cocartesian cubes in $\dg$ are also homotopy cartesian (since $\dg$ is stable).
\end{proof}

\begin{lemma}[\lbrack GIII, 1.8\rbrack]\label{P_nF n-excisive} $P_nF$ is $n$-excisive and $p_nF:F\to P_nF$ is homotopy initial.
\end{lemma}
\begin{proof}[Proof sketch]
We sketch the argument showing $P_nF$ is $n$-excisive.  That $p_nF$ is homotopy initial follows formally from Goodwillie in a similar manner.

Let $\mathcal{X}$ be a strongly cocartesian $(n+1)$-cube in $\catM$.  To show that $P_nF(\mathcal{X})$ is cartesian, we construct a cartesian cube which the map of cubes
$$F(\mathcal{X}) \xrightarrow{t_n\!F(\mathcal{X})} T_n\!F(\mathcal{X})$$
factors through.  Then we have $P_nF(\mathcal{X})$ weakly equivalent to a sequential $\hocolim_\catN$ of cartesian $(n+1)$-cubes.  Since $P_n$ commutes with sequential $\hocolim_\catN$, $P_nF(\mathcal{X})$ is cartesian as well.

We rely on the preservation of canonical homotopy limits by restriction to cofinal subindexing categories (see \S\ref{S:cofinal}) in simplicially enriched model categories.
Recall that if $\mathscr{D}:\catI \to \catC$ is a diagram in some simplicially enriched model category $\catC$ and $\catJ$ is a cofinal subcategory of the index category $\catI$, then there is a weak equivalence between the canonical homotopy limit of $\mathscr{D}|_\catJ$ and the canonical homotopy limit of the restricted diagram $\mathscr{D}$
$$ \holim_\catI (\mathscr{D})=  \holim(\mathscr{D})
  \xrightarrow{\ \weq\ } \holim \bigl(\mathscr{D}|_\catJ\bigr) = \holim_\catJ (\mathscr{D})$$ 
For any other homotopy limit functors on $\catC$ this yields a zig-zag of natural weak equivalences between $\displaystyle \holim_\catI$ and $\displaystyle \holim_\catJ$.
Recall that all of our categories $\dg$, $\dg_r$, $\dgl_{r-1}$, and $\dgc_r$ are simplicially enriched, so we have such a zig-zag of natural weak equivalences.  

The proof is therefore completed by constructing a large cube $\hat{\mathcal{X}}:\mathscr{P}(\underline{n+1})\times\catI \to \catM$ along with subcategories $\catJ, \varepsilon\subset \catI$ satisfying the following four properties:
\begin{align}
&\mathcal{X}(T) \xrightarrow{\weq} \operatorname*{holim}_{T\times
\mathcal{I}} \hat{\mathcal{X}} \label{property 1} \\
&\text{The cube } \operatorname*{holim}_{T\times
\varepsilon} F(\hat{\mathcal{X}})\ \text{is cartesian} \label{property 2}\\
& \varepsilon \text{ is left cofinal in } (\catJ \cup \varepsilon) \label{property 2.5}\\
& \operatorname*{holim}_{T\times\mathcal{J}} \hat{\mathcal{X}} \to
\operatorname*{holim}_{U \in
\mathcal{P}_0(\underline{n})}\mathcal{X}(T)\ast U \label{property 3}
\end{align}
(By $(\catJ \cup \varepsilon)$ we mean the full subcategory of $\catI$ generated by objects of $\catJ$ and $\varepsilon$.)

Our choices follow formally from Goodwillie's construction in [GIII 1.9].  Take for $\catI$ the category $\catI = \mathcal{P}(\underline{n+1})^{n+1}$ and for $\catJ$ the category $\catJ = Diagonal(\mathcal{P}_0(\underline{n+1})^{n+1})$.  The diagram 
$$\hat{\mathcal{X}}:\mathcal{P}(\underline{n+1})\times\mathcal{P}(\underline{n+1})^{n+1}\longrightarrow \catN$$
is now given by 
$$\hat{\mathcal{X}}(T,U_1,\dots,U_{n+1}) = \hocolim\biggl(\mathcal{X}(T)
 \xleftarrow{\ \ \triangledown\ \ }\!\!\!\!\coprod_{1\le s \le n+1}\ \coprod_{U_s} \mathcal{X}(T)
 \xrightarrow{\,\coprod\!\mathcal{X}(i)\,}\!\!\!\!\coprod_{1\le s\le n+1}\ \coprod_{U_s} \mathcal{X}(T\cup\{s\})
 \biggr)$$
Where $\triangledown$ is the fold map and $\coprod\!\mathcal{X}(i)$ is induced by the inclusions $i:T\cofibr(T\cup\{s\})$.  This may also be described as the colimit of the diagram mapping $\mathcal{X}(T)$ to each of the objects $\mathcal{X}(T)\ast_{\mathcal{X}(T\cup\{s\})}U_s$.

After choosing the correct $\varepsilon$ the desired factorization is given by:
$$\xymatrix{
F(\mathcal{X}) \ar[rrrr]^{t_nF(\mathcal{X})}\ar[dr]^{\weq}_{(1)} &&&& T_nF(\mathcal{X}) \\
& {} \displaystyle \phantom{F\biggl(\holim_{T\times\catI} \hat{\mathcal{X}}\biggr)} \save[]-<0pt,2pt>
  *{\displaystyle F\biggl(\holim_{T\times\catI} \hat{\mathcal{X}}\biggr)}\restore \ar[r]_{(2)}
& {} \displaystyle \phantom{F\biggl(\holim_{T\times(\catJ \cup \varepsilon)} \hat{\mathcal{X}}\biggr)}
  \save[]-<0pt,2pt> *{\displaystyle F\biggl(\holim_{T\times(\catJ \cup \varepsilon)} \hat{\mathcal{X}}\biggr)} 
  \restore \ar[r]_{(4)} \ar@{{}.{}}[d]^{\weq}_{(3)}
& {} \displaystyle \phantom{F\biggl(\holim_{T\times\catJ} \hat{\mathcal{X}}\biggr)} \save[]-<0pt,2pt>
  *{\displaystyle F\biggl(\holim_{T\times\catJ} \hat{\mathcal{X}}\biggr)} \restore \ar[ur]_{(5)} 
& \\
& & {} \displaystyle \phantom{F\biggl(\holim_{T\times\varepsilon} \hat{\mathcal{X}}\biggr)} \save-<0pt,2pt>
  *{\displaystyle F\biggl(\holim_{T\times\varepsilon} \hat{\mathcal{X}}\biggr)} \restore  & & 
}$$
The ``$\cdots$'' in (3) is a zig-zag of natural weak equivalences coming from property (\ref{property 2.5}) of the categories $\catI$, $\catJ$, and $\varepsilon$.
Maps $(1)$ and $(5)$ follow from properties (\ref{property 1}) and (\ref{property 3}).  Maps $(2)$ and $(4)$ are induced by inclusion maps of indexing categories.  Also property (\ref{property 2}) ensures that the cube $\displaystyle F\bigl(\holim_{T\times\varepsilon} \hat{\mathcal{X}}\bigr)$ is cartesian, and thus so is the cube $\displaystyle F\bigl(\holim_{T\times(\catJ\cup\varepsilon)}\hat{\mathcal{X}}\bigr)$.

In Calculus III, Goodwillie uses $\varepsilon_{\mathcal{G}} = \bigl\{\vec S
\in \mathcal{P}_0(\underline{n+1})^{n+1}\ |\ \ S_j = \{j\}, \text{\ some }
j\bigr\}$ -- the largest $\varepsilon$ for which his argument works.
The smallest $\varepsilon$ which can be chosen is
$$\varepsilon_\mathcal{W} =
\bigg\{\vec S \in \mathcal{P}_0(\underline{n+1})^{n+1}\ \ \Big|\ \ \emptyset \neq A \subset B \subset
\underline{n+1},\ \ S_j =
\begin{cases} B & j \notin A \\
 \{j\} & j \in A
\end{cases}\ \  \bigg\}$$
\end{proof}

The following now formally follow from their corresponding proofs by Goodwillie:

\begin{cor}[\lbrack GIII, 1.11\rbrack] If $0 \le m\le n$ then the map 
$$P_mF \xrightarrow{P_m(p_nF)} P_mP_nF$$ is an equivalence.
\end{cor}

\begin{lemma}[\lbrack GIII, 1.17\rbrack] $D_nF$ is $n$-homogeneous.
\end{lemma}

Note that from Lemma~\ref{G1.7} we also get that:
\begin{enumerate}
\item $D_n$ commutes with finite $\holim_\catN$.
\item $D_n$ commutes with $\hofib_\catN$.
\item $D_n$ commutes with sequential $\hocolim_\catN$.
\item for $\dg$-valued functors, $D_n$ commutes with arbitrary $\hocolim_\dg$.
\end{enumerate}

\section{A Comparison Theorem for Calculi of Functors}\label{calculus comparison}

Before moving on to Theorems II and III, we note a few formal lemmas about preservation of homogeneousness and Taylor towers by good left and right adjoint functors:

In this section let $\catM$, $\catM_1$, $\catM_2$, $\catN$, $\catN_1$ and $\catN_2$ be any pointed model categories.  We assume that our categories are such that a homotopy calculus of functors may be constructed (that is Theorem I may be proven) for functors between any combination of the above categories.  As [Kuhn] points out it, in order to construct a homotopy calculus of functors it is sufficient (though not necessary) for the categories in question to be simplicially enriched, proper model categories where very small homotopy limits commute with filtered homotopy colimits.  

Write $\mathscr{H}_n(\catM, \catN)$ for the category of $n$-homogeneous functors $\catM\to\catN$ and natural transformations.  

\begin{lemma}\label{comp}
Let $R:\catN_1\to\catN_2$ be a functor which preserves weak equivalences and homotopy pullbacks and commutes with sequential homotopy colimit.  Post-composition with $R$ commutes with $P_n$ and gives a functor $$R\circ (-) :\mathscr{H}_n(\catM,\catN_1)\to\mathscr{H}_n(\catM,\catN_2)$$

Let $L:\catM_1\to\catM_2$ be a functor which presserves weak equivalences and homotopy pushouts.  Pre-composition with $L$ commutes with $P_n$ and gives a functor $$(-)\circ L:\mathscr{H}_n(\catM_2,\catN)\to\mathscr{H}_n(\catM_1,\catN)$$
\end{lemma}
\begin{proof}
By assumption $R$ comutes with homotopy pullbacks, so $T_nRF \weq RT_nF$.  Since $R$ also commutes with sequential homotopy colimit, 
$$P_nRF = \hocolim T^i_n RF \weq \hocolim RT^i_nF \weq R\hocolim T^i_nF \weq RP_nF$$

For the second statement, note that since $L$ preserves homotopy pushouts $L(X*S) \weq L(X)*S$.  Since homotopy limits preserve weak equivalence $T_n(FL) \weq (T_nF)L$.  If $L$ preserves homotopy pushouts then it also preserves sequential homotopy colimits, so $P_n(FL) \weq (P_nF)L$.
\end{proof}

In particular, if $L:\catM_1\rightleftarrows \catN_1: R$ is a Quillen adjoint pair where $L$ and $R$ preserve \underline{all} weak equivalences and $R$ also commutes with sequential homotopy colimits, then $L$ and $R$ satisfy \ref{comp}.  Thus, for $F:\catN_1 \to \catN_2$ and $G:\catM_2\to\catM_1$ homotopy functors there are zig-zags of (level-wise) natural weak equivalences of towers exhibiting
\begin{align*}
(\text{\it Taylor tower of }F)\circ L &\weq \text{\it Taylor tower of }(FL) \\
R\circ (\text{\it Taylor tower of G}) &\weq \text{\it Taylor tower of }(RG) 
\end{align*}
In other words:

\begin{cor}\label{weak tower preservation}
If $L$ is a Quillen left adjoint functor which preserves all weak equivalences, then pre-composition with $L$ preserves Taylor towers.

If $R$ is a Quillen right adjoint functor which preserves all weak equivalences and commutes with sequential homotopy colimits, then post-composition with $R$ preserves Taylor towers.
\end{cor}

Recall that if $L:\catM\rightleftarrows\catN:R$ is a Quillen equivalence of model categories then there are natural weak equivalences $\Id_\catM \xrightarrow{\,\weq\,} R\, L$ and $L\, R \xrightarrow{\,\weq\,}\Id_\catN$.  This implies the following statement (which we call a ``theorem'' due to its importance rather than its difficulty to prove):

\begin{thm}[Comparison Theorem]\label{tower preservation}
Let $L:\catM\rightleftarrows\catN:R$ be a Quillen equivalence of model categories where $L$ and $R$ preserve all weak equivalences and $R$ also commutes with sequential homotopy colimits.  Then both pre-composition and post-composition with either $L$ or $R$ preserves Taylor towers.
\end{thm}
\begin{proof}

We already know this for pre-composition with $L$ and post-composition with $R$. 

Suppose $F$ is a homotopy functor mapping from $\catM$.  There is a (level-wise) natural weak equivalence of towers
$$\text{\it Taylor tower of }F \xrightarrow{\,\weq\,} \text{\it Taylor tower of }(F\circ RL)$$
induced by the natural weak equivalence $F \xrightarrow{\,\weq\,} F\circ RL$.  However, precomposition with $L$ preserves towers, so there are   zig-zags of (level-wise) natural weak equivalences
\begin{align*}
\bigl(\text{\it Taylor tower of }F\bigr)\circ R 
 &\weq \bigl(\text{\it Taylor tower of }(F\circ RL)\bigr) \circ R \\
 &\weq \bigl(\text{\it Taylor tower of }(F\circ R)\bigr) \circ (LR) \\
 &\weq \text{\it Taylor tower of }(F\circ R)
\end{align*}

Similarly, for $G$ a homotopy functor mapping to $\catM$ there are  zig-zags of (level-wise) natural weak equivalences
\begin{align*}
L\circ\bigl(\text{\it Taylor tower of }G\bigr)
 &\weq L\circ\bigl(\text{\it Taylor tower of }(RL\circ G)\bigr) \\
 &\weq (LR)\circ\bigl(\text{\it Taylor tower of } (L\circ G)\bigr) \\
 &\weq \text{\it Taylor tower of }(L\circ G)
\end{align*}
\end{proof}

By inspection, our right adjoint Quillen functors all commute with sequential homotopy colimits:
\begin{lemma}\label{Linf C seq hocolim}
The functors $\Omega^\infty_\dgc$ and $\C$ commute with filtered homotopy colimits.  The functor $\Omega^\infty_\dgl$ commutes with filtered homotopy colimits after fibrant replacement.
\end{lemma}

\section{Delooping Homogeneous Functors and Theorem II}

Our next step is to analyze the homogeneous fibers of the tower constructed in Theorem I. 
Let $\catM$ be $\dg$, $\dg_r$, $\dgl_{r-1}$, or $\dgc_r$, and $\catN$ be $\dg_r$, $\dgl_{r-1}$, or $\dgc_r$.  
Theorem~II will imply that $n$-homogeneous functors $H_n:\catM \to \catN$ (up to a zig-zag of natural weak equivalences) have the form $\Linf\mathbb{H}_n$ where $\mathbb{H}_n:\catM \to \dg$ is also $n$-homogeneous.  

Recall from Lemma~\ref{Linf C seq hocolim} and Corollary~\ref{weak tower preservation} that the functor $\Linf:\dg \to \catN$ preserves homotopy limits and commutes with $P_n$.  Therefore post-composition with $\Linf$ preserves $n$-homogeneity.  In particular, it gives a functor from $\mathscr{H}_n(\catM, \dg)$ to $\mathscr{H}_n(\catM, \catN)$.

\begin{note}If $\catN = \dg_r$ then by $\Linf$ we mean the $n$-reduction functor.\end{note}

\begin{thmII} The functor $\mathscr{H}_n(\catM, \dg) \xrightarrow{\Linf} \mathscr{H}_n(\catM, \catN)$ has an inverse up to a zig-zag of natural weak equivalences.
\end{thmII}

The theorem is proven by first constructing a natural single delooping (up to a zig-zag of natural weak equivalences) of homogeneous functors.  This is then used to naturally construct homogeneous $\dg$-valued functors whose infinite deloopings are (up to a zig-zag of natural weak equivalences) the given $\catN$-valued homogeneous functors.

Single deloopings follow directly from [GIII]:

\begin{lemma}[\lbrack GIII, 2.2\rbrack]\label{L:thm 2}Given $F:\catM \to \catN$ a reduced homotopy functor, there is a natural diagram of homotopy functors $\catM \to \catN$ given by
 
$$\xymatrix@C=7.5pt@R=7.5pt{
P_nF \ar[rr]^{q_nF} && P_{n-1}F \\ && \\
\hat P_nF\ar[uu]^{\weq} \ar[rr] \ar[dd] && \tilde P_{n-1}F \ar[uu]^{\weq} \ar[dd] \\ 
&  \ar@{=>}[ul]|(.05)h & \\
K_nF \save[]-<20pt,0pt>*{0\weq}\restore \ar[rr] && R_nF
}$$
where $R_nF$ is $n$-homogeneous and (as indicated) $K_nF$ is (naturally objectwise) contractible and the lower square is (naturally objectwise) cartesian.
\end{lemma}
\begin{proof}[Proof sketch]
The proof which Goodwillie gives of the corresponding lemma ([GIII 2.2]) is completely formal.  It relies entirely on taking $\holim$'s of a series of different diagram functors in $\catN$.  All maps in Goodwillie's proof are induced by inclusions of indexing categories and all weak equivalences follow from left cofinality -- similar to the proof of Lemma~\ref{P_nF n-excisive} in the previous section.  Thus Goodwillie's proof of [GIII 2.2] perfectly transports to this framework without any necessary modification.

We do not transcribe the entire proof here since it is rather long.
\end{proof}

Applying the lemma to $F$ an $n$-homogeneous functor, we get $P_{n-1}F \weq \tilde P_{n-1}F$ contractible.  In particular, there is a natural zig-zag:
$$F = P_nF \xleftarrow{\ \weq\ } \hat P_n F \xrightarrow[(1)]{\ \weq\ } 
\holim{\!}_\catN\left(\begin{aligned}
 \xymatrix@R=15pt{\tilde P_{n-1}F \ar[d] \\ R_nF \\ K_nF\ar[u]}
 \end{aligned}\right) \xrightarrow[(2)]{\ \weq\ }
\holim{\!}_\catN\left(\begin{aligned}
 \xymatrix@R=15pt{0 \ar[d] \\ R_nF \\ 0 \ar[u]}
 \end{aligned}\right) = \Omega R_nF$$
where map $(1)$ is due to the lower square in Lemma~\ref{L:thm 2} being cartesian and map $(2)$ is induced by the weak equivalence of diagrams given by the contractions of $K_nF$ and $\tilde P_{n-1}F$ to $0$.

The remainder of the proof of Theorem II deviates slightly from Goodwillie, since our category of rational spectra is $\dg$ rather than $\dgl$-{\it spectra} or $\dgc$-{\it spectra} (the stabilizations of the categories $\dgl_{r-1}$ and $\dgc_r$ respectively).

\begin{proof}[Proof of Theorem II]
Thus far, we have constructed a functor 
$R_n:\mathscr{H}_n(\catM, \catN) \to \mathscr{H}_n(\catM, \catN)$ along with a zig-zag of natural weak equivalences $ \Id_{\mathscr{H}_n} \xleftarrow{\weq} \cdot \xrightarrow{\weq} \Omega R_n$.  We use this to construct functors $B^k:\mathscr{H}_n(\catM, \catN) \to \mathscr{H}_n(\catM, \dg_{r-k})$.  The desired inverse (up to a zig-zag) for $\Linf$ is then given by $\displaystyle \lim_{k\to\infty} B^k = B^\infty$.
We will, however, need to make some sense of this expression.

Recall that $\catN$ is one of $\dg_r$, $\dgl_{r-1}$, or $\dgc_r$.  Note that the functor $\Omega:\catN\to\catN$ factors through $\dg_{r-1}$ as\footnote{By $\Linf:\dg \to \dg_r$ we mean the $r$-reduction functor $\Linf=\mathrm{red}_r:\dg\to\dg_r$.  \\ Recall $\Linf:\dg\to\dgl_{r-1}$  by $V\mapsto [s^{-1}\,\mathrm{red}_{r}(V)]_\mDGL$ and $\Linf:\dg\to\dgc_r$ by $V\mapsto\Lambda\,\mathrm{red}_r(V)$.}  

$$\Omega:\kern -20pt
\xymatrix@C=30pt@R=5pt{
\dg_r \ar[r]^(.55){\Idm}  &  \dg_{r} \ar[r]^(.45){s^{-1}}  &  \dg_{r-1} \ar[r]^{\mathrm{incl}}  
  &   \dg \ar[r]^(.4){\Linf}  &  \dg_r \phantom{Xx!}\\
{}\phantom{Xx!}\dgl_{r-1} \ar[r]^(.6){s[-]_\mDG}  &  \dg_r \ar[r]^(.45){s^{-1}} & \dg_{r-1} \ar[r]^{\mathrm{incl}}  
  &   \dg \ar[r]^(.4){\Linf}  & \dgl_{r-1} \\
{}\phantom{x}\dgc_r \ar[r]^(.6){(-)^\pr}  &  \dg_r \ar[r]^(.45){s^{-1}}  &  \dg_{r-1} \ar[r]^{\mathrm{incl}}  
  &   \dg \ar[r]^(.4){\Linf}  &  \dgc_r \phantom{xx}
}$$

For convenience let us write this factorization as 
$$\Omega:\catN \xrightarrow{\ \,f\ \,} \dg_r \xrightarrow{\ \,s^{-1}\,} \dg_{r-1} 
  \xrightarrow {\ \iota_{r-1} \ } \dg \xrightarrow{\ \,\Linf\,}\catN$$
The functors $f$ and $s^{-1}$ above preserve weak equivalences and homotopy limits and therefore also cartesian cubes.
Thus $(s^{-1}f\,R_n)$ gives a functor to $\mathscr{H}_n(\catM,\dg_{r-1})$.  Let $$B=(s^{-1}f\,R_n):\mathscr{H}_n(\catM,\catN) \to \mathscr{H}_n(\catM,\dg_{r-1})$$
and note that there is a natural zig-zag $(\Linf\, \iota_{r-1} B) = \Omega R_n\xleftarrow{\weq} \cdot \xrightarrow{\weq} \Id_{\mathscr{H}_n}$.  Iterating this process yields a sequence of functors $\bigl\{B^k:\mathscr{H}_n(\catM,\catN) \to \mathscr{H}_n(\catM,\dg_{r-k})\bigr\}$ related for $k>0$ by natural zig-zags  
$(\mathrm{red}_{r-k} B^{k+1})\xleftarrow{\,\weq\,} \cdot \xrightarrow{\,\weq\,} B^k$.  

We will make an argument which is essentially as follows:
Given a functor $F\in\mathscr{H}_n(\catM,\catN)$ we now have zig-zags of natural weak equivalences exhibiting $s^kB^kF\weq \Omega s^{k+1}B^{k+1}F$.~\footnote{Recall that $\Omega:\dg_r\to\dg_r$ by $V\mapsto \mathrm{red}_r(s^{-1} V)$ and $\mathrm{red}_r s^k = s^k\,\mathrm{red}_{r-k}$.}
Thus, given a functor $F\in\mathscr{H}_n(\catM,\catN)$ the sequence\footnote{Note that we have shifted notation: what Goodwillie calls $B^p$ is analogous to our $s^kB^k$.} $\{s^kB^kF:\catM\to\dg_{r}\}_{k>0}$ naturally determines a functor $B^\infty F$ mapping to $\dg_{r}$-spectra.  Of course, a $\dg_{r}$-spectrum is merely a \DG.  The functor $B^\infty F$ is $n$-homogeneous because each of the component functors $s^kB^kF$ are $n$-homogeneous.  By construction, $\Linf B^\infty \weq \Id_{\mathscr{H}_n}$ and an observation using bispectra along the lines of Goodwillie [GIII 2.1] shows that $B^\infty \Linf \weq \Id_{\mathscr{H}_n}$ as well.

To avoid possible unpleasantness involved in the relationship between the categories of $\dg_r$-spectra and $\dg$ we make our explicit argument without ever referring to $\dg_r$-spectra:  

Recall that $(\Linf\, \iota_{r-1} B) \xleftarrow{\weq} \cdot \xrightarrow{\weq} \Id_{\mathscr{H}_n}$, so the functor $(\iota_{r-1} B)$ is a right inverse for $\Linf$ up to a zig-zag of weak equivalences.  Unfortunately $(\iota_{r-1} B)$ is not a functor to $\mathscr{H}_n(\catM,\,\dg)$.  In particular, given $F\in\mathscr{H}_n(\catM,\catN)$, the map $(\iota_{r-1} B)F$ is in general not homogeneous, since the inclusion of categories functor $\iota_{r-1}:\dg_{r-1}\to\dg$ preserves weak equivalences but not homotopy limits.  
However $\iota_{r-1}$ isn't too bad in the sense that it does at least preserve $n$-dimensional homotopy pullbacks of objects which happen to all be at least $(n+r-2)$-reduced (the pullbacks of such diagrams in $\dg_{r-1}$ and $\dg$ coincide).
Furthermore the zig-zags $(\mathrm{red}_{r-k} B^{k+1})\xleftarrow{\,\weq\,} \cdot \xrightarrow{\,\weq\,} B^k$ combine to give natural zig-zags of weak equivalences between $\Linf (\iota_{r-k} B^k) = \Linf (\iota_{r-1}\mathrm{red}_{r-1} B^k)$ and $\Id_{\mathscr{H}_n}$.  The functor $(\iota_{r-k}B^k)$ is a clear improvement over $(\iota_{r-1}B)$ since $\iota_{r-k}$ preserves $n$-dimensional homotopy pullbacks of objects which are all at least $(n+r-k-1)$-reduced.

We combine the $B^k$'s using a construction analogous to the method of spectrification used to make an $\Omega$-spectrum out of a pre-spectrum.  In the following we view all of our functors as maps to $\dg$ rather than $\dg_{r-k}$; however, in order to simplify our notation we neglect to explicitly write the inclusion of categories functors $\iota_{j}:\dg_{j}\to \dg$.  For example, we write only $B^k$ rather than $(\iota_{r-k}B^k)$.

Recall that for $V\in\dg_{i}$ and $i<j$ there is a natural map $\mathrm{red}_j(V) \to V$ given by mapping by the identity above grading $j$ and by the inclusion of a kernel in grading $j$ (this is the map given by the composition of the adjoint pair $(\iota_{j}\, \mathrm{red}_j) \to \Id_{\dg_i}$).  Thus there are natural maps $\mathrm{red}_{r-k}B^{k+1} \to B^{k+1}$.  Let $B^\infty$ be the functor defined by taking the colimit of the $B^k$'s:

\begin{equation}\label{Binf}
B^\infty F = \colim \left(\begin{aligned}
\xymatrix@R=10pt{
B F & B^2 F & B^3 F&  \\
\cdot \ar[u]^(.4){\weq} \ar[d]_(.35){\weq} & 
 \cdot \ar[u]^(.4){\weq} \ar[d]_(.35){\weq} & 
 \cdot \ar[u]^(.4){\weq} \ar[d]_(.35){\weq} &  \cdots \\
\mathrm{red}_{r-1}B^2 F\ar[uur] 
 & \mathrm{red}_{r-2}B^3 F\ar[uur] 
 & \mathrm{red}_{r-3}B^4 F\ar[ur] &
}
\end{aligned}\right)
\end{equation}

Note that the diagram in (\ref{Binf}) is cofibrant (since the diagonal maps are all cofibrations) so its $\colim$ is weakly equivalent to its $\hocolim$.  Thus since each functor in (\ref{Binf}) preserves weak equivalences, so must $B^\infty F$.  Note as well that $\mathrm{red}_jB^\infty F:\catN\to\dg_j$ is $n$-homogeneous for all $j$ because there are weak equivalences $\mathrm{red}_jB^kF\weq \mathrm{red}_jB^{k+1}F$ for $k\ge r-j$ and $\mathrm{red}_jB^{r-j}F = B^{r-j}F$ is $n$-homogeneous.  Therefore $B^\infty F$ itself is $n$-homogeneous.  Furthermore, since for each $k$ there are natural zig-zags of weak equivalences exhibiting $\Linf B^k \weq \Id_{\mathscr{H}_n}$ there is an induced natural zig-zag of weak equivalences giving $\Linf B^\infty \weq \Id_{\mathscr{H}_n}$.

It remains to show that $B^\infty \Linf \weq \Id_{\mathscr{H}_n}$.  However, this follows from the fact that, given a functor $F\in\mathscr{H}_n(\catM,\dg)$, there are clear weak equivalences $B^k\Linf F \weq \mathrm{red}_{r-k}F$.

\end{proof}

\section{Symmetric Multilinear Functors and Theorem III}

We continue our analysis of homogeneous functors.  Let $\catM$ any of the categories $\dg_r$, $\dgl_{r-1}$, or $\dgc_r$ and let $\catN$ be any of the categories $\dg$, $\dg_r$, $\dgl_{r-1}$, or $\dgc_r$.  The main goal of this subsection is to show that certain good $n$-homogeneous functors $H_n:\catN \to \dg$ up to a zig-zag of natural weak equivalences all have the form $H_n(X) \weq \bigl(A\otimes(\Sinf X)^{\otimes n}\bigr)_{\Sigma_n}$ where $A$ is a $\dg$ with $\Sigma_n$-action (by $\Sinf$ we mean the inclusion of categories functor if $\catN=\dg_r$ and the identity functor if $\catN=\dg$).
Our proof of this is given in a slightly different manner than Goodwillie's approach in Calculus III.  In particular, we organize our work differently in order to make use of certain functors which are available to us in our algebraic setting which were not available to Goodwillie in the topological setting and in order to avoid problems introduced by the fact that our category of rational spectra is not defined as the stabilizations of either $\dgl$ or $\dgc$.

The result is proven in three steps which we state as Theorem~III, Lemma~\ref{thmIII part2}, and Lemma~\ref{thmIII part3}.  The first step is to note an equivalence between certain symmetric $n$-multilinear functors and $n$-homogeneous functors.  We write $\mathscr{L}_n(\catN, \dg)$ to denote the category consisting of $n$-multivariate functors $L:\catN\times\cdots\times\catN \to \dg$ which are $1$-homogeneous in each variable (i.e. $n$-multilinear functors) and natural transformations between them; and $\mathscr{L}^\Sigma_n(\catN, \dg)$ for the category of naturally equivariant functors in $\mathscr{L}_n(\catN,\dg^{\boldsymbol{\Sigma}_n})$ (where $\catN\times\cdots\times\catN$ has $\Sigma_n$-action permuting indices).  
More explicitly $\mathscr{L}_n^\Sigma(\catN,\dg)$ consists of $n$-multilinear functors $L$ 
equipped with natural isomorphisms 
$$\pi_L:L(X_1,\dots,X_n) \xrightarrow{\ \cong\ } L(X_{\pi^{-1}(1)},\dots,X_{\pi^{-1}(n)})\ \ \text{for each }\pi\in\Sigma_n$$
satisfying $(\pi\circ\sigma)_L = \pi_L \circ \sigma_L$.  Maps in $\mathscr{L}_n^\Sigma(\catN,\dg)$ are natural transformations between the functors $L$ which commute with the structure maps $\pi_L$.
Formally translating Goodwillie's methods from [GIII], there is a pair of functors $cr_n$ and $\Delta^\Sigma_n$ such that

\begin{thmIII}
The functors $cr_n:\mathscr{H}_n(\catN, \dg)\to \mathscr{L}^\Sigma_n(\catN, \dg)$ and $\Delta^\Sigma_n:\mathscr{L}^\Sigma_n(\catN, \dg)\to\mathscr{H}_n(\catN, \dg)$ are inverse up to a zig-zag of natural weak equivalences.\end{thmIII}

The utility of this theorem follows from the fact that it allows us to transform problems about homogeneous functors into problems about multilinear functors.  In practice, multilinear functors aren't much harder to work with than linear functors since we may work by induction one variable at a time.  The remainder of the work which we do resides in the world of multilinear functors and is considerably easier to prove than Theorem III.

At this point we diverge from Goodwillie in our organization.
Note that $\Sinf:\catM \to \dg$ preserves weak equivalences and homotopy pushouts, so pre-composition with $\Sinf$ yields a functor $\Sinf:\mathscr{H}_n(\dg, \dg) \longrightarrow \mathscr{H}_n(\catM, \dg)$.  
We would like to create an inverse (up to zig-zag) for this functor;  
however, we begin with something slightly simpler.  

Note that $\Sinf$ actually maps $\Sinf:\catM\longrightarrow\dg_r$, so precomposition with it also gives a map 
$\Sigma^\infty_r:\mathscr{H}_n(\dg_r, \dg) \longrightarrow \mathscr{H}_n(\catM, \dg)$. 
Step two of our solution consists of finding an inverse of $\Sigma^\infty_r$ up to a zig-zag of natural weak equivalences.  Actually, we invert the induced transformation of symmetric multilinear functors since that is much easier (and is equivalent by Theorem III).

\begin{lemma}\label{thmIII part2}
The functor $(\Sigma^\infty_r)^{\times n}:\mathscr{L}^\Sigma_n(\dg_r, \dg) \longrightarrow \mathscr{L}^\Sigma_n(\catM, \dg)$ has an inverse up to a zig-zag of natural weak equivalences.
\end{lemma}

With this, we have a chain of equivalences:
$$\mathscr{H}_n(\catM, \dg)\weq\mathscr{L}^\Sigma_n(\catM, \dg) \weq
\mathscr{L}^\Sigma_n(\dg_r, \dg)
$$

The third step in our analysis consists of a classification of symmetric multilinear functors $\mathscr{L}^\Sigma_n(\dg_r, \dg)$ which satisfy the colimit axiom:

\begin{lemma}\label{thmIII part3}
Let $L$ be a symmetric $n$-multilinear functor 
$L \in \mathscr{L}^\Sigma_n(\dg, \dg)$ or $L \in \mathscr{L}^\Sigma_n(\dg_r, \dg)$.  There is a zig-zag of natural weak equivalences exhibiting 
$$L(V_1,\dots,V_n) \weq A\otimes(V_1\otimes\cdots\otimes V_n)$$ 
where $A$ is a \DG\ with $\Sigma_n$-action if either $L$ satisfies the colimit axiom or if the complexes $V_i$ are all finite.
\end{lemma}
Recall that a homotopy functor $F:\catM_1\to \catM_2$ {\em satisfies the colimit axiom}\footnote{Other expressions sometimes used to denote this property are ``$F$ is {\em finitary}'' and ``$F$ is {\em continuous}.''}  if $F$ preserves filtered homotopy colimits.  
That is, for $\mathscr{D}:\catI\to\catM_1$ any filtered diagram in $\catM_1$ the canonical map
$$F\bigl(\hocolim{\!}_{\catM_1}(\mathscr{D})\bigr) \longrightarrow \hocolim{\!}_{\catM_2}F(\mathscr{D})$$
is a weak equivalence (see e.g. [GIII, 5.10]).  Lemma~\ref{thmIII part3} is proven by appealing to [Kuhn] (since the categories $\dg$ and $\dg_r$ are both simplicially enriched).

The desired result follows:

\begin{cor}\label{C:thm 3}
If $H\in\mathscr{H}_n(\catN, \dg)$ satisfies the colimit axiom then there is a zig-zag of natural weak equivalences exhibiting $H(X) \weq (A \otimes (\Sinf X)^{\otimes n})_{\Sigma_n}$ where $A$ is a \DG\ with $\Sigma_n$-action (where $\Sinf$ is the inclusion of categories functor if $\catN=\dg_r$ and the identity if $\catN=\dg$).
\end{cor}

\subsection{Cross Effect, Diagonalization, and Theorem III.}

The results of this subsection formally follow from Goodwillie's constructions in [GIII \S3].  Thus most proofs are omitted, and we give only a rough outline of the constructions necessary.  Next to each lemma or theorem we indicate which lemma or theorem of Goodwillie is analogous.

Recall that $\catN$ is one of $\dg$, $\dgl$, or $\dgc$ and an $n$-cube in $\catN$ is a functor $\mathcal{X}:\mathscr{P}(\underline{n}) \to \catN$, where $\underline{n}$ is the set $\underline{n} = \{1,\dots, n\}$ and for $S$ a set $\mathscr{P}(S)$ is the poset of subsets of $S$ and inclusion maps.  We write 0 for the initial object in $\catN$.  Let $\vec X = (X_1,\dots, X_n)\in \catN^{\times n}$ and $\mathcal{X}(\vec X)$ be the $n$-cube in $\catN$ given by $\mathcal{X}(\vec X):\mathscr{P}(\underline{n}) \to \catN$ by $$\mathcal{X}(\vec X):S \mapsto \!\! \coprod_{s\in (\underline{n}-S)} \!\!\!\! X_s$$
with maps induced by the maps $X_i \to 0$.

\begin{ex}If $\catN = \dg$ and $n = 3$ then $\mathcal{X}(\vec X)$ is the cube:
$$
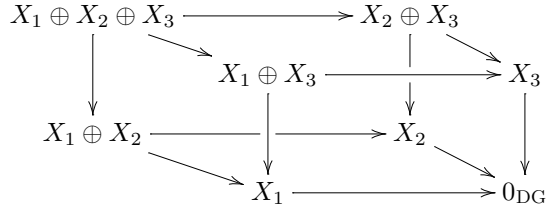
\xymatrix@C=10pt@R=8pt{
X_1 \oplus X_2 \oplus X_3 \ar[rr] \ar[dd] \ar[dr]& & X_2 \oplus X_3 \ar[dr] \ar'[d][dd] & \\
 &  X_1 \oplus X_3 \ar[rr] \ar[dd] & & X_3 \ar[dd] \\
X_1 \oplus X_2 \ar'[r][rr] \ar[dr] & & X_2 \ar[dr] & \\
 &  X_1 \ar[rr] & & 0_\mDG
}$$

If $\catN = \dgl$ and $n=2$ then $\mathcal{X}(\vec X)$ is the square:
$$\xymatrix{
 X_1 \liesum X_2 \ar[r] \ar[d] & X_2 \ar[d] \\
 X_1 \ar[r] & 0_\mDGL 
}$$
\end{ex}

\begin{defn}
Given an $n$-cube $\mathcal{Y}:\mathscr{P}(\underline{n}) \to \catN$ the total homotopy fiber of $\mathcal{Y}$  is  $$\thfib(\mathcal{Y})= \hofib \big(\mathcal{Y}(\emptyset) \longrightarrow \holim(\mathcal{Y}_0)\big)$$ where $\mathcal{Y}_0$ is the restriction of the diagram $\mathcal{Y}$ to the sub-indexing category $\mathscr{P}_0(\underline{n})$,\footnote{Recall that this is the full subcategory of $\mathscr{P}(\underline{n})$ consisting of all non-trivial subsets of $\underline{n}$ -- i.e. it is $\mathscr{P}(\underline{n})$ minus the initial object $\emptyset$.} and the map $\mathcal{Y}(\emptyset) \longrightarrow \holim(\mathcal{Y}_0)$ is the composition of the natural maps
$$\mathcal{Y}(\emptyset) \longrightarrow \lim(\mathcal{Y}_0)\longrightarrow \holim(\mathcal{Y}_0)$$  
\end{defn}

\begin{note} $\thfib(\mathscr{Y}) \weq 0$ if the cube $\mathscr{Y}$ is homotopy cartesian.  There is a dual definition for total homotopy cofiber of $\mathscr{Y}$, and dually $\thcof(\mathscr{Y}) \weq 0$ if $\mathscr{Y}$ is homotopy cocartesian.
\end{note}

It follows from definitions that the total homotopy fiber of the cube $\mathcal{Y}$ is equal to the homotopy limit of the larger diagram 
$$\tilde {\mathcal{Y}}_0:\mathscr{P}_0(\underline{n+1}) \to \catN\mathrm{\ \ by\ \ }\tilde {\mathcal{Y}}_0(S) = \begin{cases}0 &\mathrm{if\ } (n+1)\notin S \\ \mathcal{Y}(S-\{n+1\}) &\mathrm{otherwise}\end{cases}$$
It is also equal to taking the homotopy fiber of the map of total homotopy fibers induced by dividing $\mathcal{Y}$ into two $(n-1)$-cubes $\mathcal{Y} = \bigl(\mathcal{Y}^0 \to \mathcal{Y}^1\bigr)$.  From the previous statement, it is clear that total homotopy fibers may be computed by taking iterated homotopy fibers of homotopy fibers as in the example:

\begin{ex}Let $\mathcal{Y}:\mathscr{P}(\underline{3}) \to \catN$ be the $3$-cube given by $\mathcal{Y}:S \mapsto Y_S$.  Then the total homotopy fiber of $\mathcal{Y}$ may be computed by
$$\xymatrix@C=7.5pt@R=7.5pt{
{}\phantom{XXX)} \save[]+<35pt,0pt> *{\hfib(h)=\thfib(\mathcal{Y})}\restore \ar@{-->}[dr] &&& &&& \\
&\hfib(g_1) \ar[dr]^{h} \ar@{-->}[rr] 
&& \hfib(f_1) \ar[rr]^{g_1} \ar[dr] \ar@{-->}[dd] & & \hfib(f_3) \ar[dr] \ar@{-->}'[d][dd] \\
&& \hfib(g_2) \ar@{-->}[rr]
&  & \hfib(f_2) \ar[rr]^(.38){g_2} \ar@{-->}[dd] & & \hfib(f_4) \ar@{-->}[dd] \\
&&& 
 Y_\emptyset \ar'[r][rr] \ar[dd]^{f_1} \ar[dr] & & Y_{\{3\}} \ar'[d]^(.6){f_3}[dd] \ar[dr]\\
&&{}\phantom{XX}\mathcal{Y}\phantom{XX} \ar@{=}[r] &  {}\phantom{XXXXX}
 & Y_{\{2\}} \ar[rr] \ar[dd]^(.32){f_2} & & Y_{\{2,3\}} \ar[dd]^{f_4} \\
&&& 
 Y_{\{1\}} \ar'[r][rr] \ar[dr] & & Y_{\{1,3\}} \ar[dr] \\
&&&  
 & Y_{\{1,2\}} \ar[rr] & & Y_{\{1,2,3\}}
}$$
\end{ex}

\begin{defn}[Cross Effect]
Define the $n^\mathrm{th}$-cross effect functor $cr_n$ to be $$cr_n:\Funct(\catN_1, \catN_2) \to \Funct(\catN_1^{\times n}, \catN_2)\mathrm{\ \ by\ \ }cr_nF(\vec X) = \thfib\big(F\mathcal{X}(\vec X)\big)$$
where $\catN_1$ and $\catN_2$ are each one of the categories $\dg$, $\dgl_{r-1}$, or $\dgc_r$ and $F:\catN_1 \to \catN_2$.
\end{defn}

Note that the $n^\mathrm{th}$ cross-effect of the functor $F$ is naturally invariant under the action of $\Sigma_n$ permuting the inputs, since such a permutation is equivalent to rotating the cube $\mathcal{X}(\vec X)$, which changes its homotopy fiber only by a natural isomorphism.  The following immediately follow from the formality of Goodwillie's proofs of the corresponding results in the topological setting:

\begin{lemma}[GIII 3.3]If $F:\catN_1 \to \catN_2$ is $n$-excisive, then for $0\le m \le n$ the functor $cr_{m+1}F:\catN_1^{\times (m+1)} \to \catN_2$ is $(n-m)$-excisive in each variable.\end{lemma}

\begin{cor}The $n^\mathrm{th}$ cross effect gives a homotopy functor $cr_n:\mathscr{H}_n(\catN_1, \catN_2) \to \mathscr{L}^\Sigma_n(\catN_1,\catN_2)$ (weak equivalences in $\mathscr{H}_n(\catN_1, \catN_2)$ and $\mathscr{L}^\Sigma_n(\catN_1,\catN_2)$ are the natural weak equivalences).\end{cor}

We would now like to construct a homotopy inverse for $cr_n$:

\begin{defn}[Diagonalization]
Define the diagonalization functors $\Delta_n$ and $\Delta^\Sigma_n$ by\footnote{Goodwillie calls these functors $\Delta$ and $\Delta_n$.}
\begin{itemize}
\item $\Delta_n:\Funct(\catN_1^{\times n}, \catN_2) \to \Funct^{\Sigma_n}(\catN_1, \catN_2)$ by $\Delta_n G (X) = G(X, \dots, X)$
\item $\Delta^\Sigma_n:\Funct(\catN_1^{\times n}, \catN_2) \to \Funct(\catN_1, \catN_2)$ by $\Delta^\Sigma_n G (X) = G(X, \dots, X)_{\Sigma_n}$
\end{itemize}
\end{defn}

\begin{note}
 In Goodwillie's setting, he uses $\Delta^\Sigma_n G(X) = G(X,\dots,X)_{\hS_n}$.  Since orbits is already a homotopy functor on $\dg^{\boldsymbol{\Sigma}_n}$, it is weakly equivalent to the homotopy orbits (see \S\ref{symm dg}).  Thus we may use the more simply defined orbits rather than homotopy orbits in our rational $\Delta^\Sigma_n$.
\end{note}

\begin{lemma}[\lbrack GII, 3.4\rbrack] 
If $G:\catN_1^{\times n}\to \catN_2$ is $(x_1,\dots,x_n)$-excisive, then $\Delta_nG$ is $(\sum x_i)$-excisive.
\end{lemma}

\begin{lemma}[\lbrack GIII, 3.1\rbrack] 
If $G:\catN_1^{\times n}\to \catN_2$ is $(1,\dots,1)$-reduced, then $\Delta_nG$ is $n$-reduced.
\end{lemma}

\begin{cor}
The functor $\Delta_n$ is a homotopy functor $\Delta_n:\mathscr{L}^\Sigma_n(\catN_1,\catN_2) \to \mathscr{H}_n(\catN_1, \catN_2)$ and $\Delta^\Sigma_n$ is a homotopy functor $\Delta^\Sigma_n:\mathscr{L}^\Sigma_n(\catN_1,\dg) \to \mathscr{H}_n(\catN_1, \dg)$.
\end{cor}

We are now ready to state Theorem III:

\begin{thmIII}[\lbrack GIII, 3.5\rbrack]
The functors $cr_n:\mathscr{H}_n(\catN, \dg)\to \mathscr{L}^\Sigma_n(\catN, \dg)$ and $\Delta^\Sigma_n:\mathscr{L}^\Sigma_n(\catN, \dg)\to\mathscr{H}_n(\catN, \dg)$ are inverse up to a zig-zag of natural weak equivalences.\end{thmIII}

The proof of this is omitted since it follows directly from [GIII].

\subsection{The Proofs of Lemmas~\ref{thmIII part2} and \ref{thmIII part3}.}

Recall that $\Sigma^\infty_r$ is the functor $\Sigma^\infty_r:\mathscr{H}_n(\dg_r, \dg) \longrightarrow \mathscr{H}_n(\catM, \dg)$ given by precomposition with $\Sinf:\catM\to\dg_r$.
We begin by proving Lemma~\ref{thmIII part2}.

\begin{note} Lemma~\ref{thmIII part2} is analogous to [GIII 3.8]: If $F$ is $n$-homogeneous, then $F$ is determined by $F\circ \Sigma$.\end{note}

\begin{thmIII part2}
The functor $(\Sigma^\infty_r)^{\times n}:\mathscr{L}^\Sigma_n(\dg_r, \dg) \longrightarrow \mathscr{L}^\Sigma_n(\catM, \dg)$ has an inverse up to a zig-zag of natural weak equivalences.
\end{thmIII part2}
\begin{proof}
Note that the functor $\Sigma:\catM \to \catM$ factors through $\dg_r$ as
$$\Sigma:\kern -3pt
\xymatrix@C=30pt@R=5pt{
\dgl_{r-1} \ar[r]^(.56){\Sigma^\infty}  &  \dg_{r} \ar[r]^(.45){s}  &  
   \dg_{r} \ar[r]^(.45){\freeL_{s^{-1}(-)}}  &  \dgl_{r-1}\\
{}\phantom{{}_{-1}}\dgc_{r} \ar[r]^(.56){\Sigma^\infty}  &  \dg_r \ar[r]^(.45){s} & 
   \dg_{r} \ar[r]^(.45){[-]_\mDGC}  & \dgc_{r}\phantom{{}_{-1}}
}$$
For convenience, write this as
$$\Sigma:\kern -3pt
\xymatrix@C=30pt{
\catM\ar[r]^(.46){\Sigma^\infty}  &  \dg_{r} \ar[r]^(.45){s}  &  
   \dg_{r} \ar[r]^(.45){k}  &  \catM
}$$
We use this factorization in order to prove the lemma for the case $n=1$ similar to the way in which we used the factorization of $\Omega$ in our proof of Theorem II.  For general $n$, the lemma reduces to this case.

Let $L\in \mathscr{L}_1(\catM, \dg)$ be a linear functor.  Given $X\in \catM$, consider the diagrams
$$\xymatrix@R=7.5pt@C=7pt{
 X \ar[rr] \ar[dd] && cX \ar[dd]      & &   L(X) \ar[rr] \ar[dd] && L(cX) \ar[dd] \\
   &\ar@{=>}[dr]|(.05)h&  {\quad} \ \ar@{|~|>}[rr]^{L} & & \ \ \quad  &\ar@{=>}[ul]|(.05)h& \\
 0 \ar[rr]         && \Sigma X        & &   L(0) \ar[rr]         && L(\Sigma X) }$$

The left square is homotopy cocartesian and $L$ is linear, so the right square is homotopy cartesian with $L(0) \weq 0 \weq L(cX)$.  Thus there are natural weak equivalences 
$$L(X) \xrightarrow{\ \weq\ } \holim{\!}_\dg\left(\kern -2pt\begin{aligned}
 \xymatrix@R=15pt{L(cX) \ar[d] \\ L(\Sigma X) \\ L(0)\ar[u]}
 \end{aligned}\kern -2pt\right) \xrightarrow{\ \weq\ } \holim{\!}_\dg\left(\kern -5pt\begin{aligned}
 \xymatrix@R=15pt{0\vphantom{L(} \ar[d] \\ L(\Sigma X) \\ 0\vphantom{L(}\ar[u]}
 \end{aligned}\kern -5pt\right) = \Omega L(\Sigma X)$$
Recall that for $V\in \dg$, we have $\Omega V = s^{-1} V$.  Let $K L:\dg_r \to \dg$ be the functor given by $K L(V) = s^{-1}(L \circ k) (sV)$. 
The functor $k:\dg_r \to \catM$ preserves homotopy cocartesian squares and weak equivalences; therefore $K L:\dg_r \to \dg$ is linear.  By construction there is a natural weak equivalence
$$L(X) \xrightarrow{\ \weq\ } \Omega L \bigl((k s \Sinf) X\bigr) = (KL)(\Sinf X) = (\Sigma^\infty_r K L) (X)$$ 
In other words, $K:\mathscr{L}^\Sigma_n(\catM, \dg) \longrightarrow \mathscr{L}^\Sigma_n(\dg_r, \dg)$ is a right inverse (up to natural weak equivalence) of the functor $\Sigma^\infty_r:\mathscr{L}^\Sigma_n(\dg_r, \dg) \longrightarrow \mathscr{L}^\Sigma_n(\catN, \dg)$.

Now let $\hat L\in\mathscr{L}_1(\dg_r, \dg)$ be a linear functor.  We wish to exhibit a zig-zag of natural weak equivalences between $\hat L(V)$ and 
$(K\Sigma^\infty_r \hat L)(V) = K(\hat L \Sinf)(V) = s^{-1} ((\hat L \Sinf) \circ k )(s V) = s^{-1} \hat L (s V)$.  However, since $\hat L$ is linear, this follows from the diagrams
$$\xymatrix@R=7.5pt@C=7pt{
 V \ar[rr] \ar[dd] && \mathrm{c}V \ar[dd]  & &   \hat L(V) \ar[rr] \ar[dd] && \hat L(\mathrm{c}V) \ar[dd] \\
&\ar@{=>}[dr]|(.05)h&  {\quad} \ \ar@{|~|>}[rr]^{\hat L} & & \ \quad  &\ar@{=>}[ul]|(.05)h& \\
 0 \ar[rr]         && \Sigma V   & &   \hat L(0) \ar[rr]   && \hat L(\Sigma V) }$$
as above.\footnote{Note $\Omega = s^{-1}$ and $\Sigma = s$ in $\dg$.}

For the general case where $L\in\mathscr{L}_n^\Sigma(\catM,\dg)$ the result follows from applying the above construction to each of the variables of $L$ in turn.

\end{proof}

Recall Lemma~\ref{thmIII part3}:

\begin{thmIII part3}
Let $L$ be a symmetric $n$-multilinear functor 
$L \in \mathscr{L}^\Sigma_n(\dg, \dg)$ or $L \in \mathscr{L}^\Sigma_n(\dg_r, \dg)$.  There is a zig-zag of natural weak equivalences exhibiting 
$$L(V_1,\dots,V_n) \weq A\otimes(V_1\otimes\cdots\otimes V_n)$$ 
where $A$ is a \DG\ with $\Sigma_n$-action if either $L$ satisfies the colimit axiom or if the complexes $V_i$ are all finite.
\end{thmIII part3}

Rather than give a proof of this theorem using our models for homotopy limits and colimits, we note that it follows from [Kuhn 2.5] since the categories $\dg$ and $\dg_r$ are both simplicially enriched and proper.

\begin{note}
Let $L:\dg\to\dg$ be a linear functor.  There is a very simple way to make an explicit (non-natural) weak equivalence $L(V) \weq L(1_\mDG)\otimes V$.~\footnote{Recall that $1_\mDG$ is the \DG\ which has $\Q$ in degree 0 and 0 elsewhere.}  Begin with the diagram 
$$\xymatrix@R=7.5pt@C=7pt{
1_\mDG \ar[rr] \ar[dd] && \mathrm{c}1_\mDG \ar[dd]  
  & &   L(1_\mDG) \ar[rr] \ar[dd] && L(\mathrm{c}1_\mDG) \ar[dd] \\
&\ar@{=>}[dr]|(.05)h&  {\quad} \ \ar@{|~|>}[rr]^{L} & & \ \quad  &\ar@{=>}[ul]|(.05)h& \\
 0_\mDG \ar[rr] && s\otimes 1_\mDG   
  & &   L(0_\mDG) \ar[rr]   && L(s\otimes 1_\mDG) 
}$$
Since $\dg$ is a stable model category, the right diagram above is also a homotopy cocartesian diagram.  Thus there are natural weak equivalences 
$$L(s) = L(s\otimes 1_\mDG) \xleftarrow{\ \weq\ } 
\hocolim\left(\vcenter{\xymatrix@R=10pt{
 L(\mathrm{c}1_\mDG) \\ L(1_\mDG) \ar[d] \ar[u] \\ L(0_\mDG)}}\right)
\xleftarrow{\ \weq\ } 
\hocolim\left(\vcenter{\xymatrix@R=10pt{
 L(0_\mDG) \\ L(1_\mDG) \ar[d] \ar[u] \\ L(0_\mDG)}}\right) =
s\otimes L(1_\mDG)$$
Inducting on $k$ in the diagrams 
$$\xymatrix@R=7.5pt@C=6pt{
s^k \ar[rr] \ar[dd] && \mathrm{c}s^k \ar[dd]  
  & &   L(s^k) \ar[rr] \ar[dd] && L(\mathrm{c}s^k) \ar[dd] \\
&\ar@{=>}[dr]|(.05)h&  {\quad} \ \ar@{|~|>}[rr]^{L} & & \ \quad  &\ar@{=>}[ul]|(.05)h& \\
 0_\mDG \ar[rr] && s\otimes s^k  
  & &   L(0_\mDG) \ar[rr]   && L(s\otimes s^k) 
}$$
we extend to natural weak equivalences $L(s^k) \xleftarrow{\,\weq\,} \bigl(s^k\otimes L(1_\mDG)\bigr)$ for $k$ positive and negative.

If $V$ is a \DG, then by choosing basis elements for $V$ we may write it as a (possibly infinite) sequence of homotopy pushouts of \DG\,s $s^k$ analogous to the way that we may build a cellular complex inductively by adding cells of higher and higher dimension (except that now we inductively add cells of greater and greater positive and negative dimensions).\footnote{For example if $V = \bigl(\Q v_1 \oplus \Q v_2,\ d(v_2) = v_1\bigr)$ then $V$ is the homotopy pushout of 
$s \xleftarrow{\,\Idm\,} s \xrightarrow{\,\ \,} 0_\mDG$. 
}  Since $L$ commutes with filtered homotopy colimits, this yields zig-zags of natural weak equivalences $L(V) \weq \bigl(V\otimes L(1_\mDG)\bigr)$.

If $L:\dg_r\to\dg$ is a linear functor.  We may work similarly, beginning with the diagram 
$$\xymatrix@R=7.5pt@C=6pt{
s^r \ar[rr] \ar[dd] && \mathrm{c}s^r \ar[dd]  
  & &   L(s^r) \ar[rr] \ar[dd] && L(\mathrm{c}s^r) \ar[dd] \\
&\ar@{=>}[dr]|(.05)h&  {\quad} \ \ar@{|~|>}[rr]^{L} & & \ \quad  &\ar@{=>}[ul]|(.05)h& \\
 0_\mDG \ar[rr] && s\otimes s^r  
  & &   L(0_\mDG) \ar[rr]   && L(s\otimes s^r) 
}$$
which gives $L(s\otimes s^r)\xleftarrow{\,\weq\,}\bigl(s\otimes L(s^r)\bigr)$.  In particular, for $k\ge r$ we may induct to get 
$$L(s^k) \xleftarrow{\,\weq\,} \bigl(s^k\otimes s^{-r}L(s^r)\bigr)$$
The colimit axiom then gives $L(V) \weq \bigl(V\otimes s^{-r}L(s^r)\bigr)$ for all $V\in\dg_r$.
\end{note}

\chapter{Derivatives of the Rational Identity Functor and a Comparison Theorem}\label{C:d Id}
\markboth{Ben Walter}{II.8 Derivatives of the Identity Functor}

\section{DGL to DGL}\label{DGL to DGL}

Consider the identity functor $\Id_\dgl:\dgl_{r-1}\to\dgl_{r-1}$.  By our constructions of the previous section, there is a universal approximating tower of fibrations given by 
$$\cdots \to P_n\Id_\dgl \to P_{n-1}\Id_\dgl \to \cdots \to P_0\Id_\dgl$$ 
with fibers $D_n\Id_\dgl \to P_n\Id_\dgl \to P_{n-1}\Id_\dgl$.  In particular given $L$ any free \DGL, the rational Taylor tower of the identity functor evaluated at $L$ yields a tower of fibrations in $\dgl$ converging to $L$ with fibers (naturally) of the form
\begin{align*}
D_n\Id_\dgl (L)
&\weq \Omega^\infty_\dgl(\partial_n\!\Id_\dgl \otimes (\Sigma^\infty_\dgl L)^{\otimes n})_{\Sigma_n} \\
& \weq \Big[s^{-1} \Big(\partial_n\!\Id_\dgl \otimes \big(s(L)^\ab\big)^{\otimes n}\Big)_{\Sigma_n}\Big]_\mDGL \\
& \cong \Big[\Big(s^{n-1} \partial_n\!\Id_\dgl \otimes \big((L)^\ab\big)^{\otimes n}\Big)_{\Sigma_n}\Big]_\mDGL
\end{align*}
 where the $\partial_n\!\Id_\dgl$ are \DG\,s with ${\Sigma_n}$-action and $\Sigma_n$ acts on $\big[s(L)^\ab\big]_\mDG^{\otimes n}$ and $\big[(L)^\ab\big]_\mDG^{\otimes n}$ by permutation of elements with signs according to the Koszul convention, and on $s^{n-1}$ by multiplication by $(-1)^{{\mathrm sgn}(\sigma)}$.  [The $\Sigma_n$-equivariant isomorphism on the last line is given by the desuspension ($s^{-1}$) of the $\Sigma_n$-equivariant map $(s\,l_1\otimes \cdots \otimes s\,l_n) \mapsto (-1)^{k} s^n\otimes (l_1\otimes \cdots \otimes l_n)$, where $(-1)^k$ is the sign incurred under the Koszul convention by moving all of the $s$'s to the beginning of the expression ($k=\sum_{j=1}^{n-1}\sum_{i=1}^j |l_i|$) and $\Sigma_n$ acts diagonally on the right with an action on $s^n$ given by multiplication by $(-1)^{{\mathrm sgn}(\sigma)}$.]

Another natural tower of fibrations associated to any free \DGL\ is the one induced by its lower central series.  Recall that the lower central series of a \DGL\ is the (natural) sequence 
$$\xymatrix@C=0pt@R=7pt{
 L  & \supset & [L,L] & \supset 
   & \bigl[L,[L,L]\bigr] & \supset & \cdots \\
 \Gamma^1(L) \ar@{=}[u] & & \Gamma^2(L) \ar@{=}[u] 
   & & \Gamma^3(L) \ar@{=}[u] & &
}
$$ 
This induces a (natural) tower
$$\xymatrix@R=10pt@C=6pt{
{}\vphantom{|}\save[] *{\vdots} \restore   
      & && & {}\vphantom{|}\save[] *{\vdots} \restore \ar[d]  & \\
\Gamma^{n+1}(L) \ar[u]     & && &   
  \raisebox{2pt}{$L$}\!\diagup\,\raisebox{-2pt}{$\Gamma^{n+1}(L)$} \ar[d]  &  
   = B_{n}(L)  \\
\Gamma^{n}(L) \ar[u] & && &   
  \raisebox{2pt}{$L$}\!\diagup\,\raisebox{-2pt}{$\Gamma^{n}(L)$} \ar[d]  & 
  = B_{n-1}(L)\kern-9pt  \\
{}\vphantom{|}\save[] *{\vdots} \restore \ar[u] 
     & \ar@{|~|>}[rr] &{}\phantom{XXXX}& & 
     {}\vphantom{|}\save[] *{\vdots} \restore \ar[d] \\
\Gamma^3(L) \ar[u]     & && &   
    \raisebox{2pt}{$L$}\!\diagup\,\raisebox{-2pt}{$\Gamma^3(L)$} \ar[d] & = B_2(L) \\
\Gamma^2(L) \ar[u]     & && &   
    \raisebox{2pt}{$L$}\!\diagup\,\raisebox{-2pt}{$\Gamma^2(L)$}  & = B_1(L) 
}$$
with natural maps $L\to B_n(L)$ given by the quotient maps.
For $L\in \Fdgl_{r-1}$ (with $r\ge 2$) this is a tower of fibrations in $\dgl$ whose limit is $L$ because the objects $\Gamma^n(L)$ are increasingly reduced -- $\Gamma^n(L)$ is $n(r-1)$-reduced since $L$ is $(r-1)$-reduced.  When $L$ is free, the tower of objects $B_n(L)$ is called the {\em bracket-length filtration} of $L$ since for $L = (\freeL_V,\, d)$ we have 
$$B_n(L) = \bigl((T^{\le n} V) \cap \freeL_V,\  d= d_0+\cdots+d_n)$$
That is $B_n(L)$ consists of the elements of $\freeL_V$ with bracket-length $\le n$.\footnote{By ``bracket length $n$'' we mean a nested bracket expression of $n$ elements -- so bracket length 1 means there is no bracket at all.}

The fibers of this tower for $L$ free are well known -- they are $H_n(L) \to B_n(L) \to B_{n-1}(L)$ where $H_n(L)$ is ``all bracket expressions in $L$ of length exactly $n$.''  More precisely, $H_n(L)$ is naturally (equivariantly) isomorphic to the \DGL\ with trivial bracket
$$H_n(L) \cong \Bigl[\Bigl(Lie(n) \otimes \bigl((L)^\ab\bigr)^{\otimes n}\,\Bigr)_{\Sigma_n}\Bigr]_\mDGL$$ 
where $Lie(n)$ is the \DG\ concentrated in degree 0 generated by all abstract bracket expressions of $n$ elements (e.g. $[x_1,[x_2,\cdots, [x_{n-1},x_n]\cdots]]$) modulo anti-symmetry and the Jacobi identity; $\Sigma_n$ acts on $\bigl((L)^\ab\bigr)^{\otimes n}$ by permutation of terms with signs according to the Koszul convention, and the action of $\Sigma_n$ on $Lie(n)$ is generated by the permutation of the elements in a bracket expression (e.g. 
$[x_1,\cdots, [x_{n-1},x_n]] \mapsto [x_{\sigma(1)},\cdots, [x_{\sigma(n-1)}, x_{\sigma_n}]]$ with \underline{no} negative signs).

Note that $H_n$ is an $n$-homogeneous functor $H_n:\Fdgl_{r-1}\to \dgl_{n(r-1)}$.  Thus each $B_n$ is an $n$-excisive functor.\footnote{This statement uses the assumption that $r\ge2$.}  So this tower is an approximating tower of fibrations of $n$-excisive functors converging to the identity functor in the strong sense that the maps $L \to H_n(L)$ are vector space isomorphisms up to degree $n(r-1)$.  In particular the functors $\Id_\mDGL$ and $H_n$ agree up to degree $n$.  A standard theorem of homotopy calculus [GIII 1.6] implies that $P_n\Id_\mDGL \weq P_nH_n = H_n$.  Thus the bracket-length filtration is the rational Taylor tower of the identity functor evaluated at the free \DGL\ $L$.

Comparing $H_n(L)$ with $D_n\Id_\dgl(L)$ (for $L$ free) we get the following theorem:

\begin{thm}
The rational derivatives of the identity functor $\Id_\dgl:\dgl_{r-1} \longrightarrow \dgl_{r-1}$ are $Lie(n)$ graded in degree $(1-n)$ with $\Sigma_n$-action twisted by the sign of permutations.
\end{thm}

Note that it is enough to consider only $L$ free because the identity functor is a homotopy functor and every \DGL\ is quasi-isomorphic to its cofibrant replacement (which is a free \DGL).

\section{DGC to DGL}\label{DGC to DGL}

We investigate the rational Taylor tower of the functor $\L:\dgc_r \to \dgl_{r-1}$ (for $r\ge 2$) in the same manner.  The rational Taylor tower of $\L$ evaluated at $C$ a \DGC\ yields a tower of fibrations in $\dgl$ with fibers (naturally) weakly equivalent to:
\begin{align*}
D_n\L(C) 
&\weq \Omega^\infty_\dgl\bigl(\partial_n\!\L \otimes (\Sigma^\infty_\dgc C)^{\otimes n}\bigr)_{\Sigma_n} \\
& \cong \Bigl[s^{-1} \bigl(\partial_n\!\L \otimes [C]_\mDG^{\otimes n}\bigr)_{\Sigma_n}\,\Bigr]_\mDGL 
\end{align*}
Where $\Sigma_n$ acts on $[C]_\mDG^{\otimes n}$ by permutation of elements with signs according the the Koszul convention and $\Sigma_n$ acts trivially on $s^{-1}$.

Recall that $\L(C)$ is defined to be the free \DGL\ $\L(C) = (\freeL_{s^{-1}[C]_\mG},\, d_{\freeL} + d_{\Delta})$ so it has another approximating tower of fibrations in $\dgl$ given by its bracket-length filtration.  The fibers of this tower are $H_n(C) \to B_n(C) \to B_{n-1}(C)$ naturally (equivariantly) isomorphic to 
\begin{align*}
H_n(C) 
&= \Bigl[\Bigl(Lie(n) \otimes \bigl[\bigl(\L(C)\bigr)^\ab\bigr]^{\otimes n}_\mDG\Bigr)_{\Sigma_n}\,\Bigr]_\mDGL  \\
&\cong \Bigl[\Bigl(Lie(n)\otimes (s^{-1}[C]_\mDG)^{\otimes n}\Bigr)_{\Sigma_n}\,\Bigr]_\mDGL  \\
&\cong \Bigl[\Bigl(s^{-n}Lie(n)\otimes [C]_\mDG^{\otimes n}\Bigr)_{\Sigma_n}\,\Bigr]_\mDGL
\end{align*}

Again, $H_n$ is an $n$-homogeneous functor to increasingly reduced \DGL\,s, so the $B_n$ are $n$-excisive.  Since the $B_n$ form an approximating tower of fibrations of $n$-excisive functors which converges to $\L$ in the strong sense that the maps $\L C \to H_n(C)$ are vector space isomorphisms up to degree $n(r-1)$.  In particular $\L$ and $H_n$ agree up to order $n$, so $P_n\L \weq P_n H_n = H_n$ and thus the $H_n$ give the rational Taylor tower of $\L$.

Comparing $H_n(C)$ and $D_n\L(C)$ while keeping track of the $\Sigma_n$-actions,\footnote{This time, there is a twist by sign of permutation on the $\Sigma_n$-action of the $s^{-n}$ in $H_n(C)$ since the $s^{-n}$ came from $\bigl(s^{-1}[C]_\mDG\bigr)^{\otimes n}$ which has signs on its $\Sigma_n$-action according to the Koszul convention.} we have now shown: 

\begin{thm} The rational derivatives of the functor $\L:\dgc_r \to \dgl_{r-1}$ are $Lie(n)$ graded in degree $(1-n)$ with $\Sigma_n$-action twisted by the sign of permutations.
\end{thm}

\section{DGC to DGC}\label{DGC to DGC}

We outline a method for computing the derivatives of functors $\dgc_r\to\dgc_r$ using the derivatives of functors $\dgc_r\to\dgl_{r-1}$.
Let $F:\dgc_r\to\dgl_{r-1}$ be a homotopy functor.

Recall that the functor $\C:\dgl_{r-1}\to\dgc_r$ is a right Quillen adjoint which preserves all weak equivalences and commutes with sequential homotopy colimits (Lemma~\ref{Linf C seq hocolim}).  Thus post-composition with $\C$ preserves Taylor towers of homotopy functors by Theorem~\ref{tower preservation}.  In particular $D_n(\C F) \weq \C (D_nF)$.

%

Applying this to the functor $\L:\dgc_r\to\dgl_{r-1}$, it follows that there are zig-zags of natural weak equivalences
\begin{align*}
D_n(\C\L)(C) &\weq \C(D_n\L)(C) \\
 & \weq \C \, \Omega^\infty_\dgl \Bigl(\partial_n\!\L
            \otimes(\Sigma^\infty_\dgc C)^{\otimes n}\Bigr)_{\Sigma_n} \\
 & \weq \Omega^\infty_\dgc\bigl(\partial_n\!\L\otimes (\Sigma^\infty_\dgc C)^{\otimes n}\bigr)_{\Sigma_n} \notag
\end{align*}
Thus the derivatives of $\C\L$ are the same as those of $\L$.  Recall that there is a natural weak equivalence $\Id_\dgc \xrightarrow{\,\weq\,} \C\L$.  This implies that the towers of $\Id_\dgc$ and $\C\L$ are the same (or at least weakly equivalent).  In particular, $\Id_\dgc$ and $\C\L$ have the same derivatives.

\begin{thm}
The rational derivatives of the identity functor $\Id_\dgc:\dgc_r \to \dgc_{r}$ are $Lie(n)$ graded in degree $(1-n)$ with $\Sigma_n$-action twisted by the sign of permutations.
\end{thm}

\begin{note}
 We could also have attempted to analyze the rational Taylor tower of $\Id_\dgc:\dgc_r \to \dgc_r$ similar to the way in which we anaylzed the tower of $\Id_\dgl:\dgl_{r-1}\to\dgl_{r-1}$.  Associated to every cofree \DGC\ is a word-length filtration:
$$\cdots \xleftarrow{\ } \Lambda^{\le n} (C)^\pr \xleftarrow{\ } \Lambda^{\le (n-1)} (C)^\pr \xleftarrow{\ } \cdots \xleftarrow{\ } \Lambda^1 (C)^\pr \to \Q$$
which converges to the cofree \DGC\ in the strong sense that the maps $\Lambda^{\le n}(C)^\pr \to C$ are vector space isomorphisms for an increasing range of degrees.  Unfortunately, the arrows in this tower map in the wrong direction for it to possibly be a rational Taylor tower.  It could perhaps be useful in the analysis of a dual homotopy calculus of functors, however.
\end{note}

More generally, if we are able to compute the towers of all functors $\dgc_r \to \dgl_{r-1}$ then the tower of any functor $F:\dgc_r \to \dgc_r$ is given by $\C\bigl(\text{\it Tower of }(\L F:\dgc_r \to \dgl_{r-1})\bigr)$ since there is a natural weak equivalence $F\xrightarrow{\,\weq\,}\C\L F$.  In particular, an argument as above shows that the derivatives of $F$ are the same as those of $\L F$.

\section{DGL to DGC}

We analyse the tower of the functor $\C:\dgl_{r-1}\to \dgc_r$ similar to the way in which we analysed $\Id:\dgc_r\to\dgc_r$ in the previous section.

Post-composition with $\C$ also preserves Taylor towers of functors $F:\dgl_{r-1} \to \dgl_{r-1}$.  In particular it preserves homogeneous layers.
Applying this to $\Id_\dgl:\dgl_{r-1}\to\dgl_{r-1}$ yields zig-zags of natural weak equivalences
\begin{align*}
D_n\C(L) &= D_n(\C\circ\Id_\dgl)(L) \\
 &\weq \C (D_n\Id_\dgl)(L) \\
 &\weq \C\,\Omega^\infty_\dgl\bigl(\partial_n\!\Id_\dgl \otimes (\Sigma^\infty_\dgl L)^{\otimes n}\bigr)_{\Sigma_n} \\
 &\weq \Omega^\infty_\dgc \bigl(\partial_n\!\Id_\dgl \otimes (\Sigma^\infty_\dgl L)^{\otimes n}\bigr)_{\Sigma_n}
\end{align*}
So the derivatives of $\C$ are the same as those of $\Id_\dgl$.

\begin{thm}
The rational derivatives of the functor $\C:\dgl_{r-1} \to \dgc_{r}$ are $Lie(n)$ graded in degree $(1-n)$ with $\Sigma_n$-action twisted by the sign of permutations.
\end{thm}

\section{A Comparison Theorem for Derivatives}\label{Id alternate}

We could also have computed the derivatives of $\L$ and $\Id_\dgc$ from the derivatives of $\Id_\dgl$ and $\C$ by applying pre-composition with $\L:\dgc_r\to\dgl_{r-1}$.  
Since $L$ is a Quillen left adjoint which preserves weak equivalences, pre-composition with it preserves rational Taylor towers by Theorem~\ref{tower preservation}.
In particular $D_n(F\L) \weq (D_nF)\L$ for either $F:\dgl_{r-1} \to \dgc_r$ or $F:\dgl_{r-1} \to \dgl_{r-1}$.

Applying this to the functor $\C:\dgl_{r-1}\to\dgc_r$, it follows that there are zig-zags of natural weak equivalences
\begin{align*}
D_n\Id_\dgc(C) &\weq D_n(\C\L)(C) \\ 
 &\weq (D_n\C)\L(C) \\
 &\weq \Omega^\infty_\dgc \Bigl(\partial_n\!\C\otimes\bigl(\Sigma^\infty_\dgl \L C\bigr)^{\otimes n}\Bigr)_{\Sigma_n} \\
 &\weq\Omega^\infty_\dgc\Bigl(\partial_n\!\C\otimes\bigl(\Sigma^\infty_\dgc C\bigr)^{\otimes n}\Bigr)_{\Sigma_n}
\end{align*}
So the derivatives of $\Id_\dgc$ are the same as those of $\C$.

Applying this to the functor $\Id_\dgl:\dgl_{r-1}\to\dgl_{r-1}$, it follows that
\begin{align*}
D_n(\L)(C) &= D_n(\Id_\dgl \L)(C) \\
 &\weq (D_n\Id_\dgl)\L(C) \\
 &\weq \Omega^\infty_\dgl \Bigl(\partial_n\!\Id_\dgl \otimes \bigl(\Sigma^\infty_\dgl(\L C)\bigr)^{\otimes n}\Bigr)_{\Sigma_n} \\
 &\weq \Omega^\infty_\dgl \bigl(\partial_n\!\Id_\dgl\otimes(\Sigma^\infty_\dgc C)^{\otimes n}\bigr)_{\Sigma_n}
\end{align*}
So the derivatives of $\L$ are the same as those of $\Id_\dgl$.

More generally, we have the following theorem (which has an obvious generalization to any pair of Quillen equivalent model categories $F:\catB\rightleftarrows\catC:U$):

\begin{thm}Let $\catM$ be either $\dgc_r$ or $\dgl_{r-1}$.
\begin{enumerate}
\item\label{CF = F} 
 If $F:\catM \to \dgl_{r-1}$ is a homotopy functor, then post-composition with $\C$ preserves rational Taylor towers.  In particular the derivatives of $F$ are the same as those of $\C F:\catM \to \dgc_r$.
\item\label{LF = F} 
 If $F:\catM \to \dgc_{r}$ is a homotopy functor, then post-composition with $\L$ preserves rational Taylor towers.  In particular the derivatives of $F$ are the same as those of $\L F:\catM\to \dgl_{r-1}$.
\item\label{FC = F}
 If $F:\dgc_{r} \to \catM$ is a homotopy functor, then pre-composition with $\C$ preserves rational Taylor towers.  In particular the derivatives of $F$ are the same as those of $F \C:\dgl_{r-1} \to \catM$.
\item\label{FL = F}
 If $F:\dgl_{r-1} \to \catM$ is a homotopy functor, then pre-composition with $\L$ preserves rational Taylor towers.  In particular the derivatives of $F$ are the same as those of $F \L:\dgc_r \to \catM$.
\end{enumerate}
\end{thm}
\begin{proof}[Proof Sketch]
(\ref{CF = F}) follows from the facts that $\C$ is a Quillen right adjoint which preserves all weak equivalences and also $\C \Omega^\infty_\dgl \weq \Omega^\infty_\dgc$.

(\ref{LF = F}) follows from (\ref{CF = F}) and the fact that the tower of $F$ is the same as that of $\C \L F$.

(\ref{FC = F}) follows from (\ref{FL = F}) and the fact that the tower of $F$ is the same as that of $F\C\L$.

(\ref{FL = F}) follows from the facts that $\L$ is a Quillen left adjoint which preserves all weak equivalences and also $\Sigma^\infty_\dgl \L \weq \Sigma^\infty_\dgc$.
\end{proof}

\chapter{Jets:  A Preview}
\markboth{Ben Walter}{II.9  Jets}

Recall from Chapter 1 that the ``jet'' of a functor is the symmetric sequence of spectra given by the functor's derivatives along with all of the necessary structure maps required to recover the approximating Taylor tower of the functor.

\section{Motivating Example}\label{sec:motivate}

Consider a functor $F:\top_*\to\sp$  which has only two nontrivial, homogeneous layers, say $D_nF$ and $D_mF$, for $n<m$.\footnote{Bactrians are bigger than dromedaries.}  The approximating Taylor tower of this functor consists of the fibration sequence 
$$\xymatrix@C=10pt@R=10pt{
  D_mF \ar[r] & P_mF \ar[d] \ar@{=}[r] & F\phantom{D_n} \\ 
  & P_nF \ar@{=}[r] & D_nF}$$  
The sequence may be extended to the right yielding 
$$\xymatrix@C=10pt@R=10pt{
  F \ar@{=}[r] &  P_mF \ar[d] & \\ 
  & D_nF \ar[r] & \Omega^{-1} D_mF}$$
  Given only the homogeneous two layers, the map $D_nF \to \Omega^{-1} D_mF$ tells how to put them back together and recover $F$ -- the map's homotopy fiber is naturally weakly equivalent to $F$.  We show that maps between homogeneous functors to spectra are determined by maps between their coefficients half smashed with certain surjection sets -- in this case the map $D_nF \to \Omega^{-1} D_mF$ is determined by a $\Sigma_n$-equivariant map, $\partial_{n}F \to \big(\Omega^{-1}\,\partial_{m}F\wedge \text{Sur}(\underline{m},\,\underline{n})_+\big)_{\hS_m}$, which can be explicitly constructed using cross effects.  This map is a structure map of the type which we are interested in -- knowing the two nontrivial derivatives of $F$ (the coefficients of the homogeneous layers) along with their symmetric group actions and this structure map between them, we may recover the functor $F$.

\begin{note} If we instead had $F:\top\to\top$ mapping to spaces rather than spectra, then we could still use a similar argument to figure out the additional structure required to recover $F$ from $\partial_{n}F$ and $\partial_{m}F$.  Homogeneous functors can be delooped, so the fibration sequence still extends to the right to give $P_mF\to D_nF\to \Omega^{-1}D_mF$.  In order to analyze the map 
$$\Linf(\partial_{n}F\wedge\Sinf X^{\wedge n})_{\hS_n} \to 
\Linf(\Omega^{-1}\partial_{m}F\wedge\Sinf X^{\wedge m})_{\hS_m}$$
 we may consider the adjoint map 
$$\Sinf\Linf(\partial_{n}F\wedge\Sinf X^{\wedge n})_{\hS_n} \to 
(\Omega^{-1}\partial_{m}F\wedge\Sinf X^{\wedge m})_{\hS_m}$$
of functors $\top\to\sp$.  The left-hand side is a composition of two functors to $\sp$, one of them homogeneous and the other -- $\Sinf\Linf$ -- somewhat well understood (see work e.g. of Kuhn) and rationally completely understood (recall Corollary~\ref{Q Snaith}). 
\end{note}

\section{Some Structure of Rational Taylor Towers}\label{sec:dgc to dg}

Let $\catM$ and $\catN$ be any of the categories $\dg$, $\dg_r$, $\dgl_{r-1}$, or $\dgc_r$ and
suppose $F:\catM \to \catN$ is a rational homotopy map.  We would like to analyze $P_nF$ in the Taylor tower of $F$.  Our approach is
to perform the analysis in stages.  

\begin{lemma}[G-structure]\label{g-structure}
 As a graded vector space, $P_nF$ is given by: $$[P_nF]_{\rm G} = \bigoplus_{k=1}^n
[D_kF]_{\rm G}$$
\end{lemma}
\begin{proof}
In the
category $\g$ all fibrations are split.  Thus, (assuming that $F$ is
reduced) the fibration $[D_2F]_\mG \xrightarrow{i} [P_2F]_\mG \xrightarrow{p} [D_1F]_\mG$
implies that $[P_2F]_\mG = [D_1F]_\mG \oplus [D_2F]_\mG$.

Induct.
\end{proof}

\begin{lemma}[DG-structure]\label{dg-structure}
 As a differential graded vector space, $P_nF$ is given by: $$[P_nF]_{\rm DG} = \left([D_1F]_{\rm G}
\oplus \dots \oplus [D_nF]_{\rm G},\ d_{P_nF} = \left(\begin{matrix} 
 d_{11}  &    0    &  0  \\      
  \vdots & \ddots  &  0  \\   
 d_{n1}  & \cdots  & d_{nn}
\end{matrix}\right)\ \right)$$
where $d_{ij}:[D_jF]_\mG \to [D_iF]_\mG$ is a degree -1 \G-map.
Furthermore the diagonal maps are precisely the differentials of the appropriate homogeneous layers (i.e. $D_iF = ([D_iF]_\mG,\, d_{ii})$).
\end{lemma}
\begin{proof}
$[P_2F]_\mDG = ([P_2F]_\mG, d_{P_2F}) = ([D_1F]_\mG \oplus [D_2F]_\mG, 
d_{P_2F})$.
The degree -1 map, $d_{P_2F}:[D_1F]_\mG \oplus [D_2F]_\mG \to [D_1F]_\mG \oplus [D_2F]_\mG$,
can be expressed as a $2\times 2$ matrix of maps $$d_{P_2F} = \left(\begin{matrix}d_{11} &
d_{12} \\ d_{21} & d_{22}\end{matrix}\right)$$ where the $d_{ij}$ are degree -1 maps, $d_{ij}: [D_jF]_\mG \to [D_iF]_\mG$.  Since the fiber inclusion $[D_2F]_\mDG \xrightarrow{i} [P_2F]_\mDG$ is a \DG-map, we 
must have $i\circ d_{D_2F} = d_{P_2F}\circ i = (d_{12} + d_{22})\circ i$.  Thus, 
$d_{12} = 0$ and $d_{22} = d_{D_2F}$.  A similar argument applied to the map 
$[P_2F]_\dgv \xrightarrow{p} [D_1F]_\dgv$ forces $d_{11} = d_{D_1F}$.

Induct.
\end{proof}

If $F$ is a functor $\catM\to\dg$ then this gives a classification of $F$ by the collection of its homogeneous layers and the structure maps $d_{ij}:[D_jF]_\mG \to [D_iF]_\mG$.  We would prefer, however to have structure maps go between $\dg$-objects rather than just $\g$-objects.  We will analyse the structure maps $d_{ij}$ more thoroughly.  

Note that a degree -1 map of \G-spaces $f:A \to B$ can be viewed as a degree 0 map $\hat f:A \to sB$.  We use the \,$\hat{\ }$ symbol to denote degree -1 maps viewed as degree 0 maps to suspensions.

\begin{cor} The maps $d_{i,j}$ combine to give the following structure:
\begin{itemize}
\item
The maps $\hat d_{i+1, i}$ are \DG-maps 
 $$\hat d_{i+1,i}:[D_iF]_\mDG \to s[D_{i+1}F]_\mDG$$
\item
The maps $\hat d_{i+1,i-1}$ are null homotopies of the compositions of \DG-maps 
 $$(s\hat d_{i+1,i}\circ \hat d_{i,i-1}):[D_{i-1}F]_\mDG\to s^2[D_{i+1}F]_\mDG$$
\item
The maps $\hat d_{i+2,i-1}$ are homotopies of the two null homotopies of the compositions
 $$(s^2\hat d_{i+2,i+1} \circ s\hat d_{i+1,i} \circ \hat d_{i,i-1}):
  [D_{i-1}F]_\mDG \to s^3[D_{i+2}F]_\mDG$$
\item
Etc.
\end{itemize}
\end{cor}
\begin{proof}
All of these statements come from expanding $d_{P_nF}\,d_{P_nF} = 0$ from Lemma~\ref{dg-structure}.

The first statement follows from the equation 
$$d_{i+1,i+1}\,d_{i+1,i} + d_{i+1,i}\,d_{i,i} = 0$$ 
(Note that $d_{sD_{i+1}} = -sd_{D_{i+1}} = -sd_{i+1,i+1}$.)

The second statement follows from 
$$d_{i+1,i+1}\,d_{i+1,i-1} + d_{i+1,i}\,d_{i,i-1} + d_{i+1,i-1}\,d_{i-1,i-1} = 0$$
(Note that $d_{s^2D_{i+1}} = s^2d_{D_{i+1}} = s^2d_{i+1,i+1}$.)

The third statement follows from
$$d_{i+2,i+2}\,d_{i+2,i-1} + d_{i+2,i+1)}\,d_{i+1,i-1} + d_{i+2,i}\,d_{i,i-1} +
 d_{i+2,i-1}\,d_{i-1,i-1} = 0$$

Etc.
\end{proof}

The standard way to write the above information is as a series of commuting $n$-cubes of maps of suspensions $s^k[D_nF]_\mDG$ and iterated cones $\mathrm{c}^k[D_nF]_\mDG$. 

\begin{ex} The maps between $[D_1F]_\mDG$ and $[D_4F]_\mDG$ are as follows (where $i$ denotes the inclusion of a \DG\ into its cone):

$$\xymatrix@=20pt {
[D_1F]_\mDG \ar[rr]^{\hat d_{21}} \ar[dd]^{i} \ar[dr]^{i} 
  & & s[D_2F]_\mDG \ar'[d][dd]^(.4)i \ar[dr]^{s\hat d_{32}} & \\
 & \mathrm{c}[D_1F]_\mDG \ar[dd]_(.30){i} \ar[rr]^(.35){H_{31}} 
  & & s^2[D_3F]_\mDG \ar[dd]^{s^2\hat d_{43}} \\
\mathrm{c}[D_1F]_\mDG \ar'[r][rr]_(.3){\mathrm{c}\hat d_{21}} \ar[dr]^{\mathrm{c}i} 
  & & \mathrm{c}s[D_2F]_\mDG \ar[dr]^{sH_{42}} &  \\
& \mathrm{cc}[D_1F]_\mDG \ar[rr]^{H_{41}} 
  & & s^3[D_4F]_\mDG  \\
}$$

\begin{itemize}
\item $H_{31} = \hat d_{31} + s\hat d_{32} \circ \hat d_{21}$
\item $H_{42} = \hat d_{42} + s\hat d_{43} \circ \hat d_{32}$
\item $H_{41} = \hat d_{41} + s^2\hat d_{43} \circ \hat d_{31} + s\hat d_{42}\circ \hat d_{21} + s^2\hat d_{43}\circ s\hat d_{32}\circ\hat d_{21}$
\end{itemize}

\begin{note} The above information may also be encoded dually as a $3$-cube of maps of desuspensions $s^{-k}[D_nF]_\mDG$ and iterated paths $\mathrm{p}^k[D_nF]_\mDG$.
\end{note}

\end{ex}

Furthermore any $n$-excisive rational homotopy functor $F:\catM\to\dg$ 
determines a series of commuting diagrams of maps as above.  Given such a diagram, we may then recover the $P_nF$ by either extracting the structure maps $d_{ij}$ from the maps in the diagram, or by merely taking the total homotopy fiber of the $(n-1)$-cube associated to $P_nF$.

\begin{thm}\label{Q homog classification}
This gives a classification of rational $n$-excisive functors $\catM\to\dg$ by the collection of their homogeneous layers and structure maps between their suspensions, along with a system of null homotopies of compositions.
\end{thm}

Further structure maps are required to give a classification of functors which do not map to $\dg$.

\section{Rational Jets}

We write $\Sigma^m_n$ for the set of surjections $\{1,\dots,m\} \to \{1,\dots,n\}$.
Note that $\Sigma^m_n$ has an action on the right by $\Sigma_m$ and on the left by $\Sigma_n$.  Furthermore, $\Sigma^n_n \cong \Sigma_n$ as a $\Sigma_n$ bi-module.

Suppose that we are given a set of \DG-maps $\{f_\sigma:A_n\to A_m\}_{\sigma\in\sur^m_n}$ invariant in the sense that for $\pi \in \Sigma_n$ and $\eta \in \Sigma_{m}$ we have $f_\sigma \pi = f_{\pi \sigma}$ and $\eta f_\sigma = f_{ \sigma \eta}$.  If $C$ is a \DGC\ then its coproduct map $C\to C\otimes C$ induces another similarly invariant set of maps $\{\Delta_\sigma:[C]_\mDG^{\otimes n} \to [C]_\mDG^{\otimes m}\}_{\sigma \in \sur^m_n}$.  
The sets of maps $\{f_\sigma\}$ and $\{\Delta_\sigma\}$ combine to define a map
$$\Bigl({\frac{1}{|\Sigma^m_n|}}\,\sum_{\sigma\in \sur^m_n} 
 f_\sigma \otimes \Delta_\sigma \Bigr) : A_n \otimes [C]_\mDG^{\otimes n} \longrightarrow A_m \otimes [C]_\mDG^{\otimes m}$$
with the proper equivariance properties to induce a map of homogeneous functors
$$\bigl(A_n\otimes [C]^{\otimes n}_\mDG\bigr)_{\Sigma_n} \,\cong\,
 \bigl(A_n\otimes [C]^{\otimes n}_\mDG\bigr)^{\Sigma_n} \longrightarrow
 \bigl(A_m\otimes [C]^{\otimes m}_\mDG\bigr)_{\Sigma_n}$$

\begin{defn}[Rational Jets of Functors to DG]\label{Q jet}
 A rational jet consists of a symmetric sequence of \DG\,s, $\{V_k \in \dg^{\boldsymbol{\Sigma}_k}\}_{k\ge 1}$, along with structure maps
$$f_\sigma: A_k \to s\,A_{k+1}$$ for every $k$ and $\sigma \in \sur^{k+1}_k$
invariant under $\Sigma_k$ and $\Sigma_{k+1}$ in the sense that for $\pi \in \Sigma_k$ and $\eta \in \Sigma_{k+1}$ we have $f_\sigma \pi = f_{\pi \sigma}$ and $\eta f_\sigma = f_{ \sigma \eta}$, and a system of coherent homotopies consisting of:
\begin{itemize}
\item for every $\tau \in \sur^{k+2}_k$, a null homotopy $H_\tau$ of 
$$f_\tau = \sum_{\begin{subarray}{c}
  \sigma \rho = \tau, \\
  \rho \in \sur^{k+2}_{k+1}\\ \sigma \in \sur^{k+1}_k
 \end{subarray}}
f_\rho f_\sigma
$$
(the sum of $f_\rho f_\sigma$ over all of the ways of expressing $\tau$ as a composition $\tau = \sigma \rho$ of elements $\rho \in \sur^{k+2}_{k+1}$, $\sigma \in \sur^{k+1}_k$),
invariant under $\Sigma_k$ and $\Sigma_{k+2}$ in the sense above.
\item for every $\upsilon \in \sur^{k+3}_k$, a homotopy $H_\upsilon$ of the two null homotopies of 
$$f_\upsilon = \sum_{\begin{subarray}{c}
 \sigma \rho \phi = \upsilon, \\
 \phi \in \sur^{k+3}_{k+2} \\
 \rho \in \sur^{k+2}_{k+1} \\
 \sigma \in \sur^{k+1}_k \end{subarray}} f_\phi f_\rho f_\sigma 
= \sum_{\begin{subarray}{c}
 \tau \phi = \upsilon, \\
 \phi \in \sur^{k+3}_{k+2} \\
 \tau \in \sur^{k+2}_k \end{subarray}} f_\phi f_\tau
= \sum_{\begin{subarray}{c}
 \sigma \gamma = \upsilon, \\
 \gamma \in \sur^{k+3}_{k+1} \\
 \sigma \in \sur^{k+1}_k \end{subarray}} f_\gamma f_\sigma $$
invariant under $\Sigma_k$ and $\Sigma_{k+3}$ in the sense above.
\item etc...
\end{itemize}
\end{defn}

With some patience it is possible to show that the obove collection of maps and homotopies is enough to determine (using the method outlined at the beginning of this section) a collection of maps and homotopies of homogeneous functors as in Theorem~\ref{Q homog classification}.  It is somewhat harder to show that a collection of maps and homotopies as in Theorem~\ref{Q homog classification} determines a collection of maps and homotopies as in Definition~\ref{Q jet}, and that the operations of moving between these are inverse up to weak equivalence.

Recent work of Michael Ching suggests (though we have not explicitly verified) that the structure given in Definition~\ref{Q jet} is precisely the structure given by a symmetric sequence of \DG\,s being a right module over the $Lie$ operad (see [MSS 1.13 and 1.28] for the definition of the $Lie$ operad).

\section{Some Structure of Non-Rational Taylor Towers}

We may classify functors $\top_*\to\sp$ in terms of their homogeneous layers and maps between them in much the same way as we did functors to $\dg$ in Section~\ref{sec:dgc to dg}.  Further argument is required to analyze maps between homogeneous functors and their compositions or to classify functors to $\top$.

The following argument constructs a cube of $k$-invariants given a finite tower of fibrations in any stable model category.  If $P_n \to \dots \to P_3 \to P_2 \to P_1 = D_1$ is a finite tower of fibrations with fibers $D_k \to P_k \to P_{k-1}$ in the stable model category $\mathcal{C}$, we may construct an $(n-1)$-cube of maps,
$\mathcal{X}:\mathcal{P}(\underline{n-1}) \to \mathcal{C}$ as follows:  
\begin{equation}
\mathcal{X}: T \mapsto 
 \begin{cases}
   P_n      &\text{if $T = \emptyset$} \\
   P_{\min(T)} &\text{if $T \neq \emptyset$} 
 \end{cases} \notag
\end{equation} 

\begin{ex} The case ${n=4}$:  In this case the cube $\mathcal{X}$ is 
$$
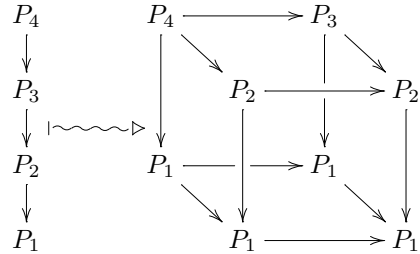
\xymatrix@=14pt{
 {\save [].[ddd]!C="b1"*[]\frm{}\restore} P_4  \ar[d] && 
 {\save [].[dddrrr]!C="b2"*[]\frm{}\restore} 
  P_4 \ar[rr] \ar[dr] \ar[dd] & & P_3 \ar'[d][dd] \ar[dr] & \\
 P_3 \ar[d] && & P_2 \ar[dd] \ar[rr] & & P_2 \ar[dd] \\
 P_2 \ar[d] && P_1 \ar'[r][rr] \ar[dr] & & P_1 \ar[dr] & \\
 P_1 && & P_1 \ar[rr] & & P_1
 \ar@{|~|>}"b1";"b2" }$$
\end{ex}

Now we can create new cube by taking the total homotopy cofibers of faces of this cube.  Define the $(n-1)$-cube $\hat{\mathcal{X}}$ by  
$$\hat{\mathcal{X}}: T \mapsto {\rm thocof}\big({\rm subcube\ under\ } \mathcal{X}(\underline {n-1}\smallsetminus T)\big)$$ 
Note that the objects in this cube are all either weakly equivalent to an
iterated delooping of a fiber in the tower or else contractible.  

\begin{ex} The case ${n=4}$.  In this case the cube $\hat{\mathcal{X}}$ is
$$
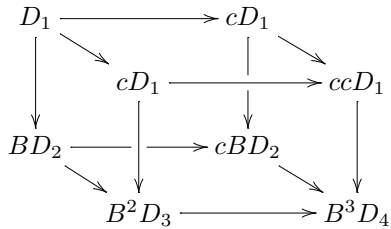
\xymatrix@=10pt{
 D_1 \ar[dd] \ar[rr] \ar[dr] && cD_1 \ar'[d][dd]\ar[dr] & \\
 & cD_1 \ar[rr]\ar[dd] && ccD_1 \ar[dd] \\
 BD_2 \ar[dr]\ar'[r][rr] && cBD_2 \ar[dr] & \\
 & B^2D_3 \ar[rr] && B^3D_4 \\
 }$$
which is weakly equivalent to the cube
$$
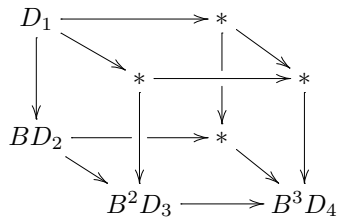
\xymatrix@=10pt{
D_1 \ar[dd] \ar[rr] \ar[dr] && \ast \ar'[d][dd]\ar[dr] & \\
 & \ast \ar[rr]\ar[dd] && \ast \ar[dd] \\
 BD_2 \ar[dr]\ar'[r][rr] && \ast \ar[dr] & \\
 & B^2D_3 \ar[rr] && B^3D_4 \\
 }$$
(using the notation $cA$ for the cone on $A$ and $BA$ for the delooping of $A$).
\end{ex}

We recover the top space of our initial tower of fibrations (up to homotopy equivalence) by taking the total homotopy fiber of this cube.

\begin{thm}\label{non-Q homog classification}
This gives a classification of $n$-excisive functors $\top_*\to\sp$ in terms of homogeneous layers and structure maps between their deloopings along with a system of null homotopies of their compositions.
\end{thm}

\section{Non-Rational Jets of Functors to Spectra}

Note that given two spectra with symmetric group actions $A_n$ and $A_m$ (where $A_i$ has a $\Sigma_i$-action) along with a $\Sigma_n$-equivariant map of spectra $A_n\to(A_m\wedge{\sur^m_n}_+)_{\hS_m}$, we can construct a map of homogeneous functors using the $(\Sigma_m\times\Sigma_n)$-equivariant map $X^{\wedge n} \wedge {\sur^m_n}_+ \to X^{\wedge m}$ induced by the diagonal $X\to X\wedge X$ as follows:  Consider the composition
$$\xymatrix@=15pt{
A_n\wedge X^{\wedge n} \ar[r] & (A_m\wedge{\sur^m_n}_+)_{\hS_m} \wedge X^{\wedge n} & \\
 & (A_m\wedge X^{\wedge n}\wedge{\sur^m_n}_+)_{\hS_m} \ar@{=}[u] \ar[r] &
   (A_m\wedge X^{\wedge m})_{\hS_m} 
}$$
Each of the above maps are $\Sigma_n$-equivariant, so their composition induces a map of homogeneous functors $(A_n\wedge X^{\wedge n})_{\hS_n} \longrightarrow  (A_m\wedge X^{\wedge m})_{\hS_m}$.

If $F$ is a functor $\top_*\to\sp$ then the jet of $F$ consists of the following data:

\begin{defn}[Jets of Functors to $\sp$]\label{non-Q jet}
A\, $jet$ is a symmetric sequence of spectra $\{A_k\}_{k\ge 1}$ (referred to as ``coefficients''), along with $\Sigma_k$-equivariant maps (referred to as ``structure maps'')
$$f_k: A_k \to \big(\Omega^{-1}\,A_{k+1} \wedge {\sur^{k+1}_k}_+\big)_{\hS_{k+1}}$$
for all $k$ and a system of coherent homotopies consisting of:
\begin{itemize}
\item null homotopies $H^{k+2}_k$ of the compositions 
$$
\xymatrix@=20pt { 
A_k  \ar[d]^{f_k} & \big(\Omega^{-2}\,A_{k+2} \wedge {\sur^{k+2}_k}_+\big)_{\hS_{k+2}} \\
\big(\Omega^{-1}\,A_{k+1} \wedge {\sur^{k+1}_k}_+\big)_{\hS_{k+1}} \ar[r]^(.40){f_{k+1}} &
\big(\Omega^{-2}\,A_{k+2} \wedge {\sur^{k+2}_{k+1}}_+ \wedge {\sur^{k+1}_k}_+\big)_{\hS_{k+1}\times\Sigma_{k+2}} \ar[u]^{({\it fold})}
}$$
\item a homotopy $H^{k+3}_k$ between the two null homotopies of compositions $$A_k \longrightarrow \big(\Omega^{-3}\,A_{k+3} \wedge {\sur^{k+3}_k}_+\big)_{\hS_{k+3}}$$
\item etc...
\end{itemize}
\end{defn}

With some patience it is possible to show that the obove collection of maps and homotopies is enough to determine (using the method outlined at the beginning of this section) a collection of maps and homotopies of homogeneous functors as in Theorem~\ref{non-Q homog classification}.  It is somewhat harder to show that a collection of maps and homotopies as in Theorem~\ref{non-Q homog classification} determines a collection of maps and homotopies as in Definition~\ref{non-Q jet}, and that the operations of moving between these are inverse up to weak equivalence.



 

\providecommand{\bysame}{\leavevmode\hbox to3em{\hrulefill}\thinspace}
\providecommand{\MR}{\relax\ifhmode\unskip\space\fi MR }
\providecommand{\MRhref}[2]{%
  \href{http://www.ams.org/mathscinet-getitem?mr=#1}{#2}
}
\providecommand{\href}[2]{#2}

\end{document}